\RequirePackage{fix-cm}
\documentclass[11pt]{article}
\usepackage[T1]{fontenc}
\usepackage{amsmath,amssymb,amsthm,mathtools}
\usepackage{erewhon}
\usepackage[erewhon]{newtxmath}
\usepackage{fullpage}
\usepackage{xurl}
\usepackage{microtype}
\usepackage{booktabs,array,longtable}
\usepackage{xcolor}
\usepackage{enumitem}
\usepackage[colorlinks=true,linkcolor=blue!45!black,citecolor=blue!45!black,urlcolor=blue!55!black]{hyperref}
\usepackage[nameinlink,noabbrev]{cleveref}
\hypersetup{pdftitle={Fisher--Hartwig asymptotics and critical moments for the circular beta ensemble},pdfauthor={Reda Chhaibi, Pierre-Loïc Méliot, Joseph Najnudel, Ashkan Nikeghbali},pdfsubject={General beta Fisher--Hartwig formulas, counting statistics, smooth backgrounds and critical moments}}
\numberwithin{equation}{section}
\newtheorem{mainthm}{Theorem}

\crefname{mainthm}{Theorem}{Theorems}
\newtheorem{theorem}{Theorem}[section]
\newtheorem{proposition}[theorem]{Proposition}
\newtheorem{lemma}[theorem]{Lemma}
\newtheorem{corollary}[theorem]{Corollary}

\theoremstyle{remark}
\newtheorem{remark}[theorem]{Remark}

\crefname{theorem}{Theorem}{Theorems}
\crefname{proposition}{Proposition}{Propositions}
\crefname{lemma}{Lemma}{Lemmas}
\crefname{corollary}{Corollary}{Corollaries}
\crefname{conjecture}{Conjecture}{Conjectures}
\crefname{assumption}{Assumption}{Assumptions}
\crefname{definition}{Definition}{Definitions}
\crefname{remark}{Remark}{Remarks}
\crefname{problem}{Problem}{Problems}
\crefname{section}{Section}{Sections}
\crefname{subsection}{Section}{Sections}
\crefname{subsubsection}{Section}{Sections}
\newcommand{\E}{\mathbb E}
\newcommand{\Prob}{\mathbb P}
\newcommand{\T}{\mathbb T}
\newcommand{\D}{\mathbb D}

\newcommand{\R}{\mathbb R}
\newcommand{\C}{\mathbb C}
\newcommand{\Z}{\mathbb Z}

\newcommand{\dd}{\,\mathrm d}

\newcommand{\Var}{\operatorname{Var}}

\newcommand{\one}{\mathbf 1}

\setlist{itemsep=3pt,topsep=5pt}
\title{Fisher--Hartwig asymptotics and critical\\ moments for the circular beta ensemble}
\author{Reda Chhaibi \quad Pierre-Loïc Méliot\\[3pt]
Joseph Najnudel \quad Ashkan Nikeghbali}
\date{1 November 2026}

\begin{document}

\maketitle
\begin{abstract}
We determine the joint moments of real powers of the absolute value of
characteristic polynomials in the circular beta ensemble, for every
$\beta>0$. At fixed distinct prescribed points, the Fisher--Hartwig
formula holds for every exponent greater than $-1$, with explicit
constants and a relative error bounded by a negative power of the
matrix size $N$. For approaching points in the same integrable range,
the error is bounded by a negative power of $N\delta_N$, provided that $N\delta_N\to\infty$
and the sum of the absolute values of the negative exponents is less
than two. Here $\delta_N$ is the minimum distance between the
prescribed points; this range includes every integrable pair of root
exponents.

For $q>0$ at the critical relation $pq^2=2\beta$, we prove the
Keating--Wong asymptotic for the $p$th moment of the integral of the
$q$th absolute power of the characteristic polynomial, for every
real $p>1$. The leading asymptotic extends to every complex spatial
weight in $L^p$. The proofs combine the recursion for orthogonal
polynomials on the unit circle with conditional expectation
identities for Gaussian multiplicative chaos and quantitative
estimates for approaching singularities.
\end{abstract}

\section{Introduction and main results}
\label{sec:intro}

The Fisher--Hartwig problem concerns the effect of prescribed
singularities on a partition function. In the Circular Unitary Ensemble,
the partition function is a Toeplitz determinant, and its asymptotics
belong to the theory initiated by Fisher and Hartwig
\cite{FisherHartwig1968}. The formula contains a power of the matrix
size, a constant at each singularity, and an interaction between pairs
of singularities.

In \cref{thm:main-root}, we prove the Fisher--Hartwig formula for
every $\beta>0$ and arbitrary integrable real root exponents at fixed
distinct points. We also obtain quantitative estimates when the
points approach one another, in the exponent domains specified
below. Our second principal result, \cref{thm:main-critical}, gives
the critical moment asymptotic conjectured by Keating and Wong for
every real order $p>1$.

The two results are connected by the contribution of configurations
in which prescribed points approach one another. At criticality,
contributions from successive scales accumulate to produce a
logarithm. For integer moments these contributions can be expressed
through joint moments of the characteristic polynomial. The
noninteger argument requires a consistent change of measure with
several marked points and quantitative estimates for nearby marked
points. We explain this distinction before \cref{thm:main-critical}
and the common method in \cref{sec:intro-reductions}.

\subsection{Fisher--Hartwig asymptotics for every \texorpdfstring{$\beta>0$}{beta > 0}}
\label{sec:intro-FH}

Fix $\beta>0$, put $b=\beta/2$, and let
$\sigma(d\theta)=d\theta/(2\pi)$ be normalized Haar measure on
$\T=\R/(2\pi\Z)$. The circular beta ensemble with $N$ particles has law
\begin{equation}\label{eq:ensemble}
 \frac1{Z_{N,\beta}}
 \prod_{1\le r<s\le N}|e^{i\theta_r}-e^{i\theta_s}|^\beta
 \prod_{r=1}^N\sigma(d\theta_r),
 \qquad Z_{N,\beta}=\frac{\Gamma(1+Nb)}{\Gamma(1+b)^N}.
\end{equation}
We normalize the characteristic polynomial to equal one at the origin:
\begin{equation}\label{eq:PN}
 P_N(z)=\prod_{r=1}^N(1-ze^{-i\theta_r}),\qquad P_N(0)=1.
\end{equation}
For prescribed points $\boldsymbol\vartheta$ and exponents $\mathbf q$,
we study the joint moments
\begin{equation}\label{eq:FN}
 F_N(\mathbf q,\boldsymbol\vartheta)
 =\E_\beta\prod_{j=1}^m|P_N(e^{i\vartheta_j})|^{q_j},
 \qquad M_N(q)=F_N(q,0).
\end{equation}
Here $\vartheta_j$ are the prescribed locations of the root singularities,
and $q_j$ their exponents. In the ensemble formulas, $\theta_1,\ldots,\theta_N$
denote the particles, that is, the eigenangles for random matrix models, and
$\vartheta_1,\ldots,\vartheta_m$ the prescribed points; the weight
$|P_N(e^{i\vartheta_j})|^{q_j}=\prod_r|e^{i\vartheta_j}-e^{i\theta_r}|^{q_j}$
is a power of the product of the distances from the prescribed point to the particles.
In integrals of deterministic test functions or random measures on
$\T$, the letter $\theta$ also denotes a variable of integration.
For fixed distinct particle locations and real $q_j$, this weight
is integrable over the prescribed point $\vartheta_j$ exactly when
$q_j>-1$. Recall that
$M_N(q)=\E_\beta|P_N(1)|^q$. 

Its asymptotic normalization is
\[
 C_\beta(q)=\lim_{N\to\infty}N^{-q^2/(2\beta)}M_N(q)>0.
\]

The Morris product gives this constant explicitly; see
\eqref{eq:constantseries} and \eqref{eq:doublegamma-constant}.
The leading factor separates the one point constants from the interaction
between pairs of prescribed points:
\begin{equation}\label{eq:complex-root-leading-factor}
 \mathcal A_\beta(\mathbf q,\boldsymbol\vartheta)
 =\prod_jC_\beta(q_j)
   \exp\left\{-\frac1\beta\sum_{j<k}q_jq_k
                   \log|e^{i\vartheta_j}-e^{i\vartheta_k}|\right\}.
\end{equation}
For fixed distinct angles, \eqref{eq:complex-root-leading-factor}
defines a holomorphic function of $\mathbf q$ on the domain 
$$\{\mathbf q\in\C^m:\Re q_j>-1\text{ for }1\le j\le m\}.$$
Here each logarithm is the ordinary real logarithm of the positive
distance $|e^{i\vartheta_j}-e^{i\vartheta_k}|$; analyticity of
$C_\beta$ on this half-plane is proved in \cref{prop:onepoint}.
The total negative exponent determines the range in which the
constants in the error estimate can be chosen uniformly as the
prescribed points approach one another. The admissible configurations and the possible exponents of
the error are specified separately below. 
With $x_+=\max(x,0)$, we introduce the two corresponding domains.
We also write $\delta(\boldsymbol\vartheta)$ for the minimum chord
distance between the prescribed points $e^{i\vartheta_j}$, capped at one, and $\mathcal S_\delta$ for configurations
whose distinct points are at least $\delta$ apart:
\begin{align}\label{eq:intro-domain}
 \mathcal D_m&=\left\{\mathbf q\in\R^m:
                              \sum_j(-q_j)_+<1\right\},\notag\\
 \mathcal E_m&=\left\{\mathbf q\in(-1,\infty)^m:
                              \sum_j(-q_j)_+<2\right\},\notag\\
 \delta(\boldsymbol\vartheta)
   &=\min\{1,\min_{j<k}|e^{i\vartheta_j}-e^{i\vartheta_k}|\},
 \qquad
 \mathcal S_\delta=\{\boldsymbol\vartheta:
         |e^{i\vartheta_j}-e^{i\vartheta_k}|\ge\delta\ (j\ne k)\}.
\end{align}
For $m=1$ put $\delta(\boldsymbol\vartheta)=1$ and omit separation
conditions. The domain $\mathcal D_m$ is contained in $\mathcal E_m$
and contains all nonnegative exponents.

The shape of the formula can be anticipated from the Gaussian
structure of $\log|P_N|$. The one point moment
$M_N(q)\approx C_\beta(q)N^{q^2/(2\beta)}$ is, up to the constant,
the moment generating function of a centered Gaussian variable of
variance $\beta^{-1}\log N$, and $\log|P_N(1)|$ does satisfy a central
limit theorem with this variance \cite{KeatingSnaith2000,BHNY2008}.
More precisely, at distances larger than the mean spacing $1/N$, the
random function $\vartheta\mapsto\log|P_N(e^{i\vartheta})|$ behaves
like a centered Gaussian field $X$ on the circle with covariance
\[
 \E X(\vartheta)X(\varphi)
 =-\frac1\beta\log|e^{i\vartheta}-e^{i\varphi}|,
\]
truncated at variance $\beta^{-1}\log N$ on the diagonal. For such a
field,
\[
 \E\exp\Bigl\{\sum_jq_jX(\vartheta_j)\Bigr\}
 =\exp\Bigl\{\frac12\Var\sum_jq_jX(\vartheta_j)\Bigr\},
\]
and the diagonal terms of the variance give the power
$N^{\sum_jq_j^2/(2\beta)}$, and the cross terms give the pair factors
$|e^{i\vartheta_j}-e^{i\vartheta_k}|^{-q_jq_k/\beta}$ in
\eqref{eq:complex-root-leading-factor}. The constants $C_\beta(q_j)$
are the corrections to the Gaussian picture at each prescribed point.
They are exactly the one point constants. Thus, at leading order in
the regimes of \cref{thm:main-root}, the interaction between the
prescribed points is given by the Gaussian pair factor. The proof
makes this heuristic precise through the multiplicative chaos of the
limiting field; see \cref{sec:gmc-identification}.

\begin{mainthm}[The Fisher--Hartwig formula at separated and moving points]
\label{thm:main-root}
Let $\beta>0$ and let $K\subset(-1,\infty)^m$ be compact.
For every fixed $\delta>0$, there is $\kappa>0$, depending only on
$\beta,m,K$, such that, uniformly on $K\times\mathcal S_\delta$,
\begin{equation}\label{eq:intro-FH}
 F_N(\mathbf q,\boldsymbol\vartheta)
 =N^{\frac1{2\beta}\sum_jq_j^2}
       \mathcal A_\beta(\mathbf q,\boldsymbol\vartheta)
       \bigl(1+O(N^{-\kappa})\bigr).
\end{equation}
For approaching points, if $K\subset\mathcal E_m$,
\begin{equation}\label{eq:intro-moving-signed}
 F_N(\mathbf q,\boldsymbol\vartheta)
 =N^{\frac1{2\beta}\sum_jq_j^2}
       \mathcal A_\beta(\mathbf q,\boldsymbol\vartheta)
       \bigl(1+O((N\delta(\boldsymbol\vartheta))^{-\kappa})\bigr),
\end{equation}
for some $\kappa>0$ and $R\ge1$ depending only on $\beta,m,K$, uniformly
in the exponents and points whenever $N\delta(\boldsymbol\vartheta)\ge R$.
For $m=2$, this covers every $q_1,q_2>-1$.
For one point the relative error is $O(N^{-1})$.
The implicit constants may depend on $\beta,m,K$ and the chosen
$\kappa$; at fixed separation they may also depend on $\delta$.
The constants in the moving point estimates are independent of the
separation and the prescribed points.
\end{mainthm}

At $\beta=2$ the expectation in \eqref{eq:FN} is a Toeplitz
determinant, by Heine's identity, and \cref{thm:main-root} takes a
classical form. Write $G$ for the Barnes $G$-function.

\begin{corollary}[Toeplitz determinants with root singularities]
\label{cor:widom}
Let $\beta=2$ and $q_j=2\alpha_j$ with $\alpha_j>-1/2$. For fixed
distinct points $\vartheta_1,\ldots,\vartheta_m$,
\[
 \E_2\prod_{j=1}^m|P_N(e^{i\vartheta_j})|^{2\alpha_j}
 =N^{\sum_j\alpha_j^2}
  \prod_{j=1}^m\frac{G(1+\alpha_j)^2}{G(1+2\alpha_j)}
  \prod_{j<k}|e^{i\vartheta_j}-e^{i\vartheta_k}|^{-2\alpha_j\alpha_k}
  \bigl(1+O(N^{-\kappa})\bigr),
\]
uniformly for $(\alpha_1,\ldots,\alpha_m)$ in a compact subset of
$(-1/2,\infty)^m$ and for points whose separation is bounded below by
a fixed positive constant. For approaching points, the conclusions of
\cref{thm:main-root} apply with $q_j=2\alpha_j$, with the corresponding
exponent domains and relative errors.
\end{corollary}

The left-hand side is the Toeplitz determinant with symbol
$\prod_j|z-e^{i\vartheta_j}|^{2\alpha_j}$, and the formula at fixed
points is Widom's theorem \cite{Widom1973}. For two merging singularities, Claeys and Krasovsky
\cite[Theorem~1.11]{ClaeysKrasovsky2015} obtained uniform asymptotics.
Fahs \cite[Theorem~1.1 and Remark~1.2]{Fahs2021} obtained estimates
for any fixed number of singularities, with an error bounded uniformly
in the locations at the level of the logarithm. For every
$\beta>0$, \cref{thm:main-root} establishes the formula throughout
the integrable root range at fixed separation, together with the
stated uniformity when the points approach one another.

The two separation regimes in \cref{thm:main-root} are summarized below.
\begin{center}
\begin{tabular}{@{}lll@{}}
\toprule
Locations & Exponents in a compact set $K$ & Relative error\\
\midrule
$\delta\ge\delta_0>0$ & $K\subset(-1,\infty)^m$ & $O(N^{-\kappa})$\\
$N\delta\to\infty$ & $K\subset\mathcal E_m$ & $O((N\delta)^{-\kappa})$\\
\bottomrule
\end{tabular}
\end{center}
The second row includes $\mathcal D_m\subset\mathcal E_m$.
The smaller domain $\mathcal D_m$ has a separate quantitative role:
direct Fourier summation gives the explicit admissible powers in
\eqref{eq:intro-rate} and the refinements in
\cref{cor:terminal-subset-rate}. On the larger domain $\mathcal E_m$,
\cref{thm:moving-signed-two} gives some positive power. For example,
$(-3/4,-3/4)\in\mathcal E_2\setminus\mathcal D_2$.
The rates are compared in \cref{rem:admissible-exponents}.
In the first row, the error constant may depend on $\delta_0$;
the second row is uniform once $N\delta\ge R$.
The one point error $O(N^{-1})$ follows from the classical Morris
product; its first correction is given in
\cref{prop:onepoint-first-correction}.

\subsection{Critical moments and the Keating--Wong conjecture}
\label{sec:intro-critical}

For a spatial weight $f$, define
\[
 \mathcal Z_N(q;f)=\int_\T f(\vartheta)|P_N(e^{i\vartheta})|^q
                            \,\sigma(d\vartheta),
 \qquad \mathcal Z_N(q)=\mathcal Z_N(q;1).
\]
All $L^p$ spaces and norms on the circle are taken with respect to the Haar measure 
$\sigma$, unless another measure is specified.
For $q>0$, the relation $pq^2=2\beta$ is the threshold at which the
$p$th moment of $\mathcal Z_N(q)$ ceases to be of order $M_N(q)^p$.
For integer $p\ge2$, this moment is the integral of
$F_N(q,\ldots,q;\boldsymbol\vartheta)$ over $\T^p$. The following
heuristic restricts to configurations of diameter between $r$ and
$2r$, with all pair distances comparable to $r$. On mesoscopic scales
$N^{-1}\ll r\ll1$, \cref{thm:main-root} gives a contribution to
$\E\mathcal Z_N(q)^p/M_N(q)^p$ of order
\[
 r^{p-1}\cdot r^{-p(p-1)q^2/(2\beta)}=r^{(p-1)(1-pq^2/(2\beta))}.
\]
The factors are the volume and the $\binom p2$ pair contributions
$r^{-q^2/\beta}$. Below the threshold, the exponent is positive,
suggesting that separated configurations dominate. At criticality,
it vanishes: successive scales give comparable contributions, and
integrating $dr/r$ from $1/N$ to $1$ produces $\log N$. The power of
$N$ is then $N^{pq^2/(2\beta)}=N$.

The bounds in \cref{lem:moments-global-bound,lem:critical-kernel}
control smaller internal gaps and the microscopic cutoff. For integer
$p$, the Fisher--Hartwig kernel determines the leading coefficient;
the microscopic profile enters the finite term of order $N$.
For noninteger $p$, we instead use the exact mass law
\eqref{eq:fractional-exact-mass}. This mass is the limit in law of
$\mathcal Z_N(q)/M_N(q)$: it is the total mass of a subcritical circular
Gaussian multiplicative chaos measure, with infinite $p$th moment.

\begin{mainthm}[Critical moments]
\label{thm:main-critical}
Let $\beta>0$, $p>1$, $q=\sqrt{2\beta/p}$ and
$A_p=(p-1)/\Gamma(1-1/p)^p$. Then
\begin{equation}\label{eq:main-critical-all-real}
 \mathbb E\mathcal Z_N(q)^p
 \sim A_pC_\beta(q)^pN\log N.
\end{equation}
For every complex $f\in L^p(\T,\sigma)$,
\[
 \lim_{N\to\infty}\frac{\E|\mathcal Z_N(q;f)|^p}{N\log N}
   =A_pC_\beta(q)^p\int_\T |f|^p\,d\sigma.
\]

For an integer $p=k\ge2$, if $f\ge0$ belongs to $L^k(\T,\sigma)$ and
satisfies
\[
 \int_0^1\sup_{|h|\le r}\|f(\cdot+h)-f\|_k\,\frac{dr}{r}<\infty,
\]
then
\[
 \E\mathcal Z_N(q;f)^k
 =A_kC_\beta(q)^k\left(\int_\T f^k\,d\sigma\right)N\log N
                  +B_{\beta,k}[f]N+o(N).
\]
The coefficient $B_{\beta,k}[f]$ is given by the absolutely convergent
formula \eqref{eq:weighted-critical-B}.
\end{mainthm}

At $\beta=2$ the constants are explicit.

\begin{corollary}[The unitary case]
\label{cor:critical-unitary}
Let $\beta=2$ and $p>1$, so that $q=2/\sqrt p$. Then
\[
 \E\mathcal Z_N(2/\sqrt p)^p
 \sim\frac{p-1}{\Gamma(1-1/p)^p}
     \left(\frac{G(1+1/\sqrt p)^2}{G(1+2/\sqrt p)}\right)^{p}N\log N.
\]
In particular, for $p=2$,
\[
 \E\left(\int_\T|P_N(e^{i\vartheta})|^{\sqrt2}\,\sigma(d\vartheta)\right)^2
 \sim\frac1\pi\,\frac{G(1+1/\sqrt2)^4}{G(1+\sqrt2)^2}\,N\log N.
\]
\end{corollary}

For integer $p$ this is the theorem of Keating and Wong
\cite[Theorem~2]{KW2022}, and for $p=2$ the logarithm with its
coefficient also follows from the merging asymptotics of Claeys and
Krasovsky \cite[Theorem~1.15]{ClaeysKrasovsky2015}. The noninteger
orders are new also in the unitary case.

Equation \eqref{eq:main-critical-all-real} proves Keating--Wong's
Conjecture~2 for every $\beta>0$ and every real $p>1$, and
\cref{cor:critical-unitary} completes their Conjecture~1.
In their notation $k=p$ and $2s=q$, so their critical relation
$k=\beta/(2s^2)$ is exactly $pq^2=2\beta$.
Their one point constant, written using $\Upsilon^\beta$ in
\cite[equation~(17)]{KW2022}, is $C_\beta(q)$: the expression is
obtained there from the same normalized Morris product. Thus both
the critical relation and the conjectured coefficient agree.
The locators here refer to the published article.
For regularized subcritical Gaussian chaos, Keating and Wong
\cite[Theorem~1 and Corollary~1]{KW2022} proved the analogous
critical moment asymptotic for nonnegative continuous weights.
The leading weighted asymptotic in \cref{thm:main-critical}
holds for the C$\beta$E, every real $p>1$, and complex weights in $L^p$.
Here $pq^2=2\beta$ is criticality of the $p$th moment, whereas
$q^2<2\beta$ means that the limiting Gaussian multiplicative chaos
itself is subcritical.

\subsection{The main reductions in the proofs}
\label{sec:intro-reductions}

We use orthogonal polynomials on the unit circle (OPUC), generated by
the Szeg\H{o} recursion from independent Verblunsky coefficients.
The reversed polynomial $\Phi_n^*$ becomes a characteristic polynomial
with $N=n+1$ particles after adding a final independent coefficient on
the unit circle; the exact relation is \eqref{eq:terminal-polynomial}.
The proof of \cref{thm:main-root} has three steps. We first prove
convergence of the normalized moments of $\Phi_n^*$. We then identify
their limit using the common Gaussian field and its multiplicative
chaos. Finally, we integrate the last coefficient to recover the
characteristic polynomial moment $F_N$. 
The first two steps are
\cref{thm:holomorphic-limit,thm:opuc-full}; the last is proved in
\cref{sec:terminal,sec:full-negative-transfer}. 
The exact identity relating the two moments, together with the
normalized limits that suffice to prove the theorem, is given at the
start of \cref{sec:cmv-unified}.

The proof requires asymptotic moments with their constants.
A correction involving finitely many oscillatory functions of the
auxiliary phases makes the conditional error in the recursion summable and proves holomorphic convergence of the normalized
moments. The conditional chaos identities of \cite{CN2026,LN2026}
identify this limit for small real exponents; holomorphic continuation
then gives the larger domain. More precisely, for $q>-1$ with
$q^2<2\beta$ and every bounded continuous test function $f$,
\[
 \int_\T f(\theta)
       \frac{|\Phi_n^*(e^{i\theta})|^q}{\nu_n(q)}\,\sigma(d\theta)
 =\E[U^{(q)}(f)\mid\mathcal F_n].
\]
Here $\nu_n(q)=\E|\Phi_n^*(1)|^q$, $\mathcal F_n$ is generated by
the first $n$ coefficients, and $U^{(q)}$ is the limiting chaos
measure. Thus the normalized polynomial measures are conditional
expectations of that random measure; see \cref{prop:chaos-input}.

The final step integrates the common independent angle in
\eqref{eq:terminal-polynomial}. Negative exponents make the resulting
conditional integral singular when auxiliary phases become close.
On $\mathcal D_m$, this integral is bounded uniformly in the phase
variables and its Fourier coefficients satisfy a weighted summability
estimate; see \cref{lem:terminal-holder,lem:terminal-weighted-fourier}.
Direct Fourier summation gives the explicit powers in
\eqref{eq:intro-rate} and the refinements of
\cref{cor:terminal-subset-rate}. At fixed separation, a bound for the
density of the weighted phase law treats all integrable root
exponents, as proved in \cref{thm:full-root-FH}.

On the larger domain $\mathcal E_m$, the conditional integral can be
unbounded as phases coincide. For moving points, the density bound
for differences of two phases in \cref{prop:moving-pair-density}
controls these singularities after averaging. We discard a small
neighborhood of the coincidences, truncate the singular factors and
estimate a finite Fourier approximation. Balancing these
approximations gives the power error in \cref{thm:moving-signed-two}.
Throughout $\mathcal E_m$, the condition $N\delta_N\to\infty$
suffices, with no prescribed rate of divergence. Outside
$\mathcal E_m$, \cref{thm:full-signed-polynomial} gives an explicit
range of polynomial merging scales. The density estimates and the
justification of integration under the singular weight occupy
\cref{sec:full-negative-transfer}.

For \cref{thm:main-critical}, the main distinction is between integer
and noninteger moment orders. At integer orders, the circular Jacobi
limit \cite{LiValko2022,Liu2017,AN2026} and the uniform formula for
separated points describe the two parts of a single moment integral.

For noninteger $p$, write $p=k+s$ with $k=\lfloor p\rfloor$ and
$0<s<1$. We mark $k$ points and apply conditional concavity only to
the remaining power $s$. A consistent Palm measure and estimates for
nearby marked points retain the full logarithmic coefficient. The argument uses
the circular mass law predicted in \cite{FyodorovBouchaud2008} and
proved in \cite{Remy2020}, together with conditional comparison for
the orthogonal polynomials and the factor arising from the final coefficient. The tail
method belongs to the GMC framework developed in
\cite{RhodesVargas2019,Wong2020}; the exact circular law fixes the
constant used here. Weak convergence of the normalized mass by itself
does not determine this divergent moment.
The two arguments are given in
\cref{sec:fractional-critical,sec:critical-moments}.

\subsection{Further estimates and extensions}
\label{sec:intro-additional}

\begin{remark}[Admissible error exponents]
\label{rem:admissible-exponents}
We can make the admissible error exponents explicit. The quantities
$A_K$ and $\eta_K$ below measure the total negative exponent and the
distance from the integrability threshold; $g_\beta$ is the factor
depending only on $\beta$ in the rate. Put
\[
 A_K=\sup_{\mathbf q\in K}\sum_{j=1}^m(-q_j)_+,
 \qquad
 \eta_K=\min\{1,\inf_{\mathbf q\in K,\,1\le j\le m}(1+q_j)\},
 \qquad g_\beta=\min\{4/\beta,1\}.
\]
At fixed separation, \eqref{eq:intro-FH} holds for every
\[
 0<\kappa<\frac{\eta_K g_\beta}{2+(A_K-1)_++\eta_K}.
\]
For $K\subset\mathcal D_m$, the rate uniform in the separation holds
for every
\begin{equation}\label{eq:intro-rate}
 \varepsilon_K=1-A_K,
 \qquad 0<\kappa<\varepsilon_K\min\{2/\beta,1/2\}.
\end{equation}
Combining the two bounds, when $K\subset\mathcal D_m$ and the
separation is fixed one may take any
\[
 0<\kappa<g_\beta\max\left\{
       \frac{\eta_K}{2+(A_K-1)_++\eta_K},
       \frac{1-A_K}{2}\right\}.
\]
The subset sums of the individual exponents can give a larger
admissible rate, as specified in \cref{cor:terminal-subset-rate}.
The estimate on $\mathcal E_m$ in \cref{thm:moving-signed-two}
provides another admissible power. On $\mathcal D_m$ these estimates
coexist, and one may use the better rate for the chosen compact set.
\end{remark}

\paragraph{Integer refinements and microscopic matching.}
For integer $p$, the Fisher--Hartwig formula and the kernel estimates
give the coefficient of $N\log N$ in \cref{thm:critical-moments}.
Matching the microscopic profile to that formula then gives the
finite part of order $N$ in \cref{thm:critical-second-term}.
The Dini condition in \cref{thm:main-critical} includes indicators
of finite unions of arcs.
The matching formula in \cref{thm:positive-merging} also determines
the joint positive moments of the Hua--Pickrell stochastic zeta
function at widely separated arguments, with explicit constants
and a uniform power error.
\Cref{thm:critical-window} also determines the transition through
the critical exponent and the next term uniformly in its bounded window.

\paragraph{Subcritical moments.}
For $q>0$, $p>1$ and $pq^2<2\beta$, the same normalizations give
\[
 \E\mathcal Z_N(q)^p\sim C_\beta(q)^pN^{pq^2/(2\beta)}
 \frac{\Gamma(1-pq^2/(2\beta))}
      {\Gamma(1-q^2/(2\beta))^p}.
\]
For integer $p$, this follows from \cref{thm:main-root}, the global
bound of \cref{lem:moments-global-bound}, and the Dyson integral
\eqref{eq:dyson-negative-power}. For general real $p$, the convergence
in \cite[Theorem~1]{LN2026} identifies the limit in law of
$\mathcal Z_N(q)/M_N(q)$ as the total mass of circular chaos.
Write $U^{(q)}$ for this limiting measure, and let
$\Phi_n^*$ be the reversed orthogonal polynomial introduced in
\cref{sec:coupling}, with $\nu_n(q)=\E|\Phi_n^*(1)|^q$.
To pass to moments, choose $r>p$ with $rq^2<2\beta$. The conditional
identity \eqref{eq:gmc-exact-conditional} and Jensen's inequality give
\[
 \sup_n\E\left(\int_\T
       \frac{|\Phi_n^*(e^{i\theta})|^q}{\nu_n(q)}\,\sigma(d\theta)
                      \right)^r
 \le\E U^{(q)}(\T)^r<\infty.
\]
The terminal representation \eqref{eq:terminal-polynomial}, with
$|1-e^{iu}|^q\le2^q$ and $M_N(q)=c(q)\nu_{N-1}(q)$, transfers this
bound to $\mathcal Z_N(q)/M_N(q)$. Its $p$th powers are therefore
uniformly integrable. The circular chaos mass formula of
\cite[Theorem~1.1]{Remy2020} and \eqref{eq:onepointasymp} give the
displayed asymptotic. Thus this consequence also follows from the
existing chaos convergence and moment identities.

\paragraph{Argument weights and analytic backgrounds.}
\Cref{thm:mixed-FH} includes
real exponential weights $e^{s_j\Im\log P_N}$, with arbitrary
real $s_j$, and local complex continuation. Genuine Fisher--Hartwig
jumps have $s_j=i\tau_j$ with $\tau_j$ real. Our continuation treats
these only in a neighborhood of zero. It gives a local version of the root and jump formula in
\cite[Conjecture~2]{FF2004}, without background.
For eigenvalue counts in a fixed arc, it gives the limiting transform
and higher cumulants, and extends the Gaussian moderate deviation
conclusion for CUE in the unpublished manuscript
\cite[Corollary~4.12 and Section~8]{KLMN2026}
to every $\beta>0$; see
\cref{cor:fixed-arc-counting,cor:fixed-arc-moderate-deviations}.
The Gaussian limit of the count itself is classical: Wieand
\cite{Wieand2002} proved it for fixed arcs in the CUE, and Killip
\cite{Killip2008} for every $\beta>0$. The earlier theorem of Costin
and Lebowitz \cite{CostinLebowitz1995} concerns counts in long
intervals of the limiting bulk process.
In \cref{thm:analytic-background} we also multiply the integrand by
$\exp\{\sum_{r=1}^N V(\theta_r)\}$, where $V$ is real analytic on the
circle. We refer to this factor as an analytic background; to obtain
a probability law one also divides by its expectation. It establishes the corresponding background case of
Conjecture~1; its full proposed smoothness range is a further question.

\subsection{Background and related work}
\label{sec:intro-background}

In the unitary case, Widom \cite{Widom1973} proved the
Fisher--Hartwig formula for real root singularities, Basor \cite{Basor1978} for purely imaginary jump
parameters, B\"ottcher and Silbermann \cite{BottcherSilbermann1985}
for root and jump parameters of real part less than $1/2$ in modulus,
and Ehrhardt \cite{Ehrhardt2001} in the whole range where a single
representation holds. Basor and Tracy \cite{BasorTracy1991} formulated
the sum over competing representations needed for general jumps.
Deift, Its and Krasovsky \cite{DIK2011} established
the asymptotics with a smooth background in the unitary setting, and
surveyed the history in \cite{DIK2013};
Claeys and Krasovsky \cite{ClaeysKrasovsky2015} treated two
merging singularities, while Fahs \cite{Fahs2021} obtained estimates
uniform in the locations of any fixed number of singularities. Forrester and Frankel \cite[Conjectures~1 and~2]{FF2004}
proposed corresponding formulas for circular beta ensembles. Exact
Jack polynomial expressions are available for certain products of
characteristic polynomials at general $\beta$
\cite{Forrester2010,Matsumoto2007,DesrosiersLiu2014}. These polynomial
identities do not directly give the joint asymptotics for arbitrary
real powers of the moduli considered here.
For positive even integer powers in the Gaussian, Laguerre and Jacobi
beta ensembles, Shen and Forrester
\cite[Propositions~2.1 and~2.2]{ShenForrester2025}
give quantitative remainder estimates and, for rational beta,
identify the associated constants with circular ensemble expressions.
Our argument treats
arbitrary $\beta$ without relying on the standard determinantal
structure at $\beta=2$ or the Pfaffian structures at $\beta=1,4$.
At $\beta=2$, Bourgade and Falconet \cite{BourgadeFalconet2022}
also obtained Fisher--Hartwig asymptotics by a probabilistic argument,
including joint moments of characteristic polynomials evaluated at
several times along the same unitary Brownian motion. 
The general beta joint moment problem remains explicitly listed among the open
questions in Lambert and Najnudel \cite[Section~1.3]{LN2026}.

For Hermitian ensembles, Charlier \cite{Charlier2019} gives
Hankel determinant asymptotics with several fixed root and jump
singularities for a one cut regular potential.
For the stationary Hermitian Ornstein--Uhlenbeck process, Keles
\cite[Theorem~1.4 and Corollary~1.6]{Keles2025}
obtains quantitative Fisher--Hartwig asymptotics with nonnegative
root powers and real exponential weights of the arguments.
His Corollary~1.6 permits singularities in the GUE bulk to
approach one another on mesoscopic scales.
Merging singularities lead to a different transition: Claeys and Fahs
\cite{ClaeysFahs2016} describe the Painlev\'e~V transition in the
model with weight $$|\det(M^2-tI)|^\alpha \,e^{-n\operatorname{Tr}V(M)},$$
where $\alpha>-1/2$ and $t$ is of order $n^{-2}$.

The fixed separation formula \eqref{eq:intro-FH} reaches the finite
ensemble integrability threshold $q_j>-1$. It proves the real root case without background of Forrester--Frankel's Conjecture~1 and answers the corresponding
joint moment question recalled in \cite{LN2026}. The theorem holds for
every real $\beta>0$, whereas \cite[Conjecture~1]{FF2004} states
the formula at rational beta. Their parameter is $q_{\rm FF}=q/\beta$. Their exact one point
product \cite[equation~(3.17)]{FF2004}, with their jump
parameter $b=0$, is $M_n(\beta q_{\rm FF})$ in our notation. Dividing it by
$n^{\beta q_{\rm FF}^2/2}$ identifies their constant
$A_{q_{\rm FF}}$ with $C_\beta(\beta q_{\rm FF})$.
In the unitary case,
$C_2(q)=G(1+q/2)^2/G(1+q)$, as in \cref{cor:widom}.
The one point constant is classical; the result concerns its
multipoint interaction, the full integrable range of exponents and the
quantitative uniformity. For nonnegative exponents, Assiotis and
Najnudel \cite[Proposition~3.5]{AN2026} already prove the global
upper bound valid even when points coincide. Our result identifies the
asymptotic constant and interaction and treats integrable negative
charges. The bound uniform over compact sets of nonnegative exponents is
proved in \cref{lem:moments-global-bound}.

The moments of the spatial integrals $\mathcal Z_N(q;f)$ arise in the
freezing picture of
Fyodorov, Hiary and Keating \cite{FHK2012} and Fyodorov and Keating
\cite{FyodorovKeating2014}. Bailey and Keating \cite{BaileyKeating2019}
determined the polynomial structure of the unweighted moments
for integer parameters in the CUE. Assiotis and Keating
\cite[Theorem~1.1]{AssiotisKeating2021} gave a combinatorial proof
of their asymptotics and expressed the leading coefficient as a volume. At the critical relation
$pq^2=2\beta$, the $p$th moment acquires
an additional logarithmic factor. Assiotis \cite{Assiotis2022}
proved asymptotics in parts of the supercritical regime and left
the critical line for further work. The general supercritical
results of Assiotis and Najnudel \cite{AN2026} concern a different
regime and do not supply the critical coefficient.

The Szeg\H{o} recursion supplies the orthogonal polynomials from
independent coefficients. It has already been used to prove Gaussian
fluctuations for phase counts \cite{Killip2008}. For distinct fixed
points on the circle, Ryckman \cite[Theorem~1.3]{Ryckman2008} proves
joint Gaussian convergence of the logarithms of characteristic
polynomial values after division by $\sqrt{\log N}$.

Hughes, Keating and O'Connell
\cite{HKO2001} obtained the Gaussian limit of $\log P_N$ in the CUE
as a random generalized function on the unit circle.
Convergence of the corresponding CUE measures to multiplicative chaos
was established in \cite{Webb2015,NSW2020}. For regularized counting
statistics, the analogous convergence was proved in \cite{LOS2018}.

In the Hermitian setting, Berestycki, Webb and Wong \cite{BWW2018}
prove convergence to Gaussian multiplicative chaos in the $L^2$ phase
for one cut regular unitarily invariant ensembles, using asymptotics
of Hankel determinants with Fisher--Hartwig singularities.

The probabilistic starting point is the representation of circular
ensembles through independent Verblunsky coefficients
\cite{KN2004,BNR2010}; for the theory of orthogonal polynomials on
the unit circle we refer to Simon \cite{Simon2005}. The present work belongs
to a broader program connecting random matrices and analytic number
theory through random analytic functions. Characteristic polynomial
moments and their arithmetic analogues
\cite{KeatingSnaith2000,CFKRS2005}, the probabilistic
product representation \cite{BHNY2008}, circular Jacobi
ensembles \cite{BNR2010}, virtual isometries \cite{BNN2013}, and the
stochastic zeta function \cite{CNN2017} are earlier parts of this
program. Here singular changes of probability measure connect the
macroscopic field to its microscopic zeros. The proofs use the
conditional chaos identities of \cite{CN2026,LN2026}, the circular mass
law of \cite{Remy2020}, the circular Jacobi microscopic limit of
\cite{LiValko2022,AN2026}, and the polynomial sublevel estimates of
\cite{CarberyWright2001}; their precise forms are given where they
enter the argument.

The large distance expansions of $\mathsf{Sine}_\beta$ correlation
functions obtained from these moment estimates, including the additional
Fourier corrections and the connected three and four point terms, are
developed in the unpublished companion manuscript \cite{XIIIB}.
The proofs here are independent of those correlation expansions:
the critical moment arguments use the Fisher--Hartwig estimates of
\cref{thm:main-root} and the microscopic matching of
\cref{thm:positive-merging}.

\paragraph{Organization.}
The proof of \cref{thm:main-root} occupies
\cref{sec:onepoint,sec:coupling,sec:cmv-unified,sec:gmc-identification,sec:terminal,sec:full-negative-transfer}.
\Cref{sec:positive-merging} proves the global moment bound at
arbitrary points and the microscopic limit, which supply the inputs
for the critical moments.
\Cref{sec:fractional-critical,sec:critical-moments} prove the critical
moment statements for noninteger and integer orders, respectively.
The two remaining sections are extensions of \cref{thm:main-root}
which the critical moment proofs do not use:
\cref{sec:mixed-root-jump} treats modulus and argument weights,
including fixed arc counting, and \cref{sec:smooth} adds a real
analytic background. \Cref{sec:extensions} describes the questions
left open. An index of notation is given at the end of the paper.

\section*{AI disclosure}
\addcontentsline{toc}{section}{AI disclosure}

This paper is the result of almost four years of joint work.
During its revision, we used ChatGPT 6 Astra (OpenAI) and, to a lesser
extent, Claude Fable 5.1 (Anthropic) for mathematical and editorial
assistance. This included reviewing arguments and references,
identifying gaps and missing hypotheses, and proposing and helping
develop corrections and extensions. In particular, the mathematical
assistance concerned quantitative signed Fisher--Hartwig estimates and complex parameter
arguments used in the mixed Fisher--Hartwig theorem. The tools also
assisted with exposition and manuscript preparation. Responsibility
for the final mathematical content, references and wording remains
with the authors.

\clearpage
\setcounter{tocdepth}{2}
\tableofcontents
\clearpage

\section{One point constants}
\label{sec:onepoint}

The exact one point integral fixes the normalization of every asymptotic
statement below. We state the result for modulus and argument together,
then recover the modulus moments by specialization. Set $b=\beta/2$
and $x_k=1+bk$.

The polynomial $P_N$ has no zero in $\D$. We use its analytic logarithm
which vanishes at zero; equivalently,
\[
 \log P_N(z)=\sum_{r=1}^N\log(1-ze^{-i\theta_r}),\qquad z\in\D,
\]
where each logarithm on the right also vanishes at $z=0$.
For $\theta\notin\{\theta_1,\ldots,\theta_N\}$ modulo $2\pi$, set
\[
 L_N(\theta)=\lim_{r\uparrow1}\log P_N(re^{i\theta}),\qquad
 X_N(\theta)=\Re L_N(\theta),\quad J_N(\theta)=\Im L_N(\theta).
\]
This finite boundary value is the sum of the boundary values of the
individual logarithms. In particular, $J_N$ is not reduced modulo
$2\pi$. At an eigenangle the real part tends to $-\infty$, so we exclude
these angles from the definition of $L_N$. For every fixed deterministic
$\theta$, the excluded event has probability zero by \eqref{eq:ensemble}.

Morris's integral \cite[equations~(1.17) and~(1.18)]{ForresterWarnaar2008}, or the
independent product representation in \cite[Proposition~4.3]{BNR2010},
gives, initially for $(a,c)$ near $(0,0)$,
\begin{equation}\label{eq:mixedmorris}
 \E e^{aL_N(0)+c\overline{L_N(0)}}
 =\prod_{k=0}^{N-1}
       \frac{\Gamma(x_k)\Gamma(x_k+a+c)}
            {\Gamma(x_k+a)\Gamma(x_k+c)}.
\end{equation}
The proof of \cref{prop:onepoint} below extends this identity to
$\Re(a+c)>-1$. In particular, setting $a=c=q/2$ gives, for
$\Re q>-1$,
\begin{equation}\label{eq:morris}
 M_N(q)=\E |P_N(1)|^q
 =\prod_{k=0}^{N-1}\frac{\Gamma(x_k)\Gamma(x_k+q)}{\Gamma(x_k+q/2)^2}.
\end{equation}
For a complex exponent, $|P_N(1)|^q$ means
$\exp(q\log|P_N(1)|)$. We use the analytic logarithm of the gamma
function on $\Re z>0$ that is real for positive $z$.

\begin{proposition}\label{prop:onepoint}
On the domain
\[
 \mathcal U_{\rm mix}=\{(a,c)\in\mathbb C^2:
         \Re(a+c)>-1,\ \Re a>-1,\ \Re c>-1\},
\]
the series
\begin{align}\label{eq:mixedconstantseries}
 \log\mathcal C_\beta(a,c)
 ={}&\frac{2ac}{\beta}\gamma_{\rm E}\notag\\
 &+\sum_{k=0}^\infty\left\{
 \log\Gamma(x_k)+\log\Gamma(x_k+a+c)
 -\log\Gamma(x_k+a)-\log\Gamma(x_k+c)
 -\frac{2ac}{\beta(k+1)}\right\}
\end{align}
converges locally uniformly and defines a holomorphic function.
Its exponential $\mathcal C_\beta$ is nonzero, and locally uniformly
on this domain,
\begin{equation}\label{eq:mixedoneasymp}
 \E e^{aL_N(0)+c\overline{L_N(0)}}
 =\mathcal C_\beta(a,c)N^{2ac/\beta}(1+O(N^{-1})).
\end{equation}
In particular, $C_\beta(q)=\mathcal C_\beta(q/2,q/2)$ is nonzero and
holomorphic for $\Re q>-1$, with
\begin{align}\label{eq:constantseries}
 \log C_\beta(q)
 &=\frac{q^2\gamma_{\rm E}}{2\beta}
 +\sum_{k=0}^{\infty}\left\{
 \log\Gamma(x_k)+\log\Gamma(x_k+q)
 -2\log\Gamma(x_k+q/2)-\frac{q^2}{2\beta(k+1)}
 \right\},
\end{align}
and, locally uniformly on that domain,
\begin{equation}\label{eq:onepointasymp}
 M_N(q)=C_\beta(q)N^{q^2/(2\beta)}\bigl(1+O(N^{-1})\bigr).
\end{equation}
For real $q>-1$, one has $C_\beta(q)>0$. At $\beta=2$,
\begin{equation}\label{eq:beta2constant}
 C_2(q)=\frac{G(1+q/2)^2}{G(1+q)},
\end{equation}
where $G$ is the Barnes function, with $G(1)=1$ and $G(z+1)=\Gamma(z)G(z)$.
\end{proposition}

\begin{proof}[Proof of \cref{prop:onepoint}]
We first extend the finite product identity to its full domain, then
subtract the harmonic term from its logarithm to obtain a convergent
series and the error estimate.
The finite identity \eqref{eq:mixedmorris} follows from the mixed
Morris formula near $(a,c)=(0,0)$. Its left side is holomorphic on
$\Re(a+c)>-1$: this condition controls the root singularity, and the
argument $J_N(0)$ is bounded for fixed $N$. The right side is
holomorphic there as well, since $1/\Gamma$ is entire and
$\Re(x_k+a+c)>0$. The identity theorem for holomorphic functions of
$(a,c)$ on this connected domain proves the stated finite identity.

Fix a compact set $K\subset\mathcal U_{\rm mix}$. All gamma arguments
have positive real part. In this proof, error constants may depend on
the fixed $\beta$ and on $K$. With $\psi=(\log\Gamma)'$ and
$\psi_1=\psi'$, Taylor expansion as $x\to+\infty$ gives, uniformly on $K$,
\begin{align*}
 g(x;a,c)
 &:=\log\Gamma(x)+\log\Gamma(x+a+c)
             -\log\Gamma(x+a)-\log\Gamma(x+c)\\
 &=ac\,\psi_1(x)+O_K(x^{-2})
 =\frac{ac}{x}+O_K(x^{-2}).
\end{align*}
Since $x_k^{-1}=2/(\beta(k+1))+O((k+1)^{-2})$, the summand
in \eqref{eq:mixedconstantseries} is $O_K((k+1)^{-2})$.
This proves local uniform convergence and holomorphy.

Put $\lambda=2ac/\beta$ and $H_N=\sum_{j=1}^N j^{-1}$.
Taking logarithms of the finite product gives
\begin{align*}
 \log\prod_{k=0}^{N-1}
       \frac{\Gamma(x_k)\Gamma(x_k+a+c)}
            {\Gamma(x_k+a)\Gamma(x_k+c)}-\lambda\log N
 ={}&\sum_{k=0}^{N-1}
       \left(g(x_k;a,c)-\frac{\lambda}{k+1}\right)\\
 &+\lambda(H_N-\log N).
\end{align*}
The omitted series tail is $O_K(N^{-1})$, and
$H_N-\log N=\gamma_{\rm E}+O(N^{-1})$. Exponentiation proves
\eqref{eq:mixedoneasymp}. The root statements follow on setting
$a=c=q/2$, and positivity follows from the real logarithmic series.
The restrictions defining $\mathcal U_{\rm mix}$ are used for this
logarithmic series and its nonzero exponential; the finite identity
requires only $\Re(a+c)>-1$.

When $\beta=2$, the functional equation of $G$ rewrites
\eqref{eq:morris} as
\[
 M_N(q)=\frac{G(N+1)G(N+1+q)}{G(N+1+q/2)^2}
 \frac{G(1+q/2)^2}{G(1+q)}.
\]
The first quotient is $N^{q^2/4}(1+O_K(N^{-1}))$ locally uniformly
for $\Re q>-1$, by the classical asymptotic expansion of $G$.
This proves \eqref{eq:beta2constant}.
\end{proof}

\begin{remark}
The one point product at $q=2$ telescopes for every $\beta>0$. With
$b=\beta/2$,
\begin{equation}\label{eq:check-q2}
 M_N(2)=\prod_{k=0}^{N-1}\frac{bk+2}{bk+1}
 =\frac{\Gamma(1/b)}{\Gamma(2/b)}
       \frac{\Gamma(N+2/b)}{\Gamma(N+1/b)}.
\end{equation}
It follows that $C_\beta(2)=\Gamma(1/b)/\Gamma(2/b)$, in agreement
with the series \eqref{eq:constantseries}, with its Barnes form
\eqref{eq:beta2constant} at $\beta=2$, and with the double gamma form
\eqref{eq:doublegamma-constant}.
\end{remark}

\subsection{The Barnes double gamma form of the constant}
\label{sec:doublegamma}

We use the zeta regularized convention of \cite{Ruijsenaars2000}:
for $b>0$, $\Re z>0$ and initially $\Re s>2$,
\[
 \zeta_2(s,z\mid1,b)=\sum_{r,k\ge0}(z+r+bk)^{-s},\qquad
 \log\Gamma_2(z\mid1,b)
   =\left.\partial_s\zeta_2(s,z\mid1,b)\right|_{s=0}.
\]
The second expression uses meromorphic continuation and specifies the
logarithm, so no multiplicative normalization is left implicit.

\begin{proposition}\label{prop:doublegamma}
With $b=\beta/2$, for $\Re q>-1$,
\begin{equation}\label{eq:doublegamma-constant}
 C_\beta(q)=b^{q^2/(4b)}
 \frac{\Gamma_2(1\mid1,b)\Gamma_2(1+q\mid1,b)}
      {\Gamma_2(1+q/2\mid1,b)^2}.
\end{equation}
The power of $b$ uses its real logarithm.
\end{proposition}

\begin{proof}[Proof of \cref{prop:doublegamma}]
We use the shift relation to express the finite Morris product through
the double gamma function, then estimate the remaining quotient at large
arguments. Throughout this proof, $b$ is fixed and its dependence is
suppressed: write $g(z)=\log\Gamma_2(z\mid1,b)$.
We use the Hurwitz zeta function
\[
 \zeta_H(s,z)=\sum_{r\ge0}(z+r)^{-s},\qquad \Re s>1,
\]
continued meromorphically in $s$, and write
$\zeta_H'(0,z)=\left.\partial_s\zeta_H(s,z)\right|_{s=0}$.
Removing the $k=0$ terms from the double zeta series gives
\[
 \zeta_2(s,z+b\mid1,b)-\zeta_2(s,z\mid1,b)=-\zeta_H(s,z).
\]
Differentiate at $s=0$ and use
$\zeta_H'(0,z)=\log\Gamma(z)-\tfrac12\log(2\pi)$ to obtain
\begin{equation}\label{eq:doublegamma-shift}
 g(z+b)-g(z)=-\log\Gamma(z)+\tfrac12\log(2\pi).
\end{equation}
To check the asymptotic with this convention, write
\[
 A=\frac{1+b}{2b},\qquad
 R(t)=\frac1{(1-e^{-t})(1-e^{-bt})}
                   -\frac1{bt^2}-\frac{A}{t}.
\]
Taylor expansion at zero and the limit at infinity show that $R$
is bounded on $(0,\infty)$. For $\Re s>2$, termwise integration gives
the Mellin representation
\[
 \zeta_2(s,z\mid1,b)=\frac1{\Gamma(s)}
 \int_0^\infty
   \frac{t^{s-1}e^{-zt}}{(1-e^{-t})(1-e^{-bt})}\,dt.
\]
After two derivatives in $z$, subtracting the two displayed singular
terms gives the continuation
\[
 \partial_z^2\zeta_2(s,z\mid1,b)
 =\frac{z^{-s}}b+A s z^{-s-1}
   +\frac1{\Gamma(s)}\int_0^\infty t^{s+1}e^{-zt}R(t)\,dt.
\]
The integral is holomorphic for $\Re s>-2$, locally uniformly in
$\Re z>0$. Differentiating at $s=0$, where
$1/\Gamma(s)=s+O(s^2)$, yields
\begin{equation}\label{eq:doublegamma-second-derivative}
 g''(z)=-\frac1b\log z+\frac{A}{z}
                  +\int_0^\infty t e^{-zt}R(t)\,dt.
\end{equation}
The integral is bounded by $\|R\|_\infty/(\Re z)^2$. Thus
\begin{equation}\label{eq:doublegamma-asymptotic}
 g''(x+w)=-\frac1b\log x+O(x^{-1}),
 \qquad x\longrightarrow+\infty,
\end{equation}
uniformly for $w$ in compact subsets of $\mathbb C$.

Let $B(q)=\exp\{g(1)+g(1+q)-2g(1+q/2)\}$,
and define $B_N(q)$ by replacing $1$ by $1+Nb$.
The shift relation \eqref{eq:doublegamma-shift}, applied to the
Morris product \eqref{eq:morris}, gives
\begin{equation}\label{eq:doublegamma-finite-product}
 M_N(q)=B(q)/B_N(q);
\end{equation}
the constant $\tfrac12\log(2\pi)$ cancels in the second difference.
Integration twice along complex line segments gives
\[
 \log B_N(q)=\frac{q^2}{4}
   \int_0^1\!\int_0^1
    g''\left(1+Nb+\frac{q(u+v)}2\right)\,du\,dv
 =-\frac{q^2}{4b}\log(Nb)+O(N^{-1}).
\]
This is locally uniform for $\Re q>-1$. Comparing
\eqref{eq:doublegamma-finite-product} with \cref{prop:onepoint}
proves \eqref{eq:doublegamma-constant}.
\end{proof}

For comparison with the ordinary Barnes convention,
\begin{equation}\label{eq:doublegamma-barnes-G}
 \Gamma_2(z\mid1,1)
 =e^{\zeta'(-1)}(2\pi)^{(z-1)/2}G(z)^{-1},\qquad \Re z>0.
\end{equation}
Indeed, the logarithmic ratio is periodic by the shift relations.
Its second derivative tends to zero under integer translation, by
\eqref{eq:doublegamma-second-derivative} and the Barnes asymptotic
\cite[equation~(36) and Theorem~2]{AlexanianKuznetsov2023};
periodicity makes that derivative zero and then makes the ratio
constant. The identity $\zeta_2(s,1\mid1,1)=\zeta(s-1)$ fixes the
constant. The constant and linear factors cancel in
\eqref{eq:doublegamma-constant}, recovering \eqref{eq:beta2constant}.

The precise deviation consequences of this normalization for
$\log|P_N(1)|$ at general beta are studied in
\cite{MeliotNikeghbali2022} and in the unpublished manuscript
\cite{KLMN2026}.

\subsection{One point modulus and argument}
The mixed transform in \cref{prop:onepoint} has two useful
specializations. Write
$L=L_N(0)$, $X=X_N(0)$ and $J=J_N(0)$. For real $s,t$,
\begin{align*}
 a=\frac{s-it}{2},\quad c=\frac{s+it}{2}
 &\quad\Longrightarrow\quad
 aL+c\bar L=sX+tJ,\qquad \frac{2ac}{\beta}=\frac{s^2+t^2}{2\beta},\\
 a=\frac{s+t}{2},\quad c=\frac{s-t}{2}
 &\quad\Longrightarrow\quad
 aL+c\bar L=sX+itJ,\qquad \frac{2ac}{\beta}=\frac{s^2-t^2}{2\beta}.
\end{align*}
The first choice gives the joint Laplace transform of log modulus and
argument, used for the real argument weights in
\cref{sec:mixed-root-jump}. The second gives a Laplace transform
in $X$ and a Fourier transform in $J$; at $s=0$ it gives the
characteristic function of the argument. The asymptotic
\eqref{eq:mixedoneasymp} applies when the chosen pair belongs to
$\mathcal U_{\rm mix}$, that is, when
$\Re(a+c)>-1$, $\Re a>-1$ and $\Re c>-1$.
For the finite identity \eqref{eq:mixedmorris}, only $s>-1$ is needed.

\begin{proposition}
\label{prop:onepoint-first-correction}
Locally uniformly for $\Re q>-1$,
\begin{equation}\label{eq:onepoint-first-correction}
 M_N(q)=C_\beta(q)N^{q^2/(2\beta)}
 \left[1+\frac{q^2}{2\beta N}
            \left(\frac{q+1}{\beta}-\frac12\right)+O(N^{-2})\right].
\end{equation}
\end{proposition}

\begin{proof}[Proof of \cref{prop:onepoint-first-correction}]
If $r_k=\Gamma(1+bk)\Gamma(1+bk+q)/\Gamma(1+bk+q/2)^2$,
the gamma expansion to the next order gives, uniformly on compact
parameter sets,
\[
 \log r_k=\frac{q^2}{4bk}
             -\frac{q^2(q+1)}{8b^2k^2}+O(k^{-3}),\qquad k\ge1.
\]
Write $a=q^2/(4b)$ and $d=-q^2(q+1)/(8b^2)$.
Regroup the defining series for $C_\beta$, absorbing the $k=0$
factor and replacing $1/(k+1)$ by $1/k$ for $k\ge1$.
Comparison with the product up to $N-1$ yields
\[
 \log M_N-\log C_\beta-a\log N
 =a\{H_{N-1}-\log N-\gamma_{\rm E}\}
       -\sum_{k=N}^{\infty}(d/k^2+O(k^{-3}))
 =\frac{-a/2-d}{N}+O(N^{-2}).
\]
Here $H_{N-1}=\log N+\gamma_{\rm E}-1/(2N)+O(N^{-2})$ and
$\sum_{k=N}^\infty k^{-2}=N^{-1}+O(N^{-2})$.
Exponentiating and substituting $b=\beta/2$ proves the formula.
\end{proof}

\section{Coupling with orthogonal polynomials}
\label{sec:coupling}

We introduce a coupling of the circular ensembles through orthogonal
polynomials on the unit circle. The same initial sequence of independent
coefficients is used for every matrix size, so that conditional
expectations can compare successive polynomials. A final independent
coefficient on the circle then gives the characteristic polynomial.
We also establish the finite moment bounds needed for negative root
exponents.

Recall that $b=\beta/2$ and
$\sigma(d\theta)=d\theta/(2\pi)$ on $\T=\R/(2\pi\mathbb Z)$.
Let $(\alpha_k)_{k\geq0}$ be independent random variables in the open unit
disk $\D$, with densities
\begin{equation}\label{eq:coefficient-law}
 \frac{b(k+1)}{\pi}(1-|z|^2)^{b(k+1)-1}\mathbf1_{\{|z|<1\}}\,d^2z.
\end{equation}
Thus $|\alpha_k|^2$ has distribution
$\operatorname{Beta}(1,b(k+1))$, independently of its uniform argument.
We use the filtration $\mathcal F_n=\sigma(\alpha_0,\ldots,\alpha_{n-1})$
generated by the first $n$ coefficients, and the sigma-field
$\mathcal F_\infty=\sigma(\alpha_k:k\geq0)$ generated by all of them.
Starting from $\Phi_0=\Phi_0^*=1$, define the monic polynomials by
\begin{equation}\label{eq:szego-recursion}
 \begin{pmatrix}\Phi_{k+1}(z)\\ \Phi_{k+1}^*(z)\end{pmatrix}
 =A(\alpha_k,z)
   \begin{pmatrix}\Phi_k(z)\\ \Phi_k^*(z)\end{pmatrix},
 \qquad
 A(\alpha,z)=
 \begin{pmatrix}z&-\bar\alpha\\-\alpha z&1\end{pmatrix}.
\end{equation}
Here $\Phi_k^*(z)=z^k\overline{\Phi_k(1/\bar z)}$ is the reversed
polynomial. Define also
\begin{equation}\label{eq:blaschke-definition}
 B_k(z)=\frac{z\Phi_k(z)}{\Phi_k^*(z)}.
\end{equation}
Induction in \eqref{eq:szego-recursion} gives
\[
 \Phi_{k+1}^*=\Phi_k^*(1-\alpha_k B_k),
 \qquad
 B_{k+1}(z)=z\frac{B_k(z)-\bar\alpha_k}{1-\alpha_k B_k(z)}.
\]
Since $B_0(z)=z$, these identities show that $B_k$ is a finite Blaschke
product: it maps $\D$ into itself and has modulus one on $\partial\D$.
They also show inductively that $\Phi_k^*$ has no zero in
$\overline\D$. Consequently its analytic logarithm is uniquely determined
by $\log\Phi_k^*(0)=0$.

The following coupling is the reason for choosing the coefficient laws
\eqref{eq:coefficient-law}. Its first $n$ coefficients can be used at
all larger matrix sizes.

\begin{proposition}
\label{prop:coupled-ensemble}
Let $\eta$ be uniform on $\T$, independently of
$\mathcal F_\infty$, and put $t=e^{i\eta}$. For every $n\geq0$, the
polynomial
\begin{equation}\label{eq:terminal-polynomial}
 \mathcal P_{n+1}(z)
 =\Phi_n^*(z)-tz\Phi_n(z)
 =\Phi_n^*(z)(1-tB_n(z))
\end{equation}
has the law of the characteristic polynomial $P_{n+1}$ for $n+1$ particles
defined in \eqref{eq:PN}, as a random polynomial. We use this coupling
and henceforth write $P_{n+1}=\mathcal P_{n+1}$.
\end{proposition}

\subsection{Reversing the order of the coefficients}

\begin{lemma}\label{lem:reversal}
Let $N\geq2$, $\alpha_0,\ldots,\alpha_{N-2}\in\D$, and let the last
coefficient be $\alpha_{N-1}=t\in\partial\D$. Put
\[
 \alpha'_j=-\overline{\alpha_{N-2-j}}\,t\quad(0\leq j\leq N-2),
 \qquad \alpha'_{N-1}=t.
\]
The monic polynomial $\Phi_N$ obtained by the recursion, including the
last coefficient $t$, is unchanged by this transformation.
\end{lemma}

\begin{proof}[Proof of \cref{lem:reversal}]
Put $r=(z,-\bar t)$, $c=(1,1)^T$, and
$P=A(\alpha_{N-2},z)\cdots A(\alpha_0,z)$. The last step of the
recursion gives $\Phi_N(z)=rPc$. For $z\ne0$, let
$D=\operatorname{diag}(-\bar t,z)$. Direct multiplication gives
\[
 A(-\bar\alpha t,z)=D A(\alpha,z)^T D^{-1}.
\]
The reversed order of the transformed coefficients therefore gives
$P'=DP^TD^{-1}$. Moreover,
\[
 rD=-z\bar t\,c^T,
 \qquad D^{-1}c=-\frac{t}{z}r^T.
\]
Since $|t|=1$, it follows that
\[
 \Phi_N(z;\boldsymbol\alpha')
 =rDP^TD^{-1}c=c^TP^Tr^T=rPc
 =\Phi_N(z;\boldsymbol\alpha).
\]
The equality at $z=0$ follows from polynomial continuity.
\end{proof}

\begin{proof}[Proof of \cref{prop:coupled-ensemble}]
The Killip--Nenciu matrix model \cite[Theorem~1]{KN2004} identifies
the zeros of the final monic polynomial with the circular beta
ensemble when its interior coefficients are independent with squared
moduli $\operatorname{Beta}(1,b(N-1-j))$, $0\leq j\leq N-2$, and
the last coefficient is independent and uniform on the circle.
Condition on that last coefficient $t$. A rotation invariant disk
law is invariant under conjugation and multiplication by $-t$.
The transformed coefficients of \cref{lem:reversal} are therefore
conditionally independent and have the original laws in reverse order.
These conditional laws do not depend on $t$. Their squared moduli
consequently have parameters $b(j+1)$, exactly as in
\eqref{eq:coefficient-law}, and they remain independent of $t$.

The transformed reversal preserves the final monic polynomial, not
only its zeros. With $N=n+1$, its constant term is $-\bar t$, so
\[
 \frac{\Phi_{n+1}(z)}{\Phi_{n+1}(0)}
 =-t\bigl(z\Phi_n(z)-\bar t\Phi_n^*(z)\bigr)
 =\Phi_n^*(z)-tz\Phi_n(z).
\]
If its zeros are $e^{i\theta_j}$, this normalization is precisely
$\prod_j(1-ze^{-i\theta_j})$. This proves the assertion for $N\geq2$;
for $N=1$, the polynomial is $1-tz$, which has one uniform zero.
\end{proof}

\subsection{Phases and logarithms}

For a real lift of $\theta\in\T$, define the lifted phase
\begin{equation}\label{eq:phase-definition}
 \psi_n(\theta)=(n+1)\theta
                 -2\Im\log\Phi_n^*(e^{i\theta}).
\end{equation}
The identity $\Phi_n(e^{i\theta})
=e^{in\theta}\overline{\Phi_n^*(e^{i\theta})}$ gives
$B_n(e^{i\theta})=e^{i\psi_n(\theta)}$. Thus, writing each logarithm
with its analytic branch normalized at zero,
\begin{align}
 \log\Phi_{n+1}^*(e^{i\theta})
 &=\log\Phi_n^*(e^{i\theta})
      +\log(1-\alpha_n e^{i\psi_n(\theta)}),
       \label{eq:log-recursion}\\
 \psi_{n+1}(\theta)
 &=\psi_n(\theta)+\theta
      -2\Im\log(1-\alpha_n e^{i\psi_n(\theta)}),
 \qquad\psi_0(\theta)=\theta.
       \label{eq:phase-recursion}
\end{align}
In particular,
\begin{equation}\label{eq:phase-log}
 \Im\log\Phi_n^*(e^{i\theta})
      =\frac{(n+1)\theta-\psi_n(\theta)}2.
\end{equation}
The polynomial identity holds jointly in all angles:
\begin{equation}\label{eq:terminal-process}
 P_{n+1}(e^{i\theta})
 =\Phi_n^*(e^{i\theta})
    \bigl(1-e^{i(\eta+\psi_n(\theta))}\bigr).
\end{equation}
There is also the exact analytic identity in $\D$
\begin{equation}\label{eq:terminal-log}
 \log P_{n+1}(z)
 =\log\Phi_n^*(z)+\log(1-tB_n(z)).
\end{equation}
Indeed $|B_n(z)|<1$ in $\D$, and both sides exponentiate to the same
nonvanishing function and equal zero at $z=0$. There is therefore no
undetermined integer multiple of $2\pi i$. Radial limits give the same
identity at boundary points that are not zeros. In particular, the
argument of the last factor is the periodic function
\[
 g(u)=\arg(1-e^{iu})=\frac{u-\pi}{2},\qquad 0<u<2\pi.
\]

\subsection{Finite moments and the one point normalization}

The product of polynomial moduli is the random weight whose expectation
we shall estimate. For $\mathbf q\in\C^m$ and angles
$\boldsymbol\vartheta$, put
\begin{equation}\label{eq:root-weight}
 Y_n(\mathbf q;\boldsymbol\vartheta)
 =\exp\!\left(\sum_{j=1}^m
        q_j\log|\Phi_n^*(e^{i\vartheta_j})|\right),
 \qquad a(\mathbf q)=\frac1{2\beta}\sum_{j=1}^m q_j^2.
\end{equation}
We omit $\boldsymbol\vartheta$ from $Y_n$ when the configuration is
fixed. For real $\mathbf q$ this weight is positive, even when some
exponents are negative. For complex exponents its modulus satisfies
$|Y_n(\mathbf q)|=Y_n(\Re\mathbf q)$.

Averaging the final coefficient in \eqref{eq:terminal-process}
introduces the following function of the phases, for real $\mathbf u$:
\begin{equation}\label{eq:terminal-W}
 W_{\mathbf u}(\boldsymbol\varphi)
 =\int_{\T}\prod_{j=1}^m
       |1-e^{i(\eta+\varphi_j)}|^{u_j}\,\sigma(d\eta),
\end{equation}
allowing the value $+\infty$. Its appearance keeps track of the fact
that all points share the same last coefficient.

Write $\boldsymbol\psi_n=(\psi_n(\vartheta_j))_{j=1}^m$.
Conditioning on $\mathcal F_n$ in \eqref{eq:terminal-process} gives,
for real $q_j>-1$ and distinct angles,
\begin{equation}\label{eq:terminal-expectation}
 F_{n+1}(\mathbf q,\boldsymbol\vartheta)
 =\E[Y_n(\mathbf q)W_{\mathbf q}(\boldsymbol\psi_n)]
 \ge\E Y_n(\mathbf q).
\end{equation}
The equality follows from Tonelli's theorem. For the inequality, apply
Jensen's inequality in $\eta$ and use
$\int_\T\log|1-e^{it}|\,\sigma(dt)=0$, which gives
$W_{\mathbf q}\ge1$. Thus the characteristic polynomial moment is
reduced to an OPUC moment and an average depending only on the phases.
We next verify the finite moments needed to use this identity.

\begin{lemma}
\label{lem:negative-initial}
Fix an integer $n\ge1$. For distinct points, the real moment
$\E Y_n(\mathbf q;\boldsymbol\vartheta)$ is finite if and only if
$q_j>-1-b$ for $1\le j\le m$. It is locally uniformly bounded on this
domain of exponents and on separated angle sets. The expectation is
holomorphic in $\mathbf q$ on the domain $\Re q_j>-1-b$ and jointly continuous in the exponents
and distinct angles.

Let $K$ be a compact subset of
$\{\mathbf q\in\mathbb C^m:\Re q_j>-1-b\ (1\le j\le m)\}$ and fix
$\delta>0$, with $\mathcal S_\delta$ as in \eqref{eq:intro-domain}.
There is $p>1$, depending on $K$, such that
\begin{equation}\label{eq:initial-uniform-p}
 \sup_{\substack{\mathbf q\in K\\
                  \boldsymbol\vartheta\in\mathcal S_\delta}}
 \E|Y_n(\mathbf q;\boldsymbol\vartheta)|^p<\infty.
\end{equation}
The bound may depend on $n$, $K$ and $\delta$.
For $n=0$ the moment equals one for every exponent, independently
of the preceding integrability criterion for $n\ge1$.
In the smaller range $q_j>-1$,
\begin{equation}\label{eq:initial-moment-bound}
 \E Y_n(\mathbf q;\boldsymbol\vartheta)
 \le F_{n+1}(\mathbf q,\boldsymbol\vartheta)<\infty.
\end{equation}
\end{lemma}

\begin{proof}[Proof of \cref{lem:negative-initial}]
We first reduce integrability to a single zero approaching a prescribed
boundary point. The same local integral will give necessity, sufficiency
and a uniform moment of order greater than one.
The zeros $z_1,\ldots,z_n$ of $\Phi_n$ have joint density, up to
normalization,
\begin{equation}\label{eq:truncated-zero-density}
 \prod_{r=1}^n(1-|z_r|^2)^{b-1}
 \prod_{r<s}|z_r-z_s|^2|1-z_r\bar z_s|^{2b-2}
 \prod_r d^2z_r.
\end{equation}
This is \cite[Theorem~6.1]{KK2017}; its independent coefficient laws
are exactly those in the present coupling, in the increasing order
\eqref{eq:coefficient-law}, rather than the reversed order in the
usual finite circular ensemble representation. In their notation,
$\Theta_\nu$ has disk density
$(\nu-1)(1-|z|^2)^{(\nu-3)/2}/(2\pi)$ for $\nu>1$;
here $\nu=\beta(k+1)+1$.
For $b\ge1$ each pair factor in \eqref{eq:truncated-zero-density}
is bounded on the bidisk. For $0<b<1$, the identity
\[
 |1-z\bar w|^2-|z-w|^2=(1-|z|^2)(1-|w|^2)\ge0
\]
gives
$|z-w|^2|1-z\bar w|^{2b-2}\le |1-z\bar w|^{2b}\le2^{2b}$.
Since $\Phi_n^*(e^{i\vartheta})=
\prod_r(1-\bar z_re^{i\vartheta})$, the desired moment is bounded
by a constant times the $n$th power of
\[
 \int_\D(1-|z|^2)^{b-1}
       \prod_j|1-\bar ze^{i\vartheta_j}|^{q_j}\,d^2z.
\]
Near $e^{i\vartheta_j}$ write $z=(1-t)e^{i(\vartheta_j+y)}$.
The singular part of this integral is comparable to
\[
 \int_0^\varepsilon\int_{-\varepsilon}^{\varepsilon}
 t^{b-1}(t^2+y^2)^{q_j/2}\,dy\,dt.
\]
Polar coordinates in the half-plane $t>0$ separate this into an
angular integral of $(\cos\phi)^{b-1}$, which is finite since $b>0$,
and a radial integral of $\rho^{b+q_j}\,d\rho$. The latter is finite
exactly when $q_j>-1-b$. Separation bounds the remaining factors above and below. This proves sufficiency, uniformly on the
stated compact sets. For necessity, restrict the other zeros to
pairwise disjoint compact disks in the interior, away from the
specified boundary point. Every pair factor is then bounded below by
a positive constant while the remaining zero approaches that point.
The same integral diverges when $q_j\le-1-b$.

Choose $p>1$ close enough to one that $p\Re q_j>-1-b$ on a given
compact set. The identity $|Y_n(\mathbf q)|^p=Y_n(p\Re\mathbf q)$
proves \eqref{eq:initial-uniform-p}. It also proves continuity: the weights converge pointwise as their
parameters and distinct angles converge, and the uniform $p$th moment
bound \eqref{eq:initial-uniform-p} gives uniform integrability. For holomorphy, enclose the real parts of a
small polydisk in a box of admissible real exponents. The sum of the
weights at its finitely many vertices dominates the supremum of the
complex weights. Dominated holomorphy applies.

Finally, \eqref{eq:terminal-expectation} gives the comparison in
\eqref{eq:initial-moment-bound}. The characteristic polynomial moment
is finite for $q_j>-1$: bound the Vandermonde and integrate each particle
near the distinct points.
\end{proof}

The OPUC one point moment omits the last coefficient of modulus one.
Its relation to $M_N$ and $C_\beta$ is given below.

\begin{lemma}\label{lem:opuc-onepoint}
For $\Re q>-1-b$, the angle independent one point moment is
\begin{equation}\label{eq:opuc-onepoint-exact}
 \nu_n(q)=\E|\Phi_n^*(1)|^q
 =\prod_{k=1}^n
 \frac{\Gamma(1+bk)\Gamma(1+bk+q)}{\Gamma(1+bk+q/2)^2}.
\end{equation}
There is a holomorphic nonzero function $D_\beta$ on this half-plane
such that, locally uniformly,
\begin{equation}\label{eq:opuc-onepoint-asymptotic}
 \nu_n(q)=D_\beta(q)\,n^{q^2/(2\beta)}\,(1+O(n^{-1})).
\end{equation}
For $\Re q>-1$ it satisfies
\begin{equation}\label{eq:terminal-c}
 \begin{gathered}
 c(q)=\int_\T|1-e^{i\eta}|^q\,\sigma(d\eta)
      =\frac{\Gamma(1+q)}{\Gamma(1+q/2)^2},\\
 D_\beta(q)=\frac{C_\beta(q)}{c(q)},\qquad
 \nu_n(q)=\frac{M_{n+1}(q)}{c(q)}.
 \end{gathered}
\end{equation}
\end{lemma}

\begin{proof}[Proof of \cref{lem:opuc-onepoint}]
Rotation invariance of each independent coefficient makes the
one point multiplicative increment independent of the incoming phase.
Its disk beta integral gives the $k$th factor of
\eqref{eq:opuc-onepoint-exact}. Alternatively, conditioning on the last coefficient
and the Morris formula give this product for $\Re q>-1$;
\cref{lem:negative-initial} and dominated holomorphy extend it to
$\Re q>-1-b$. All numerator and denominator arguments in this last
product have positive real parts. If $r_k(q)$ denotes its $k$th
factor and $a=q^2/(2\beta)$, the gamma expansion gives
$\log r_k(q)=a/k+O(k^{-2})$, locally uniformly. Hence
\[
 D_\beta(q)=\exp\left\{a\gamma_{\rm E}
       +\sum_{k=1}^{\infty}(\log r_k(q)-a/k)\right\}
\]
is holomorphic and nonzero. The tail of this series and
$\sum_{k=1}^nk^{-1}=\log n+\gamma_{\rm E}+O(n^{-1})$
prove the asymptotic. The quotient in \eqref{eq:terminal-c} is used
only where both $C_\beta(q)$ and $c(q)$ are defined; beyond that range,
the product starting at $k=1$ defines $D_\beta$.
\end{proof}

\begin{remark}
The one point threshold $-1-\beta/2$ also appears in
\cite[Lemma~21]{LN2026}. The zero density argument in \cref{lem:negative-initial} supplies the
joint integrability and local uniformity needed here. This larger
OPUC domain does not alter the characteristic polynomial threshold
$q_j>-1$.
\end{remark}

\paragraph{Two polynomials and their one point normalizations.}
The two polynomials in \cref{sec:onepoint,sec:coupling} have different
indices and one point normalizations. The characteristic polynomial $P_N$ of \eqref{eq:PN}
has $N$ particles; its one point moment is $M_N(q)=\E|P_N(1)|^q$ in
\eqref{eq:morris}, with the constant
$C_\beta(q)=\lim_NN^{-q^2/(2\beta)}M_N(q)$ of \eqref{eq:constantseries}.
The reversed orthogonal polynomial $\Phi_n^*$ of \eqref{eq:szego-recursion}
is built from the $n$ coefficients $\alpha_0,\ldots,\alpha_{n-1}$; its
one point moment is $\nu_n(q)=\E|\Phi_n^*(1)|^q$ in
\eqref{eq:opuc-onepoint-exact}, with the constant $D_\beta(q)$ of
\eqref{eq:opuc-onepoint-asymptotic}. The two polynomials are coupled by
\eqref{eq:terminal-polynomial}, $P_{n+1}=\Phi_n^*(1-tB_n)$ with an
independent uniform phase $t$, so that $N=n+1$ in every transfer from
$\Phi_n^*$ to $P_N$. Their moments and limiting constants differ
by the factor $c(q)$ according to \eqref{eq:terminal-c}. Statements about $\Phi_n^*$ are indexed by
$n$ and statements about $P_N$ by $N$.

\section{A Fourier correction and holomorphic convergence}
\label{sec:cmv-unified}

The characteristic polynomial moment is related to $Y_n$ by
\eqref{eq:terminal-expectation}. More precisely, for real $q_j>-1$ and distinct angles,
\[
 n^{-a(\mathbf q)}F_{n+1}(\mathbf q,\boldsymbol\vartheta)
 =\bigl(n^{-a(\mathbf q)}\E Y_n(\mathbf q;\boldsymbol\vartheta)\bigr)
   \frac{\E[Y_n(\mathbf q)W_{\mathbf q}(\boldsymbol\psi_n)]}
        {\E Y_n(\mathbf q)}.
\]
We will prove that the first factor converges to a holomorphic function
$H(\mathbf q;\boldsymbol\vartheta)$, and that the quotient converges to
$\prod_jc(q_j)$. We first establish the existence and holomorphy of
$H$ on the tube \eqref{eq:tube}, using estimates from the recursion
alone. The conditional chaos identities in
\cref{sec:gmc-identification} compute this limit on a smaller real
domain. Uniqueness of holomorphic continuation then identifies it
throughout the larger domain. This separates the bounds needed to
construct the limit from the argument that computes its value.
\Cref{sec:terminal,sec:full-negative-transfer} justify
the second convergence, including when $W_{\mathbf q}$ is unbounded.
Since $D_\beta(q)c(q)=C_\beta(q)$, these two limits give precisely the
constant in the Fisher--Hartwig formula.

This section supplies the first convergence and the phase estimates
needed for the second. The existence of $H$ requires
$\Re q_j>-1-b$ for integrability at a fixed initial index.
The conditional estimates at sufficiently large indices are valid on
arbitrary compact sets of exponents. They apply after conditioning
on a finite initial segment even when the corresponding unconditional
moment does not exist.

Throughout this section, $\beta$ and $m$ are fixed; constants indexed by
a compact set may also depend on these two parameters. We use the
coupling and the weight $Y_n(\mathbf q)$ from \cref{sec:coupling},
suppressing the angles when no confusion is possible. Recall the exponent
\[
 a(\mathbf q)=\frac1{2\beta}\sum_{j=1}^m q_j^2.
\]

Recall that $b=\beta/2$. The initial integrability condition and the
summability condition for the errors lead to the open tube
\begin{equation}\label{eq:tube}
 \mathcal U=\left\{\mathbf q\in\C^m:
       \Re q_j>-1-b\ (1\le j\le m),\quad a(\Im\mathbf q)<1\right\}.
\end{equation}

\begin{theorem}
\label{thm:holomorphic-limit}
There is a function $H(\mathbf q;\boldsymbol\vartheta)$, holomorphic in
$\mathbf q\in\mathcal U$ and continuous jointly in the exponents and
distinct angles, such that
\begin{equation}\label{eq:holomorphic-limit}
 n^{-a(\mathbf q)}\E Y_n(\mathbf q;\boldsymbol\vartheta)
       \longrightarrow H(\mathbf q;\boldsymbol\vartheta).
\end{equation}
For each compact $K\subset\mathcal U$ and each $\delta>0$, the convergence
is uniform on $K\times\mathcal S_\delta$, with error
\begin{equation}\label{eq:holomorphic-limit-rate}
 O_{K,\delta}(n^{-\rho_K}),\qquad
 \rho_K=1-\sup_{\mathbf q\in K}a(\Im\mathbf q)>0.
\end{equation}
In particular the error is $O_{K_0,\delta}(n^{-1})$ on a real compact
set $K_0\subset(-1-b,\infty)^m$. The limit is strictly positive for real
exponents in this domain. Its explicit value is identified in the next
section.
\end{theorem}

To study the limit and the phases together, we multiply the weight by
a Fourier character of the phase vector. For $\boldsymbol\ell\in\mathbb Z^m$,
we call $\boldsymbol\ell$ a Fourier mode and
$e^{i\boldsymbol\ell\cdot\boldsymbol\psi_n}$ its character.
These modes are distinct from Fourier coefficients of deterministic
functions of an angle. The corresponding observables are
\begin{align}
 T_n(\mathbf q,\boldsymbol\ell)
  &=Y_n(\mathbf q)\exp\!\left(i\sum_{j=1}^m
       \ell_j\psi_n(\vartheta_j)\right),\qquad
       \boldsymbol\ell\in\mathbb Z^m,\label{eq:Tdef}\\
 d(\boldsymbol\ell)&=\frac2\beta\sum_{j=1}^m\ell_j^2,
 \qquad \lambda=a(\mathbf q)-d(\boldsymbol\ell),
 \qquad \omega=\sum_{j=1}^m\ell_j\vartheta_j,
 \qquad L=1+\|\boldsymbol\ell\|_1.\label{eq:exponents}
\end{align}
Here $d(\boldsymbol\ell)$ is the damping exponent of the mode,
$\omega$ is its deterministic rotation, and $L$ measures its size.
We suppress the arguments of $T_n$ and $\lambda$ when they are fixed.
In particular,
\begin{equation}\label{eq:modT}
 |T_n(\mathbf q,\boldsymbol\ell)|=Y_n(\Re\mathbf q).
\end{equation}
Unless a subscript specifies otherwise, vector norms in estimates below are
Euclidean; changing this convention only changes constants depending on $m$.
Powers of positive real numbers
are defined using their real logarithms.
When a conditional statement is made for arbitrary exponents, its conditional
expectation means integration over the independent future coefficients,
with the coefficients in $\mathcal F_n$ fixed.
\Cref{lem:onestep,lem:conditional-weights} prove absolute integrability
of this transition kernel integral at the stated starting indices. Thus the
statement also makes sense when its initial weight has infinite
unconditional expectation. In the applications with $\Re q_j>-1-b$,
\cref{lem:negative-initial} makes these the usual integrable conditional
expectations.

To average $W_{\mathbf q}$ later, we also need its nonconstant Fourier
modes to vanish under the probability measure proportional to $Y_n$.
The following estimate controls their dependence on both the mode and
the separation.

For $d>0$ and $t\ge1$ define
\begin{equation}\label{eq:phase-rate-function}
 \mathfrak r_d(t)=
 \begin{cases}
 t^{-d},&0<d<1,\\
 (1+\log t)t^{-1},&d=1,\\
 t^{-1},&d>1.
 \end{cases}
\end{equation}
The logarithm at $d=1$ comes from summing a harmonic series in the
inhomogeneous recursion \eqref{eq:normalized-corrected-mode}.

\begin{proposition}
\label{prop:quantitative-phase-mixing}
Let $K_0\subset(-1-b,\infty)^m$ be compact. For arbitrary nonzero integer
modes put $d_*=2/\beta$. If the modes are restricted by
$\sum_j\ell_j=0$, put instead $d_*=4/\beta$. There is a constant
$C=C(K_0,\beta,m)$ such that, for $0<\delta\le1$,
$\boldsymbol\vartheta\in\mathcal S_\delta$, $\mathbf q\in K_0$,
$n\delta\ge1$, and every mode in the indicated class,
\begin{equation}\label{eq:relative-phase-mixing}
 \frac{|\E T_n(\mathbf q,\boldsymbol\ell)|}
      {\E Y_n(\mathbf q)}
 \le C(1+\|\boldsymbol\ell\|_1)^2\,
                  \mathfrak r_{d_*}(n\delta).
\end{equation}
In particular, for fixed $\delta>0$ and every
$0<\gamma<\min\{d_*,1\}$,
\begin{equation}\label{eq:quantitative-phase-mixing}
 \left|n^{-a(\mathbf q)}\E T_n(\mathbf q,\boldsymbol\ell)\right|
 \le C_{K_0,\delta,\gamma}
       (1+\|\boldsymbol\ell\|_1)^2n^{-\gamma},\qquad n\ge1.
\end{equation}
More generally, for any fixed class of nonzero integer modes satisfying
$d(\boldsymbol\ell)\ge d_0>0$, the same estimate holds with
$\mathfrak r_{d_0}$ in place of $\mathfrak r_{d_*}$ and a constant
also depending on $d_0$.
For modes with zero sum the endpoint rate is $(n\delta)^{-4/\beta}$
when $\beta>4$, $(1+\log(n\delta))/(n\delta)$ when $\beta=4$, and
$(n\delta)^{-1}$ when $\beta<4$.
\end{proposition}

The proofs use the conditional estimate in \cref{prop:unifiedblock}.
At zero Fourier mode it compares moments at indices $n\le M\le2n$;
summing its errors over successive doublings proves convergence to $H$.
For nonzero modes, the corrected recursion
\eqref{eq:corrected-recursion} gives the decay in
\cref{prop:quantitative-phase-mixing}. We prove the conditional
estimate first, then these two consequences.

\begin{proposition}\label{prop:unifiedblock}
Let $K\subset\C^m$ be compact. There are constants $C_K,c_K$,
independent of $\delta$ and the Fourier mode, such that for
$\mathbf q\in K$, $0<\delta\le1$,
$\boldsymbol\vartheta\in\mathcal S_\delta$, and integers $n\le M\le2n$
satisfying $n\ge c_KL^4$ and $n\delta\ge1$,
\begin{align}
 &\left|\E[T_M\mid\mathcal F_n]
       -e^{i(M-n)\omega}\left(\frac Mn\right)^\lambda T_n\right|
 \le \frac{C_KL^4}{n\delta}\,Y_n(\Re\mathbf q).
 \label{eq:unifiedblock}
\end{align}
Thus a fixed Fourier mode has conditional error $O_K((n\delta)^{-1})$
relative to the positive quantity $Y_n(\Re\mathbf q)$.
\end{proposition}

\subsection{One step estimates and dependence on the Fourier mode}

The one step recursion separates into a deterministic rotation and a
random multiplier. For a deterministic phase vector $\boldsymbol\varphi$,
write this multiplier, the coefficient variance, and the quadratic term
in its mean as follows:
\begin{align}
 M_k(\boldsymbol\varphi)
 &=\prod_{j=1}^m(1-\alpha_ke^{i\varphi_j})^{q_j/2-\ell_j}
                   (1-\bar\alpha_ke^{-i\varphi_j})^{q_j/2+\ell_j},
       \label{eq:multiplier}\\
 v_k&=\E|\alpha_k|^2=\frac1{1+b(k+1)},\label{eq:vk}\\
 D(\boldsymbol\varphi)
 &=\left(\sum_j(q_j/2-\ell_j)e^{i\varphi_j}\right)
   \left(\sum_j(q_j/2+\ell_j)e^{-i\varphi_j}\right).
       \label{eq:Ddef}
\end{align}
The recursions \eqref{eq:log-recursion} and \eqref{eq:phase-recursion}
give exactly
\begin{equation}\label{eq:exact-T-step}
 T_{k+1}=e^{i\omega}T_kM_k(\boldsymbol\psi_k).
\end{equation}
The powers in \eqref{eq:multiplier} use the analytic logarithm on the disk.
Since $\boldsymbol\ell$ is real, their product satisfies the useful identity
\begin{equation}\label{eq:multiplier-modulus}
 |M_k(\boldsymbol\varphi)|
   =\prod_j|1-\alpha_ke^{i\varphi_j}|^{\Re q_j}.
\end{equation}

\begin{lemma}
\label{lem:exact-mode-transition}
Fix $\mathbf q\in\mathbb C^m$ and $\boldsymbol\ell\in\mathbb R^m$,
using the real lifts of the phases in \eqref{eq:Tdef}. In this lemma
write $p_j=q_j/2$ and $t_k=b(k+1)$. The radial moments are, for integers $d\ge0$,
\[
 \mu_d(k)=\frac{d!}{(t_k+1)\cdots(t_k+d)},\qquad \mu_0(k)=1.
\]
For $\mathbf u\in\mathbb Z^m$, collect the terms producing the mode
shift $\mathbf u$ in
\begin{equation}\label{eq:exact-mode-coefficients}
 G_d^{\boldsymbol\ell}(\mathbf u)
 =\sum_{\substack{\mathbf r,\mathbf s\in\mathbb Z_{\ge0}^m\\
             |\mathbf r|=|\mathbf s|=d,\ \mathbf r-\mathbf s=\mathbf u}}
       \prod_{j=1}^m
          \binom{p_j-\ell_j}{r_j}\binom{p_j+\ell_j}{s_j}.
\end{equation}
Here $|\mathbf r|=\sum_{j=1}^m r_j$ and an empty sum is zero.
For all sufficiently large $k$,
\begin{equation}\label{eq:exact-mode-transition}
 \E[T_{k+1}(\mathbf q,\boldsymbol\ell)\mid\mathcal F_k]
 =e^{i\boldsymbol\ell\cdot\boldsymbol\vartheta}
   \sum_{d=0}^\infty\mu_d(k)
       \sum_{\mathbf u\in\mathbb Z^m}G_d^{\boldsymbol\ell}(\mathbf u)\,
                     T_k(\mathbf q,\boldsymbol\ell+\mathbf u).
\end{equation}
The inner sum is finite, and its nonzero terms satisfy
$\sum_{j=1}^m u_j=0$. The series converges absolutely, uniformly in the
phase vector. More precisely, for every compact set of parameters
$(\mathbf q,\boldsymbol\ell)$ there are $k_0$ and $C$ such that
\[
 \sum_{d=0}^{\infty}\mu_d(k)
       \sum_{\mathbf u\in\mathbb Z^m}|G_d^{\boldsymbol\ell}(\mathbf u)|
 \le C,\qquad k\ge k_0.
\]
Since $|T_k(\mathbf q,\boldsymbol\ell+\mathbf u)|
=|T_k(\mathbf q,\boldsymbol\ell)|$, this bounds the absolute sum in
\eqref{eq:exact-mode-transition} by $C|T_k(\mathbf q,\boldsymbol\ell)|$.
\end{lemma}

\begin{proof}[Proof of \cref{lem:exact-mode-transition}]
We first expand at a fixed radius and average over the argument of the
coefficient. A uniform bound on the coefficients will then justify
integration of the series over the radius.
For $|\alpha_k|<1$, expand the two products in
\eqref{eq:multiplier} by the binomial series. At each fixed radius
these series converge absolutely. Averaging over the uniform argument
of $\alpha_k$ removes terms whose total holomorphic and
antiholomorphic degrees differ. The remaining two signs are
$(-1)^d(-1)^d=1$, and
\[
 \E|\alpha_k|^{2d}
 =t_k\int_0^1u^d(1-u)^{t_k-1}\,du
 =\frac{\Gamma(t_k+1)\Gamma(d+1)}{\Gamma(t_k+d+1)}
 =\mu_d(k).
\]
The phase factor of a surviving term is
$\exp(i(\mathbf r-\mathbf s)\cdot\boldsymbol\psi_k)$.
Multiplication by $T_k(\mathbf q,\boldsymbol\ell)$ gives
\eqref{eq:exact-mode-transition}, provided the radial integration can be interchanged
with the remaining series.

We justify that interchange uniformly on compact parameter sets.
For exponents $z$ in a fixed compact subset of $\mathbb C$, the
product definition of the binomial coefficient gives constants $A,C$
such that
\[
 \left|\binom{z}{a}\right|\le C(a+1)^A,\qquad a\ge0.
\]
There are $\binom{d+m-1}{m-1}$ compositions of $d$ into $m$ parts.
Consequently, for some finite $D$,
\[
 \sum_{\mathbf u\in\mathbb Z^m}|G_d^{\boldsymbol\ell}(\mathbf u)|
 \le\sum_{|\mathbf r|=|\mathbf s|=d}
       \prod_j\left|\binom{p_j-\ell_j}{r_j}
                         \binom{p_j+\ell_j}{s_j}\right|
 \le C(d+1)^D.
\]
Choose $t_0>D+1$ and then $k_0$ with $t_{k_0}\ge t_0$.
The moments $\mu_d(k)$ decrease as $t_k$ increases, and the gamma
ratio at $t_0$ is $O((d+1)^{-t_0})$. Thus
\[
 \sum_{d=0}^{\infty}\mu_d(k)
          \sum_{\mathbf u\in\mathbb Z^m}|G_d^{\boldsymbol\ell}(\mathbf u)|
 \le C\sum_{d=0}^{\infty}(d+1)^{D-t_0}<\infty,
 \qquad k\ge k_0.
\]
This proves the required interchange and the uniform absolute
convergence. Finally, because $\mathbf u$ is real,
$|T_k(\mathbf q,\boldsymbol\ell+\mathbf u)|
 =|T_k(\mathbf q,\boldsymbol\ell)|$.
\end{proof}

The terms of degree one are the diagonal multiplier
$\sum_j(q_j^2/4-\ell_j^2)$ and the shifts
$\mathbf e_j-\mathbf e_l$. The Fourier correction
\eqref{eq:fourier-corrector} will remove
these off-diagonal shifts. The degree two terms enter the second
correction in the unpublished companion manuscript \cite{XIIIB}.
The exact identity of
\cref{lem:exact-mode-transition},
\eqref{eq:exact-mode-transition} organizes this algebra; uniform asymptotic expansions still require
estimates on the accumulated remainder and on the sum over the Fourier modes of $W_{\mathbf q}$.

\begin{lemma}\label{lem:onestep}
Let $K$ be a compact subset of $\mathbb C^m$. There are constants $C_K$ and
$k_K$ such that, for $\mathbf q\in K$, every integer vector
$\boldsymbol\ell$, every phase vector, and $k\ge k_K$,
\begin{equation}\label{eq:onestep}
 \E M_k(\boldsymbol\varphi)
  =1+v_kD(\boldsymbol\varphi)+R_k(\boldsymbol\varphi),
 \qquad |R_k(\boldsymbol\varphi)|\le \frac{C_KL^4}{(k+1)^2}.
\end{equation}
The function $D$ and its phase derivatives of orders one and two are bounded
by $C_KL^2$.
\end{lemma}

\begin{proof}[Proof of \cref{lem:onestep}]
We estimate the annulus near the boundary directly and use Taylor's
formula on the inner disk. Keeping the product in
\eqref{eq:multiplier-modulus} intact preserves polynomial dependence
on the Fourier mode.
Let $t=b(k+1)$, $u_j=\Re q_j$, $A=\sum_j(-u_j)_+$, and
$B=\sum_j(u_j)_+$. By \eqref{eq:multiplier-modulus}, the modulus of the
integrand is at most $2^B(1-|\alpha_k|)^{-A}$. The variable
$s=|\alpha_k|^2$ has density $t(1-s)^{t-1}$ on $(0,1)$. For $t>A$,
\begin{align}
 \E[(1-|\alpha_k|)^{-A}\mathbf1_{\{|\alpha_k|>1/2\}}]
 &\le 2^A t\int_{1/4}^1(1-s)^{t-A-1}\,ds\notag\\
 &=2^A\frac{t}{t-A}(3/4)^{t-A}.\label{eq:signed-annulus}
\end{align}
On a compact set $K$, this is exponentially small for $k\ge k_K$,
uniformly in the phases and in $\boldsymbol\ell$.

To apply Taylor's formula, replace the random coefficient by a
deterministic variable and put
\[
 \mathcal M_{\boldsymbol\varphi}(z)
 =\prod_{j=1}^m(1-ze^{i\varphi_j})^{q_j/2-\ell_j}
                   (1-\bar z e^{-i\varphi_j})^{q_j/2+\ell_j},
 \qquad |z|\le\tfrac12.
\]
Thus $M_k(\boldsymbol\varphi)=\mathcal M_{\boldsymbol\varphi}(\alpha_k)$
on $\{|\alpha_k|\le1/2\}$. The modulus identity
\eqref{eq:multiplier-modulus} bounds $|\mathcal M_{\boldsymbol\varphi}(z)|$
uniformly in the mode. We differentiate with respect to the two real
coordinates of $z$. Every derivative of $\log\mathcal M_{\boldsymbol\varphi}$
of order between one and four is bounded by $C_KL$: the coefficients
have total absolute value at most $C_KL$, and the denominators obtained
from $\log(1-ze^{i\varphi_j})$ and $\log(1-\bar z e^{-i\varphi_j})$
have modulus at least $1/2$. The product rule then gives
\[
 |\partial^\gamma\mathcal M_{\boldsymbol\varphi}(z)|
 \le C_KL^{|\gamma|},\qquad
 \gamma\in\mathbb Z_{\ge0}^2,\quad |\gamma|\le4,
\]
where the right-hand side is $C_K$ when $|\gamma|=0$.
Taylor's formula through total degree three therefore has remainder
bounded by $C_KL^4|z|^4$. After substituting $z=\alpha_k$, rotation
invariance eliminates every monomial $z^r\bar z^s$ with $r\ne s$,
including all odd degree monomials. The constant coefficient is one
and the coefficient of $z\bar z$ is $D(\boldsymbol\varphi)$.
Furthermore,
\[
 \E|\alpha_k|^4=\frac2{(1+t)(2+t)}.
\]
The integrals over the omitted annulus of the Taylor polynomial of degree
three have size
at most $C_KL^3(3/4)^t$, and that of the original product is bounded by
\eqref{eq:signed-annulus}. This proves \eqref{eq:onestep}. Applying a bound
to each analytic factor separately would lose the stated polynomial
dependence on $L$; the modulus identity for their product is the reason
that the preceding argument has uniform constants.
For the phase derivatives, expand
\[
 D(\boldsymbol\varphi)
 =\sum_{j,l=1}^m
   \left(\frac{q_j}{2}-\ell_j\right)
   \left(\frac{q_l}{2}+\ell_l\right)
      e^{i(\varphi_j-\varphi_l)}.
\]
Each first or second coordinate derivative multiplies a summand by
a factor of modulus at most one. Moreover,
\[
 \sum_{j,l=1}^m
   \left|\frac{q_j}{2}-\ell_j\right|
   \left|\frac{q_l}{2}+\ell_l\right|
 \le C_K\left(1+\sum_{j=1}^m|\ell_j|\right)^2=C_KL^2.
\]
This proves the stated bounds for $D$ and its phase derivatives.
\end{proof}

The random part of the phase increment in \eqref{eq:phase-recursion} is
\[
 \xi_{k,j}=-2\Im\log(1-\alpha_ke^{i\psi_k(\vartheta_j)}),
 \qquad \boldsymbol\xi_k=(\xi_{k,j})_{j=1}^m.
\]

\begin{lemma}
\label{lem:phase-moments}
For $\mathbf q$ in a compact subset $K$ of $\mathbb C^m$ and
$k\ge k_K$,
\begin{equation}\label{eq:mixed-one-step}
 \left|\E[M_k(\boldsymbol\psi_k)\,\xi_{k,j}
                    \mid\mathcal F_k]\right|
       \le \frac{C_KL}{k+1}.
\end{equation}
We also have the mode independent bound
\begin{equation}\label{eq:weighted-phase-second}
 \E[|M_k(\boldsymbol\psi_k)|\,|\boldsymbol\xi_k|^2
                     \mid\mathcal F_k]\le \frac{C_K}{k+1}.
\end{equation}
\end{lemma}

\begin{proof}[Proof of \cref{lem:phase-moments}]
Conditional rotation invariance gives
$\E[\xi_{k,j}\mid\mathcal F_k]=0$.
Moreover $|\xi_{k,j}|\le\pi$, because
$1-\alpha_ke^{i\psi_k(\vartheta_j)}$ has positive real part, and
$|\boldsymbol\xi_k|\le C|\alpha_k|$ when $|\alpha_k|\le1/2$.
On this disk the first derivative bound in the proof of
\cref{lem:onestep} gives $|M_k-1|\le C_KL|\alpha_k|$.
Subtracting the zero conditional mean of $\xi_{k,j}$, we therefore
bound the inner disk contribution to the mixed moment by
$C_KL\E|\alpha_k|^2\le C_KL/(k+1)$.
On the complementary annulus, use the boundedness of $\xi_{k,j}$,
\eqref{eq:multiplier-modulus}, and \eqref{eq:signed-annulus}.
Its contribution is exponentially small, which proves
\eqref{eq:mixed-one-step} with an initial index independent of the mode.
For \eqref{eq:weighted-phase-second}, the inner disk integrand is at
most $C_K|\alpha_k|^2$. The annulus estimate \eqref{eq:signed-annulus} controls the complement.
\end{proof}

\begin{lemma}\label{lem:conditional-weights}
Let $K_0$ be a compact subset of $\mathbb R^m$, and fix $p\ge1$. There
are $n_0$ and $C$, depending on $K_0,p,\beta,m$, such that
\begin{equation}\label{eq:Yconditionalmoments}
 \E[Y_{n+r}(\mathbf u)^p\mid\mathcal F_n]
       \le C Y_n(\mathbf u)^p,
 \qquad \mathbf u\in K_0,\quad n\ge n_0,\quad 0\le r\le n.
\end{equation}
No lower bound on the coordinates of $p\mathbf u$ is required here.
\end{lemma}

\begin{proof}[Proof of \cref{lem:conditional-weights}]
Apply \cref{lem:onestep} with exponents $p\mathbf u$ and zero Fourier mode.
For every $k\ge n_0$ the nonnegative one step multiplier has conditional
mean at most $1+C/(k+1)$. Its conditional expectation is finite by
\eqref{eq:signed-annulus}, after increasing $n_0$. Iteration gives
\[
 \E[Y_{n+r}(\mathbf u)^p\mid\mathcal F_n]
 \le Y_n(\mathbf u)^p
          \prod_{k=n}^{n+r-1}\left(1+\frac C{k+1}\right)
 \le C'Y_n(\mathbf u)^p.
\]
These are conditional statements for a finite, almost surely positive
initial weight, whose unconditional $p$th moment need not be finite.
\end{proof}

\subsection{Proof of the conditional estimate}

The diagonal part of $D$ gives the exponent $\lambda$. Its other Fourier
terms oscillate under the deterministic rotation
$\boldsymbol\varphi\mapsto\boldsymbol\varphi+\boldsymbol\vartheta$.
We incorporate their discrete antiderivative into the observable. This
leaves an error of order $k^{-2}\delta^{-1}$ at time $k$, whose sum from
$n$ to infinity is of order $(n\delta)^{-1}$.

More precisely, for a fixed mode and fixed points,
\eqref{eq:exact-T-step} and \cref{lem:onestep} give
\[
 \E[T_{k+1}\mid\mathcal F_k]
 =e^{i\omega}T_k\bigl(1+v_kD_0+v_kD_{\rm off}(\boldsymbol\psi_k)\bigr)
  +O(k^{-2})|T_k|,
\]
with $D_0$ and $D_{\rm off}$ defined below. Since $v_k\sim1/(bk)$ and
$D_0=b\lambda$, the product of the factors $1+v_kD_0$ is of order
$k^{\lambda}$. The term $v_kD_{\rm off}(\boldsymbol\psi_k)$ is also of
order $1/k$ and cannot simply be discarded; but $D_{\rm off}$ has no
constant Fourier term, and the phases $\boldsymbol\psi_k$ advance by
$\boldsymbol\vartheta$ at each step up to a small random correction.
Replacing $T_k$ by $T_k(1+v_kQ_{\boldsymbol\ell}(\boldsymbol\psi_k))$,
where $Q_{\boldsymbol\ell}$ solves \eqref{eq:corrector-equation},
cancels this term and leaves a remainder of order $k^{-2}\delta^{-1}$,
as proved in \eqref{eq:corrected-remainder}.

Write
\begin{align}
 D_0&=\sum_j(q_j^2/4-\ell_j^2)=b\lambda,\notag\\
 D_{\rm off}(\boldsymbol\varphi)
 &=\sum_{j\ne r}c_{jr}e^{i(\varphi_j-\varphi_r)},
 \qquad c_{jr}=(q_j/2-\ell_j)(q_r/2+\ell_r),
 \label{eq:offdiagonal-drift}
\end{align}
so that $D=D_0+D_{\rm off}$. For distinct points put
\begin{equation}\label{eq:fourier-corrector}
 Q_{\boldsymbol\ell}(\boldsymbol\varphi)
 =-\sum_{j\ne r}
 \frac{c_{jr}e^{i(\varphi_j-\varphi_r)}}
      {e^{i(\vartheta_j-\vartheta_r)}-1}.
\end{equation}
The dependence on the exponents and the points is suppressed
in the notation. Direct substitution gives
\begin{equation}\label{eq:corrector-equation}
 Q_{\boldsymbol\ell}(\boldsymbol\varphi+\boldsymbol\vartheta)
       -Q_{\boldsymbol\ell}(\boldsymbol\varphi)
       =-D_{\rm off}(\boldsymbol\varphi).
\end{equation}
For $\boldsymbol\vartheta\in\mathcal S_\delta$, $0<\delta\le1$, the
finite Fourier sum also gives
\begin{equation}\label{eq:corrector-derivatives}
 \|Q_{\boldsymbol\ell}\|_{C^2(\T^m)}
       \le C_K L^2\delta^{-1}.
\end{equation}
Only pairwise differences of angles occur in the denominators.
No condition on sums of several angles is needed. If the exponents
are real and $\boldsymbol\ell=0$, the terms with indices $(j,r)$ and
$(r,j)$ are complex conjugates; hence $Q_0$ is real valued.

\begin{lemma}
\label{lem:corrected-recursion}
Let $K\subset\C^m$ be compact and define
\begin{equation}\label{eq:corrected-observable}
 \widetilde T_k=T_k
        \bigl(1+v_kQ_{\boldsymbol\ell}(\boldsymbol\psi_k)\bigr).
\end{equation}
There are constants $C_K,c_K$, independent of the Fourier mode and of
$\delta$, such that for $k\ge c_K L^2$ and
$\boldsymbol\vartheta\in\mathcal S_\delta$, $0<\delta\le1$,
\begin{align}
 \E[\widetilde T_{k+1}\mid\mathcal F_k]
 &=e^{i\omega}(1+v_kD_0)\widetilde T_k+\mathcal R_k,
       \label{eq:corrected-recursion}\\
 |\mathcal R_k|&\le
       \frac{C_K L^4}{(k+1)^2\delta}\,Y_k(\Re\mathbf q).
       \label{eq:corrected-remainder}
\end{align}
No lower bound on the real parts of the exponents is required for this
conditional statement.
\end{lemma}

\begin{proof}[Proof of \cref{lem:corrected-recursion}]
We expand the correction at the deterministically advanced phase,
then use \eqref{eq:corrector-equation} to cancel the oscillatory drift.
All quantities depending on $\boldsymbol\psi_k$ are fixed under the
conditioning. Abbreviate
$Q_+=Q_{\boldsymbol\ell}(\boldsymbol\psi_k+\boldsymbol\vartheta)$.
The phase recursion gives
$\boldsymbol\psi_{k+1}=\boldsymbol\psi_k+\boldsymbol\vartheta
 +\boldsymbol\xi_k$. Taylor's formula in these real variables yields
\[
 Q_{\boldsymbol\ell}(\boldsymbol\psi_{k+1})
 =Q_++\nabla Q_{\boldsymbol\ell}
       (\boldsymbol\psi_k+\boldsymbol\vartheta)\boldsymbol\xi_k
       +O_K(L^2\delta^{-1}|\boldsymbol\xi_k|^2).
\]
Since $k\ge c_KL^2$, \cref{lem:onestep} gives
\[
 \bigl|\E[M_k\mid\mathcal F_k]-1\bigr|
 \le \frac{C_KL^2}{k+1}+\frac{C_KL^4}{(k+1)^2}
 \le \frac{C'_KL^2}{k+1}.
\]
Multiplying the Taylor expansion by $M_k$ and using
\eqref{eq:mixed-one-step}, \eqref{eq:weighted-phase-second} and
\eqref{eq:corrector-derivatives}, we therefore obtain
\begin{align}
 &\left|\E[M_kQ_{\boldsymbol\ell}(\boldsymbol\psi_{k+1})
                          \mid\mathcal F_k]-Q_+\right|\notag\\
 &\quad\le
 \frac{C_KL^4}{(k+1)\delta}
 +\frac{C_KL^3}{(k+1)\delta}
 +\frac{C_KL^2}{(k+1)\delta}
 \le\frac{C'_KL^4}{(k+1)\delta}.
 \label{eq:corrector-mixed-expansion}
\end{align}
The signed mixed first moment is used in the second term. Estimating it
by $\E[|M_k||\boldsymbol\xi_k|\mid\mathcal F_k]$ would give a larger
error and would lose summability at the required order.

Using \eqref{eq:exact-T-step}, we now have
\begin{align*}
 \E[\widetilde T_{k+1}\mid\mathcal F_k]
 =e^{i\omega}T_k\biggl\{
  1+v_kD(\boldsymbol\psi_k)+v_{k+1}Q_+
  +O_K\left(\frac{L^4}{(k+1)^2\delta}\right)\biggr\}.
\end{align*}
Here the remainder of \cref{lem:onestep} is absorbed because
$\delta\le1$. Moreover
$|v_{k+1}-v_k|\le C_\beta(k+1)^{-2}$, so replacing $v_{k+1}$ by $v_k$
in front of $Q_+$ costs at most $C_KL^2/((k+1)^2\delta)$.
Equation \eqref{eq:corrector-equation} then gives
\[
 v_kD(\boldsymbol\psi_k)+v_kQ_+
       =v_kD_0+v_kQ_{\boldsymbol\ell}(\boldsymbol\psi_k).
\]
Finally,
\begin{align*}
 &(1+v_kD_0)
       (1+v_kQ_{\boldsymbol\ell}(\boldsymbol\psi_k))\\
 &\quad=1+v_kD_0+v_kQ_{\boldsymbol\ell}(\boldsymbol\psi_k)
       +v_k^2D_0Q_{\boldsymbol\ell}(\boldsymbol\psi_k),
\end{align*}
and the last term has modulus at most
$C_KL^4/((k+1)^2\delta)$. This proves the lemma, using
$|T_k|=Y_k(\Re\mathbf q)$.
\end{proof}

\begin{proof}[Proof of \cref{prop:unifiedblock}]
We iterate the corrected recursion, bound the accumulated remainders,
and finally remove the correction at both endpoints. The products of
the deterministic multipliers will be denoted by $P_{r,s}$.
Put $\tau_\beta=1+1/b$. Since $D_0=b\lambda$,
\[
 1+v_kD_0=1+\frac{\lambda}{k+\tau_\beta}.
\]
For any subinterval of $[n,M)$ define
\[
 P_{r,s}=\prod_{k=r}^{s-1}
       \left(1+\frac{\lambda}{k+\tau_\beta}\right),
       \qquad n\le r\le s\le M.
\]
We have $|\lambda|\le C_KL^2$ and $\Re\lambda\le C_K$.
Increasing $c_K$ if necessary ensures
$|\lambda/(k+\tau_\beta)|\le1/2$ throughout the interval. The inequality
$\log|1+z|\le\Re z+C|z|^2$ therefore gives
\begin{align*}
 \log|P_{r,s}|
 &\le C_K\sum_{k=r}^{s-1}(k+\tau_\beta)^{-1}
       +C_KL^4\sum_{k=r}^{s-1}(k+\tau_\beta)^{-2}
 \le C'_K.
\end{align*}
The last bound uses $M\le2n$ and $n\ge c_KL^4$. It is independent of
the Fourier mode.

Iterate \eqref{eq:corrected-recursion} and condition on
$\mathcal F_n$. By \cref{lem:conditional-weights}, the conditional
expectation of every factor $Y_k(\Re\mathbf q)$ is bounded
by $C_KY_n(\Re\mathbf q)$. Consequently,
\begin{align*}
 \left|\E[\widetilde T_M\mid\mathcal F_n]
       -e^{i(M-n)\omega}P_{n,M}\widetilde T_n\right|
 &\le \frac{C_KL^4}{\delta}Y_n(\Re\mathbf q)
                   \sum_{k=n}^{M-1}(k+1)^{-2}\\
 &\le\frac{C'_KL^4}{n\delta}Y_n(\Re\mathbf q).
\end{align*}
The endpoint corrections satisfy
\[
 |\widetilde T_n-T_n|
       \le\frac{C_KL^2}{n\delta}Y_n(\Re\mathbf q),
 \qquad
 \E[|\widetilde T_M-T_M|\mid\mathcal F_n]
       \le\frac{C_KL^2}{n\delta}Y_n(\Re\mathbf q).
\]
Removing them leaves the same error bound with $T_M$ and $T_n$.

To evaluate the deterministic product, use the analytic logarithm near
one and the elementary sum estimate
\[
 \sum_{k=n}^{M-1}(k+\tau_\beta)^{-1}
       =\log(M/n)+O_\beta(n^{-1}).
\]
It follows that
\[
 \log P_{n,M}=\lambda\log(M/n)+O_K(L^4/n).
\]
The error in this exponential is bounded, and
$|(M/n)^\lambda|\le C_K$. Hence
\[
 |P_{n,M}-(M/n)^\lambda|\le C_KL^4/n
                              \le C_KL^4/(n\delta).
\]
This proves \eqref{eq:unifiedblock}.
All estimates above are uniform in the incoming phase vector.
Thus, for every fixed history $(\alpha_0,\ldots,\alpha_{n-1})\in\D^n$,
the same bound holds when the conditional expectation is evaluated
by integrating the independent future coefficients. Its constants do
not depend on that history; the only remaining dependence is the
factor $Y_n(\Re\mathbf q)$.
\end{proof}

\subsection{Proof of the holomorphic convergence}

We first prove the real moment bound in \cref{lem:negative-envelope},
which allows us to average \eqref{eq:unifiedblock}.
The comparison in \cref{lem:relative-moment-comparison} will also
show that the limit is strictly positive for real exponents.

\begin{lemma}
\label{lem:negative-envelope}
For every compact $K_0\subset(-1-b,\infty)^m$ and every $\delta>0$, there is
a constant $C_{K_0,\delta}$ such that
\begin{equation}\label{eq:real-envelope}
 \E Y_n(\mathbf u;\boldsymbol\vartheta)
       \le C_{K_0,\delta}\,n^{a(\mathbf u)},
 \qquad n\ge1,\quad \mathbf u\in K_0,
 \quad\boldsymbol\vartheta\in\mathcal S_\delta.
\end{equation}
\end{lemma}

\begin{proof}[Proof of \cref{lem:negative-envelope}]
We may replace $\delta$ by $\min(\delta,1)$. Take zero Fourier mode and
real exponents in \cref{prop:unifiedblock}. Choose a fixed integer
$n_*\ge c_{K_0}$ with $n_*\delta\ge1$. For
$Z_n(\mathbf u)=n^{-a(\mathbf u)}\E Y_n(\mathbf u)$, positivity gives
\[
 Z_{2n}(\mathbf u)\le Z_n(\mathbf u)
                 \{1+C_{K_0}(n\delta)^{-1}\},\qquad n\ge n_*.
\]
The initial values $Z_{n_*}$ are uniformly bounded on
$K_0\times\mathcal S_\delta$ by \cref{lem:negative-initial}. Along
$n=2^rn_*$, the product of the displayed factors is bounded, because
$\sum_{r\ge0}(2^rn_*\delta)^{-1}<\infty$. For any $N\ge n_*$,
choose $r$ with $2^rn_*\le N<2^{r+1}n_*$ and apply \eqref{eq:unifiedblock} once more. The finitely many smaller indices are covered by
\cref{lem:negative-initial}. No moment with doubled exponents is needed
at the fixed initial index.
\end{proof}

The positivity of real moments also allows a comparison whose constant
is independent of the separation. It will normalize the phase estimates
without using the value of the limiting moment.

\begin{lemma}
\label{lem:relative-moment-comparison}
Let $K_0\subset(-1-b,\infty)^m$ be compact, and put
$A_n=\E Y_n(\mathbf q;\boldsymbol\vartheta)$. There are constants
$C_{K_0},R_{K_0}$ such that, for $\mathbf q\in K_0$,
$0<\delta\le1$, $\boldsymbol\vartheta\in\mathcal S_\delta$, and
all integers $M\ge n$ with $n\delta\ge R_{K_0}$,
\begin{equation}\label{eq:relative-moment-comparison}
 \left|\left(\frac Mn\right)^{-a(\mathbf q)}
              \frac{A_M}{A_n}-1\right|
       \le\frac{C_{K_0}}{n\delta}.
\end{equation}
In particular the ratio on the left before subtracting one is bounded
above and below by positive constants independent of $M,n,\delta$ and
the points.
\end{lemma}

\begin{proof}[Proof of \cref{lem:relative-moment-comparison}]
The finite moments $A_n$ are strictly positive. For $n\le M\le2n$,
\eqref{eq:relative-moment-comparison} follows by taking expectations in
\eqref{eq:unifiedblock} with zero Fourier mode and dividing by
$(M/n)^{a(\mathbf q)}A_n$. Since $a(\mathbf q)\ge0$, this division
does not enlarge the error bound.
For a general $M$, choose consecutive indices
$n_0=n$, $n_{j+1}=\min(2n_j,M)$, stopping when $n_J=M$. Write
\[
 \left(\frac{n_{j+1}}{n_j}\right)^{-a(\mathbf q)}
       \frac{A_{n_{j+1}}}{A_{n_j}}=1+e_j,
 \qquad |e_j|\le\frac{C_{K_0}}{n_j\delta}.
\]
Choose $R_{K_0}$ large enough that every $|e_j|\le1/2$.
All these ratios are positive. Taking their real logarithms and using
$|\log(1+x)|\le2|x|$ for $|x|\le1/2$ gives
\begin{align*}
 \left|\log\left\{\left(\frac Mn\right)^{-a(\mathbf q)}
                 \frac{A_M}{A_n}\right\}\right|
 &\le2\sum_{j=0}^{J-1}|e_j|
 \le\frac{C'_{K_0}}{n\delta}\sum_{j\ge0}2^{-j}
 \le\frac{C''_{K_0}}{n\delta}.
\end{align*}
Exponentiating proves the result and the asserted uniform comparison.
\end{proof}

\begin{proof}[Proof of \cref{thm:holomorphic-limit}]
We show that the normalized moments
$Z_n(\mathbf q)=n^{-a(\mathbf q)}\E Y_n(\mathbf q)$ form a locally
uniform Cauchy sequence. The conditional estimate
\eqref{eq:unifiedblock} with zero Fourier mode, followed by
\cref{lem:negative-envelope}, gives, for $n\le M\le2n$ and
sufficiently large $n$,
\begin{align*}
 |Z_M(\mathbf q)-Z_n(\mathbf q)|
 &\le C_{K,\delta}n^{-\Re a(\mathbf q)-1}
                           \E Y_n(\Re\mathbf q)\\
 &\le C_{K,\delta}
    n^{a(\Re\mathbf q)-\Re a(\mathbf q)-1}
 =C_{K,\delta}n^{a(\Im\mathbf q)-1}.
\end{align*}
Here $a(\Re\mathbf q)-\Re a(\mathbf q)=a(\Im\mathbf q)$, by expansion
of the squares. Summing over successive doublings proves uniform
convergence along the powers of two, because $\rho_K>0$. Applying the
same bound once between a power of two and an arbitrary intermediate
integer proves convergence along all integers. Summing the geometric
tail gives \eqref{eq:holomorphic-limit-rate}.

Each $Z_n$ is holomorphic in the exponents and jointly continuous with
distinct angles by \cref{lem:negative-initial}. Local uniform
convergence therefore proves both assertions for $H$. For real exponents,
fix $n$ with $n\delta\ge R_{K_0}$ and use
\cref{lem:relative-moment-comparison}. It bounds $Z_M/Z_n$ above and
below by strictly positive constants for every $M\ge n$. Taking the
limit gives $H>0$, since $Z_n>0$.
\end{proof}

\subsection{Proof of the phase estimates}

\begin{proof}[Proof of \cref{prop:quantitative-phase-mixing}]
We first divide out the growth of the zero mode. The resulting scalar
recursion has a damping term, which we estimate separately for bounded
and large modes before removing the Fourier correction.
All exponents in the proof are real. Put $a=a(\mathbf q)\ge0$,
$d=d(\boldsymbol\ell)$, $A_n=\E Y_n(\mathbf q)$ and
$\tau_\beta=1+1/b$. In this proof, constants denoted by $C$ may depend
on $K_0,\beta,m$ and, when used, the fixed lower bound $d_0$, but not
on the mode, the points or $\delta$.
To divide out the zero mode multiplier exactly, define
\begin{equation}\label{eq:gamma-normalization}
 p_n(a)=\frac{\Gamma(n+\tau_\beta+a)}
                    {\Gamma(n+\tau_\beta)}.
\end{equation}
Then $p_n>0$ and
\[
 \frac{p_{n+1}(a)}{p_n(a)}=1+\frac{a}{n+\tau_\beta},
 \qquad p_n(a)=n^a(1+O_{K_0}(n^{-1})).
\]
Choose
\[
 n_* =\left\lceil\frac{C_{K_0}L^2}{\delta}\right\rceil
\]
with $C_{K_0}$ sufficiently large. This ensures that
\eqref{eq:corrected-recursion} applies for every $k\ge n_*$,
that the threshold $n_*\delta\ge R_{K_0}$ in
\cref{lem:relative-moment-comparison} is satisfied,
that $|v_kQ_{\boldsymbol\ell}|\le1$, and that
$0\le1-d/(k+\tau_\beta+a)\le1$ throughout this range.
Notice that
\begin{equation}\label{eq:phase-start-scale}
 1\le n_*\delta\le C'_{K_0}L^2.
\end{equation}

The constant $h=A_{n_*}/p_{n_*}(a)>0$ matches the normalization to
the moment at the starting index. By
\cref{lem:relative-moment-comparison}, $A_k/A_{n_*}$ is comparable
to $(k/n_*)^a$; the gamma ratio $p_{n_*}(a)/p_k(a)$ is comparable
to its reciprocal. Hence
\begin{equation}\label{eq:phase-normalizer-comparison}
 C^{-1}\le\frac{A_k}{p_k(a)h}\le C,
       \qquad k\ge n_*.
\end{equation}
These constants are independent of $\delta$ and $\boldsymbol\ell$.
This argument uses only finite real moments and does not require
identification of $H$.

Define the normalized corrected mode
\[
 b_k=e^{-ik\omega}\frac{\E\widetilde T_k}{p_k(a)h}.
\]
Take expectations in \eqref{eq:corrected-recursion}, then divide by
$p_{k+1}(a)h$ and remove the phase. The deterministic coefficient simplifies
exactly as follows:
\[
 \frac{1+(a-d)/(k+\tau_\beta)}
      {1+a/(k+\tau_\beta)}
       =1-\frac{d}{k+\tau_\beta+a}.
\]
Using \eqref{eq:phase-normalizer-comparison} in the remainder gives
\begin{equation}\label{eq:normalized-corrected-mode}
 b_{k+1}=\left(1-\frac{d}{k+\tau_\beta+a}\right)b_k+e_k,
 \qquad |e_k|\le\frac{CL^4}{(k+1)^2\delta},
       \qquad k\ge n_*.
\end{equation}
Moreover $|b_{n_*}|\le2$, because
$|\widetilde T_{n_*}|\le2Y_{n_*}$ by the choice of $n_*$.

A nonzero integer vector has $d\ge2/\beta$. If its coordinate sum
vanishes, at least two coordinates are nonzero, so $d\ge4/\beta$.
We prove the more general assertion for a class satisfying $d\ge d_0$.
The distinction between bounded and large modes is useful here: retaining
the actual damping of a large mode recovers two powers of $L$.

Put $c=\tau_\beta+a$, which ranges in a fixed compact subset of
$(0,\infty)$. For $n_*\le s\le n$, nonnegativity of the scalar
coefficients and comparison with an integral give
\begin{align}\label{eq:phase-actual-propagator}
 \prod_{k=s}^{n-1}\left(1-\frac d{k+c}\right)
 &\le\exp\left(-d\sum_{k=s}^{n-1}\frac1{k+c}\right)
 \le\left(\frac{s+c}{n+c}\right)^d\notag\\
 &\le\exp(cd/s)(s/n)^d\le C(s/n)^d.
\end{align}
The last constant is independent of the mode, since
$d\le C_\beta L^2$ and $s\ge n_*\ge CL^2/\delta$.
Similarly, replacing $(k+1)^d$ by $k^d$ costs at most
$\exp(d/k)\le C$ for $k\ge n_*$. Thus variation of constants in
\eqref{eq:normalized-corrected-mode} yields
\begin{equation}\label{eq:mode-forcing-sum}
 |b_n|\le C\left(\frac{n_*}{n}\right)^d
       +\frac{CL^4}{\delta n^d}
                    \sum_{k=n_*}^{n-1}k^{d-2},\qquad n\ge n_*.
\end{equation}

Suppose first that $d>2$. The function $x^{d-2}$ is increasing, so
\[
 \sum_{k=n_*}^{n-1}k^{d-2}
 \le\int_{n_*}^{n}x^{d-2}\,dx
 \le\frac{n^{d-1}}{d-1}.
\]
For nonzero integer modes,
$L\le(1+\sqrt m)\|\boldsymbol\ell\|$, so
$L^2\le C_{\beta,m}d$. Hence, with
$t=n\delta$ and $t_*=n_*\delta\le CL^2$,
\[
 |b_n|\le C(n_*/n)^d+\frac{CL^4}{(d-1)t}
       \le Cn_*/n+\frac{CL^2}{t}
       \le\frac{CL^2}{t}.
\]
Since $t^{-1}\le\mathfrak r_{d_0}(t)$ for $t\ge1$, this is the
required estimate for these modes.

The modes with $0<d\le2$ form a finite set, depending only on
$\beta,m$. In particular $L$ and $t_*$ are bounded on this set.
In \eqref{eq:phase-actual-propagator} we may replace $d$ by the
lower bound $d_0$, obtaining
\[
 |b_n|\le C(n_*/n)^{d_0}
       +\frac{C}{\delta n^{d_0}}
                   \sum_{k=n_*}^{n-1}k^{d_0-2}.
\]
The sum is bounded by $Cn_*^{d_0-1}$ if $d_0<1$, by
$C(1+\log(n/n_*))$ if $d_0=1$, and by $Cn^{d_0-1}$ if
$d_0>1$. Since $1\le t_*\le C$, these bounds give
$|b_n|\le C\mathfrak r_{d_0}(t)$.
If $d_0>2$, this finite set is empty.

It remains to remove the correction. From
\eqref{eq:corrected-observable} and \eqref{eq:corrector-derivatives},
\[
 \frac{|\E\widetilde T_n-\E T_n|}{A_n}
       \le\frac{CL^2}{n\delta}.
\]
Also $p_n(a)h/A_n$ is bounded by
\eqref{eq:phase-normalizer-comparison}. We have therefore proved
\[
 \frac{|\E T_n|}{A_n}\le CL^2\mathfrak r_{d_0}(n\delta),
       \qquad n\ge n_*.
\]
For $n<n_*$, use $|\E T_n|\le A_n$ and $1\le t\le CL^2$.
If $d_0<1$, then $L^2t^{-d_0}\ge cL^{2-2d_0}\ge c$;
if $d_0\ge1$, then $L^2\mathfrak r_{d_0}(t)\ge L^2/t\ge c$.
Increasing the constant covers this range and proves the stated
class estimate, including \eqref{eq:relative-phase-mixing}.

Finally, for a fixed $\delta>0$ replace it by $\min(\delta,1)$ if necessary.
For every $0<\gamma<\min\{d_*,1\}$,
$\mathfrak r_{d_*}(t)\le C_\gamma t^{-\gamma}$ on $t\ge1$.
Combine \eqref{eq:relative-phase-mixing} with
\cref{lem:negative-envelope} to obtain
\eqref{eq:quantitative-phase-mixing} for $n\delta\ge1$.
For the remaining indices, the same conclusion follows from
$|\E T_n|\le A_n\le C_{K_0,\delta}n^a$ and $n<\delta^{-1}$.
\end{proof}

\section{Identification of the limit by multiplicative chaos}
\label{sec:gmc-identification}

We now determine the function $H$ constructed in
\cref{thm:holomorphic-limit}. Recall that $D_\beta$ is the one point
constant in \eqref{eq:opuc-onepoint-asymptotic}.
\Cref{thm:opuc-full} identifies the limit for distinct points, and
\cref{thm:opuc-relative} retains a relative error whose constant is
independent of their separation.

\begin{theorem}
\label{thm:opuc-full}
On the domain $\mathcal U$ of \eqref{eq:tube}, the limit in
\cref{thm:holomorphic-limit} is
\begin{equation}\label{eq:opuc-limit-identified}
 H(\mathbf q;\boldsymbol\vartheta)
 =\prod_{j=1}^mD_\beta(q_j)
   \prod_{j<k}|e^{i\vartheta_j}-e^{i\vartheta_k}|^{-q_jq_k/\beta}.
\end{equation}
Consequently, for every compact $K\subset(-1-b,\infty)^m$ and
$\delta>0$, uniformly for $\mathbf q\in K$ and
$\boldsymbol\vartheta\in \mathcal S_\delta$,
\begin{align}
 n^{-a(\mathbf q)}\E Y_n(\mathbf q;\boldsymbol\vartheta)
 &=\prod_jD_\beta(q_j)
   \prod_{j<k}|e^{i\vartheta_j}-e^{i\vartheta_k}|^{-q_jq_k/\beta}
       +O_{\beta,m,K,\delta}(n^{-1}),
       \label{eq:opuc-full-asymptotic}\\
 \frac{\E Y_n(\mathbf q;\boldsymbol\vartheta)}
             {\prod_j\nu_n(q_j)}
 &=\prod_{j<k}|e^{i\vartheta_j}-e^{i\vartheta_k}|^{-q_jq_k/\beta}
       +O_{\beta,m,K,\delta}(n^{-1}).
       \label{eq:opuc-full-normalized}
\end{align}
The complex convergence is locally uniform in $\mathcal U$ and in
separated configurations, with the rate of
\cref{thm:holomorphic-limit}.
\end{theorem}

\begin{theorem}
\label{thm:opuc-relative}
For every compact $K\subset(-1-b,\infty)^m$ there are constants $C_K,R_K$
such that, for $0<\delta\le1$, $\boldsymbol\vartheta\in\mathcal S_\delta$,
$\mathbf q\in K$ and $n\delta\ge R_K$,
\begin{align}
 \E Y_n(\mathbf q;\boldsymbol\vartheta)
 &=n^{a(\mathbf q)}\prod_{j=1}^mD_\beta(q_j)
   \prod_{j<k}|e^{i\vartheta_j}-e^{i\vartheta_k}|^{-q_jq_k/\beta}
       \left(1+O_K\left(\frac1{n\delta}\right)\right),
       \label{eq:relative-opuc-asymptotic}\\
 \frac{\E Y_n(\mathbf q;\boldsymbol\vartheta)}
      {\prod_j\nu_n(q_j)}
 &=\prod_{j<k}|e^{i\vartheta_j}-e^{i\vartheta_k}|^{-q_jq_k/\beta}
       \left(1+O_K\left(\frac1{n\delta}\right)\right).
       \label{eq:relative-opuc-normalized}
\end{align}
The constants are independent of $\delta$ and the points.
Consequently these asymptotics hold for any varying configuration whose
minimum chord distance $\delta_n$ satisfies $n\delta_n\to\infty$.
\end{theorem}

The identification has two steps. For small real exponents, the common
Gaussian field in \cref{prop:chaos-input} computes the limit after
integration in the angles. The local uniform convergence already proved
in \cref{sec:cmv-unified} then identifies the density at every
prescribed distinct configuration. Holomorphic continuation extends the
identity to $\mathcal U$. The relative estimate follows by letting
$M\to\infty$ in \cref{lem:relative-moment-comparison}; this retains
constants independent of the separation.

The averaged statement used in the first step is \cref{thm:gmc-averaged}
below. Its proof uses the exact conditional expectation identity
\eqref{eq:gmc-exact-conditional} of Lambert and Najnudel and the
Gaussian exponential formula \eqref{eq:gmc-gaussian-tensor}.

For real exponents $q_j>-1$, write the normalized moment kernel and
its proposed limit as
\begin{align}
 R_n(\mathbf q;\boldsymbol\theta)
 &=\frac{\E Y_n(\mathbf q;\boldsymbol\theta)}
             {\prod_j\nu_n(q_j)},
       \label{eq:gmc-opuc-kernel}\\
 K_{\mathbf q}(\boldsymbol\theta)
 &=\prod_{j<k}|e^{i\theta_j}-e^{i\theta_k}|^{-q_jq_k/\beta}.
       \label{eq:gmc-target-kernel}
\end{align}
The first expression is a nonnegative measurable function; values on
partial diagonals may be infinite and have no effect on its Lebesgue
integrals. Its values at distinct configurations are finite by
\cref{lem:negative-initial}.

\begin{theorem}
\label{thm:gmc-averaged}
Let $m\geq2$, $q_j>-1$, and suppose
\begin{equation}\label{eq:gmc-charge-domain}
 \sum_{j=1}^m\frac{q_j^2}{2\beta}<1.
\end{equation}
Then $K_{\mathbf q}$ belongs to $L^1(\T^m,\sigma^{\otimes m})$ and,
for every $F\in C(\T^m)$,
\begin{equation}\label{eq:gmc-averaged-opuc}
 \int_{\T^m}F(\boldsymbol\theta)
              R_n(\mathbf q;\boldsymbol\theta)\,
                         \sigma^{\otimes m}(d\boldsymbol\theta)
 \longrightarrow
 \int_{\T^m}F(\boldsymbol\theta)
              K_{\mathbf q}(\boldsymbol\theta)\,
                         \sigma^{\otimes m}(d\boldsymbol\theta).
\end{equation}
\end{theorem}

\subsection{The common Gaussian field and the exact martingale identity}

The exact connection between the circular beta ensemble and Gaussian
multiplicative chaos was established by Chhaibi and
Najnudel~\cite{CN2026}. We use it in the following form, due to Lambert
and Najnudel, and state explicitly the coupling and normalization in
which it is applied.
\begin{proposition}[Lambert--Najnudel]
\label{prop:chaos-input}
In the coefficient coupling \eqref{eq:coefficient-law}, there is a
Gaussian analytic function $\varphi$, measurable with respect to
$\mathcal F_\infty$, such that, almost surely,
\[
 \log\Phi_n^*(z)\longrightarrow\varphi(z)
 \quad\text{locally uniformly in }\D.
\]
Its covariances are
\begin{equation}\label{eq:common-gaf}
 \E[\varphi(z)\overline{\varphi(w)}]
       =-\frac2\beta\log(1-z\bar w),
 \qquad \E[\varphi(z)\varphi(w)]=0.
\end{equation}
The boundary real part $X=\Re\varphi$ exists as a Gaussian generalized
function, with covariance at distinct points
\begin{equation}\label{eq:gmc-covariance}
 \E[X(\theta)X(\phi)]
       =-\frac1\beta\log|e^{i\theta}-e^{i\phi}|.
\end{equation}
For real $q>-1$ with $q^2<2\beta$, let $U^{(q)}$ be the multiplicative
chaos defined by the radial limit \eqref{eq:gmc-radial} below,
normalized to have expectation $\sigma$.

All these measures are constructed from the same $X$.
Recall that $\nu_n(q)=\E|\Phi_n^*(1)|^q$ is the angle independent
one point moment in \eqref{eq:opuc-onepoint-exact}. For bounded
continuous $f$, define
\begin{equation}\label{eq:gmc-opuc-measure}
 U_n^{(q)}(f)
 =\int_\T f(\theta)
       \frac{|\Phi_n^*(e^{i\theta})|^q}{\nu_n(q)}\,\sigma(d\theta).
\end{equation}
Then
\begin{equation}\label{eq:gmc-exact-conditional}
 U_n^{(q)}(f)=\E[U^{(q)}(f)\mid\mathcal F_n].
\end{equation}
In particular, $U_n^{(q)}$ converges almost surely weakly to $U^{(q)}$.
Moreover, with $X_\rho(\theta)=\Re\varphi(\rho e^{i\theta})$,
\begin{equation}\label{eq:gmc-radial}
 \int_\T f(\theta)
   e^{qX_\rho(\theta)-q^2\E X_\rho(\theta)^2/2}\,\sigma(d\theta)
 \longrightarrow U^{(q)}(f)
 \quad\text{in }L^1\text{ as }\rho\uparrow1.
\end{equation}
The measures $U^{(q)}$ are almost surely nonatomic.
\end{proposition}

\paragraph{References for \cref{prop:chaos-input}.}
The analytic limit and its covariance are
\cite[Proposition~2, Section~1.2, and Appendix~A]{LN2026}; in particular,
the discussion preceding its Theorem~2 identifies $\varphi$ as a
measurable function of the same Verblunsky coefficients. The chaos measures are
constructed in its Proposition~3. Its Proposition~4 is exactly
\eqref{eq:gmc-exact-conditional}, with the parameters there equal to
$\gamma=q$ and $\hat\gamma=q/\sqrt{2\beta}$; its Theorem~2 gives the
almost sure weak convergence. Equation~(22) of that paper gives
\eqref{eq:gmc-radial}, in fact in the larger $L^p$ range used below.
Nonatomicity follows from the positive exact dimension
$1-q^2/(2\beta)$ of these subcritical measures, recalled after its
Proposition~3. Only finitely many exponents are used at a time.
The almost sure convergence and nonatomicity statements just listed
therefore hold together on one event of probability one.

The restriction $q>-1$ in the proposition is enough for our purpose and
lies within the parameter domain already established for $\nu_n$; the cited
OPUC result itself holds for every real $q$ with $q^2<2\beta$.
The usual moment bound for the same circular chaos is
\begin{equation}\label{eq:gmc-moment-input}
 p q^2<2\beta,\quad p>1
 \quad\Longrightarrow\quad
 \E[U^{(q)}(\T)^p]<\infty.
\end{equation}
For a precise reference, \cite[Theorem~1.1]{Remy2020} evaluates these
moments; its parameter $\gamma=\sqrt{2/\beta}|q|$ gives exactly
\eqref{eq:gmc-moment-input} when $q\ne0$; for $q=0$ the measure is
deterministically $\sigma$. The same total mass law is obtained in
\cite[Corollary~2.5]{CN2026}. For negative $q$ one applies the construction for a positive
exponent to $-X$; this changes neither the moment threshold
nor the marginal total mass law. The joint construction still uses the
original common field with its signed exponents.

\begin{lemma}
\label{lem:gmc-measure-moments}
Let $q>-1$, $q^2<2\beta$, and $f\in C(\T)$. For every finite $r\geq1$
with $r q^2<2\beta$,
\begin{equation}\label{eq:gmc-opuc-lr}
 U_n^{(q)}(f)\longrightarrow U^{(q)}(f)\quad\text{in }L^r.
\end{equation}
The corresponding $r$th moments of the total masses are uniformly
bounded in $n$.
\end{lemma}

\begin{proof}[Proof of \cref{lem:gmc-measure-moments}]
Put $Z=U^{(q)}(f)$. By \eqref{eq:gmc-moment-input},
$Z\in L^r$; for $r=1$ this follows from the mean of the total mass.
The variable $Z$ is measurable with respect to $\mathcal F_\infty$, and
$U_n^{(q)}(f)=\E[Z\mid\mathcal F_n]$ by
\eqref{eq:gmc-exact-conditional}. The upward convergence theorem for
conditional expectations gives the asserted $L^r$ convergence.
Explicitly, approximate $Z$ in $L^r$ by a variable $Z_N$ measurable
with respect to some $\mathcal F_N$. For $n\ge N$, the contraction
property gives
\[
 \|\E[Z\mid\mathcal F_n]-Z\|_r
 \le 2\|Z-Z_N\|_r.
\]
Finally the conditional Jensen inequality, applied with $f=1$, gives
$\E[U_n^{(q)}(\T)^r]\le\E[U^{(q)}(\T)^r]<\infty$.
\end{proof}

\subsection{Gaussian Fourier approximation on the same probability space}

The products needed for \cref{thm:gmc-averaged} are the random tensor
products $\bigotimes_{j=1}^mU^{(q_j)}$ on $\T^m$, all constructed from
the same field $X$. Their expected measure is defined by
\[
 \int_{\T^m}F\,d\E\!\left[\bigotimes_{j=1}^mU^{(q_j)}\right]
 =\E\!\left[\int_{\T^m}F(\theta_1,\ldots,\theta_m)
                  \bigotimes_{j=1}^mU^{(q_j)}(d\theta_j)\right],
 \qquad F\in C(\T^m).
\]
Here $m\ge2$, the exponents are real, $q_j>-1$, and
$\sum_{j=1}^m q_j^2/(2\beta)<1$. To choose H\"older exponents, distribute the positive difference
$1-\sum_jq_j^2/(2\beta)$ among the $m$ coordinates. This gives
$r_j>1$ with $\sum_{j=1}^m1/r_j=1$ and $r_jq_j^2<2\beta$.
By H\"older's inequality and \eqref{eq:gmc-moment-input},
\[
 \E\prod_{j=1}^m U^{(q_j)}(\T)
 \le\prod_{j=1}^m\|U^{(q_j)}(\T)\|_{r_j}<\infty,
\]
so the expected measure is finite. On product tests its integral is
$\E\prod_{j=1}^m U^{(q_j)}(f_j)$, for $f_j\in C(\T)$.
No independence between the factors is assumed.

We evaluate this expected measure by truncating the common Gaussian
field. The required moment convergence follows from
\eqref{eq:gmc-radial}; no comparison of successive Fourier covariance
kernels is needed.

The Taylor coefficients of the same analytic field $\varphi$ define
independent standard real Gaussian variables $(A_k,B_k)_{k\ge1}$ by
\[
 \varphi(z)=\frac1{\sqrt\beta}\sum_{k\ge1}
                 \frac{A_k-iB_k}{\sqrt{k}}z^k,\qquad z\in\D.
\]
Consequently,
\[
 X_\rho(\theta)=\frac1{\sqrt\beta}\sum_{k\ge1}
       \frac{\rho^k}{\sqrt{k}}
       \bigl(A_k\cos(k\theta)+B_k\sin(k\theta)\bigr).
\]
Here $L$ denotes the Fourier truncation index.
Its truncated field, the sigma-field generated by its coefficients,
and its normalized exponential measure are
\begin{align}
 X_L(\theta)
 &=\frac1{\sqrt\beta}\sum_{k=1}^L
        \frac{A_k\cos(k\theta)+B_k\sin(k\theta)}{\sqrt{k}},
       \label{eq:gaussian-fourier-field}\\
 \mathcal G_L&=\sigma(A_k,B_k:1\leq k\leq L),
 \qquad
 V_L^{(q)}(d\theta)
 =e^{qX_L(\theta)-q^2\E X_L(\theta)^2/2}\sigma(d\theta).
       \label{eq:gaussian-fourier-measure}
\end{align}
These variables are the Gaussian coefficients of the same limiting
field $X$ in \cref{prop:chaos-input}.

\begin{lemma}
\label{lem:gaussian-fourier-martingale}
For $q>-1$, $q^2<2\beta$ and bounded continuous $f$,
\begin{equation}\label{eq:fourier-conditional}
 V_L^{(q)}(f)=\E[U^{(q)}(f)\mid\mathcal G_L].
\end{equation}
Consequently, whenever $r\geq1$ and $r q^2<2\beta$,
\[
 V_L^{(q)}(f)\longrightarrow U^{(q)}(f)\quad\text{in }L^r,
 \qquad
 \sup_L\E[V_L^{(q)}(\T)^r]
       \leq\E[U^{(q)}(\T)^r].
\]
\end{lemma}

\begin{proof}[Proof of \cref{lem:gaussian-fourier-martingale}]
We condition the radial approximation on its first $L$ Gaussian modes,
then pass to the boundary in $L^1$. The conditional expectation identity
will also give the moment convergence.
For $0<\rho<1$, multiply the $k$th coefficient of
\eqref{eq:gaussian-fourier-field} by $\rho^k$, and denote the resulting
field by $X_{\rho,L}$. Its limit as $L\to\infty$ is the radial field
$X_\rho$ of \eqref{eq:gmc-radial}. Let $V_\rho^{(q)}$ denote its
normalized exponential measure. The field
$X_\rho-X_{\rho,L}$ is independent of $\mathcal G_L$. Its Gaussian
exponential, after subtraction of half its variance, has expectation
one. The conditional Fubini theorem therefore gives
\begin{equation}\label{eq:fourier-radial-conditioning}
 \E[V_\rho^{(q)}(f)\mid\mathcal G_L]
 =\int_\T f(\theta)
       e^{qX_{\rho,L}(\theta)
               -q^2\E X_{\rho,L}(\theta)^2/2}\,\sigma(d\theta).
\end{equation}
For $f\geq0$ this follows from Tonelli's theorem, and for bounded signed $f$ from
integrability of the total mass, whose expectation is one.

Fix $L$ and let $\rho\uparrow1$. By \eqref{eq:gmc-radial} and the
$L^1$ contraction of conditional expectation, the left side converges in
$L^1$ to $\E[U^{(q)}(f)\mid\mathcal G_L]$. On the right,
$X_{\rho,L}\to X_L$ uniformly in $\theta$ for each realization.
Furthermore its exponential is bounded by
\[
 \exp\!\left(\frac{|q|}{\sqrt\beta}
       \sum_{k=1}^L\frac{|A_k|+|B_k|}{\sqrt{k}}\right),
\]
an integrable random variable. Dominated convergence proves convergence
of the right-hand side in $L^1$ to $V_L^{(q)}(f)$, and hence
\eqref{eq:fourier-conditional}.

The limit $U^{(q)}(f)$ is measurable with respect to all Gaussian modes,
because it is obtained by radial regularization of that field. Its
$r$th moment is finite by \eqref{eq:gmc-moment-input} and the bound by
$\|f\|_\infty U^{(q)}(\T)$; the case $r=1$ follows from its mean.
The upward convergence theorem for conditional expectations gives the
claimed $L^r$ convergence. The conditional Jensen inequality gives the
uniform bound.
\end{proof}

\subsection{Expected products of the chaos measures}

\begin{proof}[Proof of \cref{thm:gmc-averaged}]
We first pass to the limit in products of OPUC chaos integrals.
Gaussian Fourier truncation identifies the resulting expected measure
away from the diagonals; nonatomicity then excludes additional mass
on the diagonals.
Choose $w_j\in(0,1)$ with $\sum_jw_j=1$ and
$w_j>q_j^2/(2\beta)$; the strict inequality
\eqref{eq:gmc-charge-domain} makes this possible. Set $r_j=1/w_j$.
Then $r_j>1$, $\sum_jr_j^{-1}=1$, and $r_jq_j^2<2\beta$.
For $f_j\in C(\T)$, put
$Z_{n,j}=U_n^{(q_j)}(f_j)$ and $Z_j=U^{(q_j)}(f_j)$.
By \cref{lem:gmc-measure-moments},
$\|Z_{n,j}-Z_j\|_{r_j}\to0$. The identity
\[
 \prod_jZ_{n,j}-\prod_jZ_j
 =\sum_{k=1}^m (Z_{n,k}-Z_k)
              \prod_{j<k}Z_{n,j}\prod_{j>k}Z_j
\]
and H\"older's inequality show that the left side converges to zero
in $L^1$: its norm is bounded by
\[
 \sum_{k=1}^m\|Z_{n,k}-Z_k\|_{r_k}
                 \prod_{j<k}\|Z_{n,j}\|_{r_j}
                 \prod_{j>k}\|Z_j\|_{r_j}.
\]
Every factor other than the indicated difference is uniformly bounded.
Taking $f_j=1$ also gives a uniform bound on the expected products of
the total masses. Finite sums of product test functions form an algebra
containing the constants and separating the points of $\T^m$.
The Stone--Weierstrass theorem makes this algebra uniformly dense in
$C(\T^m)$. The mass bound controls the error of this approximation
uniformly in $n$, and consequently proves
\begin{equation}\label{eq:gmc-tensor-limit}
 \E\!\left[\bigotimes_{j=1}^mU_n^{(q_j)}\right]
 \Longrightarrow
 \E\!\left[\bigotimes_{j=1}^mU^{(q_j)}\right]
\end{equation}
as finite measures. No independence between the factors is used.

It remains to identify the measure on the right.
\Cref{lem:gaussian-fourier-martingale} and the same H\"older argument
give the identical expected tensor product limit when $U_n^{(q_j)}$
is replaced by $V_L^{(q_j)}$. At finite $L$, the Gaussian exponential
identity gives
\begin{equation}\label{eq:gmc-gaussian-tensor}
 \E\!\left[\bigotimes_jV_L^{(q_j)}(d\theta_j)\right]
 =\exp\!\left\{\sum_{j<k}q_jq_k C_L(\theta_j,\theta_k)\right\}
                   \prod_j\sigma(d\theta_j),
\end{equation}
where
\[
 C_L(\theta,\phi)=\frac1\beta\sum_{k=1}^L
                                  \frac{\cos(k(\theta-\phi))}{k}.
\]
On sets with $|e^{i\theta}-e^{i\phi}|\geq\delta>0$, summation by
parts bounds the tail of this Fourier series by $C/(L\delta)$.
Thus $C_L$ converges uniformly there to
$-\beta^{-1}\log|e^{i\theta}-e^{i\phi}|$.
Integrating \eqref{eq:gmc-gaussian-tensor} against continuous tests
supported in products of arcs with disjoint closures identifies the
limiting density on
\[
 \mathcal S_*=\{\boldsymbol\theta\in\T^m:
                         \theta_j\ne\theta_k\text{ for }j\ne k\}
\]
as $K_{\mathbf q}$. Such products form an open cover of
$\mathcal S_*$.

There is no additional mass on the partial diagonals. Indeed, for each
pair $j\ne k$, the nonatomicity in \cref{prop:chaos-input} gives,
almost surely,
\[
 (U^{(q_j)}\otimes U^{(q_k)})\{(\theta,\phi):\theta=\phi\}
 =\int_\T U^{(q_k)}(\{\theta\})\,U^{(q_j)}(d\theta)=0.
\]
The other total masses are finite almost surely, so the full tensor
product also gives zero mass to this diagonal. The same holds after
expectation and for the finite union of partial diagonals. Therefore
the entire limiting expected measure has density $K_{\mathbf q}$,
up to a Lebesgue null set. Its finite total mass proves the asserted
integrability. Finally the left side of \eqref{eq:gmc-tensor-limit}
has density $R_n$, by Tonelli, which proves
\eqref{eq:gmc-averaged-opuc}.
\end{proof}

\subsection{Proofs of the pointwise formulas}

\begin{proof}[Proof of \cref{thm:opuc-full}]
We first identify the pointwise limit for small positive real exponents,
then extend the identity by holomorphy. The error bounds will follow
from the convergence already proved in \cref{thm:holomorphic-limit}.
For $m=1$, the statement is \cref{lem:opuc-onepoint}. Suppose $m\geq2$,
and first take real positive exponents with $\sum_jq_j^2<2\beta$.
By \cref{thm:holomorphic-limit,lem:opuc-onepoint}, the kernel $R_n$
converges locally uniformly on $\mathcal S_*$ to the continuous function
\[
 H(\mathbf q;\boldsymbol\vartheta)\Big/\prod_jD_\beta(q_j).
\]
For any continuous test supported in a compact subset of
$\mathcal S_*$, local uniform convergence permits passage of the limit
through its integral. The averaged convergence of
\cref{thm:gmc-averaged} identifies that integral with the integral of
$K_{\mathbf q}$. Both candidate densities are continuous on
$\mathcal S_*$. If they differed at a point, a nonnegative continuous
test supported in a sufficiently small neighborhood would give
different integrals. They therefore agree at every distinct
configuration.

Now fix such a configuration. The function $H$ is holomorphic on
$\mathcal U$ by \cref{thm:holomorphic-limit}. The expression on the
right of \eqref{eq:opuc-limit-identified} is holomorphic there as well:
$D_\beta$ is holomorphic for $\Re q>-1-b$, and each chord power means
\[
 \exp\!\left\{-\frac{q_jq_k}{\beta}
                 \log|e^{i\vartheta_j}-e^{i\vartheta_k}|\right\},
\]
with the real logarithm of its positive base.
These two holomorphic functions agree on a small open real cube of
positive exponents satisfying $\sum_jq_j^2<2\beta$.
Choose a complex polydisk around a point of this cube whose real
coordinate intervals remain in the cube. Apply the one variable
identity theorem successively in each coordinate, keeping the
coordinates not yet considered real. The difference vanishes on the
whole polydisk. Since
\[
 \mathcal U=\{\Re q_j>-1-b,\ a(\Im\mathbf q)<1\}
\]
is convex, hence connected, the identity theorem extends the equality throughout
$\mathcal U$. This proves \eqref{eq:opuc-limit-identified} without
requiring chaos moments at large exponents.

For real exponents, \cref{thm:holomorphic-limit} gives the uniform error
$O(n^{-1})$ before identification of the limit; inserting the
identified expression gives \eqref{eq:opuc-full-asymptotic}.
The one point estimate \eqref{eq:opuc-onepoint-asymptotic} gives
\[
 n^{-a(\mathbf q)}\prod_j\nu_n(q_j)
       =\prod_jD_\beta(q_j)(1+O(n^{-1})).
\]
The positive constants $D_\beta(q_j)$ are bounded away from zero on
the real compact set. Division therefore proves
\eqref{eq:opuc-full-normalized} with the stated uniformity.
\end{proof}

The two convergence statements have different roles.
\Cref{thm:gmc-averaged} identifies an expected measure after integration,
whereas \cref{thm:holomorphic-limit}, proved from
\cref{prop:unifiedblock}, supplies the local uniform convergence needed
at prescribed distinct points. The latter conclusion does not follow
from weak convergence of measures alone.

\begin{proof}[Proof of \cref{thm:opuc-relative}]
Fix the exponents and a distinct configuration, and write
$A_n=\E Y_n(\mathbf q;\boldsymbol\vartheta)$ and
$Z_n=n^{-a(\mathbf q)}A_n$. By
\cref{lem:relative-moment-comparison}, for all $M\ge n$ with
$n\delta\ge R_K$,
\[
 \left|\frac{Z_M}{Z_n}-1\right|\le\frac{C_K}{n\delta}.
\]
The OPUC limit already identified in \cref{thm:opuc-full} is strictly
positive at this real vector of exponents. Letting $M\to\infty$ gives
\[
 \left|\frac{H(\mathbf q;\boldsymbol\vartheta)}{Z_n}-1\right|
       \le\frac{C_K}{n\delta}.
\]
Increase $R_K$ until the right-hand side is at most $1/2$. Inverting this
positive ratio yields $Z_n/H=1+O_K((n\delta)^{-1})$. Substitution of
\eqref{eq:opuc-limit-identified} proves
\eqref{eq:relative-opuc-asymptotic}. Although the limit was taken with
the configuration fixed, the bound before passage to the limit has a
constant independent of that configuration and of $\delta$. The same
uniformity therefore holds in the resulting estimate.

Finally, the one point asymptotic \eqref{eq:opuc-onepoint-asymptotic} gives
\[
 \prod_j\nu_n(q_j)=n^{a(\mathbf q)}\prod_jD_\beta(q_j)
                       (1+O_K(n^{-1})).
\]
Its positive leading constant is bounded away from zero on $K$.
Division proves \eqref{eq:relative-opuc-normalized}, because
$n^{-1}\le(n\delta)^{-1}$ for $\delta\le1$.
\end{proof}

\section{Quantitative transfer to the characteristic polynomial}
\label{sec:terminal}

The remaining task is to estimate the contribution of the last
independent coefficient in \eqref{eq:terminal-polynomial}.
Recall that $\boldsymbol\psi_n=(\psi_n(\vartheta_j))_{j=1}^m$
is the vector of phases from \eqref{eq:terminal-process}.
For real $q_j>-1$ and distinct angles, the phase law weighted by $Y_n$
is the probability measure on $\T^m$ defined by
\[
 \mu_{n,\mathbf q}(B)
 =\frac{\E[Y_n(\mathbf q)
                     \mathbf1_{\{\boldsymbol\psi_n\in B\}}]}
        {\E Y_n(\mathbf q)}.
\]
Here the phases are taken modulo $2\pi$. Recall from
\eqref{eq:terminal-W} that
\[
 W_{\mathbf q}(\boldsymbol\varphi)
 =\int_\T\prod_{j=1}^m
       |1-e^{i(\eta+\varphi_j)}|^{q_j}\,\sigma(d\eta).
\]
The exact identity \eqref{eq:terminal-expectation} gives
\begin{equation}\label{eq:terminal-conditioned}
 F_{n+1}(\mathbf q,\boldsymbol\vartheta)
 =\E Y_n(\mathbf q)\int W_{\mathbf q}\,d\mu_{n,\mathbf q}.
\end{equation}

All integrands before conditioning are nonnegative; Tonelli's theorem
applies, and finiteness follows from \cref{lem:negative-initial}.

By \cref{thm:opuc-relative}, the asymptotic of $\E Y_n$ is already
known. We therefore need
\[
 \int W_{\mathbf q}\,d\mu_{n,\mathbf q}
       =\prod_jc(q_j)+o(1).
\]
The target constant follows from $D_\beta(q)c(q)=C_\beta(q)$ in
\eqref{eq:terminal-c}. We first prove a quantitative version on
$\mathcal D_m$. On this domain
$W_{\mathbf q}$ is bounded, and its Fourier coefficients are summable
with a positive power of the mode. \Cref{sec:full-negative-transfer}
will remove this restriction for fixed distinct points.

For a compact set $K\subset\mathcal D_m$, let $\varepsilon_K$ measure
the distance from the threshold for the sum of negative exponents,
and choose the H\"older exponent $\alpha_K$ by
\begin{equation}\label{eq:terminal-alpha}
 \varepsilon_K
 =1-\sup_{\mathbf q\in K}\sum_j(-q_j)_+>0,
 \qquad \alpha_K=\frac{\varepsilon_K}{2m}.
\end{equation}
In particular $0<\alpha_K\leq1/2$.

\begin{theorem}
\label{thm:signed-characteristic-FH}
Let $\beta>0$, $m\ge2$, and $K\subset\mathcal D_m$ be compact. For every
\begin{equation}\label{eq:terminal-rate}
 0<\kappa<\varepsilon_K\min\left\{\frac2\beta,\frac12\right\},
\end{equation}
there are $C,R<\infty$ such that, uniformly for $0<\delta\le1$,
$\mathbf q\in K$, $\boldsymbol\vartheta\in\mathcal S_\delta$, and
$N\delta\ge R$,
\begin{align}\label{eq:signed-characteristic-FH}
 F_N(\mathbf q,\boldsymbol\vartheta)
 ={}&N^{a(\mathbf q)}\prod_jC_\beta(q_j)
       \prod_{j<l}|e^{i\vartheta_j}-e^{i\vartheta_l}|^{-q_jq_l/\beta}
       \bigl(1+O((N\delta)^{-\kappa})\bigr).
\end{align}
The constant in the error does not depend on $\delta$. In particular the
formula holds along every sequence for which $N\delta_N\to\infty$.
For $m=1$, \cref{prop:onepoint} gives relative error $O(N^{-1})$.
\end{theorem}

The proof uses two estimates on $W_{\mathbf q}$: continuity in
\cref{lem:terminal-holder}, and the weighted Fourier summability in
\cref{lem:terminal-weighted-fourier}. Combining the latter with
\eqref{eq:relative-phase-mixing} gives the error in the theorem. We
establish these analytic estimates before carrying out that summation.

\subsection{Translation estimates for singular powers}

For real $q>-1$ let $h_q(t)=|1-e^{it}|^q$, interpreted as a measurable
function on $\mathbb T$; its value at zero is immaterial for integration.
Write $\|t\|_{2\pi}=\min_{k\in\mathbb Z}|t-2\pi k|$ for angular distance to zero. All $L^p$ norms in this section use normalized Haar measure $\sigma$.

\begin{lemma}\label{lem:terminal-translation}
Fix $0<\alpha<1$, $P<\infty$, and $Q<\infty$. Uniformly for
\[
 |q|\leq Q,\qquad 1\leq p\leq P,
 \qquad q+\frac1p\geq2\alpha,
\]
the function $h_q$ satisfies
\begin{equation}\label{eq:terminal-translation}
 \|h_q\|_{L^p(\sigma)}\leq C,
 \qquad
 \|h_q(\,\cdot+t)-h_q\|_{L^p(\sigma)}
       \leq C\|t\|_{2\pi}^{\alpha}.
\end{equation}
\end{lemma}

\begin{proof}[Proof of \cref{lem:terminal-translation}]
We integrate directly near the singularity and use a derivative bound
on its complement. Write $r(x)=\|x\|_{2\pi}$. On $0<r(x)\leq\pi$, the chord length
$|1-e^{ix}|$ is comparable to $r(x)$, with absolute comparison constants.
Since $pq+1=p(q+1/p)\geq2\alpha p$, integration of $r^{pq}$ gives
the first estimate, uniformly on the stated parameter set.

It suffices to prove the translation estimate for $0<|t|\leq1/4$;
larger translations are bounded by twice the first norm. Let
$E_t=\{x:r(x)\le2|t|\text{ or }r(x+t)\le2|t|\}$ be the region
near the two possible singularities. Both distances are at most
$3|t|$ there, so
\[
 \int_{E_t}|h_q(x+t)-h_q(x)|^p\,\sigma(dx)
 \leq C\int_0^{3|t|}r^{pq}\,\mathrm dr
 =C\frac{(3|t|)^{pq+1}}{pq+1}.
\]
The denominator is uniformly positive because $pq+1\geq2\alpha p$.
Its $p$th root is bounded by
$C|t|^{q+1/p}\leq C|t|^\alpha$.

On the complement, the segment joining $x$ to $x+t$ avoids the
singularity and its distance to the singularity is comparable to $r(x)$.
Differentiation gives $|h_q'(x)|\leq C r(x)^{q-1}$ away from zero,
uniformly for $|q|\leq Q$. The mean value estimate therefore bounds the
remaining integral by
\[
 C|t|^p\int_{2|t|}^{\pi}r^{p(q-1)}\,\mathrm dr
 \leq C|t|^{\alpha p}
           \int_{2|t|}^{\pi}r^{p(q-\alpha)}\,\mathrm dr
 \leq C|t|^{\alpha p}.
\]
For the middle inequality use $|t|/r\leq1$ and $\alpha<1$.
The last integral is uniformly bounded because
$p(q-\alpha)+1\geq\alpha p\geq\alpha$. Taking $p$th roots proves
\eqref{eq:terminal-translation}. All constants remain uniform as $p$
varies in $[1,P]$.
\end{proof}

\begin{lemma}
\label{lem:terminal-holder}
Uniformly for $\mathbf q\in K$, $W_{\mathbf q}$ is bounded on
$\mathbb T^m$ and satisfies
\begin{equation}\label{eq:terminal-holder}
 |W_{\mathbf q}(\boldsymbol\varphi+\mathbf t)
       -W_{\mathbf q}(\boldsymbol\varphi)|
 \leq C_K\sum_{j=1}^m\|t_j\|_{2\pi}^{\alpha_K}.
\end{equation}
Moreover,
\begin{equation}\label{eq:terminal-mean}
 \int_{\mathbb T^m}W_{\mathbf q}\,d\sigma^{\otimes m}
       =\prod_{j=1}^m c(q_j),
\end{equation}
where, for real $q>-1$,
\[
 c(q)=\int_\T|1-e^{i\eta}|^q\,\sigma(d\eta)
     =\frac{\Gamma(1+q)}{\Gamma(1+q/2)^2}.
\]
For $\boldsymbol\ell\in\mathbb Z^m$, the Fourier coefficient
\[
 \widehat W_{\mathbf q}(\boldsymbol\ell)
 =\int_{\T^m}W_{\mathbf q}(\boldsymbol\varphi)
             e^{-i\boldsymbol\ell\cdot\boldsymbol\varphi}
                 \,\sigma^{\otimes m}(d\boldsymbol\varphi)
\]
vanishes unless $\sum_{j=1}^m\ell_j=0$.
\end{lemma}

\begin{proof}[Proof of \cref{lem:terminal-holder}]
We choose H\"older exponents by assigning weight $(-q_j)_+$ to each
singular factor and distributing the remaining weight equally.
For each $\mathbf q\in K$, put
\[
 a_j=(-q_j)_+,\qquad \rho=1-\sum_j a_j,
 \qquad \frac1{p_j}=a_j+\frac\rho m.
\]
Then $\sum_j1/p_j=1$, $1\leq p_j\leq m/\varepsilon_K$, and
\[
 q_j+\frac1{p_j}=q_j+a_j+\frac\rho m
       \geq\frac{\varepsilon_K}{m}=2\alpha_K.
\]
Consequently \cref{lem:terminal-translation} applies with constants
independent of the chosen vector of exponents, even though the exponents
$p_j$ depend on that vector. H\"older's inequality gives
\[
 W_{\mathbf q}(\boldsymbol\varphi)
 \leq\prod_j\|h_{q_j}\|_{L^{p_j}(\sigma)}\leq C_K.
\]
This also proves that its defining integral is finite at every phase
configuration, including coincident phases.

To compare $W_{\mathbf q}(\boldsymbol\varphi+\mathbf t)$ and
$W_{\mathbf q}(\boldsymbol\varphi)$, telescope the product of the $m$
factors. In its $j$th difference, apply H\"older with exponent $p_j$ to
the translated difference of $h_{q_j}$ and with the remaining exponents
to the other factors. Translation invariance of the norms and
\eqref{eq:terminal-translation} bound that term by
$C_K\|t_j\|_{2\pi}^{\alpha_K}$. Their sum proves
\eqref{eq:terminal-holder}.

Tonelli's theorem and separate translation invariance in each phase
give \eqref{eq:terminal-mean}; here each $c(q_j)$ is finite because
$q_j>-1$. Finally
$W_{\mathbf q}(\boldsymbol\varphi+t(1,\ldots,1))
=W_{\mathbf q}(\boldsymbol\varphi)$ by changing the integration variable
$\eta$. Its Fourier coefficient must therefore equal itself multiplied
by $e^{it\sum_j\ell_j}$ for every $t$. This proves the asserted
vanishing.
\end{proof}

\subsection{Exact Fourier coefficients and their summability}

The function $W_{\mathbf q}$ has more structure than a general continuous
function on $\mathbb T^m$. Its Fourier coefficients factor into
one variable coefficients and lie on the hyperplane where the sum of
the indices vanishes. We use both properties to avoid a loss depending
on the number of Fourier modes.

\begin{lemma}
\label{lem:root-fourier-coefficients}
For $q>-1$ and $k\in\mathbb Z$,
\begin{equation}\label{eq:root-fourier-coefficients}
 \widehat h_q(k)
 =\frac{(-1)^k\Gamma(1+q)}
  {\Gamma(1+q/2-k)\Gamma(1+q/2+k)}.
\end{equation}
Reciprocal gamma factors in this formula are interpreted by their entire
continuations. Uniformly for $q$ in a compact interval
$J\subset(-1,\infty)$,
\begin{equation}\label{eq:root-fourier-decay}
 |\widehat h_q(k)|\le C_J(1+|k|)^{-1-q}.
\end{equation}
\end{lemma}

\begin{proof}[Proof of \cref{lem:root-fourier-coefficients}]
Integration by parts gives a recurrence for the Fourier coefficients.
We solve it from the beta integral at zero, then use gamma ratios to
estimate its solution uniformly in $q$.
Evenness gives $\widehat h_q(-k)=\widehat h_q(k)$. The beta integral
gives the value at $k=0$, namely $c(q)$. On $(0,2\pi)$,
$h'_q(t)=\frac q2\cot(t/2)h_q(t)$. Since
$\sin(t/2)h_q(t)$ vanishes at both endpoints and its derivative is
integrable for $q>-1$, integration by parts gives, for $k\ge1$,
\begin{align*}
 0&=\int_0^{2\pi}\frac d{dt}
      \{e^{-i(k-1/2)t}\sin(t/2)h_q(t)\}\,dt\\
  &=\frac{\pi}{2}\{(q+2-2k)\widehat h_q(k-1)
                              +(q+2k)\widehat h_q(k)\}.
\end{align*}
Thus
\[
 \widehat h_q(k)
 =\frac{k-1-q/2}{k+q/2}\widehat h_q(k-1).
\]
Iteration from $\widehat h_q(0)=c(q)$ proves
\eqref{eq:root-fourier-coefficients}, including its zero coefficients
when $q$ is a nonnegative even integer.

For $k$ larger than a fixed bound depending on $J$, reflection gives
\[
 \widehat h_q(k)
 =-\frac{\Gamma(1+q)\sin(\pi q/2)}\pi
            \frac{\Gamma(k-q/2)}{\Gamma(k+1+q/2)}.
\]
Uniformly for $u,v$ in a fixed compact subset of $\C$,
\begin{equation*}\tag{6.9a}\label{eq:uniform-gamma-ratio}
 \frac{\Gamma(k+u)}{\Gamma(k+v)}
       =k^{u-v}\bigl(1+O(k^{-1})\bigr),\qquad k\to\infty.
\end{equation*}
This follows by subtracting the uniform Stirling expansions of the
two logarithms; see \href{https://dlmf.nist.gov/5.11.E13}{DLMF, equation (5.11.13)}.
Taking $u=-q/2$ and $v=1+q/2$ proves
\eqref{eq:root-fourier-decay} for these $k$. The sine factor causes no
loss near $q=0$ or a positive even integer. The finitely many remaining
coefficients are bounded uniformly by the original formula and the
continuity of reciprocal gamma.
\end{proof}

\begin{lemma}
\label{lem:terminal-weighted-fourier}
Let $K\subset\mathcal D_m$ be compact. For $m\ge2$, put
\begin{equation}\label{eq:terminal-subset-index}
 s_K=\inf_{\mathbf q\in K}\left\{
       1+\min_{\substack{J\subset\{1,\ldots,m\}\\|J|\ge2}}
                          \sum_{j\in J}q_j\right\},
\end{equation}
and set $s_K=\infty$ for $m=1$. Then $s_K\ge\varepsilon_K>0$,
and, for every $0<s<s_K$,
\begin{equation}\label{eq:terminal-weighted-fourier}
 \sup_{\mathbf q\in K}\sum_{\boldsymbol\ell\in\mathbb Z^m}
 (1+\|\boldsymbol\ell\|_1)^s
          |\widehat W_{\mathbf q}(\boldsymbol\ell)|<\infty.
\end{equation}
The coefficients are given by
\begin{equation}\label{eq:terminal-coefficients}
 \widehat W_{\mathbf q}(\boldsymbol\ell)
 =\mathbf1_{\{\sum_j\ell_j=0\}}
                         \prod_j\widehat h_{q_j}(\ell_j).
\end{equation}
In particular the estimate holds for every $0<s<\varepsilon_K$.
If $K\subset[q_*,\infty)^m$ with $q_*>0$, it holds for every
$0<s<1+2q_*$. For $m\ge2$, the sharper value in
\eqref{eq:terminal-subset-index} is one plus the infimum over $K$
of the sum of the two smallest exponents when every coordinate is positive.
\end{lemma}

\begin{proof}[Proof of \cref{lem:terminal-weighted-fourier}]
We first factor the coefficients. In each shell of modes with zero sum,
at least two coordinates are large; their decay will compensate for
the number of modes in the shell.
Write the Fourier coefficient as
\[
 \widehat W_{\mathbf q}(\boldsymbol\ell)
 =\int_\T\int_{\T^m}
       \prod_{j=1}^m h_{q_j}(\eta+\varphi_j)
       e^{-i\boldsymbol\ell\cdot\boldsymbol\varphi}
       \,\sigma^{\otimes m}(d\boldsymbol\varphi)\,\sigma(d\eta).
\]
Fubini's theorem applies because the integral of the absolute value is
$\prod_jc(q_j)<\infty$. In this integral make the substitutions
$u_j=\eta+\varphi_j$.

The remaining factor is
$e^{i\eta\sum_j\ell_j}$, whose integral vanishes unless
$\sum_j\ell_j=0$. This proves \eqref{eq:terminal-coefficients}.
When $m=1$, only the zero coefficient remains and every positive
weight exponent is allowed. We henceforth suppose $m\ge2$.

Every subset sum of the exponents is at least
$-\sum_j(-q_j)_+$, so $s_K\ge\varepsilon_K$.
The coefficient estimate \eqref{eq:root-fourier-decay} implies,
uniformly over the present compact set,
\begin{equation}\label{eq:terminal-one-factor-truncation}
 \sum_{|l|\le2R}|\widehat h_{q_j}(l)|
 \le C(1+\log R)R^{(-q_j)_+},\qquad R\ge2.
\end{equation}
Indeed, the part with $2^r\le1+|l|<2^{r+1}$ is at most
$C2^{-rq_j}$. There are $O(1+\log R)$ such parts, and each is
bounded by $CR^{(-q_j)_+}$. This argument stays uniform when an
exponent passes through zero.

Consider a shell $R\le\max_j|\ell_j|<2R$ subject to
$\sum_j\ell_j=0$. Choose $j$ attaining the maximum. There is
$h\ne j$ with $|\ell_h|\ge R/(m-1)$. There are only finitely many
choices for $(j,h)$. Eliminate $\ell_j$ by the zero sum relation.
The two selected coefficients contribute at most
$CR^{-2-q_j-q_h}$, and there are $O(R)$ values for $\ell_h$.
Sum each remaining coordinate by
\eqref{eq:terminal-one-factor-truncation}. Since the weight on the
shell is at most $CR^s$, its contribution is bounded by
\[
 C(1+\log R)^{m-2}
 R^{s-1-q_j-q_h+\sum_{l\notin\{j,h\}}(-q_l)_+}.
\]
The number
\[
 q_j+q_h-\sum_{l\notin\{j,h\}}(-q_l)_+
\]
is the sum over the subset consisting of $j,h$ and all negative
remaining exponents. That subset has at least two indices, so the
number is at least $s_K-1$. The shell is consequently bounded by
$C(1+\log R)^{m-2}R^{s-s_K}$. Sum this bound over
$R=2,4,8,\ldots$. It is summable for $s<s_K$, and the finitely
many remaining coefficients are uniformly bounded. This proves
\eqref{eq:terminal-weighted-fourier}.
If all exponents are positive, the smallest subset sum with at least
two terms is the sum of the two smallest exponents. This proves the
last assertions.
\end{proof}

By \eqref{eq:terminal-weighted-fourier}, the Fourier series is
absolutely and uniformly convergent.
By \cref{lem:terminal-holder}, $W_{\mathbf q}$ is continuous.
Uniqueness of Fourier coefficients therefore identifies the series
with $W_{\mathbf q}$ at every phase vector.

\subsection{Proof of the characteristic polynomial asymptotic}

\begin{proof}[Proof of \cref{thm:signed-characteristic-FH}]
Put $n=N-1$. By \eqref{eq:terminal-conditioned}, it remains to
estimate the integral of $W_{\mathbf q}$ against $\mu_{n,\mathbf q}$.
We interpolate the decay of each Fourier mode with the bound by one,
so that the resulting power of the mode can be summed using
\cref{lem:terminal-weighted-fourier}.

Choose $0<s<\varepsilon_K$ and $0<\gamma<\min\{4/\beta,1\}$.
For nonzero modes with $\sum_j\ell_j=0$, the normalized estimate
\eqref{eq:relative-phase-mixing} and the probability bound imply
\[
 \left|\int e^{i\boldsymbol\ell\cdot\boldsymbol\varphi}
                                    \,d\mu_{n,\mathbf q}\right|
 \le\min\{1,CL^2(n\delta)^{-\gamma}\},
 \qquad L=1+\|\boldsymbol\ell\|_1.
\]
For $0<s\le2$ and $x\ge0$, $\min(1,x)\le x^{s/2}$.
Consequently this bound is at most
$C L^s(n\delta)^{-\gamma s/2}$.
Integrate the absolutely convergent Fourier series and apply
\cref{lem:terminal-weighted-fourier} to obtain
\begin{equation}\label{eq:terminal-balanced-error}
 \int W_{\mathbf q}\,d\mu_{n,\mathbf q}
 =\prod_jc(q_j)+O((n\delta)^{-\gamma s/2}).
\end{equation}
For any $\kappa<\varepsilon_K\min\{2/\beta,1/2\}$, choose
$s<\varepsilon_K$ and $\gamma<\min\{4/\beta,1\}$ sufficiently
close to their upper bounds that $\gamma s/2>\kappa$.
The error in \eqref{eq:terminal-balanced-error} is then
$O_{\beta,m,K,\kappa}((n\delta)^{-\kappa})$.

By \cref{thm:opuc-relative},
\begin{equation}\label{eq:terminal-OPUC-input}
 \E Y_n(\mathbf q)
 =n^{a(\mathbf q)}
     \prod_j\frac{C_\beta(q_j)}{c(q_j)}
     \prod_{j<l}|e^{i\vartheta_j}-e^{i\vartheta_l}|^{-q_jq_l/\beta}
     \bigl(1+O((n\delta)^{-1})\bigr).
\end{equation}
The constants in this relative estimate are uniform in $\delta$.
Moreover $c(q)\ge1$ for real $q>-1$, by Jensen's inequality and
$\int\log|1-e^{it}|\,d\sigma(t)=0$.
Thus \eqref{eq:terminal-balanced-error} provides a relative estimate for its
positive main term. Multiplication in \eqref{eq:terminal-conditioned}
proves the formula with $n$. Replacing $n=N-1$ by $N$ costs
$1+O_K(N^{-1})$; for $N\delta\ge R$ this is absorbed by the claimed
error. This proves the theorem.
\end{proof}

\begin{corollary}
\label{cor:terminal-subset-rate}
Let $m\ge2$ and let $K\subset\mathcal D_m$ be compact. With $s_K$
from \eqref{eq:terminal-subset-index}, the relative error in
\eqref{eq:signed-characteristic-FH} is $O((N\delta)^{-\kappa})$
for every
\begin{equation}\label{eq:positive-terminal-rate}
 0<\kappa<\min\{4/\beta,1\}\min\{s_K/2,1\}.
\end{equation}
In particular, for $K\subset[q_*,\infty)^m$ with $q_*>0$, every
$\kappa<\min\{4/\beta,1\}\min\{(1+2q_*)/2,1\}$ is allowed.
If $s_K>2$ (in particular if $q_*>1/2$), it is $O(\mathfrak r_{4/\beta}(N\delta))$, where
\begin{equation}\label{eq:positive-terminal-endpoint-rate}
 \mathfrak r_{4/\beta}(t)=
 \begin{cases}
 t^{-4/\beta},&\beta>4,\\
 (1+\log t)/t,&\beta=4,\\
 t^{-1},&0<\beta<4.
 \end{cases}
\end{equation}
\end{corollary}

\begin{proof}[Proof of \cref{cor:terminal-subset-rate}]
In the proof of \cref{thm:signed-characteristic-FH}, the Fourier
summation permits any $0<s<\min\{s_K,2\}$ by
\cref{lem:terminal-weighted-fourier}. With
$0<\gamma<\min\{4/\beta,1\}$, it gives error in the integral of $W_{\mathbf q}$
$O((n\delta)^{-\gamma s/2})$. Since $\prod_jc(q_j)$ is uniformly
positive, combining this with \cref{thm:opuc-relative} yields every exponent in
\eqref{eq:positive-terminal-rate} by letting $s$ and $\gamma$
approach their upper bounds.
If $s_K>2$, use $s=2$ in
\cref{lem:terminal-weighted-fourier} and apply the endpoint mode
bound \eqref{eq:relative-phase-mixing} directly. Its summed error
is $O(\mathfrak r_{4/\beta}(n\delta))$. The OPUC error from
\cref{thm:opuc-relative} and the replacement of $n=N-1$ by $N$
are absorbed by this bound.
\end{proof}

\begin{remark}
For two points, $s_K=\inf_K(1+q_1+q_2)$. In particular,
$\mathbf q=(a,-a)$ with $a$ in a compact subset of $(-1,1)$
gives every $\kappa<\min\{2/\beta,1/2\}$. This exponent does not
deteriorate with $|a|$, although the implicit constant may deteriorate
as $|a|$ approaches one. For positive exponents, the sum of the two
smallest exponents gives the relevant Fourier threshold.
\end{remark}

\begin{remark}
When all exponents are fixed nonnegative even integers, each $h_{q_j}$
is a trigonometric polynomial. The function $W_{\mathbf q}$ then has finitely
many Fourier modes, and the endpoint error
$O(\mathfrak r_{4/\beta}(N\delta))$ follows directly, including when
some exponents vanish.
\end{remark}

\begin{remark}
The domain $\mathcal D_m$ ensures that $W_{\mathbf q}$ is bounded
and has the weighted Fourier summability needed for a constant uniform
in the separation. It is not the integrability domain of the original
separated moments. In the following section we remove this restriction
at fixed separation by controlling the weighted phase law instead.
\end{remark}

\section{The full integrable range of root exponents}
\label{sec:full-negative-transfer}

The integral $W_{\mathbf q}$ is still in $L^1(\sigma^{\otimes m})$
for every $q_j>-1$, but it need not be bounded. Consequently the
decay of nonconstant Fourier modes in
\cref{prop:quantitative-phase-mixing} does not by itself imply
\[
 \int W_{\mathbf q}\,d\mu_{n,\mathbf q}
       \longrightarrow \int W_{\mathbf q}\,d\sigma^{\otimes m}
       =\prod_jc(q_j).
\]
For fixed distinct points, the density bound in
\cref{prop:phase-law-decomposition} controls the error when
$W_{\mathbf q}$ is approximated in $L^1$ by
trigonometric polynomials. This reduction comes first, followed by the
density bound. Its explicit dependence on the separation also yields
a range of polynomial approaches of the points.
On $\mathcal E_m$, the final subsection bounds $W_{\mathbf q}$ by the
sum of singular pair powers in \eqref{eq:terminal-pair-majorant}.
The pair density bound in \cref{prop:moving-pair-density} controls
the error when each singular factor is truncated near its zero.
Fourier approximation then gives a relative error bounded by a
negative power of $N\delta$ whenever $N\delta$ is sufficiently large.

Unless otherwise stated, $m\ge2$, $K\subset(-1,\infty)^m$ is compact,
and $0<\delta\le1$. Constants may depend on $\beta,m,K$; dependence
on a fixed separation is indicated when it is used. Write
\[
 A_n(\mathbf q)=\E Y_n(\mathbf q),\qquad
 \mu_{n,\mathbf q}(B)
 =\frac{\E[Y_n(\mathbf q)\mathbf1_{\{\boldsymbol\psi_n\in B\}}]}
        {A_n(\mathbf q)}.
\]
The phase variables in these measures are taken modulo $2\pi$.
For $m=1$, \cref{prop:onepoint} already gives the characteristic
asymptotic with relative error $O(N^{-1})$; no separation condition or
phase density argument is needed.

The Vandermonde product will measure the dependence of the full density
bound on the prescribed points:
\[
 \mathcal V(\boldsymbol\vartheta)
   =\prod_{j<l}|e^{i\vartheta_j}-e^{i\vartheta_l}|.
\]
We shall balance a translation error against a truncated Fourier sum.
For a compact set $K\subset(-1,\infty)^m$, the parameters for this
balance are
\begin{align*}
 A_K&=\sup_{\mathbf q\in K}\sum_j(-q_j)_+,
 \qquad M_K=2+(A_K-1)_+,\\
 \eta_K&=\min\{1,\inf_{\mathbf q\in K,j}(1+q_j)\}>0,
 \qquad g_\beta=\min\{4/\beta,1\}.
\end{align*}
Here $A_K$ bounds the total negative exponent, $\eta_K$ bounds the
available translation regularity, and $g_\beta$ is the phase decay
threshold. The power $M_K$ combines
\eqref{eq:terminal-truncated-mass} with the quadratic dependence on the mode in \eqref{eq:relative-phase-mixing}.
\begin{theorem}
\label{thm:full-root-FH}
Fix $\beta>0$, $m\ge1$, a compact $K\subset(-1,\infty)^m$, and
$0<\delta\le1$. With the notation above, for every
\begin{equation}\label{eq:full-root-rate}
 0<\kappa<\frac{\eta_K g_\beta}{M_K+\eta_K},
\end{equation}
the Fisher--Hartwig formula
\begin{align}\label{eq:full-root-FH}
 F_N(\mathbf q,\boldsymbol\vartheta)
 ={}&N^{a(\mathbf q)}\prod_jC_\beta(q_j)
       \prod_{j<l}|e^{i\vartheta_j}-e^{i\vartheta_l}|^{-q_jq_l/\beta}
       \bigl(1+O_{\beta,m,K,\delta,\kappa}(N^{-\kappa})\bigr)
\end{align}
holds uniformly on $K\times\mathcal S_\delta$.

If $K\subset\mathcal D_m$, the same conclusion holds for every
\[
 0<\kappa<g_\beta\max\left\{
          \frac{\eta_K}{M_K+\eta_K},\frac{1-A_K}{2}\right\},
\]
by combining this bound with \cref{thm:signed-characteristic-FH} at fixed separation.
The subset sum parameter $s_K$ in \cref{cor:terminal-subset-rate}
can give a larger admissible rate than either bound displayed here. 
For $m=1$, \cref{prop:onepoint} gives the stronger error $O(N^{-1})$.
\end{theorem}

The main additional estimate is a bound for the phase law itself. The
power of the separation in this bound is harmless when the points are
fixed.

\begin{proposition}
\label{prop:phase-law-decomposition}
There is $R_K=R_K(\beta,m,K)<\infty$ such that, if
$n\delta\ge R_K$, the normalized weighted phase law is absolutely
continuous and satisfies
\begin{equation}\label{eq:phase-law-decomposition}
 \mu_{n,\mathbf q}=f_{n,\mathbf q}\sigma^{\otimes m},\qquad
 \|f_{n,\mathbf q}\|_\infty
 \le\frac{C_K\delta^{-m/2}}{\mathcal V(\boldsymbol\vartheta)}
 \le C_K\delta^{-B_m},\qquad B_m=\frac{m^2}{2}.
\end{equation}
In particular the densities are uniformly bounded on every fixed
separated angle set. The constants are uniform for $\mathbf q\in K$.
\end{proposition}

Assuming \cref{prop:phase-law-decomposition}, we can prove
\cref{thm:full-root-FH}. Its
proof also uses the Fourier estimate of \cref{lem:terminal-truncated-mass},
proved below. The density estimate will then be deduced from a finite
sequence of transitions and an adjoint identity.

\begin{proof}[Proof of \cref{thm:full-root-FH}]
We approximate $W_{\mathbf q}$ by a Fej\'er polynomial. The density
bound controls the approximation error, and the phase estimate controls
its nonconstant modes. We then balance the two errors.
Suppose $m\ge2$. Choose $0<\alpha<\eta_K$,
$0<\gamma<g_\beta$ and $\zeta>0$ such that
$\kappa<\alpha\gamma/(M_K+\zeta+\alpha)$.
We use the two estimates in the proof of
\cref{lem:terminal-translation} with $p=1$. That lemma assumes
$q+1\ge2\alpha$ in this case; the same split only requires
$q+1>\alpha$, with a uniform positive gap, to give
\[
 \|h_q(\cdot+t)-h_q\|_1\le C\|t\|_{2\pi}^{\alpha}
\]
uniformly for the coordinates of $\mathbf q\in K$. Indeed the integral
near zero is $O(|t|^{q+1})$, while away from zero it is at most
\[
 C|t|\int_{|t|}^{\pi}r^{q-1}\,dr
 \le C|t|^\alpha\int_{|t|}^{\pi}r^{q-\alpha}\,dr
 \le C|t|^\alpha,
\]
because $q+1-\alpha$ is bounded below by a positive number.
Telescoping the product defining $W_{\mathbf q}$ and integrating all
phase coordinates gives
\begin{equation}\label{eq:terminal-Lone-translation}
 \|W_{\mathbf q}(\cdot+\mathbf t)-W_{\mathbf q}\|_1
 \le C\sum_j\|t_j\|_{2\pi}^{\alpha}.
\end{equation}
Here Fubini's theorem applies since
$\|W_{\mathbf q}\|_1=\prod_jc(q_j)<\infty$.

Let $\mathcal K_L$ be the Fej\'er kernel. We approximate $W_{\mathbf q}$
by the trigonometric polynomial
$W_{\mathbf q,L}=W_{\mathbf q}*\mathcal K_L^{\otimes m}$.
The kernel bound
\[
 \mathcal K_L(t)\le C\min\{L+1,((L+1)\|t\|_{2\pi}^2)^{-1}\}
\]
gives $\int_\T\|t\|_{2\pi}^{\alpha}\mathcal K_L(t)\,d\sigma(t)
\le CL^{-\alpha}$: split the integral at $(L+1)^{-1}$ and use
$0<\alpha<1$. Combining this with
\eqref{eq:terminal-Lone-translation} yields
$\|W_{\mathbf q,L}-W_{\mathbf q}\|_1\le CL^{-\alpha}$.
The full density bound in \cref{prop:phase-law-decomposition} therefore
gives, at fixed separation,
\[
 \left|\int(W_{\mathbf q}-W_{\mathbf q,L})\,d\mu_{n,\mathbf q}\right|
 \le CL^{-\alpha}.
\]
Only modes with $|\ell_j|\le L$ and $\sum_j\ell_j=0$ occur in
$W_{\mathbf q,L}$: the restriction on the sum follows from invariance
of $W_{\mathbf q}$ under a common translation of all phases.
Its Fej\'er multipliers have modulus at most one.
By \eqref{eq:relative-phase-mixing} and
\cref{lem:terminal-truncated-mass}, the contribution of its nonconstant
modes is at most $Cn^{-\gamma}L^{M_K+\zeta}$. Here the factor $L^2$
in $L^{M_K+\zeta}$ comes from
$(1+\|\boldsymbol\ell\|_1)^2$ in
\eqref{eq:relative-phase-mixing}; the remaining power comes from
\eqref{eq:terminal-truncated-mass}.
Its constant coefficient is $\prod_jc(q_j)$. Consequently
\[
 \left|\int W_{\mathbf q}\,d\mu_{n,\mathbf q}-\prod_jc(q_j)\right|
 \le C\{L^{-\alpha}+n^{-\gamma}L^{M_K+\zeta}\}.
\]
Choose $L=\max(1,\lfloor n^{\gamma/(M_K+\zeta+\alpha)}\rfloor)$.
The right-hand side is $O(n^{-\kappa})$.

Finally combine the terminal identity \eqref{eq:terminal-conditioned},
namely $F_{n+1}=A_n\int W_{\mathbf q}\,d\mu_{n,\mathbf q}$,
with \cref{thm:opuc-relative} and the identity
$D_\beta(q)c(q)=C_\beta(q)$.
All one point factors are uniformly positive on compact sets of real
exponents, so the error in this integral is relative. Replacing $n$ by $N=n+1$
costs $1+O_K(N^{-1})$, proving \eqref{eq:full-root-FH}.
\end{proof}

\subsection{The adjoint transition and propagation of densities}

For $\alpha\in\mathbb D$, let
\[
 g_\alpha(x)=x-2\Im\log(1-\alpha e^{ix}),\qquad
 w_{\mathbf q}(\alpha,\mathbf x)
       =\prod_j|1-\alpha e^{ix_j}|^{q_j}.
\]
The map $g_\alpha$ is an orientation preserving circle diffeomorphism,
and acts separately on the coordinates of $\mathbf x$.
For bounded measurable $f$, define the positive weighted transition
\[
 (\mathcal P_k f)(\mathbf x)
 =\E\left[w_{\mathbf q}(\alpha_k,\mathbf x)
    f(\boldsymbol\vartheta+g_{\alpha_k}(\mathbf x))\right].
\]
Its adjoint $\mathcal P_k^*$ acts on densities relative to
$\sigma^{\otimes m}$. In this subsection, the radial parameter, its
shift under the adjoint, and the dual exponents are denoted by
\[
 t_k=b(k+1),\qquad S=\sum_jq_j+m,
 \qquad q'_j=-q_j-2.
\]
Write $\E_t$ for expectation under the rotational disk law with planar
density $\frac t\pi(1-|z|^2)^{t-1}\mathbf1_{\{|z|<1\}}$.

\begin{lemma}
\label{lem:adjoint-density-propagation}
For $k$ sufficiently large and $\mathbf y=\boldsymbol\varphi-
\boldsymbol\vartheta$, one has the exact identity
\begin{equation}\label{eq:adjoint-dual-identity}
 (\mathcal P_k^* f)(\boldsymbol\varphi)
 =\frac{t_k}{t_k+S}\,
   \E_{t_k+S}\left[w_{\mathbf q'}(\zeta,\mathbf y)
                         f(g_\zeta(\mathbf y))\right].
\end{equation}
Moreover, there are $R_K,C_K$ such that, for every $n\ge k$ with
$k\delta\ge R_K$ and bounded nonnegative $f$,
\begin{equation}\label{eq:adjoint-density-propagation}
 \|\mathcal P_{n-1}^*\cdots\mathcal P_k^* f\|_\infty
 \le C(n/k)^{a(\mathbf q)}\|f\|_\infty,
\end{equation}
uniformly for $\mathbf q\in K$ and
$\boldsymbol\vartheta\in\mathcal S_\delta$, with constants independent
of $\delta$.
\end{lemma}

\begin{proof}[Proof of \cref{lem:adjoint-density-propagation}]
A change of variables gives the exact adjoint identity. We then adapt
the Fourier correction to the adjoint and construct a positive
function that bounds the growth of propagated densities.
The inverse of $g_\alpha$ is $g_{-\alpha}$. Direct calculation gives
\[
 |1-\alpha e^{ig_{-\alpha}(y)}|
      =\frac{1-|\alpha|^2}{|1+\alpha e^{iy}|},
 \qquad
 g_{-\alpha}'(y)
      =\frac{1-|\alpha|^2}{|1+\alpha e^{iy}|^2}.
\]
Changing variables $\mathbf x=g_{-\alpha}(\mathbf y)$ in the adjoint
integral therefore multiplies the disk density by
\[
 (1-|\alpha|^2)^S
       \prod_j|1+\alpha e^{iy_j}|^{-q_j-2}.
\]
Absorb the first factor into the radial density and then put
$\zeta=-\alpha$. This proves \eqref{eq:adjoint-dual-identity} whenever
$t_k+S>0$.

The dual exponents $q'_j$ need not exceed $-1$. Only conditional disk
integrals at large indices are being used here, so the annulus estimate
\eqref{eq:signed-annulus} still supplies the required integrability.
Define
\[
 D'(\mathbf y)
 =\left(\sum_j\frac{q'_j}{2}e^{iy_j}\right)
  \left(\sum_j\frac{q'_j}{2}e^{-iy_j}\right)
 =D'_0+D'_{\rm off}(\mathbf y).
\]
The diagonal term satisfies
\begin{equation}\label{eq:adjoint-diagonal-cancellation}
 D'_0-S
 =\sum_j\frac{(q_j+2)^2}{4}-\sum_j(q_j+1)
 =\sum_j\frac{q_j^2}{4}=:D_0.
\end{equation}
Since $t_k/(t_k+S)=1-Sv_k+O(k^{-2})$ and
$(1+t_k+S)^{-1}=v_k+O(k^{-2})$, the expansion in
\cref{lem:onestep} gives
\begin{equation}\label{eq:adjoint-one-step}
 \mathcal P_k^*1(\boldsymbol\varphi)
 =1+v_k\{D_0+D'_{\rm off}(\boldsymbol\varphi-
          \boldsymbol\vartheta)\}+O(k^{-2}).
\end{equation}
The estimates are uniform on the specified compact sets.

To cancel the oscillatory part of \eqref{eq:adjoint-one-step}, introduce
the real trigonometric polynomial
\[
 R(\mathbf x)
 =\sum_{j\ne l}\frac{q'_jq'_l}{4}
       \frac{e^{i(x_j-x_l)}}{e^{i(\vartheta_j-\vartheta_l)}-1}.
\]
Terms indexed by $(j,l)$ and $(l,j)$ are conjugate. Thus $R$ is real,
$\|R\|_{C^2}\le C/\delta$, and
\begin{equation}\label{eq:adjoint-corrector}
 R(\mathbf x+\boldsymbol\vartheta)-R(\mathbf x)
       =D'_{\rm off}(\mathbf x).
\end{equation}
Taylor expansion in the phase increment, using the signed mixed moment
\eqref{eq:mixed-one-step} at zero Fourier mode, gives
\begin{equation}\label{eq:adjoint-R-transition}
 \mathcal P_k^*R(\boldsymbol\varphi)
 =R(\boldsymbol\varphi-\boldsymbol\vartheta)
                              +O((k\delta)^{-1}).
\end{equation}
Here the constant Taylor term uses
$\E w_{\mathbf q'}=1+O(k^{-1})$; the linear term uses
$\E[w_{\mathbf q'}\boldsymbol\xi]=O(k^{-1})$; and the quadratic
remainder uses $\E[w_{\mathbf q'}|\boldsymbol\xi|^2]=O(k^{-1})$.
The last bound is \eqref{eq:weighted-phase-second}, whose proof splits
at $|\zeta|=1/2$. The shifted radial parameter $t_k+S$ is comparable
to $k$, with $S$ in a fixed compact interval, so the same proof and
\eqref{eq:signed-annulus} give uniform constants here.

The corrected constant function is $r_k=1+v_kR$. Since $\|R\|_\infty\le C_K/\delta$,
$1/2\le r_k\le3/2$ whenever $k\delta\ge R_K$, after increasing $R_K$.
Equations \eqref{eq:adjoint-one-step} to~\eqref{eq:adjoint-R-transition}
and $v_{k+1}-v_k=O(k^{-2})$ imply
\[
 \mathcal P_k^*r_k=(1+v_kD_0)r_{k+1}+O((k^2\delta)^{-1}).
\]
Positivity consequently gives
\[
 \mathcal P_k^*r_k
 \le\left(1+v_kD_0+\frac C{k^2\delta}\right)r_{k+1}.
\]
Iterating this inequality and using $f\le2\|f\|_\infty r_k$ yields
\[
 \mathcal P_{n-1}^*\cdots\mathcal P_k^*f
 \le2\|f\|_\infty r_n
       \prod_{j=k}^{n-1}\left(1+v_jD_0+\frac C{j^2\delta}\right).
\]
Finally, $D_0/b=a(\mathbf q)$,
$\sum_{j=k}^{n-1}bv_j=\log(n/k)+O(k^{-1})$, and
$\sum_{j\ge k}(j^2\delta)^{-1}\le C/(k\delta)$.
The product is therefore bounded by $C(n/k)^{a(\mathbf q)}$, proving
\eqref{eq:adjoint-density-propagation}.
\end{proof}

\subsection{A bounded density after finitely many coefficients}
\label{sec:unconditional-phase-smoothing}

The phase map is rational in its coefficients. The proof of
\cref{lem:polynomial-inverse-gaussian} bounds the probability that a
polynomial is small and deduces an inverse moment estimate. This yields
an $L^p$ bound with $p>1$ for the entire transition density.
Applying the same argument to the adjoint
and composing a fixed number of blocks then gives a bounded density.
Throughout this subsection the exponents range over a fixed compact
set $K_0\subset\mathbb R^m$. The conditional transition estimates themselves allow
arbitrary real root exponents: the starting index is chosen sufficiently
large that all the coefficient integrals exist.

For the finite number of transitions used below, put
\[
 \mathfrak d_m=2m^2,\qquad J_m=\mathfrak d_m+2,
 \qquad r_m=mJ_m.
\]
The integer $r_m$ is a convenient choice, not an optimal number of steps.
\begin{theorem}
\label{thm:unconditional-phase-kernel}
With $r_m=m(2m^2+2)$, the entire root weighted transition over
$k,\ldots,k+r_m-1$ has a density bounded by
\begin{equation}\label{eq:unconditional-phase-kernel}
 C\frac{k^{m/2}}{\mathcal V(\boldsymbol\vartheta)},
 \qquad k\ge k_{K_0},
\end{equation}
uniformly in the incoming and outgoing phase vectors. Here $K_0$ may be any compact subset of $\mathbb R^m$.
\end{theorem}

We prove this by first obtaining an $L^p$ bound with $p>1$ for both rows
and columns of the transition density over $m$ coefficients. A fixed
number of compositions then gives an $L^\infty$ bound. The polynomial
estimate in \cref{lem:polynomial-inverse-gaussian} supplies the negative
Jacobian moment in that argument.

\begin{lemma}
\label{lem:polynomial-inverse-gaussian}
Fix $d,D\ge1$, $c>0$ and $0<\epsilon<1/D$.
For every real polynomial $P$ on $\mathbb R^d$ of degree at most $D$
with $P(0)\ne0$, and every $k\ge1$,
\begin{equation}\label{eq:polynomial-inverse-gaussian}
 \int_{\mathbb R^d}e^{-ck|z|^2}|P(z)|^{-\epsilon}\,dz
 \le C k^{-d/2}|P(0)|^{-\epsilon}.
\end{equation}
The constant depends only on $d,D,c,\epsilon$.
\end{lemma}

\begin{proof}[Proof of \cref{lem:polynomial-inverse-gaussian}]
Let $Z$ have a centered nondegenerate Gaussian law with density
proportional to $e^{-c|z|^2}$. The Carbery--Wright inequality
\cite[Theorem~8]{CarberyWright2001} gives
\[
 \Prob\{|Q(Z)|\le u\}
 \le C_D\left(\frac{u}{\E|Q(Z)|}\right)^{1/D}
\]
for a nonzero real polynomial $Q$ of degree at most $D$; the bound may
always be replaced by one if its right-hand side is larger. On the
finite dimensional space of such polynomials, the Gaussian $L^1$ norm
is a norm and evaluation at zero is a continuous linear functional.
Consequently $\E|Q(Z)|\ge c_{d,D,c}|Q(0)|$.
For $Q(0)\ne0$, put $X=|Q(Z)|/|Q(0)|$. The two bounds imply
$\Prob(X<u)\le\min\{1,Cu^{1/D}\}$. The layer cake identity gives
\[
 \E X^{-\epsilon}
 \le 1+\epsilon\int_0^1u^{-\epsilon-1}\Prob(X<u)\,du
 \le 1+\frac{C\epsilon}{1/D-\epsilon}.
\]
Thus $\E|Q(Z)|^{-\epsilon}\le C|Q(0)|^{-\epsilon}$, with a
constant independent of the coefficients of $Q$. Apply this to
$Q(z)=P(z/\sqrt{k})$ and change variables.
\end{proof}

\begin{lemma}
\label{lem:whole-block-Lp}
Let $K_0\subset\mathbb R^m$ be a compact set of real root exponents. Let $K_k(\mathbf x,\mathbf y)$ be the kernel,
relative to Haar measure in the outgoing phase, of the root weighted
transition over $k,\ldots,k+m-1$. There is $k_{K_0}$ such that this kernel
exists for $k\ge k_{K_0}$. If
\[
 1<p<1+1/\mathfrak d_m,\qquad a_p=1-1/p,
\]
then
\begin{equation}\label{eq:whole-block-Lp}
 \sup_{\mathbf x}\|K_k(\mathbf x,\cdot)\|_p
 +\sup_{\mathbf y}\|K_k(\cdot,\mathbf y)\|_p
 \le C_p\left(\frac{k^{m/2}}{\mathcal V(\boldsymbol\vartheta)}
                                                    \right)^{a_p}.
\end{equation}
The transition is taken under the full coefficient law. The constants
are independent of the distinct angles and their separation.
\end{lemma}

\begin{proof}[Proof of \cref{lem:whole-block-Lp}]
We first bound the rows. We select a nonzero Jacobian minor at zero
coefficients, control its inverse moment by
\cref{lem:polynomial-inverse-gaussian}, and bound the number of regular
preimages by B\'ezout's theorem. The area formula then gives the
$L^p$ estimate. Applying the same argument to the adjoint will bound
the columns.
Fix an incoming vector $\mathbf x$ and denote
the $2m$ real coefficient coordinates by $z$.
If $w_j=e^{ix_j}$, one step of the phase map on the circle is
\[
 w_j\longmapsto e^{i\vartheta_j}
                  \frac{w_j-\bar\alpha}{1-\alpha w_j}.
\]
After $m$ steps its $j$th coordinate has the form $P_j(z)/Q_j(z)$,
where $P_j,Q_j$ are complex polynomials in the real coefficient
coordinates, of degree at most $m$. One obtains them recursively from
\[
 P\longmapsto e^{i\vartheta_j}(P-\bar\alpha Q),
 \qquad Q\longmapsto Q-\alpha P,
 \qquad (P,Q)=(e^{ix_j},1)\text{ initially}.
\]
Inside the coefficient polydisk, $Q_j$ never vanishes,
$|P_j|=|Q_j|$, and $|Q_j|\le2^m$. The equality of the moduli
follows inductively because each displayed M\"obius map preserves the
unit circle. At zero coefficients $|Q_j|=1$.
For any local real lift $F_j$ of the outgoing angle,
\[
 \partial F_j
 =\frac{\Im\{(\partial P_j)\overline{P_j}
                  -(\partial Q_j)\overline{Q_j}\}}{|Q_j|^2}.
\]
The numerator is a real polynomial of degree at most $2m-1$.

At zero coefficients the real differential is
$2\Im[\operatorname{diag}(e^{ix_j})V]$, where
$V_{jh}=e^{ih\vartheta_j}$, $0\le h<m$.

Indeed, perturbing the coefficient at step $h$, where $0\le h<m$,
changes the final $j$th phase to first order by
\[
 2\Im\{\delta\alpha_h e^{i(x_j+h\vartheta_j)}\}.
\]
The phase entering that step is $x_j+h\vartheta_j$, and every subsequent
step with zero coefficient has phase derivative one.
For $F=(F_1,\ldots,F_m)$, we claim that the real Gram matrix at zero satisfies
\[
 \det\bigl(DF(0)DF(0)^{\mathsf T}\bigr)\ge4^m|\det V|^2,
\]
uniformly in the incoming phases. To prove this, put
$H=\operatorname{diag}(e^{ix_j})VV^*
\operatorname{diag}(e^{-ix_j})$. The real Gram matrix is $4\Re H$.
For any positive Hermitian matrix $H=A+iB$, where $A$ is real symmetric
and $B$ real skew symmetric, the eigenvalues of
$A^{-1/2}BA^{-1/2}$ occur in pairs $\pm it_l$. Positivity of $H$
implies $|t_l|<1$, and hence
\[
 \det H=\det A\prod_l(1-t_l^2)\le\det A.
\]
It follows that the determinant of the real Gram matrix is at least
$4^m\det H=4^m|\det V|^2$.
The Cauchy--Binet formula therefore selects $m$ real coordinate columns, denoted
by $u$, for which the minor at zero has absolute value at least
$c_m\mathcal V(\boldsymbol\vartheta)$. Denote the unused coordinates
by $v$. The selected Jacobian throughout the polydisk is
\begin{equation}\label{eq:rational-phase-minor}
 \det D_uF(z)=\frac{P(z)}{\prod_j|Q_j(z)|^2},\qquad
 \deg P\le m(2m-1)\le \mathfrak d_m,
 \quad |P(0)|\ge c_m\mathcal V(\boldsymbol\vartheta).
\end{equation}
The choice of columns may depend on $\mathbf x,\boldsymbol\vartheta$;
there are only finitely many choices and all estimates are uniform.

The polynomial $P$ is nonzero. Its zero set has Lebesgue measure zero.
On its complement, change of variables in $u$ is available locally.
There is also a uniform bound on its multiplicity. Fix $v$ and
$\mathbf y$. A preimage is a solution in the $m$ real variables $u$
of the $m$ polynomial equations
\[
 \Im\{e^{-iy_j}P_j(u,v)\overline{Q_j(u,v)}\}=0,
 \qquad j=1,\ldots,m,
\]
each of degree at most $2m$. At an actual preimage with nonzero selected
Jacobian, the Jacobian of these polynomial equations is nonzero too:
its rows are $|Q_j|^2$ times those of $D_uF$.
Such a solution is an isolated nonsingular complex solution after
complexification. The elementary B\'ezout bound consequently gives at
most $(2m)^m$ such solutions. Additional solutions representing the
antipodal output can only enlarge this upper bound. Critical coefficient
points have zero measure and make no contribution. The area formula on
the regular set now supplies the whole outgoing density.

Let $\rho_k(u,v)$ be the density in coefficient space, including the
$m$ incremental root weights. Its pushforward under the phase map is
the transition measure we wish to estimate. Uniformly on the coefficient
polydisk, for $k\ge k_{K_0}$,
\begin{equation}\label{eq:whole-weighted-Gaussian}
 \rho_k(u,v)\le C k^m e^{-c_0 k(|u|^2+|v|^2)}.
\end{equation}
To verify this even at the coefficient boundary, the negative root
powers at one step are bounded by
$C_{K_0}(1-|\alpha|^2)^{-A_{K_0}}$, with
$A_{K_0}=\sup_{\mathbf q\in K_0}\sum_j(-q_j)_+$. The disk density has exponent $b(k+h+1)-1$. For large $k$
it absorbs $A_{K_0}$ and leaves an exponent at least $c k$, giving
\eqref{eq:whole-weighted-Gaussian}.

To integrate over the unused coordinates $v$ by Jensen's inequality,
choose $0<c_1\le c_0$ and use the normalized Gaussian density
$\gamma_k(v)=C_{c_1}k^{m/2}e^{-c_1k|v|^2}$.
Then the product $\rho_k^p\gamma_k^{1-p}$ has Gaussian exponent
$-[pc_0-(p-1)c_1]k|v|^2$ in the unused coordinates, with
$pc_0-(p-1)c_1\ge c_0>0$. Thus the same Gaussian domination remains
available after Jensen's inequality.
Extend $\rho_k$ by zero outside the coefficient polydisk. For almost
every $v$, the area formula gives the Euclidean angle density
\[
 H_v(\mathbf y)=
 \sum_{u:F(u,v)=\mathbf y}
          \frac{\rho_k(u,v)}{|\det D_uF(u,v)|},
\]
where the sum runs over regular actual preimages. There are at most
$(2m)^m$ terms. The Euclidean version of $K_k(\mathbf x,\cdot)$
is $\int H_v\,dv$. Jensen's inequality in $\gamma_k(v)\,dv$, followed
by the finite sum inequality and the area formula, therefore gives
\begin{align*}
 \int |K_k(\mathbf x,\mathbf y)|^p\,d\sigma^m(\mathbf y)
 &\le C_p\int\gamma_k(v)^{1-p}
                          \int H_v(\mathbf y)^p\,d\mathbf y\,dv\\
 &\le C_p\int \rho_k(u,v)^p\gamma_k(v)^{1-p}
                 |\det D_uF(u,v)|^{1-p}\,du\,dv.
\end{align*}
Fubini's theorem discards the exceptional $v$, since the full
polynomial critical set has measure zero. All changes between normalized Haar
measure and Euclidean angle measure have been absorbed into $C_p$.
By \eqref{eq:rational-phase-minor}, the negative power of the
Jacobian is at most $C|P(u,v)|^{-(p-1)}$.
Since $\gamma_k$ uses precisely the unused $m$ coordinates,
\[
 \rho_k^p\gamma_k^{1-p}
 \le C_p k^{mp-m(p-1)/2}e^{-c'k(|u|^2+|v|^2)}.
\]
Apply \cref{lem:polynomial-inverse-gaussian} in dimension $2m$ with
$\epsilon=p-1<1/\mathfrak d_m$. The resulting bound is
\[
 \int |K_k(\mathbf x,\mathbf y)|^p\,d\sigma^m(\mathbf y)
 \le C_p k^{m(p-1)/2}
           \mathcal V(\boldsymbol\vartheta)^{-(p-1)}.
\]
Taking the $p$th root proves the row estimate.

For the column estimate use the exact adjoint identity
\eqref{eq:adjoint-dual-identity}. It replaces
$\mathbf q$ by $-\mathbf q-2\mathbf1$, shifts each radial parameter by
$S=\sum_jq_j+m$, and multiplies the transition by $t/(t+S)$.
Over a fixed block these scalar factors are bounded above and below.
The reversed coefficient parameters remain comparable to $k$.
The adjoint map is $\mathbf y\mapsto g_\zeta(\mathbf y-
\boldsymbol\vartheta)$. Conjugating by translation makes it the same
phase recursion with drift $-\boldsymbol\vartheta$; its Vandermonde
has the same absolute determinant. Thus the row argument applies to
the adjoint with these dual real exponents, proving the column estimate.
\end{proof}

\begin{proof}[Proof of \cref{thm:unconditional-phase-kernel}]
Each block of $m$ transitions improves the reciprocal integrability
exponent by $1/J_m$. We compose $J_m$ such blocks to reach
$L^\infty$ from $L^1$.
Set $J=J_m$ and $p=J/(J-1)$, so that $p-1<1/\mathfrak d_m$ and
$1-1/p=1/J$. Regard each kernel over $m$ steps as an operator on densities.
By \cref{lem:whole-block-Lp}, its row estimate gives $L^1\to L^p$
and its column estimate gives $L^{p'}\to L^\infty$, both with norm at most
\[
 C\left(k^{m/2}/\mathcal V(\boldsymbol\vartheta)\right)^{1/J}.
\]
Riesz--Thorin interpolation gives $L^r\to L^s$ with
$1/s=1/r-1/J$ whenever $1\le r\le J$.
Compose $J$ consecutive blocks, using successively reciprocal
exponents $1,1-1/J,\ldots,1/J,0$. Their indices are comparable to
$k$, since the total block length depends only on $m$.
Multiplication of the $J$ operator bounds gives
$Ck^{m/2}/\mathcal V$.
The first block acts on arbitrary finite measures by its uniform
row bound, so the same conclusion applies in particular to every
incoming point mass. This is \eqref{eq:unconditional-phase-kernel}.
\end{proof}

\begin{proof}[Proof of \cref{prop:phase-law-decomposition}]
Choose $k=\lceil C_K/\delta\rceil$, with the constant large enough
for \cref{thm:unconditional-phase-kernel}, the moment comparison in
\cref{lem:relative-moment-comparison}, and the density propagation
estimate in \cref{lem:adjoint-density-propagation}.
Increase $R_K$ so that $n\delta\ge R_K$ implies $n\ge k+r_m$.
The weighted incoming measure has mass $A_k=\E Y_k$.
Applying \cref{thm:unconditional-phase-kernel} gives
an unnormalized density at time $k+r_m$ bounded by
$Ck^{m/2}A_k/\mathcal V$.
The estimate \eqref{eq:adjoint-density-propagation} propagates this
density to time $n$ with factor $C(n/k)^{a(\mathbf q)}$.
Its constant is independent of separation because $k\delta$ is bounded
below. The comparison \eqref{eq:relative-moment-comparison} gives
$A_k/A_n\le C(k/n)^{a(\mathbf q)}$. Division by $A_n$
cancels these powers, proving the first bound in
\eqref{eq:phase-law-decomposition}. The second follows from
$\mathcal V\ge\delta^{m(m-1)/2}$ and $k\le C_K/\delta$.
\end{proof}

\subsection{Fourier summation and dependence on the separation}
\label{sec:terminal-sharper-rate}

We now prove \cref{lem:terminal-truncated-mass}, the Fourier mass bound
used in the Fej\'er approximation.
It remains valid when $W_{\mathbf q}$ is unbounded and keeps track of
the total negative exponent through $A_K$.

\begin{lemma}
\label{lem:terminal-truncated-mass}
For every $\zeta>0$ there is $C<\infty$ such that
\begin{equation}\label{eq:terminal-truncated-mass}
 \sum_{\substack{\sum_j\ell_j=0\\\max_j|\ell_j|\le L}}
       \prod_j|\widehat h_{q_j}(\ell_j)|
 \le C(1+L)^{(A_K-1)_++\zeta},
 \qquad \mathbf q\in K,\quad L\ge1.
\end{equation}
The same assertion holds for the mixed factors $h_{q_j,s_j}$ when
$(\mathbf q,\mathbf s)$ ranges over a compact subset $K$ of
$\{(\mathbf q,\mathbf s)\in\C^{2m}:\Re q_j>-1\text{ for all }j\}$,
with $A_K=\sup_K\sum_j(-\Re q_j)_+$.
\end{lemma}

\begin{proof}[Proof of \cref{lem:terminal-truncated-mass}]
We dominate the absolute Fourier coefficients by positive coefficients
of slightly more singular powers. Poisson regularization then expresses
the truncated sum as a one variable integral.
Choose $t>0$ so small that $r_j=(-q_j)_++t<1$ throughout $K$ and
$mt<\zeta/2$. The numbers $r_j$ belong to a compact subset of $(0,1)$.
The gamma formula in \cref{lem:root-fourier-coefficients} gives
\[
 |\widehat h_{q_j}(k)|\le C(1+|k|)^{r_j-1}.
\]
For $0<r<1$, all Fourier coefficients of $h_{-r}$ are positive, and
the same gamma formula gives
\[
 C^{-1}(1+|k|)^{r-1}\le\widehat h_{-r}(k)
                         \le C(1+|k|)^{r-1}.
\]
The constants are uniform for the present values of $r$.

Let $P_\rho$ be the Poisson kernel and take $\rho=e^{-1/L}$.
Since $\rho^{\sum_j|\ell_j|}\ge e^{-m}$ on the truncated sum,
positivity and the convolution identity imply that its left side is at most
\[
 C\int_\T\prod_{j=1}^m(P_\rho*h_{-r_j})(x)\,\sigma(dx).
\]
Writing $a=1-\rho\asymp L^{-1}$, the Poisson kernel bound gives
\begin{equation}\label{eq:poisson-root-bound}
 (P_\rho*h_{-r})(x)\le C(a+\|x\|_{2\pi})^{-r}.
\end{equation}
This can be verified directly as follows. If $\|x\|_{2\pi}\le2a$, rescale the
integral near zero by $a$ and use
$P_\rho(y)\le Ca/(a^2+\|y\|_{2\pi}^2)$; the resulting integral is
bounded uniformly for $r$ in a compact subset of $(0,1)$. If
$\|x\|_{2\pi}>2a$, split into $\|y\|_{2\pi}\le\|x\|_{2\pi}/2$
and its complement. On the first part the translated singular factor is
at most $C\|x\|_{2\pi}^{-r}$; on the second part the kernel bound and a
further split at the singularity give the same estimate.

Put $R=\sum_jr_j\le A_K+mt$. It follows that the required integral is
bounded by
\[
 C\int_0^\pi(a+x)^{-R}\,dx
 \le C(1+L)^{(A_K-1)_++\zeta}.
\]
The last estimate is uniform even when $R$ approaches one: a logarithmic
factor is absorbed by the unused part of $\zeta$. This proves
\eqref{eq:terminal-truncated-mass}. For the mixed assertion on a complex
compact set, choose $r_j=(-\Re q_j)_++t$ with the same restrictions
on $t$. The complex coefficient bound
\eqref{eq:mixed-terminal-fourier-decay} below gives
$|\widehat h_{q_j,s_j}(k)|\le C(1+|k|)^{r_j-1}$.
The positive majorants $h_{-r_j}$ and the Poisson integral estimate
\eqref{eq:poisson-root-bound} then give the same conclusion.
\end{proof}

\begin{remark}
The full range $q_j>-1$ is the finite ensemble integrability domain at
distinct points. The density bound \eqref{eq:phase-law-decomposition}
contains the factor
$\delta^{-m/2}/\mathcal V(\boldsymbol\vartheta)$. Thus
\cref{thm:full-root-FH} by itself does not give a uniform asymptotic
for every moving configuration with $N\delta_N\to\infty$.
\Cref{thm:full-signed-polynomial} gives an explicit polynomial merging range.
On $\mathcal D_m$, \cref{thm:signed-characteristic-FH} also gives
an error bound that depends on the separation only through $N\delta$.
\end{remark}

\subsection{Approaching points throughout the integrable range}

\begin{theorem}
\label{thm:full-signed-polynomial}
Let $m\ge2$ and $K\subset(-1,\infty)^m$ be compact. Define
$A_K,M_K,\eta_K,g_\beta$ as at the start of
\cref{sec:full-negative-transfer}, and put $B_m=m^2/2$. Suppose
\begin{equation}\label{eq:signed-merging-threshold}
 0\le\theta<
 \frac{\eta_Kg_\beta}{\eta_Kg_\beta+B_mM_K}.
\end{equation}
Then \eqref{eq:full-root-FH} holds uniformly for $\mathbf q\in K$
and $\delta(\boldsymbol\vartheta)\ge N^{-\theta}$, with relative
error $O(N^{-\kappa})$, for every
\begin{equation}\label{eq:signed-merging-rate}
 0<\kappa<
 \frac{\eta_Kg_\beta(1-\theta)-\theta B_mM_K}{M_K+\eta_K}.
\end{equation}
For $\theta=0$, one may instead fix any positive separation, with
constants depending on it. For $m=1$, \cref{prop:onepoint} gives
$O(N^{-1})$ without any separation parameter.
\end{theorem}

\begin{proof}[Proof of \cref{thm:full-signed-polynomial}]
We repeat the balance in the proof of \cref{thm:full-root-FH}, now
retaining the dependence of the approximation error on the separation.
Choose $\alpha<\eta_K$, $\gamma<g_\beta$ and $M>M_K$ so close to
their limiting values that \eqref{eq:signed-merging-threshold} and
\eqref{eq:signed-merging-rate} remain valid with these choices.
Put $x=n\delta$.
The density bound \eqref{eq:phase-law-decomposition}, the Fej\'er
approximation obtained from \eqref{eq:terminal-Lone-translation},
and \cref{lem:terminal-truncated-mass} give, uniformly whenever $x\ge R_K$,
\[
 \left|\int W_{\mathbf q}\,d\mu_{n,\mathbf q}-\prod_jc(q_j)\right|
 \le C\{\delta^{-B_m}L^{-\alpha}+x^{-\gamma}L^M\}.
\]
The first term is the approximation error integrated against the full
density. The second follows from the normalized mode bound
\eqref{eq:relative-phase-mixing}. Choose $L$
comparable to $(\delta^{-B_m}x^\gamma)^{1/(M+\alpha)}$. We obtain
\begin{equation}\label{eq:quantitative-signed-terminal-error}
 \left|\int W_{\mathbf q}\,d\mu_{n,\mathbf q}-\prod_jc(q_j)\right|
 \le C\delta^{-B_mM/(M+\alpha)}
                   (n\delta)^{-\alpha\gamma/(M+\alpha)}.
\end{equation}
If $\delta\ge N^{-\theta}$ and $n=N-1$, this is $O(N^{-\kappa})$
whenever
\[
 \kappa<\frac{\alpha\gamma(1-\theta)-\theta B_mM}{M+\alpha}.
\]
The stated threshold on $\theta$ is less than one, so $n\delta$
tends uniformly to infinity in this range and the density estimate
applies for all sufficiently large $N$.

Finally combine the identity \eqref{eq:terminal-conditioned} with \cref{thm:opuc-relative}.
Its relative error $O((n\delta)^{-1})$ is smaller than the bound in
\eqref{eq:quantitative-signed-terminal-error}. The means $\prod_jc(q_j)$ are uniformly positive on $K$, so this
error is relative as well. Replacing $n$ by $N=n+1$ costs $O(N^{-1})$.
Letting $\alpha,\gamma$ increase to their allowed limits and $M$
decrease to $M_K$ proves every exponent in
\eqref{eq:signed-merging-rate}.
\end{proof}

\begin{remark}
We do not determine whether this separation threshold is optimal.
The proof uses \eqref{eq:phase-law-decomposition}, which contains a
power of the separation. It therefore does not
establish the signed formula under the sole assumption
$N\delta(\boldsymbol\vartheta)\to\infty$. The next argument reaches
that conclusion on $\mathcal E_m$ by controlling only pairs of phases.
\end{remark}

\subsection{Uniform convergence for approaching points}
\label{sec:moving-pair-smoothing}

We now remove the restriction on the rate at which $N\delta$ tends to
infinity when $\sum_j(-q_j)_+<2$, and obtain an error bounded by a
negative power of $N\delta$. Only a bound
for each difference of two phases is needed. The argument is related to
the continuous phase smoothing in
\cite[Lemma~3.1]{AssiotisBKT2026}. Here we work with the finite sequence
of Verblunsky coefficients, under the weight $Y_n$, and retain constants
uniform in the moving points. The proof of
\cref{prop:moving-pair-density} uses one integration by parts. It does not give a uniform bound for the full
phase density.

\begin{theorem}
\label{thm:moving-signed-two}
Let $K\subset\mathcal E_m$ be compact. There are
$\kappa_K>0$ and $C,R<\infty$, depending only on $\beta,m,K$, such
that, uniformly for $\mathbf q\in K$ and distinct points with
$N\delta(\boldsymbol\vartheta)\ge R$,
\begin{equation}\label{eq:moving-signed-two}
 F_N(\mathbf q,\boldsymbol\vartheta)
 =N^{a(\mathbf q)}\mathcal A_\beta(\mathbf q,\boldsymbol\vartheta)
   \left(1+O\bigl((N\delta(\boldsymbol\vartheta))^{-\kappa_K}\bigr)\right).
\end{equation}
The implicit constant is at most $C$. In particular, for two points
this holds for every compact $K\subset(-1,\infty)^2$.
\end{theorem}

We need only the following bound for each difference of two phases,
rather than a bound for the full density on $\T^m$.

\begin{proposition}
\label{prop:moving-pair-density}
Let $K\subset(-1,\infty)^m$ be compact and let $j\ne l$. There are
$C,R<\infty$, depending only on $\beta,m,K$, such that the image of
$\mu_{n,\mathbf q}$ under
$\boldsymbol\varphi\mapsto\varphi_j-\varphi_l\pmod {2\pi}$
has density bounded by $C$ with respect to angular Lebesgue measure,
provided
\[
 \mathbf q\in K,\qquad
 n|e^{i\vartheta_j}-e^{i\vartheta_l}|\ge R.
\]
All points are assumed distinct. The bound is independent of
their configuration.
\end{proposition}

\Cref{prop:moving-pair-density} implies the uniform moment bound in
\cref{lem:moving-terminal-ui}. We combine its deterministic majorant
\eqref{eq:terminal-pair-majorant} with the Fourier estimate
\eqref{eq:relative-phase-mixing}: truncating each singular factor gives
a Lipschitz function, and the pair densities control the error of this
truncation uniformly in the prescribed points.

\begin{lemma}
\label{lem:moving-terminal-ui}
Let $K\subset\mathcal E_m$ be compact. There exist $p>1$ and
$C,R<\infty$ such that
\begin{equation}\label{eq:moving-terminal-ui}
 \sup_{\substack{\mathbf q\in K,\ n\delta\ge R}}
       \int W_{\mathbf q}^{p}\,d\mu_{n,\mathbf q}\le C.
\end{equation}
\end{lemma}

\begin{proof}[Proof of \cref{thm:moving-signed-two}]
We truncate each singular factor near its zero, approximate the resulting
Lipschitz function by a trigonometric polynomial, and apply the mode
estimate \eqref{eq:relative-phase-mixing}. The pair densities bound
the truncation error near coincident phases. Choosing the truncation
and approximation scales as powers of $n\delta$ gives the result.

For $m=1$, $W_{\mathbf q}$ is the constant $c(q_1)$, and
\cref{prop:onepoint} gives the assertion with $\kappa_K=1$.
Suppose $m\ge2$. Recall that
\[
 A_K=\sup_{\mathbf q\in K}\sum_{j=1}^m(-q_j)_+<2,
 \qquad
 \eta_K=\min\{1,\inf_{\mathbf q\in K,\,1\le j\le m}(1+q_j)\}>0.
\]
Choose $s$ and $p$, depending only on $K$, such that
\[
 \max\{A_K-1,0\}<s<1,\qquad 1<p<1/s,
 \qquad \zeta=1-1/p>0.
\]
These are the conditions used in the proof of
\cref{lem:moving-terminal-ui}. We put $x=n\delta$ and assume
$x\ge R$, where $R$ is at least the threshold in
\cref{prop:moving-pair-density}. In the estimates below, $H$ denotes
Haar measure $\sigma^{\otimes m}$ on the phase torus, and $\nu$ denotes
either $\mu_{n,\mathbf q}$ or $H$. All constants are uniform for
$\mathbf q\in K$ and depend only on $\beta,m,K$.

For $0<\tau<1$, truncate the distance from the singularity at level
$\tau$ by defining
\[
 h_{q,\tau}(t)=(|1-e^{it}|\vee\tau)^q,
 \qquad
 W_{\mathbf q,\tau}(\boldsymbol\varphi)
 =\int_\T\prod_{j=1}^m h_{q_j,\tau}(u+\varphi_j)\,\sigma(du).
\]
Write $\mathbf q^-=(\min\{q_j,0\})_{j=1}^m$. Truncation decreases
a negative power, while all positive powers remain bounded uniformly
on $K$. Therefore
$W_{\mathbf q}+W_{\mathbf q,\tau}\le C W_{\mathbf q^-}$.
The vectors $\mathbf q^-$ range over a compact subset of
$\mathcal E_m$ with the same bound $A_K$ on their total negative
exponent. The deterministic estimate
\eqref{eq:terminal-pair-majorant}, applied to these vectors, gives
\[
 W_{\mathbf q}(\boldsymbol\varphi)
       +W_{\mathbf q,\tau}(\boldsymbol\varphi)
 \le C\left(1+\sum_{1\le j<l\le m}
               \|\varphi_j-\varphi_l\|_{2\pi}^{-s}\right).
\]
We integrate this majorant under $\mu_{n,\mathbf q}$, using
\cref{prop:moving-pair-density} for the original exponent vector
$\mathbf q$. Since $sp<1$, the same calculation as in the proof of
\cref{lem:moving-terminal-ui} yields
\[
 \int (W_{\mathbf q}+W_{\mathbf q,\tau})^p\,d\nu\le C.
\]
For Haar measure the bound follows directly by integrating each pair
power on the circle.

Each factor $h_{q_j,\tau}$ has supremum norm at most
$C\tau^{-(-q_j)_+}$ and Lipschitz constant at most
$C\tau^{-(-q_j)_+-1}$. Indeed, where the distance exceeds $\tau$,
its derivative is bounded by $C|1-e^{it}|^{q_j-1}$; on the remaining
arc the factor is constant. Multiplying the factors and integrating
the common phase therefore gives
\[
 \|W_{\mathbf q,\tau}\|_\infty\le C\tau^{-A_K},
 \qquad
 \operatorname{Lip}(W_{\mathbf q,\tau})\le C\tau^{-A_K-1},
\]
where the Lipschitz constant is taken with respect to the sum of the
circular distances in the $m$ coordinates.

Let $B_\varepsilon$ be the set of phase vectors for which
$\|\varphi_j-\varphi_l\|_{2\pi}<\varepsilon$ for at least one pair
$j<l$, where $0<\varepsilon<1$. The pair density bound implies
$\nu(B_\varepsilon)\le C\varepsilon$. H\"older's inequality and
the preceding $p$th moment bound thus give
\[
 \int_{B_\varepsilon}|W_{\mathbf q}-W_{\mathbf q,\tau}|\,d\nu
 \le C\varepsilon^\zeta.
\]
Outside $B_\varepsilon$, assume $\tau\le\varepsilon/4$. The two
common phase integrands differ only on the arcs
$\{|1-e^{i(u+\varphi_j)}|<\tau\}$. On the $j$th such arc,
for every $l\ne j$,
\[
 |1-e^{i(u+\varphi_l)}|
 \ge |e^{i\varphi_j}-e^{i\varphi_l}|
       -|1-e^{i(u+\varphi_j)}|
 \ge \frac2\pi\varepsilon-\tau\ge\varepsilon/4.
\]
The product of the other factors is consequently bounded by
$C\varepsilon^{-A_K}$. Also, uniformly for the individual exponents
appearing in $K$,
\[
 \int_{\{|1-e^{it}|<\tau\}}
       \bigl||1-e^{it}|^q-h_{q,\tau}(t)\bigr|\,\sigma(dt)
 \le C\tau^{1+q}\le C\tau^{\eta_K}.
\]
To obtain the first bound, use $|1-e^{it}|\asymp|t|$ on this arc
and integrate $|t|^q+\tau^q$; the constants are uniform because
$q+1\ge\eta_K$. A telescoping expansion of the difference of the
products now gives
\[
 \int|W_{\mathbf q}-W_{\mathbf q,\tau}|\,d\nu
 \le C\left(\varepsilon^\zeta+
                   \tau^{\eta_K}\varepsilon^{-A_K}\right),
 \qquad \tau\le\varepsilon/4.
\]

We next convolve $W_{\mathbf q,\tau}$ with the product Fej\'er kernel
of degree $L\ge2$ in each coordinate, and denote the resulting
trigonometric polynomial by $P_L$. The one variable Fej\'er kernel
is bounded by
$C\min\{L,(L\|t\|_{2\pi}^2)^{-1}\}$, so its first circular
moment is at most $C(1+\log L)/L$. Applying the Lipschitz estimate
coordinate by coordinate gives
\[
 \|P_L-W_{\mathbf q,\tau}\|_\infty
 \le C\tau^{-A_K-1}\frac{1+\log L}{L}
 \le C\tau^{-A_K-1}L^{-1/2}.
\]
There are at most $(2L+1)^m$ nonzero Fourier coefficients of $P_L$,
each bounded in modulus by $C\tau^{-A_K}$. Moreover,
$\widehat P_L(0)=\int W_{\mathbf q,\tau}\,dH$.
The function $W_{\mathbf q,\tau}$ is invariant under a common shift
of all phases, by translation of the common integration variable.
Hence every nonzero Fourier coefficient of $P_L$ has
$\sum_j\ell_j=0$.

Fix $\gamma=g_\beta/2$, where $g_\beta=\min\{4/\beta,1\}$.
For a nonzero mode with zero sum,
\[
 \left|\int e^{i\boldsymbol\ell\cdot\boldsymbol\varphi}
                    \,d\mu_{n,\mathbf q}\right|
 =\frac{|\E T_n(\mathbf q,\boldsymbol\ell)|}{\E Y_n(\mathbf q)}
 \le C(1+\|\boldsymbol\ell\|_1)^2x^{-\gamma}
\]
by \eqref{eq:relative-phase-mixing}. Summing over the modes of
$P_L$ gives
\[
 \left|\int P_L\,d\mu_{n,\mathbf q}-\int P_L\,dH\right|
 \le C\tau^{-A_K}L^{m+2}x^{-\gamma}.
\]
Combining this estimate with the truncation and approximation bounds
under both probability measures, we obtain
\[
 \left|\int W_{\mathbf q}\,d\mu_{n,\mathbf q}
                    -\int W_{\mathbf q}\,dH\right|
 \le C\left(\varepsilon^\zeta+
       \tau^{\eta_K}\varepsilon^{-A_K}
       +\tau^{-A_K-1}L^{-1/2}
       +\tau^{-A_K}L^{m+2}x^{-\gamma}\right).
\]

We choose the three scales so that every term is a negative power of
$x$. Put
\[
 C_0=A_K+4(A_K+1)(m+2),\qquad
 v=\frac{\gamma}{2C_0},\qquad
 u=\frac{v\eta_K}{2(A_K+1)},\qquad
 w=4v(A_K+1),
\]
and take $\varepsilon=x^{-u}$, $\tau=x^{-v}$ and
$L=\lceil x^w\rceil$. Since $\eta_K\le1$, we have $u<v$;
thus $\tau\le\varepsilon/4$ after increasing $R$. The four error
terms are bounded, respectively, by
\[
 Cx^{-u\zeta},\qquad Cx^{-v\eta_K/2},\qquad
 Cx^{-v(A_K+1)},\qquad Cx^{-\gamma/2}.
\]
Each exponent is at least
\[
 \kappa_K=u\zeta
 =\frac{\gamma\eta_K\zeta}{4C_0(A_K+1)}>0.
\]
Fubini's theorem and rotation invariance give
$\int W_{\mathbf q}\,dH=\prod_{j=1}^m c(q_j)$, which is bounded
below by a positive constant on $K$. We have therefore proved
\[
 \int W_{\mathbf q}\,d\mu_{n,\mathbf q}
 =\prod_{j=1}^m c(q_j)\bigl(1+O(x^{-\kappa_K})\bigr).
\]
Insert this formula into \eqref{eq:terminal-conditioned} and apply
\cref{thm:opuc-relative}, whose relative error $O((n\delta)^{-1})$
is smaller because $\kappa_K<1$. The identity
$c(q)D_\beta(q)=C_\beta(q)$ identifies the constants. Finally,
$n=N-1$ gives $n\delta\ge N\delta/2$ for $N\ge2$, and replacing
$n^{a(\mathbf q)}$ by $N^{a(\mathbf q)}$ costs a relative
$O(N^{-1})$. Since $\delta\le1$, this is also absorbed by the
stated error, after increasing $R$ once more.
\end{proof}

The exponent $\kappa_K$ obtained in this proof need not be optimal.
For $K\subset\mathcal D_m$, the estimates of
\cref{thm:signed-characteristic-FH,cor:terminal-subset-rate}
remain available as well; one may use whichever gives the larger
admissible exponent for the chosen compact set.

We now prove the two estimates used above. The pair density bound is
the longer argument. The subsequent proof of
\cref{lem:moving-terminal-ui} supplies the deterministic comparison
\eqref{eq:terminal-pair-majorant} and integrates it against those
pair densities.

The conditional martingale inequality in
\cref{lem:conditional-martingale-maximal} is used in
\cref{lem:pair-occupation-estimate,lem:pair-derivative-estimate};
see also \cite{Burkholder1973}.

\begin{lemma}[Conditional maximal inequality]
\label{lem:conditional-martingale-maximal}
Fix $v\ge2$. Let $(\Delta_j)_{j=1}^J$ be a finite real martingale
difference sequence with respect to $(\mathcal G_j)_{j=0}^J$, with
$\Delta_j\in L^v$ and $\mathcal F_a\subset\mathcal G_0$.
Write $\E_a=\E[\,\cdot\mid\mathcal F_a]$. Denote the partial sums,
their maximum, and the predictable sum of conditional variances by
\[
 S_l=\sum_{j=1}^l\Delta_j,\qquad
 S^*=\max_{0\le l\le J}|S_l|,\qquad
 V=\sum_{j=1}^J\E[\Delta_j^2\mid\mathcal G_{j-1}].
\]
Then
\begin{equation}\label{eq:conditional-martingale-maximal}
 \E_a(S^*)^v
 \le C_v\left\{\E_a V^{v/2}
                   +\sum_{j=1}^J\E_a|\Delta_j|^v\right\}.
\end{equation}
The constant $C_v$ depends only on $v$.
\end{lemma}

\begin{proof}[Proof of \cref{lem:conditional-martingale-maximal}]
For $v>2$ the
second order Taylor bound
\[
 |x+y|^v\le |x|^v+v|x|^{v-2}xy
                 +C_v\bigl(|x|^{v-2}y^2+|y|^v\bigr)
\]
and conditional centering give, after summation,
\[
 \E_a|S_J|^v
 \le C_v\E_a[(S^*)^{v-2}V]
                  +C_v\sum_{j=1}^J\E_a|\Delta_j|^v.
\]
Conditional Doob's inequality bounds $\E_a(S^*)^v$ by a constant
times $\E_a|S_J|^v$. Conditional H\"older then gives
\[
 \E_a[(S^*)^{v-2}V]
 \le \bigl(\E_a(S^*)^v\bigr)^{1-2/v}
             \bigl(\E_a V^{v/2}\bigr)^{2/v}.
\]
Young's inequality bounds this by
$\epsilon\E_a(S^*)^v+C_{v,\epsilon}\E_a V^{v/2}$.
Choose $\epsilon$ small enough to absorb the first term into the
left side. This proves \eqref{eq:conditional-martingale-maximal}. For $v=2$, conditional
orthogonality and Doob's inequality give it directly. Apply this to
the real and imaginary parts when the increments are complex.
\end{proof}

Fix the pair $j\ne l$ and a compact set $K\subset(-1,\infty)^m$.
Put $a=\lfloor n/2\rfloor$ and $L=n-a$, so that we work on the last
half of the coefficient sequence. In the next three lemmas and the
proof of \cref{prop:moving-pair-density}, write
$\E_a=\E[\,\cdot\mid\mathcal F_a]$ and
$\|A\|_{v,a}=(\E_a|A|^v)^{1/v}$. All moment bounds are conditional
on the incoming phases at time $a$.
All differentiations below involve only
$\alpha_a,\ldots,\alpha_{n-1}$. For an index $r$, write
$\psi_{k,r}=\psi_k(\vartheta_r)$. The real relative phase
$x_k=\psi_{k,j}-\psi_{k,l}$ satisfies
\[
 x_{k+1}=x_k+\theta+\xi_k,\qquad
 \theta=\vartheta_j-\vartheta_l,
\]
where
\[
 \xi_k=-2\Im\{\log(1-\alpha_ke^{i\psi_{k,j}})
                    -\log(1-\alpha_ke^{i\psi_{k,l}})\}.
\]
The average separation of the two phases during this block is measured by
\[
 I=\frac1L\sum_{k=a}^{n-1}(1-\cos x_k).
\]

Throughout the next three lemmas, $n|1-e^{i\theta}|$ is assumed
sufficiently large. Each lemma states how this threshold depends on
the moment order; no assertion requires all moment orders at one fixed
starting index. Write $t_k=b(k+1)$.

\begin{lemma}[Occupation estimate]
\label{lem:pair-occupation-estimate}
For each integer $r\ge1$ there are $C_r,R_r<\infty$, depending only
on $\beta$ and $r$, such that, whenever $n|1-e^{i\theta}|\ge R_r$,
\begin{equation}\label{eq:pair-small-occupation}
 \Prob(I\le\varepsilon\mid\mathcal F_a)\le C_r\varepsilon^r,
 \qquad 0<\varepsilon\le1.
\end{equation}
The estimate is uniform in the phases at time $a$.
\end{lemma}

\begin{proof}[Proof of \cref{lem:pair-occupation-estimate}]
The deterministic rotation spends a positive fraction of the block
away from coincident phases. If the random occupation average is small,
the accumulated random increments must therefore be large. We bound
that event by stopping before their conditional variance becomes large.
Rotation invariance gives $\E[\xi_k\mid\mathcal F_k]=0$.
For every fixed integer $r\ge1$, and all sufficiently large $k$,
\begin{equation}\label{eq:pair-increment-moment}
 \E[|\xi_k|^{2r}\mid\mathcal F_k]
 \le C_r k^{-r}(1-\cos x_k)^r.
\end{equation}
Indeed, integration of the analytic derivative along the chord joining
$z,w\in\partial\mathbb D$ gives
\[
 |\log(1-\alpha z)-\log(1-\alpha w)|
 \le\frac{|\alpha|}{1-|\alpha|}|z-w|.
\]
With $t_k=\beta(k+1)/2$, the variable $|\alpha_k|^2$ has density
$t_k(1-u)^{t_k-1}\,du$. For $t_k>2r$,
\[
 \E\frac{|\alpha_k|^{2r}}{(1-|\alpha_k|)^{2r}}
 \le 2^{2r}t_k\int_0^1u^r(1-u)^{t_k-2r-1}\,du
 =2^{2r}t_k\frac{\Gamma(r+1)\Gamma(t_k-2r)}
                         {\Gamma(t_k-r+1)}
 \le C_r k^{-r}.
\]
This proves \eqref{eq:pair-increment-moment}, including the conditional
variance bound $Ck^{-1}(1-\cos x_k)$.

The geometric series bound gives
\[
 \frac1L\sum_{s=0}^{L-1}(1-\cos(x_a+s\theta))
 \ge1-\frac{2}{L|1-e^{i\theta}|}\ge\frac12
\]
once $n|1-e^{i\theta}|$ exceeds a fixed constant. Since cosine is
Lipschitz and $x_{a+s}=x_a+s\theta+\sum_{k=a}^{a+s-1}\xi_k$,
\[
 \frac12\le I+\max_{0\le s\le L}
                    \left|\sum_{k=a}^{a+s-1}\xi_k\right|.
\]
For $\varepsilon\le1/4$, stop the martingale predictably before the
first increment whose cost $1-\cos x_k$ would make the accumulated
cost exceed $L\varepsilon$. The cost is $\mathcal F_k$ measurable,
so the stopped increments still have conditional mean zero. On
$\{I\le\varepsilon\}$ every increment is retained. Their conditional
variances have sum at most $C\varepsilon$. Moreover,
\[
 \sum_{\text{retained }k}
   \E[|\xi_k|^{2r}\mid\mathcal F_k]
 \le C_r n^{-r}\sum_{\text{retained }k}(1-\cos x_k)^r
 \le C_r n^{-r}(L\varepsilon)^r
 \le C_r\varepsilon^r.
\]
The martingale maximal inequality \eqref{eq:conditional-martingale-maximal}
therefore bounds the
$2r$th moment of the stopped maximum by $C_r\varepsilon^r$.
On $\{I\le\varepsilon\}$ that maximum is at least $1/4$.
Markov's inequality proves \eqref{eq:pair-small-occupation}; enlarging
the constant covers $\varepsilon>1/4$.

\end{proof}

We next choose a direction in the coefficient variables whose effect
on $x_n$ is nonnegative. Regard each coefficient disk as $\mathbb R^2$
and represent real gradients by complex numbers. The quantities
$P_{k,r}$ and $u_{k,r}$ below are the phase and coefficient derivatives
of the phase map. The difference of the two normalized coefficient
gradients defines the direction $h_k$:
\begin{gather*}
 P_{k,r}=\frac{1-|\alpha_k|^2}
                   {|1-\alpha_ke^{i\psi_{k,r}}|^2},\qquad
 u_{k,r}=\frac{2i}{e^{i\psi_{k,r}}-\overline{\alpha_k}},\\
 h_k=t_k^{-1}\left(\frac{u_{k,j}}{|u_{k,j}|}
                         -\frac{u_{k,l}}{|u_{k,l}|}\right),
 \qquad D_h=\sum_{k=a}^{n-1}h_k\cdot\nabla_{\alpha_k}.
\end{gather*}
\begin{lemma}[A positive derivative direction]
\label{lem:pair-direction-estimate}
Let $D=D_hx_n$. For every real $v>0$ there are $C_v,R_v<\infty$
such that, whenever $n|1-e^{i\theta}|\ge R_v$, one has $D>0$
almost surely and
\begin{equation}\label{eq:pair-inverse-direction}
 \E_a[D^{-v}]\le C_v\qquad(v>0).
\end{equation}
The constants depend only on $\beta,v$ and are uniform in the
incoming phases.
\end{lemma}

\begin{proof}[Proof of \cref{lem:pair-direction-estimate}]
We express $D$ as a sum of nonnegative contributions, then bound it
below by the occupation average times a product of phase derivatives.
The occupation estimate \eqref{eq:pair-small-occupation} and
martingale bounds for these products will give the inverse moments.
The phase map is $g_\alpha(\phi)=\phi-2\Im\log(1-\alpha e^{i\phi})$;
its derivatives are $\partial_\phi g_\alpha=P$ and
$\nabla_\alpha g_\alpha=u$. The product
$\Pi_{k,r}=\prod_{v=k+1}^{n-1}P_{v,r}$ propagates a perturbation at
step $k$ to the final phase. Differentiating the recursion gives exactly
\begin{equation}\label{eq:pair-positive-direction}
 D:=D_hx_n
 =\sum_{k=a}^{n-1}t_k^{-1}
       (\Pi_{k,j}|u_{k,j}|+\Pi_{k,l}|u_{k,l}|)(1-\cos\omega_k),
\end{equation}
where $\omega_k$ is the angle between $u_{k,j}$ and $u_{k,l}$.
For $|z|<1$,
\[
 \frac d{d\phi}\arg(e^{i\phi}-z)
 =\frac12+\frac{1-|z|^2}{2|e^{i\phi}-z|^2}\ge\frac12.
\]
Apply this inequality on both complementary arcs joining the two
phases. The smaller angular distance between their images is at least
half their original smaller distance. Consequently
$1-\cos\omega_k\ge(1-\cos x_k)/4$. Also $|u_{k,r}|\ge1$.
Thus
\begin{equation}\label{eq:pair-direction-lower}
 D\ge c_\beta I\min_{a\le k<n,\ r\in\{j,l\}}\Pi_{k,r}.
\end{equation}

For every fixed real $v$, rotation invariance and the disk split used
in \cref{lem:onestep} imply
\[
 c_v(k):=\E[P_{k,r}^{v}\mid\mathcal F_k]=1+O_v(k^{-1}),
\]
independently of the incoming phase, once $k$ is sufficiently large.
One may verify the error by expansion on $|\alpha_k|\le1/2$;
rotation invariance cancels the linear term, and the complementary
annulus has an exponentially small integral, also for negative $v$.
To relate partial products to the suffix products in
\eqref{eq:pair-direction-lower}, write
\[
 A_{b,r}=\prod_{k=a}^b P_{k,r},\qquad a\le b<n,
 \qquad A_{a-1,r}=1.
\]
For every fixed real $v$, the normalized products
\[
 M_{b,r}^{(v)}=
 \frac{A_{b,r}^v}{\prod_{k=a}^b c_v(k)},\qquad a-1\le b<n,
\]
form a positive martingale. Iterated conditional expectation gives
\[
 \E_a\bigl[(M_{n-1,r}^{(v)})^2\bigr]
 =\prod_{k=a}^{n-1}\frac{c_{2v}(k)}{c_v(k)^2}\le C_v,
\]
because $\sum_{k=a}^{n-1}k^{-1}$ is bounded. The deterministic
normalizing products are bounded above and below by positive constants
for the same reason. Conditional Doob's inequality consequently gives
\[
 \E_a\left[\max_{a-1\le b<n}A_{b,r}^{2v}\right]\le C_v.
\]
This estimate applies to both signs of $v$.
In particular, for $v>0$, the identity
$\Pi_{k,r}=A_{n-1,r}/A_{k,r}$ gives
\[
 \left(\min_{a\le k<n}\Pi_{k,r}\right)^{-v}
 =\left(\max_{a\le k<n}A_{k,r}\right)^v A_{n-1,r}^{-v}.
\]
By H\"older's inequality, its conditional expectation is at most
\[
 \left(\E_a\max_{a\le k<n}A_{k,r}^{2v}\right)^{1/2}
 \left(\E_a A_{n-1,r}^{-2v}\right)^{1/2}\le C_v.
\]
Taking the maximum over $r\in\{j,l\}$ and bounding it by the sum
therefore yields
\[
 \E_a\left[\left(\min_{a\le k<n,\ r\in\{j,l\}}\Pi_{k,r}\right)^{-v}
                   \right]\le C_v,
 \qquad v>0.
\]
Together with \eqref{eq:pair-small-occupation} and
\eqref{eq:pair-direction-lower}, this proves
\eqref{eq:pair-inverse-direction}.
For each displayed moment one first chooses a sufficiently large
occupation exponent, and then increases the starting index $a$.

\end{proof}

For every $1\le r\le m$, denote the first and second derivatives
along $h$ by $V_{k,r}=D_h\psi_{k,r}$ and
$B_{k,r}=D_h^2\psi_{k,r}$. The weight acquired during the block is
$Z=Y_n(\mathbf q)/Y_a(\mathbf q)$.
Write $\rho_k(z)=t_k\pi^{-1}(1-|z|^2)^{t_k-1}$ for the disk density.
The negative divergence of $h$ relative to the product coefficient
law, which is the score in the integration by parts formula, is
\[
 \mathcal H=-\sum_{k=a}^{n-1}
      \{\operatorname{div}_{\alpha_k}h_k+
                              h_k\cdot\nabla\log\rho_k(\alpha_k)\}.
\]
\begin{lemma}[Moments of the derivatives and the score]
\label{lem:pair-derivative-estimate}
For each fixed $v\ge2$ there are $C_{K,v},R_{K,v}<\infty$ such that,
if $n|1-e^{i\theta}|\ge R_{K,v}$ and $\mathbf q\in K$, then
\[
 \E_a\left[
  \max_{a\le k\le n,\,1\le r\le m}(|V_{k,r}|^v+|B_{k,r}|^v)
  +Z^v+|D_h\log Z|^v+|\mathcal H|^v
       \right]\le C_{K,v},
 \qquad \E_a[Z]\ge1.
\]
The estimates are uniform in the incoming phases and the distinct
configuration of prescribed points.
\end{lemma}

\begin{proof}[Proof of \cref{lem:pair-derivative-estimate}]
We first bound the derivatives of the phase flow, then the weight and
its derivative, and finally the score. To control second derivatives
at order $v$, the phase recursion will require first derivative
moments at order $2v$.
The dependence of $h_k$ on the earlier coefficients must be retained.
Both $V_{k,r}$ and $B_{k,r}$ vanish at $k=a$. With
$F(\phi,\alpha)=g_\alpha(\phi)$, the exact recursions are
\begin{align}
 V_{k+1,r}&=P_{k,r}V_{k,r}+F_\alpha\cdot h_k,
                                                    \label{eq:pair-V-recursion}\\
 B_{k+1,r}&=P_{k,r}B_{k,r}+F_{\phi\phi}V_{k,r}^2
       +2F_{\phi\alpha}[V_{k,r},h_k]
       +F_{\alpha\alpha}[h_k,h_k]
       +F_\alpha\cdot D_hh_k.                     \label{eq:pair-B-recursion}
\end{align}
All derivatives on the right are evaluated at
$(\psi_{k,r},\alpha_k)$.

Put $z=\alpha e^{i\phi}$ and regard $v,w\in\mathbb C$
as real directions in the coefficient disk. Direct differentiation gives
\[
 \begin{aligned}
 F_\phi&=1+2\Re\frac{z}{1-z},&
 F_{\phi\phi}&=-2\Im\frac{z}{(1-z)^2},\\
 F_{\phi\alpha}[v]&=2\Re\frac{v e^{i\phi}}{(1-z)^2},&
 F_{\alpha\alpha}[v,w]&=2\Im\frac{vw e^{2i\phi}}{(1-z)^2}.
 \end{aligned}
\]
Consequently
\[
 |F_{\phi\phi}|\le
      \frac{2|\alpha|}{(1-|\alpha|)^2},\qquad
 \|F_{\phi\alpha}\|,\ \|F_{\alpha\alpha}\|
      \le\frac{2}{(1-|\alpha|)^2}.
\]
The rotational mean of $F_{\phi\phi}$ is zero, as is seen by expanding
$z/(1-z)^2=\sum_{h\ge1}h z^h$. Uniformly in the incoming phases,
\[
 \E[P_{k,r}-1\mid\mathcal F_k]
 =\E[F_{\phi\phi}\mid\mathcal F_k]=0,
 \qquad
 \|P_{k,r}-1\|_{v,k}+\|F_{\phi\phi}\|_{v,k}
 \le C_v k^{-1/2},
\]
where $\|X\|_{v,k}=(\E[|X|^v\mid\mathcal F_k])^{1/v}$.
These estimates follow from the same disk split: the functions have
zero rotational mean, vanish at $\alpha=0$, and their derivatives
are bounded by fixed powers of $(1-|\alpha|)^{-1}$. Every remaining
derivative in \eqref{eq:pair-V-recursion} and~\eqref{eq:pair-B-recursion}
has bounded conditional moments of each fixed order for large $k$.
Furthermore,
\[
 |h_k|\le2/t_k,\qquad
 \|D_hh_k\|_{v,k}
 \le C_vk^{-1}(|V_{k,j}|+|V_{k,l}|)+C_vk^{-2}.
\]
To see the last inequality, write
$w=e^{i\phi}-\bar\alpha$, so that $u/|u|=i\bar w/|w|$.
For the normalization map on $\mathbb C\setminus\{0\}$,
\[
 D\left(\frac{w}{|w|}\right)[v]
 =\frac{v}{|w|}
  -\frac{w\Re(\bar wv)}{|w|^3},
 \qquad
 \left|D\left(\frac{w}{|w|}\right)[v]\right|
 \le\frac{|v|}{|w|}.
\]
Since $|w|\ge1-|\alpha|$, its derivatives in the phase and coefficient
arguments are bounded by $C(1-|\alpha|)^{-1}$. After the factor
$t_k^{-1}$ in $h_k$, phase differentiation gives
$C_vk^{-1}(|V_{k,j}|+|V_{k,l}|)$ in conditional norm.
Differentiation in the current coefficient introduces $h_k=O(k^{-1})$
and gives the remaining $C_vk^{-2}$.

It follows that, for every fixed $v\ge2$,
\begin{equation}\label{eq:pair-flow-moments}
 \E_a\left[\max_{a\le k\le n,\,1\le r\le m}
       (|V_{k,r}|^v+|B_{k,r}|^v)\right]\le C_v.
\end{equation}
We apply \cref{lem:conditional-martingale-maximal} to
\eqref{eq:pair-V-recursion} and \eqref{eq:pair-B-recursion}.

Since $\sum_{k=a}^{n-1}k^{-1}$ is bounded, we can partition the block
into a bounded number of intervals on each of which
$\sum k^{-1}\le\eta_v$, with $\eta_v$ chosen below.
For a predictable process $U_k$ and centered increments $A_k$ with
$\|A_k\|_{v,k}\le C_vk^{-1/2}$, the maximal inequality
\eqref{eq:conditional-martingale-maximal} and Minkowski's inequality give,
on each such interval,
\[
 \left\|\max_b\left|\sum_{k<b}A_kU_k\right|\right\|_{v,a}
 \le C_v\left(\sum_k k^{-1}\|U_k\|_{v,a}^2\right)^{1/2}.
\]
The additional sum of $v$th moments in
\eqref{eq:conditional-martingale-maximal} is bounded by the same
right-hand side, since $v\ge2$. This is at most $C_v\eta_v^{1/2}$
times the maximal norm of $U$. Choosing $C_v\eta_v^{1/2}<1/2$
allows this term to be absorbed into the left side of the recursion estimate. In
\eqref{eq:pair-V-recursion}, the remaining terms have conditional norm
$O(k^{-1})$. This first gives all required moments of $V$. In
\eqref{eq:pair-B-recursion}, the term $F_{\phi\phi}V^2$ is centered
and is bounded using the moments of $V$ at order $2v$; the other forcing
terms have conditional norm at most
$C_vk^{-1}(1+\sum_r|V_{k,r}|^2)$. The same absorption controls the
term $(P-1)B$. Iterating over the bounded number of intervals proves
\eqref{eq:pair-flow-moments}, conditionally on any initial phases.

For fixed $v>1$, apply \eqref{eq:onestep} at $v\mathbf q$ and zero mode. Thus
$\E[w_{v\mathbf q}\mid\mathcal F_k]\le1+C_{K,v}/k$. Iteration gives
\begin{equation}\label{eq:pair-tilt-moments}
 \E_a[Z^v]\le C_{K,v},\qquad
 \E_a[Z]\ge1.
\end{equation}
The second inequality follows from Jensen's inequality and
$\E_a[\log Z]=0$. The contribution of one coefficient
to the logarithmic weight at point $r$ is
$G_r(\phi,\alpha)=q_r\Re\log(1-\alpha e^{i\phi})$. Thus
\[
 D_h\log Z=\sum_{k=a}^{n-1}\sum_{r=1}^m
       \{G_{r,\phi}V_{k,r}+G_{r,\alpha}\cdot h_k\}.
\]
The first terms are centered and have conditional norm
$C_{K,v}k^{-1/2}|V_{k,r}|$; the second have conditional norm
$C_{K,v}/k$. Applying \eqref{eq:conditional-martingale-maximal}
and the flow bound \eqref{eq:pair-flow-moments} yields
\begin{equation}\label{eq:pair-weight-derivative}
 \E_a[|D_h\log Z|^v]\le C_{K,v}.
\end{equation}

Each summand defining $\mathcal H$ is conditionally centered, by disk integration by parts
with the earlier phases fixed. Boundary terms vanish for large $k$.
Its conditional norm of fixed order is $O(k^{-1/2})$: the divergence
has norm $O(k^{-1})$, while
\[
 \nabla\log\rho_k(z)=-\frac{2(t_k-1)z}{1-|z|^2},
 \qquad |h_k|\le2/t_k.
\]
For a coefficient with density
$t\pi^{-1}(1-|z|^2)^{t-1}$, polar integration gives, whenever $t>v$,
\[
 \E\left(\frac{|z|}{1-|z|^2}\right)^v
 =t\,B(1+v/2,t-v)
 =\frac{t\,\Gamma(1+v/2)\Gamma(t-v)}
        {\Gamma(t+1-v/2)}
 =O_v(t^{-v/2}).
\]
Since $t_k$ is of order $k$, this bounds the radial part of the score
in conditional $L^v$ by $C_vk^{-1/2}$. Summing these centered
terms with \eqref{eq:conditional-martingale-maximal} gives
\begin{equation}\label{eq:pair-score-moments}
 \E_a[|\mathcal H|^v]\le C_v.
\end{equation}

\end{proof}

\begin{proof}[Proof of \cref{prop:moving-pair-density}]
We integrate by parts in the direction $h/D$, whose derivative of
$x_n$ is one. The moment estimates will bound the resulting score
in weighted $L^1$, and periodic primitives will then turn this bound
into a density estimate. Finally we average over the incoming phases.
We use only finitely many moment orders. Choose
the threshold $R$ so that \cref{lem:pair-occupation-estimate} applies
with $r=9$, \cref{lem:pair-direction-estimate} applies with $v=4$,
and \cref{lem:pair-derivative-estimate} applies with $v=4$.
In the proof of \cref{lem:pair-derivative-estimate}, eighth moments of $V$ control
fourth moments of $B$ through the quadratic forcing terms in
\eqref{eq:pair-B-recursion}; this requirement is included in the
threshold $R_{K,4}$.
In the proof of \cref{lem:pair-direction-estimate}, occupation order nine gives
$\E_a[I^{-8}]\le C$. The eighth inverse moment of the
minimum suffix product, together with H\"older's inequality, then gives
the required fourth inverse moment of $D$. Those suffix moments use
only the partial product estimates with exponents $16$ and $-16$ and
their second moments, hence disk moments with exponents at most $32$
in absolute value. Increase $R$ once more for these estimates and for
$\E[(1-|\alpha_k|)^{-2}]<\infty$ throughout the block.
All these choices depend only on $\beta,m,K$.

For a smooth periodic $f$, integration by parts now gives
\begin{equation}\label{eq:pair-weighted-ibp}
 \E_a[Zf'(x_n)]
 =\E_a[Zf(x_n)H],\qquad
 H=\frac{\mathcal H-D_h\log Z}{D}+\frac{D_hD}{D^2}.
\end{equation}
Indeed, apply disk integration by parts to the vector field $h/D$;
its negative weighted divergence with respect to the product
coefficient density $\rho$ is
$-\rho^{-1}\operatorname{div}(\rho h/D)
 =\mathcal H/D+(D_hD)/D^2$. Also
$D_hD=B_{n,j}-B_{n,l}$.
To justify this, first use compact cutoffs inside each coefficient
disk and replace $D^{-1}$ by $(D+u)^{-1}$, $u>0$. Before taking the
limit, the left side is
$\E_a[Zf'(x_n)D/(D+u)]$.
At each fixed $n$ there are finitely many boundary cutoffs.
For a cutoff varying where $1-|\alpha_k|$ is of order $\varepsilon$,
its derivative is bounded by $C(1-|\alpha_k|)^{-1}$ on that region.
Since $|h_k|\le2/t_k$, its error is dominated by a constant times
$ZD^{-1}(1-|\alpha_k|)^{-1}$. This is integrable by conditional
H\"older with exponents $4,4,2$. Dominated convergence therefore removes
the cutoff errors. The derivative estimates control the remaining terms.
Remove these cutoffs first and then let $u$ decrease
to zero. Equations \eqref{eq:pair-inverse-direction} and
\eqref{eq:pair-flow-moments} to~\eqref{eq:pair-score-moments} justify
convergence in $L^1$ on the right. Dominated convergence applies
on the left, with dominating variable $Z\|f'\|_\infty$.
With the conditional norms defined before \cref{lem:pair-occupation-estimate},
H\"older's inequality with exponents $4,4,2$, together with the finite
moment choices above, gives
\[
 \begin{aligned}
 \E_a[Z|H|]
 &\le \|Z\|_{4,a}
       \bigl(\|\mathcal H\|_{4,a}+\|D_h\log Z\|_{4,a}\bigr)
       \|D^{-1}\|_{2,a}\\
 &\quad+\|Z\|_{4,a}\|B_{n,j}-B_{n,l}\|_{4,a}
                         \|D^{-1}\|_{4,a}^{\,2}
 \le C_K.
 \end{aligned}
\]
Since $\E_a[Z]\ge1$, this proves the uniform bound
\[
 \frac{\E_a[Z|H|]}{\E_a[Z]}\le C_K.
\]

For an interval $J$ of length $|J|$ on the circle, a periodic primitive
of $\mathbf1_J-|J|/(2\pi)$ can be chosen with supremum norm at most
$|J|$. Smooth approximation in \eqref{eq:pair-weighted-ibp} gives
\[
 \frac{\E_a[Z\mathbf1_{\{x_n\bmod2\pi\in J\}}]}
      {\E_a[Z]}
 \le (C_K+(2\pi)^{-1})|J|.
\]
One can first use smooth majorants supported in an arbitrarily small
enlargement of $J$; thus this deduction does not assume absence of
atoms. The interval bound implies absolute continuity and a bounded
density. Finally multiply the conditional inequality by
$Y_a\E_a[Z]$ and average. Division by $\E Y_n$
gives the same bound for the relative phase under $\mu_{n,\mathbf q}$.
\end{proof}

\begin{proof}[Proof of \cref{lem:moving-terminal-ui}]
We bound the common phase integral by a sum of singular pair powers,
then integrate those powers against the densities from
\cref{prop:moving-pair-density}.
Write $b_j=(-q_j)_+$, $A=\sum_jb_j$, and $A_* =\sup_K A<2$.
The exponent $s$ will dominate the singularity of the integral, while
$sp<1$ will make its $p$th power integrable against a bounded density.
Choose $s$ and $p$ with
\[
 \max\{A_*-1,0\}<s<1,\qquad 1<p<1/s.
\]
We first claim the deterministic estimate
\begin{equation}\label{eq:terminal-pair-majorant}
 W_{\mathbf q}(\boldsymbol\varphi)
 \le C\left(1+\sum_{j<l}
                 \|\varphi_j-\varphi_l\|_{2\pi}^{-s}\right).
\end{equation}
Positive powers are uniformly bounded, so it is enough to consider
the factors with $b_j>0$. If there is at most one such factor, the common phase integral is
bounded because each $b_j$ stays below one.
Otherwise let $d$ be the minimum separation of their centers.
On the interval of radius $d/2$ about center $j$, the other distances
are at least $d/2$, and the integral is bounded by
\[
 C d^{-(A-b_j)}\int_0^{d/2}r^{-b_j}\,dr\le C d^{1-A}.
\]
On the complement, split into cells according to the nearest center.
If its distance is $r$, every other distance is at least $r$, giving
an upper bound $C\int_{d/2}^{\pi}r^{-A}\,dr$.
The sum of these bounds is at most $C(1+d^{-s})$, uniformly on $K$,
since $s>A_*-1$. This remains true when $A$ crosses one. As $d^{-s}$
is bounded by the sum of the corresponding pair powers,
\eqref{eq:terminal-pair-majorant} follows. It is interpreted with an
infinite right-hand side on a diagonal where two phases coincide.

By \cref{prop:moving-pair-density}, every pair distance under
$\mu_{n,\mathbf q}$ has bounded angular density when $n\delta\ge R$.
Hence
\[
 \int\|\varphi_j-\varphi_l\|_{2\pi}^{-sp}\,d\mu_{n,\mathbf q}
 \le C\int_0^\pi r^{-sp}\,dr<\infty.
\]
Raise \eqref{eq:terminal-pair-majorant} to the power $p$ and use the
finite sum inequality to prove \eqref{eq:moving-terminal-ui}.
The same estimate also holds under Haar measure.
\end{proof}

\begin{remark}
The restriction $\sum_j(-q_j)_+<2$ enters only in the integrability
of $W_{\mathbf q}$ against the pair density. It is automatic
for two integrable root exponents. For three or more points with
$\sum_j(-q_j)_+\ge2$, pair estimates alone do not control
all simultaneous coincidences of phases. \Cref{thm:full-signed-polynomial}
remains available throughout the full integrable domain of exponents.
\end{remark}

\begin{proof}[Proof of \cref{thm:main-root}]
For $m\ge2$ and fixed separation, \cref{thm:full-root-FH}
gives the formula and the displayed exponent on the full integrable
domain of exponents.
On $\mathcal D_m$, \cref{thm:signed-characteristic-FH} gives the
error depending only on $N\delta$, with its uniform constants;
\cref{cor:terminal-subset-rate} gives the sharper rate depending on
the individual exponents. On $\mathcal E_m$, \cref{thm:moving-signed-two} gives
\eqref{eq:intro-moving-signed}. For $m=1$, \cref{prop:onepoint} gives
the error $O(N^{-1})$. All constants are the Morris constants
identified in \cref{sec:onepoint}.
\end{proof}

\section{The microscopic limit and matching at large distances}
\label{sec:positive-merging}

This section supplies two inputs for the critical moments: a bound
for the moments with nonnegative exponents at arbitrary points,
including coincident ones, and the microscopic profile needed for the
finite part of the integer critical moment.
The matching of that profile to separated points uses the uniform
Fisher--Hartwig estimate of \cref{thm:signed-characteristic-FH}.
Local fusion expansions \cite{Fang2026,FangCritical2026,FangFourth2026} and
joint continuity of the pair density \cite{FangContinuity2026}
concern complementary properties of these microscopic limits.
The same microscopic limit appears in the supercritical regime
$pq^2>2\beta$ of \cite{AN2026}, through the growing integral recalled
in \cref{sec:extensions}; at the critical relation $pq^2=2\beta$ of
\cref{thm:main-critical}, the intermediate scales produce the leading
$N\log N$ term. Matching the microscopic profile with the asymptotic
for separated points identifies the finite term of order $N$.

\subsection{A global bound at arbitrary points}
\label{sec:global-bound}

The following bound controls the characteristic polynomial moments
with nonnegative exponents at every configuration, including
coincident points. It is uniform in the exponents. It is used for the
global bound in \cref{thm:positive-merging}, for the moment estimates
of \cref{sec:fractional-critical}, and again in
\cref{sec:critical-moments} for the transition through the critical
value. For nonnegative exponents, \cite[Proposition~3.5]{AN2026}
proves the same bound with the same induction; the estimate between
two polynomial degrees, obtained there from Propositions~3.2 to~3.4
of that paper, follows here from \cref{prop:unifiedblock}, with
constants uniform on compact sets of exponents.

In \cref{lem:moments-global-bound}, $k$ denotes the number of
prescribed points, so $\mathbf q$ and $\boldsymbol\vartheta$ have
$k$ coordinates. This notation anticipates the integer moment of order
$k$ in \cref{sec:critical-moments}. In the induction below,
$m$ denotes the current number of points, with $1\le m\le k$.

\begin{lemma}
\label{lem:moments-global-bound}
Fix $\beta>0$, an integer $k\ge1$, and $Q<\infty$. Write
$d_{ij}=|e^{i\vartheta_i}-e^{i\vartheta_j}|$ for the chord distance
between prescribed points. There is a
constant $C_{\beta,k,Q}$ such that, for $q_j\ge0$ with
$\sum_{j=1}^kq_j\le Q$, $N\ge1$, and arbitrary points,
\begin{equation}\label{eq:moments-global-bound}
 F_N(\mathbf q,\boldsymbol\vartheta)
 \le C_{\beta,k,Q} N^{\sum_jq_j^2/(2\beta)}
         \prod_{i<j}(d_{ij}\vee N^{-1})^{-q_iq_j/\beta}.
\end{equation}
Coincident points are allowed.
\end{lemma}

The proof uses induction on the number of points, with a bound on the
sum of the nonnegative exponents throughout. No critical moment
asymptotic is used here.

\begin{proof}[Proof of \cref{lem:moments-global-bound}]
We keep track of the dependence on the exponents, since it is needed
in \cref{sec:critical-window}, where $q$ varies through the critical
value.

For $m\le k$ put
\[
 \mathcal K_m(Q)=\{(q_1,\ldots,q_m):q_i\ge0,\ \sum_iq_i\le Q\},
 \qquad
 A_n(\mathbf q,\boldsymbol\vartheta)=\E Y_n(\mathbf q;\boldsymbol\vartheta),
 \qquad b_n(d)=\min\{n,d^{-1}\},
\]
where $b_n(0)=n$. These sets of exponents are compact and are preserved
when two exponents are added. We first prove the OPUC bound
\begin{equation}\label{eq:compact-positive-opuc-bound}
 A_n(\mathbf q,\boldsymbol\vartheta)
 \le C_{\beta,k,Q}n^{a(\mathbf q)}
                 \prod_{i<j}b_n(d_{ij})^{q_iq_j/\beta}
 \qquad(n\ge1,\ \mathbf q\in\mathcal K_m(Q)).
\end{equation}

Take the zero Fourier mode in \cref{prop:unifiedblock}. There are
$C_0,R<\infty$, depending only on $\beta,k,Q$, such that, for
$0<\delta\le1$, $\boldsymbol\vartheta\in\mathcal S_\delta$, and
$n\ge j$ with $j\delta\ge R$,
\begin{equation}\label{eq:compact-positive-index-bound}
 A_n(\mathbf q,\boldsymbol\vartheta)
 \le C_0(n/j)^{a(\mathbf q)}
                          A_j(\mathbf q,\boldsymbol\vartheta).
\end{equation}
The iteration is as follows. Choose $R\ge1$ greater than the finitely
many degree thresholds for the compact sets $\mathcal K_m(Q)$.
For $j_0=j$, $j_{r+1}=\min(2j_r,n)$, the block estimate and
$a(\mathbf q)\ge0$ give
\[
 A_{j_{r+1}}\le (j_{r+1}/j_r)^{a(\mathbf q)}
                    \left(1+\frac{C_1}{j_r\delta}\right)A_{j_r}.
\]
Multiply these inequalities until $j_r=n$. The power factors
multiply to $(n/j)^{a(\mathbf q)}$, whereas
\[
 \prod_r\left(1+\frac{C_1}{j_r\delta}\right)
 \le\exp\left(\frac{C_1}{j\delta}\sum_{r\ge0}2^{-r}\right)
 \le e^{2C_1/R}.
\]
This proves \eqref{eq:compact-positive-index-bound} uniformly on every
$\mathcal K_m(Q)$.

We now prove \eqref{eq:compact-positive-opuc-bound} by induction on
$m$. For $m=1$, it is the locally uniform one point estimate
\eqref{eq:opuc-onepoint-asymptotic}, with its finitely many initial
indices included in the constant. Suppose it holds for $m-1$.
Relabel the points so that $d=d_{m-1,m}$ is a smallest
pair distance. If $d>0$, put
\[
 \delta=\min(1,d),\qquad
 j=\min\{n,\lceil R/\delta\rceil\}.
\]
If $d=0$, put $j=n$. In either case
\eqref{eq:compact-positive-index-bound} applies when $j<n$; when
$j=n$ the same inequality is immediate. For $u,v\ge0$ and $x,y>0$,
\[
 x^uy^v\le\max\{x^{u+v},y^{u+v}\}
            \le x^{u+v}+y^{u+v}.
\]
Use this with $u=q_{m-1}$, $v=q_m$ and the two corresponding OPUC
absolute values at degree $j$. Each resulting term has $m-1$
points and sum of the exponents at most $Q$. The induction
hypothesis therefore applies with the same bound on their sum.

For $i\le m-2$, the choice of the closest pair gives
\[
 d_{i,m-1}\le2d_{i,m},\qquad d_{i,m}\le2d_{i,m-1},
 \qquad
 \frac12\le\frac{b_j(d_{i,m-1})}{b_j(d_{i,m})}\le2.
\]
If one of these distances is zero, both are zero and the last ratio
is one. For example, in the term in which the two exponents are
placed at $\vartheta_{m-1}$,
\begin{align*}
 b_j(d_{i,m-1})^{q_i(q_{m-1}+q_m)/\beta}
 &\le 2^{q_iq_m/\beta}
       b_j(d_{i,m-1})^{q_iq_{m-1}/\beta}
       b_j(d_{i,m})^{q_iq_m/\beta}.
\end{align*}
The other term satisfies the same bound with the indices exchanged.
All factors of two have total exponent at most $Q^2/\beta$.
Since adding the two exponents changes $a(\mathbf q)$ by
$q_{m-1}q_m/\beta$, we obtain
\[
 A_n\le C(n/j)^{a(\mathbf q)}
           j^{a(\mathbf q)+q_{m-1}q_m/\beta}
     \prod_{\substack{i<l\\(i,l)\ne(m-1,m)}}
                      b_j(d_{il})^{q_iq_l/\beta}.
\]
The choice of $j$ also gives
\[
 j\le 2(R+1)b_n(d),\qquad b_j(d_{il})\le b_n(d_{il}).
\]
Indeed $d\le2$, so $\lceil R/\min(1,d)\rceil
\le2(R+1)d^{-1}$ when $d>0$; taking the minimum with $n$ proves
the first inequality. For $d=0$ it is immediate.
These bounds prove \eqref{eq:compact-positive-opuc-bound}. The
constant stays uniform because the number of induction steps is at
most $k-1$ and every sum of exponents stays at most $Q$.

For $N\ge2$, the terminal representation gives
$|P_N(e^{i\vartheta})|\le2|\Phi_{N-1}^*(e^{i\vartheta})|$.
Thus \eqref{eq:compact-positive-opuc-bound} at $n=N-1$, followed by
$b_{N-1}(d)\le b_N(d)$ and $(N-1)^{a(\mathbf q)}\le N^{a(\mathbf q)}$,
proves \eqref{eq:moments-global-bound}. Each factor in \eqref{eq:terminal-process} has modulus at most
two, so their powers contribute at most $2^Q$. The case $N=1$ follows from $|P_1|\le2$ and $d_{ij}\le2$,
after increasing the same uniform constant.
\end{proof}

\subsection{The microscopic limit}
\label{sec:microscopic-limit}

The circular Jacobi limit gives a merging formula for arbitrary
nonnegative exponents. We state the result in the coordinate
$x=N(\vartheta-\vartheta_0)$. The unweighted rescaled circular ensemble
has intensity $1/(2\pi)$ in this coordinate, as does its
$\mathsf{Sine}_\beta$ limit.
The stochastic zeta input is the analytic circular Jacobi limit of
Li and Valk\'o \cite{LiValko2022}, in the form recalled in
\cite[Propositions~2.7 and~2.8]{AN2026}. The analogous entire functions
for the CUE and the $\mathsf{Sine}_\beta$ process are studied in
\cite{CNN2017,ValkoVirag2022}. The additional conclusion supplied
by our separation uniform asymptotic is the quantitative matching in
\eqref{eq:merging-matching} below.

For every $\beta>0$ and $\tau>0$, let $\mathcal{HP}_{\beta,\tau}$ denote
the Hua--Pickrell point process obtained as the microscopic limit, in
the coordinate $x=N\theta$ with $\theta\in(-\pi,\pi]$, of the circular
Jacobi ensemble whose density relative to the circular beta ensemble
is $|P_N(1)|^{2\tau}/M_N(2\tau)$. This convergence follows from the
coupling in \cite[Proposition~2.7]{AN2026}, based on
\cite{LiValko2022}. We use the stochastic zeta function
normalized by the symmetric principal value product
\[
 \xi^{\beta,\tau}(z)
 =\lim_{A\to\infty}
   \prod_{\substack{u\in\mathcal{HP}_{\beta,\tau}\\|u|<A}}
       (1-z/u),\qquad \xi^{\beta,\tau}(0)=1.
\]
The limit exists almost surely locally uniformly in $z\in\C$.
More precisely, \cite[Propositions~2.7 and~2.8, and Remark~2.9]{AN2026} give a coupling
under which $e^{-iz/2}P_N(e^{iz/N})/P_N(1)$ converges to this function
locally uniformly almost surely. This specifies the microscopic
coordinate and the exponential normalization of the entire function.
\begin{theorem}
\label{thm:positive-merging}
Fix $\beta>0$ and $q_1,\ldots,q_m\geq0$, with
total exponent $Q=\sum_jq_j>0$. Zero exponents can be omitted.
Choose an anchor index $a$ with
$q_a>0$, and define
\begin{equation}\label{eq:merging-function}
 \mathcal H_{\beta,\mathbf q}(\mathbf x)
 =C_\beta(Q)\,
   \E\prod_{j\ne a}
       |\xi^{\beta,Q/2}(x_j-x_a)|^{q_j},
 \qquad \mathbf x\in\mathbb R^m.
\end{equation}
The expression \eqref{eq:merging-function} is independent of the choice
of $a$. The function is finite,
strictly positive, continuous, and invariant under simultaneous
translation of the coordinates. For every compact
$K\subset\mathbb R^m$,
\begin{equation}\label{eq:merging-limit}
 \sup_{\mathbf x\in K}\left|
 N^{-Q^2/(2\beta)}
 F_N\left(\mathbf q,
     (\vartheta_0+x_j/N)_{j=1}^m\right)
 -\mathcal H_{\beta,\mathbf q}(\mathbf x)\right|
 \longrightarrow0.
\end{equation}
The convergence is uniform in $\vartheta_0$.

If $I_1,\ldots,I_r$ partition $\{1,\ldots,m\}$ and $x_j=y_\ell$
for $j\in I_\ell$, then, with $Q_\ell=\sum_{j\in I_\ell}q_j$,
\begin{equation}\label{eq:merging-collision}
 \mathcal H_{\beta,\mathbf q}(\mathbf x)
 =\mathcal H_{\beta,(Q_1,\ldots,Q_r)}(y_1,\ldots,y_r).
\end{equation}
In particular, $\mathcal H_{\beta,\mathbf q}(x,\ldots,x)=C_\beta(Q)$.
There is a constant depending only on $\beta$ and $\mathbf q$ such that
\begin{equation}\label{eq:merging-global-bound}
 0<\mathcal H_{\beta,\mathbf q}(\mathbf x)
 \leq C\prod_{i<j}(1+|x_i-x_j|)^{-q_iq_j/\beta}.
\end{equation}

If at least two exponents are positive, let
$d(\mathbf x)=\min_{i<j:q_iq_j>0}|x_i-x_j|$.
For every $0<\kappa<\min\{2/\beta,1/2\}$, as $d(\mathbf x)\to\infty$,
\begin{equation}\label{eq:merging-matching}
 \mathcal H_{\beta,\mathbf q}(\mathbf x)
 =\prod_{j=1}^m C_\beta(q_j)
    \prod_{i<j}|x_i-x_j|^{-q_iq_j/\beta}
       \left(1+O_{\beta,\mathbf q}(d(\mathbf x)^{-\kappa})\right).
\end{equation}
The error is uniform over all such real configurations.
\end{theorem}

\begin{proof}[Proof of \cref{thm:positive-merging}]
We first express the root moment under a circular Jacobi law and
prove uniform integrability of the resulting products. This yields the
microscopic limit and its continuity. A periodic moment bound then
gives the global estimate, and \cref{thm:signed-characteristic-FH}
gives the matching at large distances.
Under the circular Jacobi law with parameter $Q/2$, define the
normalized analytic function directly by
\[
 f_N(z)=e^{-iz/2}\frac{P_N(e^{iz/N})}{P_N(1)}.
\]
The change of measure defining this law has density
$|P_N(1)|^Q/M_N(Q)$ with respect to the circular beta ensemble.
Rotating the prescribed point with index $a$ to one therefore gives the exact identity
\begin{equation}\label{eq:merging-change-measure}
 F_N\left(\mathbf q,(\vartheta_0+x_j/N)_j\right)
 =M_N(Q)\,
   \E_{\mathrm{CJ}(Q/2)}
      \prod_{j\ne a}|f_N(x_j-x_a)|^{q_j}.
\end{equation}
For real arguments the factor $e^{-iz/2}$ has modulus one.

The positive anchor exponent leaves enough room for a moment of order
greater than one. The following bound supplies the required uniform
integrability. For nonnegative $r_1,\ldots,r_\ell$ with
$S=\sum_jr_j\leq Q$, and arbitrary real $y_j$, the same change of measure
and generalized H\"older inequality give
\begin{align}\label{eq:cj-holder-bound}
 \E_{\mathrm{CJ}(Q/2)}
     \prod_{j=1}^\ell|f_N(y_j)|^{r_j}
 &=\frac{\E\left[|P_N(1)|^{Q-S}
          \prod_j|P_N(e^{iy_j/N})|^{r_j}\right]}{M_N(Q)}
 \notag\\
 &\leq\frac{M_N(Q)^{(Q-S)/Q}
                 \prod_jM_N(Q)^{r_j/Q}}{M_N(Q)}=1.
\end{align}
Factors with exponent zero are omitted. The equality of the one point
moments at different angles follows from rotational invariance.

By \cite[Propositions~2.7 and~2.8]{AN2026}, these functions can be
coupled so that $f_N\to\xi^{\beta,Q/2}$ almost surely, uniformly on
compact subsets of $\mathbb C$. Put $S_a=Q-q_a<Q$. If $S_a>0$, choose
$p>1$ with $pS_a\leq Q$, which is possible because $S_a<Q$.
Bound \eqref{eq:cj-holder-bound}, applied with
$r_j=pq_j$ for $j\ne a$, gives a uniform $L^p$ bound for the products
in \eqref{eq:merging-change-measure}; if $S_a=0$, these products equal
one. Fatou's lemma gives the same moment bound for the limiting products.

This yields locally uniform convergence of their
expectations without requiring a bound on a random supremum. Suppose
$\mathbf x_N\to\mathbf x$ in a fixed compact set. Almost sure locally
uniform convergence of $f_N$ gives convergence of the corresponding
products, and their uniform $L^p$ bound gives convergence in $L^1$.
The analogous argument for $\xi^{\beta,Q/2}$ shows that the limiting
expectation is continuous in $\mathbf x$. If convergence of expectations
were not uniform on a compact set, a sequence witnessing failure would
have a convergent subsequence, contradicting these two observations.
Combining this conclusion with
the one point limit
$M_N(Q)N^{-Q^2/(2\beta)}\to C_\beta(Q)$ from \cref{prop:onepoint} proves
\eqref{eq:merging-limit}. Rotation invariance makes every expression
independent of $\vartheta_0$.

The left side of \eqref{eq:merging-change-measure} is independent of
the anchor. Passing to the limit therefore proves independence of $a$.
At coincident points the powers on that same left side add, so
\eqref{eq:merging-limit} also gives \eqref{eq:merging-collision}. The Hua--Pickrell process has no atom at
a deterministic real point: in its phase representation
\cite[Definition~1.5]{AN2026}, the final random phase is independent of the
driving Brownian motion and has an absolutely continuous law.
Consequently $\xi^{\beta,Q/2}(y)\ne0$ almost surely for each fixed
$y$, and the expectation in \eqref{eq:merging-function} is positive.

For the global bound, \cref{lem:moments-global-bound} gives
\[
 F_N(\mathbf q,\boldsymbol\theta)
 \leq C N^{\sum_jq_j^2/(2\beta)}
  \prod_{i<j}
  \min\left\{N,|e^{i\theta_i}-e^{i\theta_j}|^{-1}\right\}^{q_iq_j/\beta}.
\]
For fixed $\mathbf x$, substitute $\theta_j=x_j/N$, divide by
$N^{Q^2/(2\beta)}$, and let $N\to\infty$. Since
$N|e^{ix_i/N}-e^{ix_j/N}|\to|x_i-x_j|$, this yields
\[
 \mathcal H_{\beta,\mathbf q}(\mathbf x)
 \leq C\prod_{i<j}\min\{1,|x_i-x_j|^{-1}\}^{q_iq_j/\beta},
\]
which is equivalent to \eqref{eq:merging-global-bound} up to a constant.
This also explains why the distances in an estimate at finite $N$ must be
periodic; only after passage to the limit are they ordinary real distances.

Finally, omit the zero exponents, fix a real configuration with
$d(\mathbf x)$ sufficiently large, and apply
\cref{thm:signed-characteristic-FH} at
$\theta_j=x_j/N$. Its separation parameter $\delta_N$ satisfies
$N\delta_N\to d(\mathbf x)$. After division by $N^{Q^2/(2\beta)}$,
the powers of $N$ cancel because
\[
 \frac1{2\beta}\sum_jq_j^2+
       \frac1\beta\sum_{i<j}q_iq_j=\frac{Q^2}{2\beta}.
\]
The limiting identity is \eqref{eq:merging-matching}, with the same
uniform error constant. Since that constant is independent of the
angular separation and of the configuration, the conclusion is
uniform over all $\mathbf x$ with $d(\mathbf x)$ sufficiently large.
\end{proof}

\begin{remark}
The merging representation uses the analytic circular Jacobi limit of
Li and Valk\'o \cite{LiValko2022}, recalled in \cite{AN2026}, and the elementary bound
\eqref{eq:cj-holder-bound}. It should be distinguished from the new
matching estimate \eqref{eq:merging-matching}, which uses
\cref{thm:signed-characteristic-FH}. The latter treats configurations
in which every distance between points with positive exponents tends to infinity.
\end{remark}

\section{Critical moments of noninteger order}
\label{sec:fractional-critical}

At $q=\sqrt{2\beta/p}$, the $p$th moment of the limiting chaos
mass is infinite, whereas the corresponding polynomial moment is
finite at every index. We shall determine its logarithmic growth.
We first normalize the masses so that the coefficient is independent of $\beta$.

Fix $\beta>0$ and $p>1$, and put
\[
 q=\sqrt{\frac{2\beta}{p}},\qquad
 d_p=\Gamma(1-1/p)^{-p},\qquad A_p=(p-1)d_p.
\]
Throughout this section $\sigma(d\theta)=d\theta/(2\pi)$.
The normalized OPUC density, its random measure and its total mass are
\[
 D_n(\theta)=\frac{|\Phi_n^*(e^{i\theta})|^q}{\nu_n(q)},\qquad
 \mu_n(d\theta)=D_n(\theta)\sigma(d\theta),\qquad U_n=\mu_n(\mathbb T).
\]
For the characteristic polynomial $P_N$ with $N$ eigenvalues, we use
$\chi_N$ for the normalized density and $W_N$ for its total mass:
\[
 \chi_N(\theta)=\frac{|P_N(e^{i\theta})|^q}{M_N(q)},\qquad
 W_N=\int_{\mathbb T}\chi_N(\theta)\sigma(d\theta).
\]
We suppress the fixed parameters $p$ and $\beta$ in this notation.
Unless stated otherwise, constants may depend on these parameters;
dependencies on additional parameters remain indicated when needed.

\begin{proposition}\label{thm:fractional-critical}
For every $\beta>0$ and every noninteger real $p>1$,
\begin{equation}\label{eq:fractional-critical-normalized}
 \mathbb E U_n^p\sim A_p\log n,
 \qquad
 \mathbb E W_N^p\sim A_p\log N.
\end{equation}
\end{proposition}

By \eqref{eq:onepointasymp}, $M_N(q)^p\sim C_\beta(q)^pN$.
Thus \cref{thm:fractional-critical} proves the unweighted
fractional assertion of \cref{thm:main-critical}. The corresponding
orthogonal polynomial constant is $A_pD_\beta(q)^p$.

We first treat the OPUC mass $U_n$. Its upper bound combines the
conditional chaos identity \eqref{eq:fractional-chaos-martingale}
with a convex truncation of the exact mass law
\eqref{eq:fractional-exact-mass}. We state the bound now and prove it
after the conditional comparisons and annular estimates needed for
this truncation.

\begin{proposition}\label{prop:all-real-critical-upper}
For every $\beta>0$ and every real $p>1$,
\begin{equation}\label{eq:all-real-critical-upper}
 \limsup_{n\to\infty}\frac{\mathbb E U_n^p}{\log n}
 \le\frac{p-1}{\Gamma(1-1/p)^p},
 \qquad \mathbb E W_N^p=O_{p,\beta}(\log(2+N)).
\end{equation}
\end{proposition}

For the matching lower bound, a change of measure represents an
integer part of the moment by $k=\lfloor p\rfloor$ marked points. The
remaining power $p-k$ lies in $(0,1)$, so conditional concavity
applies to it. Three inputs connect the two bounds. The conditional
product estimates of \cref{lem:fractional-conditional-product}
compare polynomial masses at finite indices with the limiting chaos
under the change of measure. The exact circular mass law
\eqref{eq:fractional-exact-mass} determines the tail constant from
which the coefficient $A_p$ is recovered. Finally,
\cref{lem:fractional-rooted-cutoff,lem:critical-chaos-collision} show
that removing small neighborhoods of the marked points, as the
conditional comparisons require, does not change that coefficient;
for several marked points, these estimates also control
configurations in which two of the points are close. Averaging over
the independent angle in \eqref{eq:terminal-process} then gives the
lower bound for $W_N$ by Jensen's inequality, and the conditional
estimate \eqref{eq:critical-terminal-multipoint}, proved for any
fixed number of marked points, gives the upper bound.

The section is organized by the range of $p$.
\Cref{sec:fc-one-marked,sec:fc-annular,sec:fc-cutoff,sec:fc-averaging}
treat $1<p<2$ with one marked point; there we write $t=p-1\in(0,1)$
for the power remaining after marking that point.
\Cref{sec:fc-upper-all} proves the upper bound of
\cref{prop:all-real-critical-upper} for every real $p>1$.
\Cref{sec:fc-several,sec:fc-distinguished,sec:fc-nearby,sec:fc-upper-transfer}
treat noninteger $p>2$ with $k$ marked points, and
\cref{sec:fractional-critical-weights} establishes the comparison of
disjoint arcs needed for spatial weights at every real order $p>1$. The chaos identities, the conditional
comparisons and
\cref{lem:critical-terminal-multipoint} remain valid for every $p>1$.

\subsection{The change of measure and conditional comparison}
\label{sec:fc-one-marked}

By the conditional chaos identity of Lambert and Najnudel
\cite[Proposition~4]{LN2026}, there is a subcritical circular Gaussian
multiplicative chaos $\mu$ on the coefficient probability space such that
\begin{equation}\label{eq:fractional-chaos-martingale}
 \mu_n(f)=\mathbb E[\mu(f)\mid\mathcal F_n].
\end{equation}
Its Gaussian covariance, with the exponent included in the field, is
$-(2/p)\log|e^{i\theta}-e^{i\phi}|$. In particular,
\begin{equation}\label{eq:fractional-exact-mass}
 U:=\mu(\mathbb T)\ \stackrel{\mathrm d}{=}
 \frac{\mathcal E^{-1/p}}{\Gamma(1-1/p)},
 \qquad \mathcal E\sim\operatorname{Exp}(1).
\end{equation}
This is the circular chaos mass law of \cite[Theorem~1.1]{Remy2020},
with its parameter $\gamma=2/\sqrt p$. In that normalization the
Gaussian covariance is $-2\log|e^{i\theta}-e^{i\phi}|$. Multiplying
the field by $\gamma/2$ gives exactly the covariance
$-(2/p)\log|e^{i\theta}-e^{i\phi}|$.
Both total masses use the probability measure $\sigma$.

For each $\theta$, the recursion \eqref{eq:szego-recursion} and the
normalization \eqref{eq:opuc-onepoint-exact} make $D_n(\theta)$ a
positive martingale of mean one.
Let $\mathbb Q$ be the probability measure on the coefficient sequence
whose restriction to $\mathcal F_n$ has density $D_n(0)$ relative to
$\mathbb P$. When working under this possibly singular change of
measure, define the limiting random measure by the weak limit of
$\mu_n$, rather than by an arbitrary version of the Gaussian field
defined almost surely under $\mathbb P$. This limit exists under $\mathbb Q$, as shown
below, and has the rooted law
\begin{equation}\label{eq:fractional-rooted-law}
 \mu^0(d\theta)\ \overset{\rm law}=\ |1-e^{i\theta}|^{-2/p}\widetilde\mu(d\theta),
\end{equation}
where $\widetilde\mu$ is an ordinary, unrooted circular chaos on a
fresh probability space. The law on the left is taken under $\mathbb Q$.
Moreover, $\mu_n$ converges weakly to this limiting measure
almost surely under $\mathbb Q$.

We construct this rooted law without choosing a conditional version
at an exceptional marked angle. On coefficient sequences put
\[
 (\rho_\theta\omega)_k=e^{i(k+1)\theta}\alpha_k.
\]
The recursion \eqref{eq:szego-recursion} gives
$\Phi_n^*(z;\rho_\theta\omega)=\Phi_n^*(e^{i\theta}z;\omega)$.
Consequently
$D_n(0;\rho_\theta\omega)=D_n(\theta;\omega)$, and the
coefficient law is invariant under $\rho_\theta$. Write $\mu_\omega$
for the realization of the limiting chaos on the coefficient sequence
$\omega$. We mark a point with the measure $\mu_\omega$ and rotate
that point to zero. This defines
\[
 \mathbb Q(A)=\mathbb E\int_\T
             \mathbf1_A(\rho_\theta\omega)\,\mu_\omega(d\theta).
\]
This is a probability measure because $\mathbb E U=1$. For
$A\in\mathcal F_n$, the conditional chaos identity
\eqref{eq:fractional-chaos-martingale} and rotation
invariance give
\begin{align*}
 \mathbb Q(A)
 &=\int_\T\mathbb E[
       \mathbf1_A(\rho_\theta\omega)D_n(\theta;\omega)]
                                                        \,\sigma(d\theta)\\
 &=\mathbb E[\mathbf1_A(\omega)D_n(0;\omega)].
\end{align*}
It therefore has exactly the restrictions specified above.
The exceptional set for almost sure weak convergence of $\mu_n$ has
zero measure under $\mathbb P(d\omega)\mu_\omega(d\theta)$, since
it has probability zero under $\mathbb P$ and $U<\infty$ almost surely.
Rotation covariance of the entire sequence $(\mu_n)$ preserves
convergence, so it holds under $\mathbb Q$ as well. Finally the
Gaussian Cameron--Martin identity for regularized chaos, followed by
its subcritical limit, gives \eqref{eq:fractional-rooted-law}.
The same argument, applied to a function of the total mass, gives
\begin{equation}\label{eq:fractional-size-bias}
 \mathbb E_{\mathbb Q}f(\mu^0(\mathbb T))
       =\mathbb E[U f(U)]
\end{equation}
for every nonnegative Borel function $f$.

Rotation invariance also gives, for every real $r>1$,
\begin{equation}\label{eq:fractional-rooting-finite}
 \mathbb E U_n^r
 =\int_{\mathbb T}\mathbb E[D_n(\theta)U_n^{r-1}]\sigma(d\theta)
 =\mathbb E_{\mathbb Q}U_n^{r-1}.
\end{equation}

The identity \eqref{eq:fractional-rooting-finite} reduces the moment
of $U_n$ to a lower power under $\mathbb Q$. To compare that power
with the limiting mass outside a neighborhood of $0$, the following
conditional form of the positive moment comparison is needed. Set
$\delta(\theta)=\min(1,|1-e^{i\theta}|)$.

\begin{lemma}\label{lem:fractional-conditional-product}
For every fixed integer $j\ge1$ there are constants $C_j,R_j$ such that
\begin{equation}\label{eq:critical-conditional-products}
 \left|
 \frac{\mathbb E[\prod_{i=1}^jD_M(\theta_i)\mid\mathcal F_n]}
      {\prod_{i=1}^jD_n(\theta_i)}-1\right|
 \le\frac{C_j}{n\delta},\qquad M\ge n,
\end{equation}
provided $n\delta\ge R_j$, where $\delta$ is the minimum pairwise
chord separation, clipped above by $1$; for $j=1$ the equality is exact.
The estimate is uniform in the coefficients revealed by time $n$.
\end{lemma}

\begin{proof}[Proof of \cref{lem:fractional-conditional-product}]
We first compare two indices in a doubling block, where
\cref{prop:unifiedblock} applies. We then sum the relative errors over
successive blocks.
For $n\le M\le2n$, apply \cref{prop:unifiedblock} with exponents
$(q,\ldots,q)$ and zero Fourier mode. Its leading factor is
$(M/n)^{j/p}$. The gamma product \eqref{eq:opuc-onepoint-exact} gives
\[
 \left(\frac{\nu_M(q)}{\nu_n(q)}\right)^j
       =(M/n)^{j/p}(1+O_j(n^{-1}))
\]
uniformly in this block, so division by the normalizers cancels
that leading factor.
The error is pointwise relative to the initial positive product.
Iterating over successive doubling blocks is legitimate by positivity
and conditional expectation; the accumulated relative error is bounded
by a constant times $\sum_{l\ge0}(2^ln\delta)^{-1}$.
Increase $R_j$ so that all lower factors remain positive.
This proves \eqref{eq:critical-conditional-products}. It imposes no condition
$j<p$, and involves no limiting $j$th GMC moment.
\end{proof}

For two points we write $R_0=R_2$ and $C=C_2$, increasing $R_0$
when necessary.

\begin{lemma}\label{lem:fractional-rooted-conditional}
Let $A$ be an arc or a finite union of arcs at chord distance at least
$\delta>0$ from $0$. If $n\delta\ge R_0$, then
\begin{equation}\label{eq:fractional-rooted-conditional}
 (1-C/(n\delta))\mu_n(A)
 \le\mathbb E_{\mathbb Q}[\mu^0(A)\mid\mathcal F_n]
 \le(1+C/(n\delta))\mu_n(A).
\end{equation}
The same inequalities hold with $\mu_M(A)$ in place of $\mu^0(A)$,
for every $M\ge n$.
\end{lemma}

\begin{proof}[Proof of \cref{lem:fractional-rooted-conditional}]
We prove the comparison first for $\mu_M(A)$, then pass to the
limiting rooted mass in $L^1(\mathbb Q)$.
The assertion at finite $M$ follows by the change of measure formula
\[
 \mathbb E_{\mathbb Q}[D_M(\theta)\mid\mathcal F_n]
  =\frac{\mathbb E[D_M(0)D_M(\theta)\mid\mathcal F_n]}{D_n(0)},
\]
and \cref{lem:fractional-conditional-product}, followed by integration.
For fixed $A$, almost sure weak convergence and absence of atoms at its
endpoints give $\mu_M(A)\to\mu^0(A)$ almost surely under $\mathbb Q$.
Since $A$ is separated from the root, the relative two point
asymptotic of \cref{thm:opuc-relative} gives
\[
 \sup_{\theta\in A}\left|
   \mathbb E[D_M(0)D_M(\theta)]-|1-e^{i\theta}|^{-2/p}
 \right|\longrightarrow0.
\]
Consequently,
\[
 \mathbb E_{\mathbb Q}\mu_M(A)
 =\int_A\mathbb E[D_M(0)D_M(\theta)]\sigma(d\theta)
 \longrightarrow\int_A|1-e^{i\theta}|^{-2/p}\sigma(d\theta)
 =\mathbb E_{\mathbb Q}\mu^0(A).
\]
These nonnegative random variables converge almost surely and their
expectations converge to the expectation of the limit. Scheff\'e's
lemma therefore gives convergence in $L^1(\mathbb Q)$. Conditional expectation
is an $L^1$ contraction. The conditional expectations therefore
converge in $L^1$, and the inequalities pass to the limit along an
almost surely convergent subsequence. We may apply the resulting assertion to any deterministic arc
chosen separately for each $n$.
\end{proof}

\subsection{Upper bounds from annular estimates}
\label{sec:fc-annular}

\begin{lemma}\label{lem:fractional-supercritical-bound}
For $p\le r<2$, set
$\xi(r)=(r-1)(r/p-1)$. Then
\begin{equation}\label{eq:fractional-supercritical-bound}
 \mathbb E U_n^p\le C\log(2+n),\qquad
 \mathbb E U_n^r\le C_r(n+1)^{\xi(r)}\quad(p<r<2).
\end{equation}
The same upper bounds hold for $W_N$, with $n+1$ replaced by $N$.
\end{lemma}

\begin{proof}[Proof of \cref{lem:fractional-supercritical-bound}]
We estimate each annulus at the polynomial index corresponding to
its distance from the root. Subadditivity of the remaining power
then reduces the total moment to a sum over scales.
For this proof put $s=r-1\in(0,1)$. Partition the circle into an inner arc of radius
$O(R_0/n)$ and $O(\log(2+n))$ dyadic annuli $A$ whose lengths and
distances from zero are comparable to a number $a$. For each annulus
let $\delta_A$ be its minimum clipped chord distance from zero and
choose $m\asymp R_0/a$, increasing the fixed implicit constant so that
$m\delta_A\ge R_0$, and combine the bounded number of innermost pieces
until $m\le n$. The upper comparison at finite indices in
\cref{lem:fractional-rooted-conditional}, with initial index $m$, gives
$\mathbb E_{\mathbb Q}[\mu_n(A)\mid\mathcal F_m]\le C\mu_m(A)$.
Conditional concavity therefore yields
\[
 \mathbb E_{\mathbb Q}\mu_n(A)^s
 \le C_r\mathbb E_{\mathbb Q}\mu_m(A)^s.
\]
H\"older's inequality with conjugate exponents $r$ and $r/(r-1)$,
followed by Jensen's inequality in the spatial integral, gives
\begin{align*}
 \mathbb E_{\mathbb Q}\mu_m(A)^{r-1}
 &=\mathbb E[D_m(0)\mu_m(A)^{r-1}]\\
 &\le(\mathbb E D_m(0)^r)^{1/r}
       (\mathbb E\mu_m(A)^r)^{(r-1)/r}\\
 &\le\sigma(A)^{r-1}\mathbb E D_m(0)^r.
\end{align*}
Here rotation invariance makes $\mathbb E D_m(\theta)^r$ independent
of $\theta$. The one point product \eqref{eq:opuc-onepoint-exact} implies
\[
 \mathbb E D_m(0)^r=\frac{\nu_m(rq)}{\nu_m(q)^r}
 \le C_r(m+1)^{(r^2-r)/p}.
\]
Since $\sigma(A)\asymp m^{-1}$ up to fixed constants, the annular
bound is $C_r(m+1)^{\xi(r)}$. On the inner arc the same H\"older and
Jensen argument applies directly at index $n$ and gives
$C_r(n+1)^{\xi(r)}$.

Now $(\sum_j x_j)^s\le\sum_jx_j^s$ for $x_j\ge0$.
Use \eqref{eq:fractional-rooting-finite} and sum over the annuli.
When $r=p$, every annulus contributes a bounded amount. When $r>p$,
$\xi(r)>0$, and the dyadic sum is bounded by a constant multiple of
$(n+1)^{\xi(r)}$. This proves the OPUC bounds.

For $n=N-1$, averaging \eqref{eq:terminal-process} over its
independent uniform angle gives
\begin{equation}\label{eq:fractional-terminal-bounds}
 \mathbb E[W_N\mid\mathcal F_n]=U_n,
 \qquad 0\le W_N\le\frac{2^q}{c(q)}U_n,
 \quad c(q)=\int_{\mathbb T}|1-e^{iu}|^q\sigma(du).
\end{equation}
The characteristic polynomial upper bounds follow.
\end{proof}

The estimate in \cref{lem:fractional-supercritical-bound} gives
$O(\log n)$ at the critical order with an unspecified constant.
We now use the martingale identity \eqref{eq:fractional-chaos-martingale}
to identify the sharp upper
coefficient $A_p$. For $T>0$ let $F_T$ be the convex function equal to $x^p$ on
$[0,T]$ and equal to its tangent at $T$ on $[T,\infty)$. Thus
\[
 F_T(x)=pT^{p-1}x-(p-1)T^p\quad(x>T).
\]
The exact law \eqref{eq:fractional-exact-mass} gives
\begin{equation}\label{eq:fractional-convex-cutoff}
 \mathbb E F_T(U)=p d_p\log T+O(1),\qquad T\to\infty.
\end{equation}
Indeed $\mathbb P(U>x)=1-\exp(-d_px^{-p})$; integrating this tail
against $F_T'$ gives the displayed logarithm, with a bounded remainder
from both the integrable tail correction and the linear part of $F_T$.
By \eqref{eq:fractional-chaos-martingale} and the conditional Jensen
inequality, $\mathbb E F_T(U_n)\le\mathbb E F_T(U)$.
For $r>p$,
\[
 U_n^p\le F_T(U_n)+T^{p-r}U_n^r.
\]
We choose the truncation level so that the moment of order $r$ makes
a negligible contribution. More precisely, take $a>(p-1)/p$,
$T=n^a$, and $p<r<2$ sufficiently close to $p$ that $(r-1)/p<a$.
Then
\[
 \xi(r)+a(p-r)=(r-p)((r-1)/p-a)<0.
\]
Using \cref{lem:fractional-supercritical-bound}, dividing by $\log n$,
and letting $a\downarrow(p-1)/p$, we obtain
\begin{equation}\label{eq:fractional-opuc-limsup}
 \limsup_{n\to\infty}\frac{\mathbb E U_n^p}{\log n}\le A_p.
\end{equation}
\subsection{The mass outside a neighborhood of the marked point}
\label{sec:fc-cutoff}

\begin{lemma}\label{lem:fractional-rooted-cutoff}
Let $1<p<2$. Let $V=\mu^0(\mathbb T)$ and
$X_\epsilon=\mu^0(\{\theta:\epsilon\le|\theta|\le\pi\})$, where
the angular representative belongs to $[-\pi,\pi]$. Then
\begin{equation}\label{eq:fractional-rooted-cutoff-lower}
 \liminf_{\epsilon\downarrow0}
 \frac{\mathbb E_{\mathbb Q}X_\epsilon^{p-1}}{\log(1/\epsilon)}
 \ge A_p.
\end{equation}
\end{lemma}

\begin{proof}[Proof of \cref{lem:fractional-rooted-cutoff}]
We compare the tail of the mass outside the removed arc with the
tail of the full rooted mass. A moment of order below $t=p-1$ will
show that the removed mass is negligible over a growing range of
tail levels.
The size bias identity \eqref{eq:fractional-size-bias} and the
unrooted mass law \eqref{eq:fractional-exact-mass} imply
\begin{equation}\label{eq:fractional-rooted-tail}
 \mathbb Q(V>x)=\mathbb E[U\mathbf1_{\{U>x\}}]
 \sim K_p x^{-t},\qquad K_p=\frac{p d_p}{t}.
\end{equation}
The constant follows by writing the expectation as
$x\mathbb P(U>x)+\int_x^\infty\mathbb P(U>y)\,dy$.

Let $Y_\epsilon=V-X_\epsilon$ be the mass in the removed arc.
For every $0<\alpha<t$,
\begin{equation}\label{eq:fractional-rooted-smallball}
 \mathbb E_{\mathbb Q}Y_\epsilon^\alpha
 \le C_\alpha\epsilon^{\alpha(t-\alpha)/p}.
\end{equation}
To prove \eqref{eq:fractional-rooted-smallball}, rescale
$\theta=\epsilon u$ on $[-1,1]$.
Uniformly for distinct $u,v$ in this interval,
\[
 -\log|e^{i\epsilon u}-e^{i\epsilon v}|
 =\log(1/\epsilon)-\log|e^{iu}-e^{iv}|+O(1).
\]
The root weight is comparable to
$\epsilon^{-2/p}|1-e^{iu}|^{-2/p}$. Kahane's convexity inequality,
with an independent constant Gaussian variable to absorb the bounded
covariance difference, compares the moment of order $\alpha$ to a
constant multiple of the moment of the same order of
\[
 \epsilon^{1-2/p}
 \exp\{\sqrt{2/p}\,\Omega-(1/p)\log(1/\epsilon)\}
 \int_{-1}^1|1-e^{iu}|^{-2/p}\mu(du),
 \quad \Omega\sim N(0,\log(1/\epsilon)).
\]
The Gaussian variable is independent of the chaos in the integral.
Its moment of order $\alpha$, together with the deterministic prefactor, is
$\epsilon^{\alpha(t-\alpha)/p}$. The remaining moment is finite:
the integral is bounded by a full rooted circle mass, whose
moment of order $\alpha$ equals $\mathbb E U^{1+\alpha}<\infty$ by
\eqref{eq:fractional-size-bias}. Here Kahane's inequality
is applied with the root weight cut off at a fixed positive distance
from zero and with a common Gaussian regularization. The comparison
constant is independent of both cutoffs. First remove the Gaussian
regularization: the weight is bounded, subcritical chaos converges in
$L^1$. The inequality $|x^\alpha-y^\alpha|\le|x-y|^\alpha$ for
$x,y\ge0$, followed by Jensen's inequality, then gives convergence
of its moments of order $\alpha$.
Then let the root cutoff decrease to zero and use monotone convergence.
The reference moments are bounded by $\mathbb E U^{1+\alpha}$
throughout. This proves \eqref{eq:fractional-rooted-smallball} without assuming uniform
integrability at the divergent critical moment. This is
the usual covariance comparison, as stated for instance in
\cite[Section~2.1; Lemma~2.2 in the arXiv version]{KW2022}. To check
the direction for the concave function $x^\alpha$, add an independent
constant Gaussian of variance $C$ to the actual Gaussian field until
its covariance dominates that of the reference field. Kahane's
inequality decreases the moment of order $\alpha$. The added normalized
Gaussian multiplies that moment by $e^{(\alpha^2-\alpha)C/2}$;
moving this factor to the other side gives the required upper bound.

We now use \eqref{eq:fractional-rooted-smallball} to compare tails
up to level $\epsilon^{-a}$. Fix $0<a<t/p$ and choose $pa<\alpha<t$;
these inequalities make the power of $\epsilon$ below positive.
For fixed $L>0$ and
$0<\eta<1$, Markov's inequality and
\eqref{eq:fractional-rooted-smallball} imply, uniformly for
$L\le x\le\epsilon^{-a}$,
\begin{align*}
 \mathbb Q(Y_\epsilon>\eta x)
 &\le C_\alpha\eta^{-\alpha}
             \epsilon^{\alpha(t-\alpha)/p}x^{-\alpha}\\
 &\le C_\alpha\eta^{-\alpha}
       \epsilon^{(t-\alpha)(\alpha/p-a)}x^{-t}.
\end{align*}
With $L,a,\alpha,\eta$ fixed, this is $o(x^{-t})$ as
$\epsilon\downarrow0$, uniformly in the specified range. Since
$V=X_\epsilon+Y_\epsilon$,
\[
 \mathbb Q(X_\epsilon>(1-\eta)x)
 \ge\mathbb Q(V>x)-\mathbb Q(Y_\epsilon>\eta x)
 \ge (K_p-o_L(1)-o_\epsilon(1))x^{-t}.
\]
The identity $\mathbb E_{\mathbb Q}X_\epsilon^t=
t\int_0^\infty y^{t-1}\mathbb Q(X_\epsilon>y)\,dy$ allows us to
integrate this lower bound. Restricting the integral to
$y\in[(1-\eta)L,(1-\eta)\epsilon^{-a}]$ yields
\[
 \liminf_{\epsilon\downarrow0}
 \frac{\mathbb E_{\mathbb Q}X_\epsilon^t}{\log(1/\epsilon)}
 \ge t\bigl(K_p-o_L(1)\bigr)a(1-\eta)^t.
\]
After this limit in $\epsilon$, let $L\to\infty$, then
$\eta\downarrow0$ and $a\uparrow t/p$.
The resulting constant is $tK_p(t/p)=t d_p=A_p$.
\end{proof}

The truncation away from the marked point is essential: the limiting
critical moment $\E U^p$ is infinite, so the comparison must be made
with the truncated mass before any limit in $n$ is taken.
Take $A_n=\{\theta:|\theta|\ge R/n\}$, with $R$ fixed and large.
By \cref{lem:fractional-rooted-conditional}, angular and chord distance
being comparable near zero,
\[
 \mu_n(A_n)\ge(1+C/R)^{-1}
            \mathbb E_{\mathbb Q}[\mu^0(A_n)\mid\mathcal F_n].
\]
Since $0<t<1$, the conditional Jensen inequality gives
\begin{align*}
 \mathbb E U_n^p
 &=\mathbb E_{\mathbb Q}U_n^t
 \ge(1+C/R)^{-t}
 \mathbb E_{\mathbb Q}
       \bigl(\mathbb E_{\mathbb Q}[\mu^0(A_n)\mid\mathcal F_n]\bigr)^t\\
 &\ge(1+C/R)^{-t}\mathbb E_{\mathbb Q}\mu^0(A_n)^t.
\end{align*}
Apply \cref{lem:fractional-rooted-cutoff}, let $n\to\infty$, and then
$R\to\infty$. Together with \eqref{eq:fractional-opuc-limsup}, this
proves the first assertion in \eqref{eq:fractional-critical-normalized}.

\subsection{Averaging the common independent angle}
\label{sec:fc-averaging}

\begin{lemma}\label{lem:critical-terminal-multipoint}
Recall that
$\chi_N(\theta)=|P_N(e^{i\theta})|^q/M_N(q)$.
Fix an integer $j\ge1$ and $\eta>0$. There is $L_0>1$ such that,
for every fixed $L\ge L_0$, there is $R>0$ for which, with
$n=N-1$ and $m=\lfloor n/L\rfloor$, for all sufficiently large $N$,
\begin{equation}\label{eq:critical-terminal-multipoint}
 \mathbb E\left[\prod_{a=1}^j\chi_N(\theta_a)\mid\mathcal F_m\right]
 \le(1+\eta)\prod_{a=1}^jD_m(\theta_a)
\end{equation}
whenever $m\min_{a\ne b}|e^{i\theta_a}-e^{i\theta_b}|\ge R$.
The bound is uniform in the coefficients revealed by time $m$.
For $j=1$ the conditional equality is exact.
\end{lemma}
\begin{proof}[Proof of \cref{lem:critical-terminal-multipoint}]
We approximate the terminal factor by finitely many Fourier modes.
The conditional block estimate then suppresses each nonzero mode,
while the zero mode gives the product on the right of
\eqref{eq:critical-terminal-multipoint}.

The common uniform angle in \eqref{eq:terminal-process} is
independent of the coefficients. Put
$h(u)=|1-e^{iu}|^q$, so $\widehat h(0)=c(q)>0$. Since $q>0$,
$h$ is continuous. Let $h_J=h*F_J$, where $F_J$ is the Fej\'er
kernel normalized to have integral one against $\sigma$. These
trigonometric polynomials converge uniformly to $h$ and have the
same constant coefficient. Therefore
\[
 H(\boldsymbol\psi)=\int\prod_{a=1}^j h(u+\psi_a)\,d\sigma(u)
\]
is uniformly approximated by the corresponding finite sum $H_J$.
All modes of $H_J$ are neutral:
\[
 H_J(\boldsymbol\psi)
 =\sum_{\substack{\boldsymbol\ell\in\mathbb Z^j\\
                   \sum_a\ell_a=0}}
 c_{J,\boldsymbol\ell}e^{i\boldsymbol\ell\cdot\boldsymbol\psi},
 \qquad c_{J,\boldsymbol0}=c(q)^j,
\]
where only finitely many coefficients are nonzero. In particular,
no bound uniform in all Fourier modes is needed.

For a fixed retained mode, iterate \cref{prop:unifiedblock} over
the finitely many doubling intervals from $m$ to $n$. In this proof write
$Y_r=\prod_a|\Phi_r^*(e^{i\theta_a})|^q$ and
$\delta=\min(1,\min_{a\ne b}|e^{i\theta_a}-e^{i\theta_b}|)$.
The zero mode has exponent $j/p$, and a neutral mode
$\boldsymbol\ell$ has exponent $j/p-2\|\boldsymbol\ell\|^2/\beta$.
The initial mode has absolute value $Y_m$. Each intermediate block
remainder is bounded by its initial positive product, multiplied by
$C_{\boldsymbol\ell}/(r\delta)$ at that block's initial index $r$.
By \eqref{eq:critical-conditional-products}, its conditional expectation
back to time $m$ is at most
$C_{\boldsymbol\ell,L}(r/m)^{j/p}Y_m/(r\delta)$.
There are only finitely many blocks for fixed $L$, and deterministic
phase factors have modulus one. The normalizer ratios introduce a
further $O_{\boldsymbol\ell,L}(m^{-1})$, absorbed in the same error. Thus
\begin{equation}\label{eq:fractional-conditional-fixed-mode}
 \frac{\left|\mathbb E[Y_ne^{i\boldsymbol\ell\cdot\boldsymbol\psi_n}
                         \mid\mathcal F_m]\right|}{Y_m}
 \left(\frac{\nu_m(q)}{\nu_n(q)}\right)^j
 \le (n/m)^{-2\|\boldsymbol\ell\|^2/\beta}
                        +\frac{C_{\boldsymbol\ell,L}}{m\delta}.
\end{equation}
For $\boldsymbol\ell=0$, the normalized conditional expectation is
$1+O((m\delta)^{-1})$ by \eqref{eq:critical-conditional-products}.
The error from replacing $H$ by $H_J$ is at most
$\|H-H_J\|_\infty$ times that same positive conditional product.

Divide by $M_N(q)^j=(c(q)\nu_n(q))^j$, using \eqref{eq:terminal-c}.
First choose $J$ so that
the approximation error is less than $\eta/3$ after normalization.
Next choose $L_0$ so that, for every $L\ge L_0$, the finitely many
nonzero leading terms contribute less than $\eta/3$. Every nonzero
neutral integer vector has positive squared norm, so this is possible.
Fix such an $L$, and then choose $R$ so that all remainder terms,
including the zero mode error, contribute less than $\eta/3$. Make the inequalities strict before taking $N$
large, since $n/m\to L$. These choices prove
\eqref{eq:critical-terminal-multipoint}. The constants do not depend on the
revealed coefficients.
\end{proof}

\begin{proof}[Completion of the proof of
\cref{thm:fractional-critical} for $1<p<2$]
The lower bound follows by averaging the terminal angle. For the
upper bound, we compare the mass outside a small arc about the root
with its OPUC counterpart, and bound the mass of the removed arc
directly.

The first identity in \eqref{eq:fractional-terminal-bounds} and
Jensen's inequality give $\mathbb E W_N^p\ge\mathbb E U_{N-1}^p$, proving the required
lower bound for characteristic polynomials.

For the upper bound introduce the rooted law $\mathbb Q_N$ at finite $N$,
used only in this proof:
\[
 \frac{d\mathbb Q_N}{d\mathbb P}=\chi_N(0)
\]
on the coefficients and the terminal angle. Its restriction to
$\mathcal F_m$, for every $m\le N-1$, is exactly $\mathbb Q|_{\mathcal F_m}$:
terminal angle averaging gives $D_{N-1}(0)$, and the one point
martingale gives its conditional expectation $D_m(0)$.
Rotation invariance gives
$\mathbb E W_N^p=\mathbb E_{\mathbb Q_N}W_N^t$.

Fix $\eta>0$ and choose $L,R$ from
\cref{lem:critical-terminal-multipoint} for $j=2$.
With $m=\lfloor(N-1)/L\rfloor$,
let $B=\{\theta:m\delta(\theta)<R\}$ and $A=\mathbb T\setminus B$.
Write $W_N(A)=\int_A\chi_N\,d\sigma$, and similarly for $B$.
The change of measure formula and
\eqref{eq:critical-terminal-multipoint} show that
\[
 \mathbb E_{\mathbb Q_N}[W_N(A)\mid\mathcal F_m]
 \le(1+\eta)\mu_m(A).
\]
Conditional concavity and the agreement of the rooted laws on
$\mathcal F_m$ therefore yield
\[
 \mathbb E_{\mathbb Q_N}W_N(A)^t
 \le(1+\eta)^t\mathbb E_{\mathbb Q}\mu_m(A)^t
 \le(1+\eta)^t\mathbb E U_m^p.
\]
For the excluded arc, H\"older's inequality and Jensen's inequality in
the spatial integral give directly
\begin{align*}
 \mathbb E_{\mathbb Q_N}W_N(B)^t
 &=\mathbb E[\chi_N(0)W_N(B)^{p-1}]\\
 &\le\sigma(B)^{p-1}\mathbb E\chi_N(0)^p
 =\sigma(B)^t\frac{M_N(pq)}{M_N(q)^p}
 \le C_{p,\beta}(RL)^t.
\end{align*}
Indeed $\sigma(B)=O(R/m)=O(RL/N)$ and \eqref{eq:onepointasymp}
makes the one point ratio $O(N^{p-1})$.
This bound is independent of $N$, with $R,L$ fixed.
Since $(x+y)^t\le x^t+y^t$,
\[
 \mathbb E W_N^p\le(1+\eta)^t\mathbb E U_m^p+C_{p,\beta}(RL)^t.
\]
Divide by $\log N$, let $N\to\infty$, and use $\log m/\log N\to1$.
The OPUC upper bound \eqref{eq:fractional-opuc-limsup} gives a
limsup at most $(1+\eta)^t A_p$.
Finally let $\eta\downarrow0$. This proves
\eqref{eq:fractional-critical-normalized}. The one point
asymptotics \eqref{eq:onepointasymp} and \eqref{eq:opuc-onepoint-asymptotic}
give $M_N(q)^p\sim C_\beta(q)^pN$ and
$\nu_n(q)^p\sim D_\beta(q)^pn$, proving the unnormalized formulas.
\end{proof}

\begin{remark}
The restriction $p<2$ is used twice in an essential way: the rooted
power $p-1$ is concave, and a moment $r$ strictly between $p$ and $2$
is available in \cref{lem:fractional-supercritical-bound}.
For noninteger $p>2$, we replace the single marked point below by
$\lfloor p\rfloor$ marked points. The conclusion is an asymptotic
equivalence, without an $O(N)$ remainder.
\end{remark}

\subsection{Proof of the upper bound for every real critical order}
\label{sec:fc-upper-all}

The conditional comparisons and chaos martingale identity hold for
every $p>1$. We now complete the proof of
\cref{prop:all-real-critical-upper} by extending the annular estimate.

\begin{proof}[Proof of \cref{prop:all-real-critical-upper}]
Only $p\ge2$ remains. We estimate each annular mass by conditional
Jensen for a convex power, sum by Minkowski's inequality, and then
apply the convex truncation from \eqref{eq:fractional-opuc-limsup}.
Choose $p<r<p+1$ and, in this proof, set $s=r-1$, so $1<s<p$.
We use the same exponent $\xi(r)=(r-1)(r/p-1)$ in this range.
For an annulus $A$ of length and distance from the root
comparable to $a$, with $na\ge R_0$, the lower bound in
\eqref{eq:fractional-rooted-conditional} gives
\[
 \mu_n(A)\le C\mathbb E_{\mathbb Q}[\mu^0(A)\mid\mathcal F_n].
\]
The conditional Jensen inequality now applies to the convex power $s$
and yields
\[
 \mathbb E_{\mathbb Q}\mu_n(A)^s
 \le C_s\mathbb E_{\mathbb Q}\mu^0(A)^s
 \le C_s a^{s(p-1-s)/p}=C_s a^{-\xi(r)}.
\]
The last inequality is the Gaussian scaling calculation used in
\eqref{eq:fractional-rooted-smallball}, now on an annulus avoiding
the root. On the rescaled annulus the root weight is bounded, and the
unrooted GMC has a finite moment of order $s$ because $s<p$. Bounded covariance
comparison holds for this convex power as well; adding an independent
constant Gaussian to the appropriate side of the comparison gives
the same fixed multiplicative bound.

On the inner arc of length $O(n^{-1})$, the one point H\"older estimate
from \cref{lem:fractional-supercritical-bound} gives
$\mathbb E_{\mathbb Q}\mu_n(A)^s\le C_r n^{\xi(r)}$.
Minkowski's inequality over dyadic annuli therefore gives
\[
 (\mathbb E_{\mathbb Q}U_n^s)^{1/s}
 \le C_r\sum_{j\le\log_2 n}2^{j\xi(r)/s}
 \le C_r n^{\xi(r)/s}.
\]
Thus $\mathbb E U_n^r\le C_r n^{\xi(r)}$ for $p<r<p+1$.
The proof of \eqref{eq:fractional-opuc-limsup} now applies verbatim:
given $a>(p-1)/p$, take $r>p$ sufficiently close to $p$ that
$(r-1)/p<a$, and use the same convex truncation at $T=n^a$.
This proves the first assertion. The bound for $W_N/U_{N-1}$ in
\eqref{eq:fractional-terminal-bounds} proves the second.
\end{proof}

\subsection{Several marked points and the lower critical coefficient}
\label{sec:fc-several}

For the following unweighted arguments, take a noninteger $p>2$ and put
\[
 k=\lfloor p\rfloor,\qquad s=p-k\in(0,1),\qquad
 A_p=\frac{p-1}{\Gamma(1-1/p)^p}.
\]
We represent the integer part $k$ by marked points and use concavity
for the remaining power $s$. Write $\rho$ for circular angular
distance and $B(x,\epsilon)$ for the corresponding ball.
The following quantities describe the mass that remains after excluding
nearby roots and their neighborhoods, and the mass of nearby pairs:
\[
 \begin{split}
 I_\epsilon&=\E
   \int_{\min_{i\ne j}\rho(\theta_i,\theta_j)>\epsilon}
   \mu\left(\mathbb T\setminus\bigcup_{i=1}^k
                          B(\theta_i,\epsilon)\right)^s
                    \mu^{\otimes k}(d\boldsymbol\theta),\\
 P_\epsilon(\mu)&=\int_{\mathbb T}\mu(dx)\mu(B(x,\epsilon)).
 \end{split}
\]
Thus $P_\epsilon(\mu)$ is the mass, under $\mu^{\otimes2}$, of
pairs at circular distance at most $\epsilon$. The comparison under
the measure with marked points, proved in
\eqref{eq:critical-palm-lower}, gives
\[
 \E U_n^p\ge(1+C/R)^{-1-s}I_{R/n}.
\]
The deterministic cutoff estimate
\eqref{eq:critical-deterministic-cutoff} and the bound for nearby
pairs in \cref{lem:critical-chaos-collision} then give
\[
 I_\epsilon\ge
 \E[U^p;U\le\epsilon^{-a}]
 -C\epsilon^{(p-r)((r-1)/p-a)},
\]
where $0<a<(p-1)/p$ and
$\max(2,p-1)<r<p$ is chosen so that $a<(r-1)/p$.
The error tends to zero, while the exact mass law
\eqref{eq:fractional-exact-mass} makes the first term
$pa\,\Gamma(1-1/p)^{-p}\log(1/\epsilon)+O(1)$.
Letting $a$ increase to $(p-1)/p$ yields the required coefficient.
We now prove these comparisons, starting with the change of measure;
in particular, no product martingale is assumed.

\subsubsection{The measure with marked points}

We first construct the joint finite measure of the coefficients and the marked
points. We then compare its conditional density with the finite
OPUC product and apply concavity to the exterior mass.
Work on the coefficient space
$\Omega=\mathbb D^{\mathbb N_0}$ with its product Borel sigma-field,
where $\mathbb D=\{z\in\mathbb C:|z|<1\}$. For $\omega\in\Omega$,
write $\mu_\omega$ for the realization of the limiting chaos $\mu$.
The product $\prod_iD_n(\theta_i)$ is not a martingale when $k\ge2$.
Define instead the finite measure
\[
 \mathcal M(d\omega,d\boldsymbol\theta)
 =\mathbb P(d\omega)\mu_\omega^{\otimes k}(d\boldsymbol\theta).
\]
Its total mass is $\mathbb E U^k<\infty$, since $k<p$. The Gaussian
Cameron--Martin formula shows that its spatial marginal is
\[
 K_k(\boldsymbol\theta)d\sigma^{\otimes k},\qquad
 K_k(\boldsymbol\theta)
 =\prod_{i<j}|e^{i\theta_i}-e^{i\theta_j}|^{-2/p}.
\]
Both $\Omega$ and $\mathbb T^k$ are standard Borel spaces. The
disintegration theorem for finite measures therefore gives a
measurable probability kernel $\mathbb Q_{\boldsymbol\theta}$ such that
\[
 \mathcal M(d\omega,d\boldsymbol\theta)
 =K_k(\boldsymbol\theta)\mathbb Q_{\boldsymbol\theta}(d\omega)
   \sigma^{\otimes k}(d\boldsymbol\theta).
\]
This identity defines the Palm distribution used here: under
$\mathcal M$, it is the conditional law of the coefficient sequence
given the marked points. The kernel is defined for almost every
spatial configuration. Only measurable dependence on the marked
points is needed for the integrals below; no continuity is required.

Choose $k<r<p$. Since $U\in L^r$, the conditional chaos identity
gives $U_M\to U$ in $L^r$. This convergence and almost sure
weak convergence of $\mu_M$ give convergence of
$\mu_M^{\otimes k}(F)$ in $L^1$ for every bounded continuous spatial
test $F$: pointwise convergence follows from weak convergence, while
$|\mu_M^{\otimes k}(F)|\le\|F\|_\infty U_M^k$ and $r/k>1$
give uniform integrability from the uniform $r$th mass moment.
Pass \eqref{eq:critical-conditional-products}, integrated against nonnegative
tests on separated configurations, to this conditional tensor limit.
Thus the measure $\mathbb E[\mu^{\otimes k}\mid\mathcal F_n]$
has a density $H_{n,\boldsymbol\theta}$ with respect to
$\sigma^{\otimes k}$ on $n\delta\ge R$. Almost everywhere there,
\begin{equation}\label{eq:critical-palm-density-comparison}
 (1-C/R)\prod_iD_n(\theta_i)
 \le H_{n,\boldsymbol\theta}
 \le(1+C/R)\prod_iD_n(\theta_i).
\end{equation}
To identify this with the Palm density, test also against a bounded
random variable measurable with respect to $\mathcal F_n$ in the
preceding passage.
Disintegration gives
\[
 \frac{d(\mathbb Q_{\boldsymbol\theta}|\mathcal F_n)}
      {d(\mathbb P|\mathcal F_n)}
 =\frac{H_{n,\boldsymbol\theta}}{K_k(\boldsymbol\theta)}.
\]
For each $n$, first use indicators from a countable determining
class for $\mathcal F_n$, which is generated by the first $n$
coefficient coordinates. A monotone class argument then gives the
identity for every bounded test measurable with respect to
$\mathcal F_n$. Taking the countable
union of the exceptional spatial sets over these classes and over
$n$ permits all versions to be chosen outside one spatial null set.
For every distinct configuration outside that set, the separation
condition holds for all sufficiently large $n$.
Absolute continuity at an earlier index follows by restricting one
of these later densities.

Restrict the roots to angular separation greater than $R/n$, and
let $A$ be the complement of their arcs of radius $R/n$. Angular
and chord distances satisfy $2\rho/\pi\le d\le\rho$ on
$0\le\rho\le\pi$. Choosing $R$ sufficiently large therefore
permits every conditional product estimate below, with error $C/R$.
For $M\ge n$, the change of measure formula is
\[
 \mathbb E_{\mathbb Q_{\boldsymbol\theta}}
       [\mu_M(A)\mid\mathcal F_n]
 =\frac{\mathbb E[H_{M,\boldsymbol\theta}\mu_M(A)
                               \mid\mathcal F_n]}
        {H_{n,\boldsymbol\theta}}.
\]
Use the upper comparison for $H_M$ and the lower comparison for $H_n$
in \eqref{eq:critical-palm-density-comparison}, with $R$ large enough
that the latter is positive. Then apply
\eqref{eq:critical-conditional-products} at $k+1$ separated points.
After increasing the constant,
\[
 \mathbb E_{\mathbb Q_{\boldsymbol\theta}}
       [\mu_M(A)\mid\mathcal F_n]\le(1+C/R)\mu_n(A).
\]
Let $E\subset\Omega$ be the event that $\mu_M\to\mu$ weakly and
$\mu$ has no atoms. By \cref{prop:chaos-input} and the conditional
identity \eqref{eq:fractional-chaos-martingale}, it has probability
one under $\mathbb P$, so
\[
 \mathcal M(E^c\times\mathbb T^k)
 =\mathbb E[\mathbf1_{E^c}U^k]=0.
\]
Disintegration therefore gives $\mathbb Q_{\boldsymbol\theta}(E)=1$
for almost every configuration. Nonatomicity implies that the
finitely many endpoints of $A$ have zero mass, so
$\mu_M(A)\to\mu(A)$ almost surely under these Palm distributions.
The conditional Fatou lemma gives
\begin{equation}\label{eq:critical-palm-exterior}
 \mathbb E_{\mathbb Q_{\boldsymbol\theta}}
       [\mu(A)\mid\mathcal F_n]\le(1+C/R)\mu_n(A).
\end{equation}
This argument does not assume that the global $(k+1)$st GMC moment is
finite. The comparison at $k+1$ points is localized off every partial diagonal.

Expand $k$ powers of $U_n$ and restrict the roots to the separated
region $S_n=\{\min_{i\ne j}\rho(\theta_i,\theta_j)>R/n\}$.
Since $U_n\ge\mu_n(A)$, \eqref{eq:critical-palm-density-comparison}
gives
\begin{align*}
 \mathbb E U_n^p
 &\ge\int_{S_n}\mathbb E\left[\prod_iD_n(\theta_i)
                                      \mu_n(A)^s\right]d\sigma^{\otimes k}\\
 &\ge(1+C/R)^{-1}\int_{S_n}
       K_k(\boldsymbol\theta)
       \mathbb E_{\mathbb Q_{\boldsymbol\theta}}\mu_n(A)^s
                                                  d\sigma^{\otimes k}.
\end{align*}
On the other hand, \eqref{eq:critical-palm-exterior} and conditional
concavity imply
\[
 \mu_n(A)^s\ge(1+C/R)^{-s}
   \left(\mathbb E_{\mathbb Q_{\boldsymbol\theta}}
                     [\mu(A)\mid\mathcal F_n]\right)^s
 \ge(1+C/R)^{-s}\mathbb E_{\mathbb Q_{\boldsymbol\theta}}
                     [\mu(A)^s\mid\mathcal F_n].
\]
Integrating this inequality against the spatial Palm density yields
\begin{equation}\label{eq:critical-palm-lower}
 \mathbb E U_n^p\ge(1+C/R)^{-1-s}I_{R/n}.
\end{equation}

\subsubsection{The contribution of nearby points in the limiting chaos}
\begin{lemma}\label{lem:critical-chaos-collision}
For every $\max(2,p-1)<r<p$,
\begin{equation}\label{eq:critical-chaos-collision}
 \mathbb E[U^{r-2}P_\epsilon(\mu)]
 \le C_r\epsilon^{(r-1)(p-r)/p},\qquad 0<\epsilon<\epsilon_0.
\end{equation}
\end{lemma}

\begin{proof}[Proof of \cref{lem:critical-chaos-collision}]
We mark two nearby points and estimate the remaining mass by scaling
their common neighborhood. Dyadic annuli control the exterior mass,
after which we integrate over the distance between the points.
In this proof put $t=r-2>0$ for the remaining power of the mass.
The Gaussian Cameron--Martin formula at two points gives
\begin{align*}
 \mathbb E[U^t\mu(dx)\mu(dy)]
 &=|e^{ix}-e^{iy}|^{-2/p}\mathbb E[V_{x,y}^t]\,d\sigma(x)d\sigma(y),\\
 V_{x,y}
 &=\int|e^{iz}-e^{ix}|^{-2/p}|e^{iz}-e^{iy}|^{-2/p}\mu(dz).
\end{align*}
This identity can be obtained with a common Gaussian regularization
and bounded test functionals. Choose $r<r'<p$; since $r'/r>1$,
the uniform $r'$th mass moment gives uniform integrability of the
total mass functional of order $r=t+2$.
Remove Gaussian regularization using subcritical moment
convergence. Remove bounded spatial and functional cutoffs afterwards
by monotone convergence. Monotonicity of the Gaussian regularizations
themselves is not needed.

By rotation take $x=0$, and put $d=\rho(x,y)$, initially small.
On the interval $|z|\le3d$, rescale $z=du$. The two roots are then
at $0$ and $1$ or $-1$. Use as reference on $[-3,3]$ a log correlated field on a circle of
length $16$. Its covariance is
\[
 K_0(u,v)=-\frac2p\log\left|2\sin\frac{\pi(u-v)}{16}\right|.
\]
On this interval the difference
\[
 -\frac2p\log|e^{idu}-e^{idv}|
 -\frac2p\log(1/d)-K_0(u,v)
\]
is uniformly bounded, including its continuous extension on the
diagonal. This follows from $|\sin x|/|x|$ being bounded above and
below on the fixed ranges in question. The two weights contribute $d^{-4/p}$
and the spatial measure contributes $d$. Apply Kahane's covariance inequality with bounded spatial weights
and a common Gaussian regularization. For $t\ge1$, add a constant
independent Gaussian to the reference field to dominate the actual
covariance. For $0<t<1$, add it to the actual field instead and use
the reversed comparison for the concave power. The added variance is
bounded independently of $d$, so in both cases it changes the bound
by only a fixed multiplicative factor. We obtain a bound by a
constant times the $t$th moment of
\[
 d^{1-4/p}e^{\Omega-\frac12\mathbb E\Omega^2}
 \int_{-3}^3 |u|^{-2/p}|u-1|^{-2/p}\,\widetilde\mu(du),
 \qquad \mathbb E\Omega^2=\frac2p\log(1/d),
\]
Here the displayed linear weights are comparable, by constants
independent of $d$, to the circular weights of the reference field;
when $y=-d$, reflect $u$ to place the second root at $1$.
The coarse Gaussian is independent of the reference chaos.
Its contribution gives the power
\[
 \nu=t(1-4/p)-\frac{t^2-t}{p}
     =\frac{t(p-3-t)}p
     =\frac{(r-2)(p-r-1)}p<0,
\]
since $r>2$ and $r>p-1$.
The reference weighted $t$th moment is finite. Split its interval into
fixed disjoint neighborhoods of the two roots and the complement.
On either neighborhood the other root weight is bounded. After a fixed change of scale, the single root moment is bounded
by a constant times the $t$th moment of a rooted circular mass. By
the Gaussian Cameron--Martin formula and rotation invariance, the latter is
$\mathbb E U^{t+1}$ by \eqref{eq:fractional-size-bias}.
This is finite because $t+1=r-1<p$ by
\eqref{eq:fractional-exact-mass}. On the complement the weights are bounded and
$t<p$. Subadditivity for $t\le1$, or the finite sum inequality for
$t\ge1$, finishes this check.

On each exterior dyadic annulus
$\{z:a\le|z|\le\min(2a,\pi)\}$ in the circular coordinate
$z\in[-\pi,\pi]$, where $3d\le a\le\pi$, both root weights are
$O(a^{-2/p})$. The unrooted chaos mass on that annulus has $t$th
moment bounded by
$C a^{(1+1/p)t-t^2/p}$ by the same bounded covariance scaling.
Thus the weighted annular mass has moment at most $C a^\nu$.
For $t\le1$ sum the moments; for $t\ge1$ sum their $t$th roots by
Minkowski's inequality. Since $\nu<0$, the geometric sum is dominated by the
innermost scale. The macroscopic remainder is bounded using $t<p$.
We have proved, uniformly in the two points,
\[
 \mathbb E[V_{x,y}^t]\le C d^\nu.
\]
Integrating this formula over $d<\epsilon$ now gives
\[
 \mathbb E[U^{r-2}P_\epsilon(\mu)]
 \le C\int_0^\epsilon d^{-2/p+\nu}\,dd
 =C\epsilon^{(r-1)(p-r)/p},
\]
since $1-2/p+\nu=(r-1)(p-r)/p>0$.
This proves \eqref{eq:critical-chaos-collision}.
\end{proof}

\subsubsection{Recovery of the logarithmic coefficient}
The following inequality holds for every finite measure $\mu$,
without taking expectations. In its left side,
``separated'' means that the marked points have mutual distance
greater than $\epsilon$, and ``exterior'' is the complement of their
balls of radius $\epsilon$, as in the definition of $I_\epsilon$:
\begin{equation}\label{eq:critical-deterministic-cutoff}
 \int_{\rm separated}\mu(\text{exterior})^s\,d\mu^{\otimes k}
 \ge U^p-C_kU^{p-2}P_\epsilon(\mu),\qquad
 C_k=\binom{k}{2}+k.
\end{equation}
For $U>0$, sample the $k$ roots independently from $\mu/U$.
The probability that a marked pair is within distance $\epsilon$ is at most
$\binom{k}{2}P_\epsilon/U^2$.
The expected fraction of mass in the union of the $k$ root balls is
at most $kP_\epsilon/U^2$.
For a fraction $z\in[0,1]$, $(1-z)^s\ge1-z$.
These three observations give \eqref{eq:critical-deterministic-cutoff};
for $U=0$ it is immediate. In particular this deterministic estimate
must be truncated before subtracting expectations, since
$\mathbb E U^p=\infty$.

Fix $0<a<(p-1)/p$, and choose
$\max(2,p-1)<r<p$ such that $a<(r-1)/p$.
Restrict \eqref{eq:critical-deterministic-cutoff} to the event
$U\le\epsilon^{-a}$. On this event,
$U^{p-r}\le\epsilon^{-a(p-r)}$.
The integral on the left is nonnegative, so
\eqref{eq:critical-chaos-collision} gives
\begin{align*}
 I_\epsilon
 &\ge\mathbb E[U^p;U\le\epsilon^{-a}]
     -C\epsilon^{-a(p-r)}\mathbb E[U^{r-2}P_\epsilon(\mu)]\\
 &\ge\mathbb E[U^p;U\le\epsilon^{-a}]
     -C\epsilon^{(p-r)((r-1)/p-a)}.
\end{align*}
The last error tends to zero because both factors in its exponent
are positive. By \eqref{eq:fractional-exact-mass} and the definition of $d_p$,
\[
 \mathbb P(U>x)=1-e^{-d_px^{-p}},\qquad
 \mathbb E[U^p;U\le\epsilon^{-a}]
 =p d_p a\log(1/\epsilon)+O(1).
\]
Let $a\uparrow(p-1)/p$, choosing $r$ sufficiently close to $p$ for
each fixed $a$. We obtain
\[
 \liminf_{\epsilon\downarrow0}
 \frac{I_\epsilon}{\log(1/\epsilon)}\ge A_p.
\]
Now fix large $R$, send $n\to\infty$ in
\eqref{eq:critical-palm-lower}, and then send $R\to\infty$.
Together with \cref{prop:all-real-critical-upper}, this proves
\begin{equation}\label{eq:all-orders-opuc-critical}
 \mathbb E U_n^p\sim A_p\log n.
\end{equation}

\subsection{A moment bound with one distinguished point}
\label{sec:fc-distinguished}

The transfer to characteristic polynomials requires a bound uniform
when marked points are at distance $O(N^{-1})$. The global positive
moment estimate of \cref{lem:moments-global-bound} gives, for fixed
nonnegative $a_1,\ldots,a_j$,
\begin{equation}\label{eq:critical-normalized-global-envelope}
 \mathbb E\prod_{i=1}^j\chi_N(\theta_i)^{a_i}
 \le C N^{(\sum_i a_i^2-\sum_i a_i)/p}
 \prod_{i<l}(\rho(\theta_i,\theta_l)\vee N^{-1})^{-2a_i a_l/p}.
\end{equation}
The estimate allows coincident points. Its proof uses the
conditional block estimate and induction on the number of nonnegative
exponents, not the critical moment asymptotic.

\begin{lemma}\label{lem:critical-marked-moment}
Let $m\ge0$ be an integer and $w>1$ satisfy $m+w=p$.
Then
\begin{equation}\label{eq:critical-marked-moment}
 \mathbb E[\chi_N(0)^w W_N^m]=O(N^{w-1}).
\end{equation}
The same bound holds for the spatial integral of the deterministic
envelope in \eqref{eq:critical-normalized-global-envelope}, with weights $(w,1,\ldots,1)$.
\end{lemma}

\begin{proof}[Proof of \cref{lem:critical-marked-moment}]
For $m=0$, this follows from \eqref{eq:onepointasymp}.
For $m\ge1$, we expand the $m$ copies of $W_N$ and rescale their
angular variables by $N$. The power of $N$ gives the desired bound;
the remaining integral is controlled by the Dyson integral.
Apply \eqref{eq:critical-normalized-global-envelope} with
weights $(w,1,\ldots,1)$ to obtain
\begin{align*}
 \mathbb E\left[\chi_N(0)^w\prod_{i=1}^m\chi_N(\theta_i)\right]
 \le{}& C N^{(w^2-w)/p}
 \prod_{i=1}^m\left(|1-e^{i\theta_i}|\vee N^{-1}\right)^{-2w/p}\\
 &\times\prod_{i<j}
 \left(|e^{i\theta_i}-e^{i\theta_j}|\vee N^{-1}\right)^{-2/p}.
\end{align*}
The normalizing factor is $M_N(q)^{w+m}=M_N(q)^p$;
the excess one point exponent is $(w^2-w)/p$.

Cut the circle at a largest cyclic gap among the $m+1$ marked points,
including $0$. That gap has length at least $2\pi/(m+1)$.
On each of the finitely many resulting ordering domains, unwrap the
circle and choose real representatives $t_i$, with the representative
of the distinguished point equal to zero. All chord distances are
comparable to ordinary differences, with constants depending only on
$m$. Indeed all unwrapped differences are at most
$2\pi-2\pi/(m+1)$, on which $|1-e^{it}|/|t|$ has a positive lower bound.
The changes of representatives are translations by $2\pi$ and preserve
Lebesgue measure on each ordering domain.

Set $x_i=Nt_i$, and enlarge the integration domain to $\mathbb R^m$.
The total power of $N$ is
\[
 \frac{w^2-w+2mw+m(m-1)}p-m
 =\frac{p^2-p}p-m=w-1.
\]
It remains to show that
\[
 \int_{\mathbb R^m}
 \prod_i(1+|x_i|)^{-2w/p}
 \prod_{i<j}(1+|x_i-x_j|)^{-2/p}\,d\mathbf x<\infty.
\]
Replace each pair factor by $|x_i-x_j|^{-2/p}$, which is larger,
and use $(1+|x|)^{-2w/p}\le C(1+x^2)^{-w/p}$.
Under $x_i=\tan(\phi_i/2)$, the Jacobian and pair factors leave the
power
\[
 |\cos(\phi_i/2)|^{2(w+m-1)/p-2}
 =|\cos(\phi_i/2)|^{-2/p}
\]
at each variable. Thus the preceding integral is bounded by a constant times the circular
Dyson integral with $m+1$ identical exponents and pair exponent $-2/p$,
one point being fixed at $\phi=\pi$.
This integral is finite because $m+1<p$, which is exactly $w>1$.
Finiteness follows from \eqref{eq:dyson-negative-power}, taking
$t=2/p$ and $k=m+1$ in that formula. It also follows by counting
powers near each possible collision: a cluster of $j\le m+1$ points
has homogeneity $-j(j-1)/p$ in
$j-1$ relative coordinates, and $j<p$ for every such cluster.
This proves \eqref{eq:critical-marked-moment} and its envelope version.
\end{proof}

\subsection{Nearby marked points and their neighborhoods}
\label{sec:fc-nearby}

\subsubsection{Pairs of nearby marked points}

Recall that $k=\lfloor p\rfloor$ and $s=p-k\in(0,1)$.
We interpolate between two instances of
\cref{lem:critical-marked-moment}, so that the powers of $W_N$
in both expectations are integers.
Apply H\"older's inequality with weights $1-s,s$:
\begin{align*}
 \mathbb E[\chi_N(0)^2W_N^{p-2}]
 \le{}&
 \bigl(\mathbb E[\chi_N(0)^{2+s}W_N^{k-2}]\bigr)^{1-s}\\
 &\times\bigl(\mathbb E[\chi_N(0)^{1+s}W_N^{k-1}]\bigr)^s.
\end{align*}
The pointwise exponents are exactly
\[
 (2+s)(1-s)+(1+s)s=2,\qquad
 (k-2)(1-s)+(k-1)s=p-2.
\]
By \eqref{eq:critical-marked-moment}, the two expectations are respectively
$O(N^{1+s})$ and $O(N^s)$. Since $(1+s)(1-s)+s^2=1$,
\begin{equation}\label{eq:critical-collision-point}
 \mathbb E[\chi_N(0)^2W_N^{p-2}]=O(N).
\end{equation}

To measure nearby pairs, write $\lambda_N=\chi_N\sigma$ for the
random measure with density $\chi_N$, whose total mass is $W_N$, and put
\[
 P_\epsilon(\lambda_N)
 =\int\lambda_N(d\theta)\lambda_N(B(\theta,\epsilon)).
\]
The Cauchy--Schwarz inequality and rotation invariance imply
\[
 \mathbb E[W_N^{p-2}\chi_N(\theta)\chi_N(\phi)]
 \le\mathbb E[W_N^{p-2}\chi_N(0)^2].
\]
Consequently, for every fixed $C>0$,
\begin{equation}\label{eq:critical-finite-collision}
 \mathbb E[W_N^{p-2}P_{C/N}(\lambda_N)]=O_C(1).
\end{equation}
No critical $L^p$ convergence is used here.

\subsubsection{Neighborhoods of separated roots}

Fix $C>0$. Let the $k$ marked points have pairwise angular
distance strictly greater than $4C/N$, and put $B_i=B(\theta_i,C/N)$.
We claim
\begin{equation}\label{eq:critical-near-ball}
 \int_{\text{separated}}
 \mathbb E\left[\prod_{j=1}^k\chi_N(\theta_j)
                         \lambda_N(B_i)^s\right]d\sigma^{\otimes k}
 =O_C(1).
\end{equation}
Indeed, for the fixed marked points, set $\chi_j=\chi_N(\theta_j)$ and
$G=\chi_i^{1+s}\prod_{j\ne i}\chi_j$.
Write the integrand as $G(\lambda_N(B_i)/\chi_i)^s$ and normalize
the finite measure with density $G$ to a probability measure.
The ratio can be set to zero where $\chi_i=0$, since $G$ also vanishes there.
Jensen's inequality for the concave power $s$ gives
\[
 \mathbb E\left[\prod_j\chi_j\lambda_N(B_i)^s\right]
 \le(\mathbb EG)^{1-s}
 \left(\int_{B_i}\mathbb E\left[
       \chi_i^s\prod_{j\ne i}\chi_j\chi_N(x)\right]d\sigma(x)\right)^s.
\]
In this estimate let $E_N$ denote the deterministic expression on the
right of \eqref{eq:critical-normalized-global-envelope}, without its constant,
with weight $1+s$ at root $i$ and weight $1$ at the other roots.
We suppress its dependence on the fixed marked points.
Then $\mathbb EG\le C E_N$.
For $x\in B_i$, every distance from $x$ to another root is comparable
to the corresponding distance from $\theta_i$, because the roots are
more than four ball radii apart. Apply
\eqref{eq:critical-normalized-global-envelope} with separate weights
$s,1,\ldots,1,1$ to obtain
\[
 \mathbb E\left[\chi_i^s\prod_{j\ne i}\chi_j\chi_N(x)\right]
 \le C_C E_N.
\]
The exponent check is exact: the separate one point excess
$(s^2-s)/p$, plus the clipped interaction $2s/p$ between $x$ and
$\theta_i$, is $s(s+1)/p$, the excess of the merged weight $1+s$.
The interactions with each other root add to $2(1+s)/p$.
Both expectations are normalized by $M_N(q)^p$.
It follows that
\[
 \mathbb E\left[\prod_j\chi_j\lambda_N(B_i)^s\right]
 \le C_C N^{-s}E_N.
\]
Integrate over the roots and enlarge the domain. Rotation invariance
and the envelope part of \eqref{eq:critical-marked-moment}, with $w=1+s$
and $m=k-1$, give
\[
 N^{-s}\int E_N\,d\sigma^{\otimes k}=O(1).
\]
This proves \eqref{eq:critical-near-ball}. The constants remain uniform if
$C$ ranges over a fixed compact interval in $(0,\infty)$; an upper
bound for $C$ is enough when enlarging the balls.

\subsection{The upper coefficient for characteristic polynomials}
\label{sec:fc-upper-transfer}

\begin{proof}[Completion of \cref{thm:fractional-critical} for $p>2$]
We first bound the contributions of nearby roots and their small
neighborhoods. Both are $O(1)$, so only the exterior mass contributes
to the logarithmic coefficient. We compare that mass with the OPUC mass.

Fix $\eta>0$. In \cref{lem:critical-terminal-multipoint}, choose the
Fourier truncations for $j=k$ and $j=k+1$ first, then a common $L$
above both thresholds. For this fixed $L$, choose $R_0$ above the two
resulting separation thresholds.
Put $n=N-1$, and use the OPUC index $m=\lfloor n/L\rfloor$
for the conditional comparison.
Set $R=\pi R_0/2$ for the angular
cutoffs below. Then an angular distance at least $R/m$ has chord
distance at least $R_0/m$.
Expand the $k$ integral powers in $W_N^p=W_N^kW_N^s$.
Configurations of marked points containing a pair at angular distance at most
$4R/m$ contribute at most
\[
 \binom{k}{2}\,
 \mathbb E[W_N^{p-2}P_{4R/m}(\lambda_N)]=O_{L,R}(1)
\]
by \eqref{eq:critical-finite-collision}, because $N/m$ remains bounded for fixed $L$.
On the remaining configurations remove the root neighborhoods
$B_i=B(\theta_i,R/m)$ and put $A=\mathbb T\setminus\bigcup_i B_i$.
Since $0<s<1$,
\[
 W_N^s\le\lambda_N(A)^s+\sum_i\lambda_N(B_i)^s.
\]
The integrated neighborhood terms are $O_{L,R}(1)$ by \eqref{eq:critical-near-ball}. Here
$C=NR/m$ ranges over a fixed compact interval of positive values,
on which the proof of \eqref{eq:critical-near-ball} is uniform.

For the exterior term, conditional concavity gives
\begin{align*}
 &\mathbb E\left[\prod_i\chi_N(\theta_i)\lambda_N(A)^s
                         \mid\mathcal F_m\right]\\
 &\quad\le
 \left(\mathbb E\left[\prod_i\chi_N(\theta_i)\mid\mathcal F_m\right]
                                                               \right)^{1-s}
 \left(\mathbb E\left[\prod_i\chi_N(\theta_i)\lambda_N(A)
                            \mid\mathcal F_m\right]\right)^s\\
 &\quad\le(1+\eta)\prod_iD_m(\theta_i)\mu_m(A)^s.
\end{align*}
The last step uses \eqref{eq:critical-terminal-multipoint} at $k$ and $k+1$ points. Every point in $A$
is at the required chord distance from every root by the choice
of the angular cutoff $R$.
Integrating over the roots and enlarging the domain gives
\begin{equation}\label{eq:critical-terminal-transfer}
 \mathbb E W_N^p\le(1+\eta)\mathbb E U_m^p+O_{L,R}(1).
\end{equation}
Let $N\to\infty$, with $L,R$ fixed. Since $\log m/\log N\to1$,
\eqref{eq:all-orders-opuc-critical} gives
\[
 \limsup_{N\to\infty}\frac{\mathbb E W_N^p}{\log N}
 \le(1+\eta)A_p.
\]
Finally let $\eta\downarrow0$.
The exact identity $\mathbb E[W_N\mid\mathcal F_{N-1}]=U_{N-1}$
in \eqref{eq:fractional-terminal-bounds}
and Jensen's inequality give the matching lower bound. This proves
\eqref{eq:fractional-critical-normalized} for every noninteger $p>2$.
\end{proof}

The conclusion is an asymptotic equivalence. Neither an $O(N)$
unnormalized remainder nor a second coefficient at noninteger orders
is obtained by this argument.

\subsection{Spatial weights for critical moments}
\label{sec:fractional-critical-weights}

We now allow every real $p>1$. The leading critical coefficient for a
complex spatial weight $f$ is proportional to $\int_\T|f|^p\,d\sigma$.
The estimate for two disjoint arcs in
\cref{prop:fractional-critical-additivity} is the common input.
For $1<p<2$ it follows from concavity and a two point bound; for
$p\ge2$ we use the rooted conditional comparison and sum annular
moments. The resulting arc asymptotic and the passage to arbitrary
$L^p$ weights are given in
\cref{cor:weighted-fractional-critical,prop:critical-lp-weights}.

\begin{proposition}
\label{prop:fractional-critical-additivity}
Let $p>1$, $q=\sqrt{2\beta/p}$, and let
\[
 \chi_N(\theta)=\frac{|P_N(e^{i\theta})|^q}{M_N(q)},
 \qquad W_N(I)=\int_I\chi_N\,d\sigma.
\]
If $I,J$ are arcs with disjoint interiors, then
\begin{equation}\label{eq:fractional-disjoint-mixed-bound}
 \sup_N\E[W_N(I)W_N(J)^{p-1}]<\infty.
\end{equation}
Consequently, for a fixed finite partition of the circle into arcs
$I_1,\ldots,I_l$ and fixed complex numbers $c_1,\ldots,c_l$,
\begin{equation}\label{eq:fractional-additivity}
 \E\left|\sum_{j=1}^l c_jW_N(I_j)\right|^p
 =\sum_{j=1}^l |c_j|^p\E W_N(I_j)^p+O(1).
\end{equation}
The constant may depend on the partition and the coefficients.
\end{proposition}

\begin{proof}[Proof of \cref{prop:fractional-critical-additivity}]
We treat the mixed moment separately for $1<p<2$ and $p\ge2$.
An elementary inequality for complex sums then gives
\eqref{eq:fractional-additivity} in both ranges.

Suppose first that $1<p<2$, and put $t=p-1\in(0,1)$. For each fixed $\theta$, the random variable
$\chi_N(\theta)$ has mean one and defines a probability measure
$\mathbb Q_{N,\theta}$ by change of density. Fubini's theorem and
concavity give
\begin{align*}
 \E[W_N(I)W_N(J)^t]
 &=\int_I\E_{\mathbb Q_{N,\theta}}W_N(J)^t\,d\sigma(\theta)\\
 &\le\int_I\left(\int_J
      \E[\chi_N(\theta)\chi_N(\phi)]\,d\sigma(\phi)\right)^t
                                               d\sigma(\theta).
\end{align*}
The OPUC two point formula in \cref{thm:opuc-relative} and
\eqref{eq:terminal-process} give
the global estimate
\begin{equation}\label{eq:fractional-clipped-two-point}
 \E[\chi_N(\theta)\chi_N(\phi)]
 \le C\bigl(N^{-1}+|e^{i\theta}-e^{i\phi}|\bigr)^{-2/p}.
\end{equation}
The two ranges are treated as follows. If $(N-1)|e^{i\theta}-e^{i\phi}|\ge R$,
\cref{thm:opuc-relative} bounds the normalized OPUC product by a
constant times $|e^{i\theta}-e^{i\phi}|^{-2/p}$.
The factor relating $\chi_N(\theta)$ to $D_{N-1}(\theta)$
is at most $2^q/c(q)$, so the same bound holds for the normalized
characteristic product. In the complementary
range, the Cauchy--Schwarz inequality and rotation invariance give
\[
 \E[\chi_N(\theta)\chi_N(\phi)]
 \le\E\chi_N(0)^2
 =\frac{M_N(2q)}{M_N(q)^2}\le C N^{2/p}.
\]
These two estimates imply \eqref{eq:fractional-clipped-two-point},
with finitely many small values of $N$ absorbed into the constant.
For $\theta$ in the interior of $I$, put
$d(\theta)=\operatorname{dist}_{\T}(\theta,J)>0$.
Since $2/p>1$, integration of \eqref{eq:fractional-clipped-two-point}
over $J$ gives
\[
 \int_J\E[\chi_N(\theta)\chi_N(\phi)]\,d\sigma(\phi)
 \le C\bigl(1+d(\theta)^{1-2/p}\bigr).
\]
This also covers arcs crossing the origin by splitting them into
at most two intervals. Raising the right-hand side to $t$ produces at
an adjacent endpoint the singular power
$d(\theta)^{-(p-1)(2-p)/p}$. Its exponent is strictly larger than
$-1$, so it is integrable in $\theta$. This proves
\eqref{eq:fractional-disjoint-mixed-bound} even when $I$ and $J$
share an endpoint. Zero length arcs contribute zero.

Now suppose $p\ge2$, and put $t=p-1\ge1$. We first prove the
corresponding bound for the OPUC masses $\mu_n(I)$ and $\mu_n(J)$.
For each $\theta\in\T$, let $\mathbb Q_\theta$ be the law on the
coefficient sequence whose restriction to $\mathcal F_n$ has density
$D_n(\theta)$. The construction of $\mathbb Q$ and rotation covariance
in \cref{sec:fc-one-marked} give this law and the rooted limit
\[
 \mu^\theta(d\phi)\ \overset{\rm law}=
 |e^{i\theta}-e^{i\phi}|^{-2/p}\widetilde\mu(d\phi)
 \quad\hbox{under }\mathbb Q_\theta,
\]
where $\widetilde\mu$ is ordinary circular chaos. Fubini's theorem gives
\[
 \E[\mu_n(I)\mu_n(J)^t]
 =\int_I\E_{\mathbb Q_\theta}\mu_n(J)^t\,d\sigma(\theta).
\]
We shall bound the integrand by an integrable power of the logarithm
of the distance from $\theta$ to $J$.

Write $\rho$ for circular angular distance. Choose $R$ large enough
that \cref{lem:fractional-rooted-conditional}, after rotation, has
lower factor at least $1/2$ on every set at angular distance at
least $R/n$ from $\theta$. This is possible because chord distance
is at least $2/\pi$ times angular distance. For $n\ge1$, decompose
$J$ into the inner piece and the annular pieces
\[
 \begin{split}
 B&=J\cap\{\phi:\rho(\theta,\phi)<R/n\},\\
 A_j&=J\cap\{\phi:2^jR/n\le\rho(\theta,\phi)<2^{j+1}R/n\},
       \qquad j\ge0.
 \end{split}
\]
Empty pieces are omitted. For a nonempty $A_j$, put $a=2^jR/n$.
The conditional comparison gives
\[
 \mu_n(A_j)\le
 2\E_{\mathbb Q_\theta}[\mu^\theta(A_j)\mid\mathcal F_n].
\]
Conditional Jensen applies to the convex power $t$. The root weight
on $A_j$ is bounded by $Ca^{-2/p}$, and the unrooted chaos satisfies
\[
 \E\widetilde\mu(A_j)^t
 \le C a^{(1+1/p)t-t^2/p}.
\]
For small $a$, the set $A_j$ is contained in two arcs of length at
most $a$. Rescaling either arc gives covariance
$(2/p)\log(1/a)$ plus a fixed reference covariance and a uniformly
bounded error. The independent constant Gaussian with variance
$(2/p)\log(1/a)$ contributes $a^{-(t^2-t)/p}$ to the moment, while
the spatial Jacobian contributes $a^t$. Bounded covariance comparison,
as in the proof of \cref{prop:all-real-critical-upper}, changes this
by a fixed factor. The rescaled unrooted moment is finite because
$t=p-1<p$. The inequality for the sum of two nonnegative terms gives
the displayed bound on $A_j$. For $a$ bounded away from zero,
increase the constant and use $\E U^t<\infty$. Consequently,
\[
 \E_{\mathbb Q_\theta}\mu_n(A_j)^t
 \le 2^t\E_{\mathbb Q_\theta}\mu^\theta(A_j)^t
 \le C a^{-2t/p+(1+1/p)t-t^2/p}
 =C a^{t(p-1-t)/p}=C.
\]
The constant is independent of $j,n$ and $\theta$. In particular,
the annular estimate remains valid at $t=p-1$, where its exponent
vanishes.

The inner piece has $\sigma(B)\le CR/n$. Jensen's inequality in the
spatial variable and H\"older's inequality in probability give
\begin{align*}
 \E_{\mathbb Q_\theta}\mu_n(B)^t
 &=\E[D_n(\theta)\mu_n(B)^t]\\
 &\le\sigma(B)^{t-1}\int_B
          \E[D_n(\theta)D_n(\phi)^t]\,d\sigma(\phi)\\
 &\le\sigma(B)^t\E D_n(0)^p\le C_R.
\end{align*}
Indeed, H\"older uses the conjugate exponents $p$ and $p/t$, and
rotation invariance makes the two moments of order $p$ equal.
The one point formula \eqref{eq:opuc-onepoint-exact} gives
\[
 \E D_n(0)^p=\frac{\nu_n(pq)}{\nu_n(q)^p}
       \le C(n+1)^{p-1},
\]
which cancels the factor $n^{-t}$ from $\sigma(B)^t$.
At $t=1$, spatial Jensen is an equality. Empty inner pieces
contribute zero.

For $\theta$ in the interior of $I$, let
$d_\theta=\inf_{\phi\in J}\rho(\theta,\phi)>0$.
The number of nonempty pieces is at most
$C_R(1+\log_+(1/(d_\theta\vee n^{-1})))$. Minkowski's inequality
in $L^t(\mathbb Q_\theta)$ therefore gives
\[
 \E_{\mathbb Q_\theta}\mu_n(J)^t
 \le C_R\left(1+\log_+\frac1{d_\theta\vee n^{-1}}\right)^t
 \le C_R(1+|\log d_\theta|)^t.
\]
This bound is integrable over $I$, including when $I$ and $J$ share
an endpoint: near such an endpoint $d_\theta$ is the angular distance
to that endpoint, and every fixed power of its logarithm is
integrable. It follows that
\[
 \sup_{n\ge0}\E[\mu_n(I)\mu_n(J)^{p-1}]<\infty.
\]
The index $n=0$ is immediate because $\mu_0=\sigma$; arcs of zero
length also contribute zero. The terminal representation
\eqref{eq:terminal-process}, the bound $|1-e^{iu}|^q\le2^q$, and
$M_N(q)=c(q)\nu_{N-1}(q)$ from \eqref{eq:terminal-c} give pointwise
\[
 W_N(I)\le\frac{2^q}{c(q)}\mu_{N-1}(I),\qquad
 W_N(J)\le\frac{2^q}{c(q)}\mu_{N-1}(J).
\]
These inequalities prove \eqref{eq:fractional-disjoint-mixed-bound}
for $p\ge2$.

Finally, for every $p>1$ and $z_1,\ldots,z_l\in\C$,
\[
 \left|\left|\sum_{j=1}^l z_j\right|^p
                  -\sum_{j=1}^l|z_j|^p\right|
 \le C_{p,l}\sum_{i\ne j}|z_i|\,|z_j|^{p-1}.
\]
To see this, relabel so that $|z_1|$ is maximal. The mean value
bound for $|z|^p$ gives
\[
 \left|\left|\sum_{j=1}^l z_j\right|^p-|z_1|^p\right|
 \le p l^{p-1}|z_1|^{p-1}\sum_{j=2}^l|z_j|,
 \qquad
 \sum_{j=2}^l|z_j|^p
 \le |z_1|^{p-1}\sum_{j=2}^l|z_j|.
\]
Both right sides occur in the ordered pair sum, up to its
constant. Apply the resulting inequality to $z_j=c_jW_N(I_j)$
and use \eqref{eq:fractional-disjoint-mixed-bound} to obtain
\eqref{eq:fractional-additivity}.
\end{proof}

\begin{corollary}
\label{cor:weighted-fractional-critical}
For every $\beta>0$, $p>1$, $q=\sqrt{2\beta/p}$, and fixed
nonnegative Riemann integrable function $f$ on the circle,
\begin{equation}\label{eq:weighted-fractional-critical}
 \lim_{N\to\infty}
 \frac{\E\mathcal Z_N(q;f)^p}{N\log N}
 =\frac{p-1}{\Gamma(1-1/p)^p}
     C_\beta(q)^p\int_\T f^p\,d\sigma.
\end{equation}
In particular, every nonnegative H\"older function is allowed.
\end{corollary}

\begin{proof}[Proof of \cref{cor:weighted-fractional-critical}]
We first find the coefficient for equal arcs, then for arbitrary
arcs by monotonicity. Step function approximation gives the result
for $f$.

For noninteger $p$, \cref{thm:fractional-critical} gives
\[
 \frac{\E W_N(\T)^p}{\log N}\longrightarrow A_p.
\]
For integer $p$, the same limit follows from
\cref{thm:critical-moments} with $f=1$, proved below directly from
the Fisher--Hartwig estimate and the deterministic kernel bounds.
Thus the integer case of this corollary uses that later unweighted
result, whose proof is independent of the present spatial comparison.
Partition the circle into $l$ equal arcs. Rotation invariance makes
the $p$th moments of their masses equal, and
\eqref{eq:fractional-additivity} shows that each divided by $\log N$
tends to $A_p/l$. Apply \eqref{eq:fractional-additivity} to unions of consecutive
arcs. Every arc whose normalized length is rational therefore has
limit $A_p\sigma(I)$. For an arbitrary arc, inscribe and circumscribe
arcs of rational normalized lengths. Monotonicity of their masses
and then passage of those lengths to $\sigma(I)$ give
\[
 \frac{\E W_N(I)^p}{\log N}\longrightarrow A_p\sigma(I).
\]
Using \eqref{eq:fractional-additivity} once more proves the weighted
assertion for every nonnegative step function constant on finitely
many arcs.

A nonnegative Riemann integrable $f$ is bounded and has upper and lower
step approximations. Uniform continuity of $x\mapsto x^p$ on the
bounded range of these approximations implies that their integrals
of the $p$th powers tend to
$\int f^p\,d\sigma$. The masses and their $p$th moments are ordered
between the two approximations. Their limits thus give
\[
 \frac1{\log N}\E\left(\int f\chi_N\,d\sigma\right)^p
 \longrightarrow A_p\int f^p\,d\sigma.
\]
Finally \eqref{eq:onepointasymp} gives $M_N(q)^p/N\to C_\beta(q)^p$, proving
\eqref{eq:weighted-fractional-critical}.
\end{proof}

\section{Critical moments of integer order}
\label{sec:critical-moments}

We now prove the integer assertions of \cref{thm:main-critical}.
Configurations in which all $k$ points approach one another produce
the term $N\log N$. The uniform Fisher--Hartwig formula determines
its coefficient away from distances of order $N^{-1}$. The
microscopic profile determines the remaining term of order $N$;
its contribution is separated explicitly from the circular chord
correction in \cref{thm:critical-second-term}.

We first prove the leading formula for H\"older weights from two
estimates: an integral estimate for the singular kernel,
\cref{lem:critical-kernel}, proved in \cref{sec:critical-integral-estimates},
and the moment bound of \cref{lem:moments-global-bound}, valid even at
coincident points. After identifying the
next term, \cref{cor:weighted-critical-second-term} extends that
formula to the Dini condition in the main theorem.

\begin{proposition}
\label{thm:critical-moments}
Let $\beta>0$, let $k\ge2$ be an integer, and set
\begin{equation}\label{eq:critical-charge-constant}
 q_k=\sqrt{\frac{2\beta}{k}},\qquad
 A_k=\frac{k-1}{\Gamma(1-1/k)^k}.
\end{equation}
If $f:\T\to[0,\infty)$ is H\"older continuous, then
\begin{equation}\label{eq:weighted-critical-moment}
 \E_\beta\bigl[\mathcal Z_N(q_k;f)^k\bigr]
 =A_k C_\beta(q_k)^k
       \left(\int_{\T}f(\vartheta)^k\,\sigma(d\vartheta)\right)
       N\log N+O_{\beta,k,f}(N).
\end{equation}
For a fixed H\"older exponent, the error is uniform when the supremum
norm and the H\"older seminorm of $f$ are bounded.
\end{proposition}

\subsection{Reduction to two estimates}
\label{sec:critical-reduction}

The integer moment is the integral of $F_N$ over $k$ points. We
separate configurations whose minimum spacing exceeds $R/N$ from
its complement. For
$\boldsymbol\vartheta\in\T^k$, write
\[
 d_{ij}=|e^{i\vartheta_i}-e^{i\vartheta_j}|,\qquad
 \delta=\min_{i<j}d_{ij},\qquad
 K_t(\boldsymbol\vartheta)=\prod_{i<j}d_{ij}^{-t}.
\]
At the critical value $t=2/k$, abbreviate $K=K_{2/k}$ and define
\[
 K_\epsilon(\boldsymbol\vartheta)
       =\prod_{i<j}(d_{ij}\vee\epsilon)^{-2/k},
       \qquad 0<\epsilon<1.
\]
The exponents below concern ordinary Lebesgue integrability; values on
coincidence sets may be assigned arbitrarily.

The proof of \cref{thm:critical-moments} uses two estimates. The
first evaluates the logarithmic integral away from coincident points
and bounds the part within distance $O(N^{-1})$.

\begin{lemma}
\label{lem:critical-kernel}
For every fixed $R>0$ and $\kappa>0$, as $\epsilon\downarrow0$,
\begin{align}
 \int_{\{\delta\ge R\epsilon\}}K\,d\sigma^k
      &=A_k\log(1/\epsilon)+O_{k,R}(1),
                                       \label{eq:critical-kernel-log}\\
 \int_{\{\delta<R\epsilon\}}K_\epsilon\,d\sigma^k
      &=O_{k,R}(1),                     \label{eq:critical-kernel-small}\\
 \int_{\{\delta\ge R\epsilon\}}K
                    (\epsilon/\delta)^\kappa\,d\sigma^k
      &=O_{k,R,\kappa}(1).              \label{eq:critical-kernel-error}
\end{align}
With the additional factor $\prod_j f(\vartheta_j)$, where $f$ is
H\"older continuous, the coefficient in
\eqref{eq:critical-kernel-log} becomes $A_k\int f^k\,d\sigma$.
The absolute values of the two other weighted integrals are bounded
by the corresponding bounds times $\|f\|_\infty^k$.
\end{lemma}

The second estimate is the global bound of
\cref{lem:moments-global-bound}, proved in \cref{sec:global-bound}.
It controls the characteristic polynomial moments also when some
points coincide, and it is uniform in the exponents.

Together, \cref{lem:critical-kernel,lem:moments-global-bound} give
the leading coefficient from the Fisher--Hartwig formula in
\cref{thm:signed-characteristic-FH}.

\begin{proof}[Proof of \cref{thm:critical-moments}]
Since $k$ is an integer, Tonelli's theorem gives
\[
 \E_\beta[\mathcal Z_N(q_k;f)^k]
 =\int_{\T^k}\prod_{j=1}^k f(\vartheta_j)
              F_N(q_k\mathbf1,\boldsymbol\vartheta)\,d\sigma^k.
\]
Choose an admissible $\kappa>0$ in
\cref{thm:signed-characteristic-FH}, and choose its constant $R$ once
and for all. On $\{\delta\ge R/N\}$ that theorem gives
\[
 F_N(q_k\mathbf1,\boldsymbol\vartheta)
       =N C_\beta(q_k)^k K(\boldsymbol\vartheta)
                              \bigl(1+O((N\delta)^{-\kappa})\bigr),
\]
since $kq_k^2/(2\beta)=1$. By
\eqref{eq:critical-kernel-error}, the integrated error is $O(N)$.
The weighted form of \eqref{eq:critical-kernel-log} gives the leading
term in \eqref{eq:weighted-critical-moment}.
On the complementary region, \cref{lem:moments-global-bound} gives
$F_N\le C N K_{1/N}$. Its integral is $O(N)$ by
\eqref{eq:critical-kernel-small}. This proves the theorem and the
stated uniformity in $f$.
\end{proof}

\subsection{Proof of the integral estimates}
\label{sec:critical-integral-estimates}

To prove \cref{lem:critical-kernel}, we separate the common scale of
a configuration from the ratios of its consecutive gaps. The radial
integral produces the logarithm, for the following reason. The kernel
$K=\prod_{i<j}d_{ij}^{-2/k}$ is the pair product of the Fisher--Hartwig
formula at the exponents $q_j=q_k$, since $q_k^2/\beta=2/k$. After one point is fixed by rotation, there are $k-1$ relative
coordinates. Replacing circular chords by linear distances gives a
homogeneous kernel of degree $-\binom k2\cdot\frac2k=-(k-1)$,
exactly minus the dimension. Writing the configuration as a common
scale $r$ times a shape therefore produces the radial factor $dr/r$.
The chord correction is $1+O(r^2)$, as proved in
\eqref{eq:critical-chord-comparison}, and changes the integral by a
bounded amount. The resulting logarithm is $\log(1/\epsilon)$. A proper subset of $2\le j<k$ points carries
only $\binom j2\cdot\frac2k=j(j-1)/k<j-1$ powers, so the collapse of a
subcluster is integrable: only the collapse of all $k$ points produces
the logarithm. \Cref{lem:critical-shape-integrals} makes this precise
and justifies the separation of scale and shape even when only some of
the points become close.

\begin{lemma}
\label{lem:critical-shape-integrals}
Let
\[
 \Delta_{k-2}=\left\{\mathbf u\in(0,1)^{k-1}:
                                  \sum_{j=1}^{k-1}u_j=1\right\},
 \qquad
 t_1(\mathbf u)=0,\quad
 t_j(\mathbf u)=\sum_{h<j}u_h\quad(2\le j\le k).
\]
On this simplex, $d\mathbf u$ denotes Lebesgue measure in the first
$k-2$ coordinates. For $k=2$, the simplex consists of the point $u_1=1$
and its measure is one. Put
\begin{equation}\label{eq:critical-shape}
 v(\mathbf u)=\min_j u_j,\qquad
 V_k(\mathbf u)=\prod_{i<j}
                   (t_j(\mathbf u)-t_i(\mathbf u))^{-2/k}.
\end{equation}
If $k>2$, then
\begin{equation}\label{eq:shape-negative-moment}
 \int_{\Delta_{k-2}} V_k(\mathbf u)v(\mathbf u)^{-\eta}
                       \,d\mathbf u<\infty,
 \qquad 0\le\eta<\frac{k-2}{k}.
\end{equation}
For $k=2$, this integral equals one for every $\eta\ge0$.
More generally, for $k>2$ and $t,\eta\ge0$,
\[
 \int_{\Delta_{k-2}}
       \prod_{i<j}(t_j(\mathbf u)-t_i(\mathbf u))^{-t}
       v(\mathbf u)^{-\eta}\,d\mathbf u<\infty
\]
whenever
\[
 \eta<\min_{1\le l\le k-2}
                  l\left(1-\frac{(l+1)t}{2}\right).
\]
The bound is uniform on compact subsets of this strict parameter
range. The same conclusion holds with a fixed product of nonnegative
powers of $1+|\log u_j|$ and $1+|\log(t_j-t_i)|$ inserted. For $k=2$, the shape integral is one
for every $t,\eta\ge0$.
For the circular integral, if $0\le t<2/k$ and
\begin{equation}\label{eq:subcritical-collision-margin}
 0\le\eta<\lambda_k(t):=
             \min_{2\le j\le k}(j-1)\left(1-\frac{jt}{2}\right),
\end{equation}
then
\begin{equation}\label{eq:subcritical-delta-moment}
 \int_{\T^k}K_t(\boldsymbol\vartheta)
             \delta(\boldsymbol\vartheta)^{-\eta}
                         \,d\sigma^k(\boldsymbol\vartheta)<\infty.
\end{equation}
\end{lemma}

\begin{proof}[Proof of \cref{lem:critical-shape-integrals}]
We order the gaps and use their successive ratios to separate the
possible collisions. Integrability then reduces to checking the power
of each ratio; a final radial integral gives the corresponding statement
on the circle. Divide the simplex
into the finitely many regions on which the gaps have a prescribed
order, and write the ordered gaps as
$u_{\pi(1)}\le\cdots\le u_{\pi(n)}$, where $n=k-1$. Introduce
\[
 z_l=\frac{u_{\pi(l)}}{u_{\pi(l+1)}}\quad(1\le l<n),\qquad
 w_n=1,\quad w_j=\prod_{h=j}^{n-1}z_h,
 \qquad S=\sum_{j=1}^n w_j.
\]
Then $u_{\pi(j)}=w_j/S$, $1\le S\le n$, and the Jacobian is
\begin{equation}\label{eq:gap-ratio-jacobian}
 d\mathbf u=S^{-n}\prod_{l=1}^{n-1}z_l^{l-1}
                                      \,dz_1\cdots dz_{n-1}.
\end{equation}
Every difference $t_j-t_i$ is a sum of consecutive gaps. It is bounded
below by its largest gap and above by $n$ times that gap. Consequently,
on the chosen region, the product
$\prod_{i<j}(t_j-t_i)^{-t}$ is comparable to the product obtained by
replacing each difference by its largest constituent gap.

Let $P_l$ be the number of pairs joined using only the $l$ smallest
gaps. These gaps form disjoint consecutive components. If their
numbers of vertices are $c_1,\ldots,c_a$, then
\[
 \sum_{h=1}^a(c_h-1)=l,\qquad
 P_l=\sum_{h=1}^a\frac{c_h(c_h-1)}2
                            \le\frac{l(l+1)}2.
\]
The last inequality follows by successively joining components: for
nonnegative integers $d_h=c_h-1$ with sum $l$, one has
$\sum_h d_h(d_h+1)\le l(l+1)$.
The exponent of $z_l$ in the resulting integrand, including the
Jacobian and $v^{-\eta}$, is $l-1-tP_l-\eta$. Thus its integral is
finite whenever
\begin{equation}\label{eq:gap-integrability-margin}
 l-tP_l-\eta>0\qquad(1\le l\le k-2).
\end{equation}
Since $P_l\le l(l+1)/2$, the general condition in the statement
implies \eqref{eq:gap-integrability-margin} on every gap ordering.
On a compact subset of that strict range, the exponents have a
common positive margin. Products of fixed powers of $|\log z_l|$
remain integrable against these powers of $z_l$. Each gap is a
product of some $z_l$, divided by $S\in[1,k-1]$, and every difference
$t_j-t_i$ is comparable to its largest constituent gap. Thus the same
margin also controls logarithmic factors of the gaps and differences.
This proves the general assertion and its local uniformity.
For $t=2/k$,
\[
 l-\frac2kP_l\ge\frac{l(k-l-1)}k
                         \ge\frac{k-2}{k},
\]
which proves \eqref{eq:shape-negative-moment}. When $k=2$ there are no
ratio variables, so the stated convention gives the assertion directly.

For the circular integral, cut the configuration at a largest cyclic
gap and order the remaining points. Apart from a null set, this gives
finitely many coordinate regions with consecutive gaps
$x_1,\ldots,x_{k-1}>0$. Write $x_j=ru_j$. The removed gap is $2\pi-r$,
and being largest means
\begin{equation}\label{eq:gap-chart-radius}
 0<r<b(\mathbf u):=\frac{2\pi}{1+\max_j u_j}
                 \le\frac{2\pi(k-1)}k.
\end{equation}
In particular, every chord is comparable, with constants depending
only on $k$, to its linear difference $r(t_j-t_i)$, and the minimum
chord is comparable to $rv$. The radial integral in
\eqref{eq:subcritical-delta-moment} is therefore bounded by a constant
times
\[
 \int_0^{2\pi(k-1)/k}
                 r^{k-2-tk(k-1)/2-\eta}\,dr.
\]
It is finite under the $j=k$ inequality in
\eqref{eq:subcritical-collision-margin}. For the shape integral,
\eqref{eq:gap-integrability-margin} follows from its inequalities
with $j=l+1$. This proves \eqref{eq:subcritical-delta-moment}.
\end{proof}

The same coordinates give an exact integration formula. If a function
of $k$ points is symmetric and invariant under common rotation, then
its integral against $\sigma^k$ is
\begin{equation}\label{eq:exact-gap-integration}
 c_k\int_{\Delta_{k-2}}\int_0^{b(\mathbf u)}
   G\bigl(r\mathbf t(\mathbf u)\bigr)
                  r^{k-2}\,dr\,d\mathbf u,
 \qquad c_k=\frac{k!}{(2\pi)^{k-1}}.
\end{equation}
Indeed, after fixing the first point by rotation there are $(k-1)!$
orders for the remaining labels. Cyclically shifting the gap vector
partitions the gap simplex into $k$ equal regions according to the
position of its largest gap. The change $\mathbf x=r\mathbf u$ then
has Jacobian $r^{k-2}$.
For a common weight $f$, its contribution in these coordinates is
\begin{equation}\label{eq:critical-averaged-weight}
 w_f(r,\mathbf u)=\int_{\T}
       \prod_{j=1}^k f\bigl(\vartheta+rt_j(\mathbf u)\bigr)
                                             \,\sigma(d\vartheta).
\end{equation}
If $f$ has H\"older exponent $\alpha>0$, then uniformly in $\mathbf u$,
\begin{equation}\label{eq:critical-weight-diagonal}
 w_f(r,\mathbf u)=\int_{\T}f^k\,d\sigma+O_f(r^\alpha).
\end{equation}
This follows by replacing the factors one at a time and using
$0\le t_j\le1$.

\begin{proof}[Proof of \cref{lem:critical-kernel}]
We first bound the contribution of small gaps and the integrated error,
then identify the logarithmic coefficient using Dyson's integral.
In the gap coordinates of \eqref{eq:exact-gap-integration}, the
correction from linear distances to circular chords is
\begin{equation}\label{eq:critical-chord-ratio}
 D_k(r,\mathbf u)=\prod_{i<j}
  \left\{\frac{2\sin\bigl(r(t_j-t_i)/2\bigr)}
                         {r(t_j-t_i)}\right\}^{-2/k}.
\end{equation}
The ratio in braces is interpreted as one when its argument is zero.
By \eqref{eq:gap-chart-radius}, it is bounded above and below by
positive constants depending only on $k$. Moreover,
\begin{equation}\label{eq:critical-chord-comparison}
 K(r\mathbf t)=r^{-(k-1)}V_k(\mathbf u)D_k(r,\mathbf u),\qquad
 |D_k(r,\mathbf u)-1|\le C_k r^2.
\end{equation}
The smallest circular gap is $rv$, since the removed gap is largest.
As the smallest gap is at most $2\pi/k$, we have the exact formula
\begin{equation}\label{eq:critical-minimum-chord}
 \delta(r\mathbf t)=2\sin(rv/2)\asymp_k rv.
\end{equation}

For $k>2$, choose $0<\eta<(k-2)/k$. By
\cref{lem:critical-shape-integrals},
\begin{equation}\label{eq:shape-tail-bound}
 \int_{\{v<h\}}V_k\,d\mathbf u\le C_k h^\eta,
 \qquad
 \int_{\Delta_{k-2}}V_k|\log v|\,d\mathbf u<\infty.
\end{equation}
For $k=2$, the same conclusions hold with any $\eta>0$, after changing
the constant, because $v=1$.

If $r<\epsilon$, then $K_\epsilon\le\epsilon^{-(k-1)}$, whereas
this region has volume $O_k(\epsilon^{k-1})$. For $r\ge\epsilon$,
we use $K_\epsilon\le K$ and
\eqref{eq:critical-minimum-chord}. The part where
$\delta<R\epsilon$ is bounded by
\[
 C_k\int_\epsilon^{2\pi(k-1)/k}\frac{dr}{r}
        \int_{\{v<C_kR\epsilon/r\}}V_k\,d\mathbf u
 \le C_{k,R}\int_\epsilon^{2\pi(k-1)/k}
                   \left(\frac\epsilon r\right)^\eta\frac{dr}{r}
 \le C_{k,R}.
\]
This proves \eqref{eq:critical-kernel-small}, including regions where several pairs of points are simultaneously close.

On $\{\delta\ge R\epsilon\}$, the radial variable satisfies
$r\ge c_k'R\epsilon/v$. Hence the left side of
\eqref{eq:critical-kernel-error} is at most
\[
 C_k\int_{\Delta_{k-2}} V_k(\mathbf u)
       \int_{c_k'R\epsilon/v}^{b(\mathbf u)}
             \left(\frac{C_k\epsilon}{rv}\right)^\kappa
                               \frac{dr}{r}\,d\mathbf u
 \le C_{k,R,\kappa}\int_{\Delta_{k-2}}V_k\,d\mathbf u.
\]
An inner integral is set to zero when its lower endpoint exceeds its
upper endpoint. In particular, no inverse moment of $v$ of order
$\kappa$ is needed.

To evaluate the logarithmic coefficient, put
$a_\epsilon=2\arcsin(R\epsilon/2)$, for small enough $\epsilon$.
The condition $\delta\ge R\epsilon$ is exactly $rv\ge a_\epsilon$.
Replacing $D_k$ by one changes the integral by $O_k(1)$, by
\eqref{eq:critical-chord-comparison}. Thus the integral equals
\[
 c_k\int_{\Delta_{k-2}}V_k(\mathbf u)
             \left[\log\frac{b(\mathbf u)v}{a_\epsilon}\right]_+
                     d\mathbf u+O_k(1).
\]
The function $b$ is bounded above and below by positive constants,
$a_\epsilon=R\epsilon(1+O_R(\epsilon^2))$, and
\eqref{eq:shape-tail-bound} controls $\log v$. The coefficient of
$\log(1/\epsilon)$ is consequently
$c_k\int V_k\,d\mathbf u$.

We identify this integral from Dyson's circular integral, in its
integrable range; see
\cite[equations~(1.10) and~(1.12)]{ForresterWarnaar2008}. For $0\le t<2/k$,
\begin{equation}\label{eq:dyson-negative-power}
 I_k(t):=\int_{\T^k}K_t\,d\sigma^k
       =\frac{\Gamma(1-kt/2)}{\Gamma(1-t/2)^k}.
\end{equation}
Write $M=k(k-1)/2$ and $t=2/k-h$. The radial exponent in the gap
coordinates is $-1+hM$. The shape integrals converge to
$\int V_k$ as $h\downarrow0$, by the strict inequalities in
\eqref{eq:gap-integrability-margin}; bounded factors and their
logarithms are controlled by a slightly smaller positive margin.
The chord correction contributes a bounded quantity uniformly for
small $h$. It follows that
\[
 \lim_{h\downarrow0}h I_k(2/k-h)
          =\frac{c_k}{M}\int_{\Delta_{k-2}}V_k\,d\mathbf u.
\]
On the other hand, the gamma quotient in
\eqref{eq:dyson-negative-power} gives
\[
 \lim_{h\downarrow0}h I_k(2/k-h)
       =\frac{2}{k\Gamma(1-1/k)^k}.
\]
Therefore
\begin{equation}\label{eq:critical-shape-normalization}
 c_k\int_{\Delta_{k-2}}V_k\,d\mathbf u
                    =\frac{k-1}{\Gamma(1-1/k)^k}=A_k,
\end{equation}
which proves \eqref{eq:critical-kernel-log}. Finally,
\eqref{eq:critical-weight-diagonal} changes the weighted radial
integral by an integrable term $O_f(r^{\alpha-1})$ and gives the
weighted coefficient.
\end{proof}

\subsection{The next term at the critical moment order}
\label{sec:critical-second-term}

The $O(N)$ term in \cref{thm:critical-moments} contains the contribution
of configurations at the mean spacing. The compact convergence and
large distance matching of the microscopic profile identify its limit.
In the following statement, all shape integrals use the conventions
of \cref{lem:critical-shape-integrals}.

\begin{theorem}
\label{thm:critical-second-term}
Let $\beta>0$, let $k\ge2$ be an integer, and let $q_k,A_k$ be as in
\eqref{eq:critical-charge-constant}. Define
\[
 C_k=C_\beta(q_k)^k,\qquad
 H_{\beta,k}(y,\mathbf u)
       =\mathcal H_{\beta,q_k\mathbf1}
                           \bigl(y\mathbf t(\mathbf u)\bigr),
\]
where $\mathcal H$ is the microscopic profile of
\cref{thm:positive-merging}. Put
\begin{align}
 B_{\beta,k}^{\mathrm{mic}}
  =c_k\int_{\Delta_{k-2}}\bigg\{
       &\int_0^1 y^{k-2}H_{\beta,k}(y,\mathbf u)\,dy\notag\\
       &+\int_1^\infty
          \left(y^{k-2}H_{\beta,k}(y,\mathbf u)
                         -\frac{C_kV_k(\mathbf u)}y\right)dy
                         \bigg\}\,d\mathbf u,
                                      \label{eq:critical-mic-constant}\\
 B_{\beta,k}^{\mathrm{circ}}
  =c_k C_k\int_{\Delta_{k-2}}V_k(\mathbf u)
       &\left\{\log b(\mathbf u)
          +\int_0^{b(\mathbf u)}
                    \frac{D_k(r,\mathbf u)-1}{r}\,dr\right\}
                                     d\mathbf u.
                                      \label{eq:critical-circle-constant}
\end{align}
All the integrals in these definitions converge absolutely. With
$B_{\beta,k}=B_{\beta,k}^{\mathrm{mic}}+
B_{\beta,k}^{\mathrm{circ}}$, one has
\begin{equation}\label{eq:critical-next-term}
 \E_\beta[\mathcal Z_N(q_k;1)^k]
       =A_k C_\beta(q_k)^k N\log N+B_{\beta,k}N+o(N).
\end{equation}
\end{theorem}

The subtraction in \eqref{eq:critical-mic-constant} removes precisely
the nonintegrable $y^{-1}$ term of the profile at large distances.
The second integral accounts for the difference between linear distances
and circular chords. Their sum is the coefficient of $N$ and is
independent of the auxiliary division at $y=1$, provided both
subtractions are changed consistently: division at $y=a>0$ adds
$A_kC_k\log a$ to the microscopic part and subtracts the same quantity
from the circular part.

For $k=2$, the circular term can be evaluated directly. Here
$c_2=1/\pi$, $b=\pi$, $V_2=1$ and
$D_2(r)=r/(2\sin(r/2))$. Hence
\[
 \log\pi+\int_0^\pi
       \left(\frac1{2\sin(r/2)}-\frac1r\right)dr=\log4,
 \qquad
 B_{\beta,2}^{\mathrm{circ}}
       =\frac{C_\beta(\sqrt\beta)^2\log4}{\pi}.
\]
This also checks the normalization of the circular subtraction.

\begin{proof}[Proof of \cref{thm:critical-second-term}]
We subtract the common $1/y$ tail before passing to the microscopic
limit. An integrable bound uniform in the shape justifies this limit;
the subtracted term can then be integrated explicitly.
Within this proof, write $H=H_{\beta,k}$, since $\beta$ and $k$ are
fixed. The integrability of $V_k$ in \cref{lem:critical-shape-integrals}
and the chord estimate \eqref{eq:critical-chord-comparison} imply
absolute convergence of \eqref{eq:critical-circle-constant}.
To control the profile and the finite $N$ kernels with the same
bounds, put
\[
 T_N(y,\mathbf u)=N^{-k}
       F_N\left(q_k\mathbf1,\frac yN\mathbf t(\mathbf u)\right),
             \qquad 0<y<Nb(\mathbf u).
\]
The normalization is $N^k$ because the combined exponent is $kq_k$ and
$(kq_k)^2/(2\beta)=k$. By \cref{thm:positive-merging},
\begin{equation}\label{eq:critical-profile-compact}
 T_N(y,\mathbf u)\longrightarrow H(y,\mathbf u)
\end{equation}
uniformly when $y$ ranges over a compact interval and $\mathbf u$
ranges over the closed simplex.

The moment bound \eqref{eq:moments-global-bound} and the chord
comparison \eqref{eq:critical-chord-comparison} give, uniformly on the
whole integration domain,
\begin{equation}\label{eq:critical-scaled-global}
 T_N(y,\mathbf u)
       \le C\prod_{i<j}\bigl(1\vee y(t_j-t_i)\bigr)^{-2/k}.
\end{equation}
In particular, $T_N$ is bounded for $0<y\le1$, while for $y\ge1$,
\begin{equation}\label{eq:critical-scaled-dominance}
 y^{k-2}T_N(y,\mathbf u)\le\frac{C V_k(\mathbf u)}y.
\end{equation}
The same bounds hold for $H$ by passage to the limit.

Choose $\kappa>0$ from \cref{thm:signed-characteristic-FH}.
By \eqref{eq:critical-minimum-chord}, there is $R_1<\infty$ such that
$yv\ge R_1$ implies that this theorem applies, uniformly in $N$.
On that region,
\begin{equation}\label{eq:critical-scaled-tail-error}
 \left|y^{k-2}T_N(y,\mathbf u)
     -\frac{C_kV_k(\mathbf u)D_k(y/N,\mathbf u)}y\right|
 \le\frac{C V_k(\mathbf u)}y(yv)^{-\kappa}.
\end{equation}
On $yv<R_1$, use \eqref{eq:critical-scaled-dominance} and boundedness
of $D_k$ instead. Combining the two bounds, and increasing the
constant on the fixed transition region, gives for $y\ge1$
\begin{equation}\label{eq:critical-tail-envelope}
 \left|y^{k-2}T_N(y,\mathbf u)
     -\frac{C_kV_k(\mathbf u)D_k(y/N,\mathbf u)}y\right|
 \le\frac{C V_k(\mathbf u)}y
                           \min\{1,(yv)^{-\kappa}\}.
\end{equation}
For $k>2$, choose
$0<\gamma<\min\{\kappa,(k-2)/k\}$. For $k=2$, choose
$0<\gamma<\kappa$. Since
$\min\{1,x^{-\kappa}\}\le x^{-\gamma}$ for $x>0$,
\cref{lem:critical-shape-integrals} implies
\begin{equation}\label{eq:critical-integrated-tail}
 \int_{\Delta_{k-2}}
      \frac{V_k(\mathbf u)}y\min\{1,(yv)^{-\kappa}\}
                                        \,d\mathbf u
                   \le C y^{-1-\gamma}\qquad(y\ge1).
\end{equation}
This estimate also covers shapes with a small internal gap. Taking
$N\to\infty$ in \eqref{eq:critical-tail-envelope} proves absolute
convergence of the second integral in
\eqref{eq:critical-mic-constant}. Its first integral converges by the
bounded small $y$ estimate.

By \eqref{eq:exact-gap-integration} and the change of variables $y=Nr$,
\[
 \frac1N\E_\beta[\mathcal Z_N(q_k;1)^k]
 =c_k\int_{\Delta_{k-2}}\int_0^{Nb(\mathbf u)}
                     y^{k-2}T_N(y,\mathbf u)\,dy\,d\mathbf u.
\]
Subtract inside this integral
$C_kV_k(\mathbf u)D_k(y/N,\mathbf u)\mathbf1_{\{y\ge1\}}/y$.
On each bounded $y$ interval, the subtracted integrand converges by
\eqref{eq:critical-profile-compact} and $D_k(y/N,\mathbf u)\to1$.
Equations \eqref{eq:critical-scaled-global} and
\eqref{eq:critical-integrated-tail} give an integrable majorant,
independent of $N$, for this convergence. Extending the integrand by
zero when $y\ge Nb(\mathbf u)$, dominated convergence therefore gives
$B_{\beta,k}^{\mathrm{mic}}$.

The integral of the subtracted term can be calculated directly:
\begin{align*}
 &c_k C_k\int_{\Delta_{k-2}}V_k(\mathbf u)
          \int_1^{Nb(\mathbf u)}D_k(y/N,\mathbf u)\frac{dy}{y}
                                                       \,d\mathbf u\\
 &\quad=c_k C_k\int_{\Delta_{k-2}}V_k(\mathbf u)
       \left\{\log N+\log b(\mathbf u)
             +\int_{1/N}^{b(\mathbf u)}
                      \frac{D_k(r,\mathbf u)-1}{r}\,dr\right\}
                                                        d\mathbf u\\
 &\quad=A_k C_k\log N+B_{\beta,k}^{\mathrm{circ}}+O(N^{-2}),
\end{align*}
where we used \eqref{eq:critical-shape-normalization} and the uniform
bound $D_k-1=O(r^2)$. Combining the two limits proves
\eqref{eq:critical-next-term}.
\end{proof}

The same subtraction identifies the next coefficient for spatial
weights with substantially less regularity. This includes sharp
restrictions to arcs. The microscopic part depends only on the value
of the product of weights on the diagonal; the circular part also
accounts for the variation of the weight away from it.
For $f\in L^k(\T,\sigma)$, define its translation modulus by
\[
 \omega_k(f,r)=\sup_{|h|\le r}
       \|f(\cdot+h)-f\|_{L^k(\sigma)},\qquad 0<r\le1.
\]

\begin{corollary}
\label{cor:weighted-critical-second-term}
Let $\beta>0$ and let $k\ge2$ be an integer. Suppose $f\ge0$ belongs
to $L^k(\T,\sigma)$ and satisfies
\begin{equation}\label{eq:critical-dini-modulus}
 \int_0^1\omega_k(f,r)\,\frac{dr}{r}<\infty.
\end{equation}
Put $w_0=\int_\T f^k\,d\sigma$ and let $w_f$ be defined by
\eqref{eq:critical-averaged-weight}. With the notation of
\cref{thm:critical-second-term}, define
\begin{align}
 B_{\beta,k}[f]
 ={}&w_0 B_{\beta,k}^{\mathrm{mic}}\notag\\
 &+c_kC_k\int_{\Delta_{k-2}}V_k(\mathbf u)
   \left\{w_0\log b(\mathbf u)
       +\int_0^{b(\mathbf u)}
          \frac{D_k(r,\mathbf u)w_f(r,\mathbf u)-w_0}{r}\,dr
   \right\}\,d\mathbf u.
                                      \label{eq:weighted-critical-B}
\end{align}
The integrals converge absolutely, and
\begin{equation}\label{eq:weighted-critical-second-term}
 \E_\beta[\mathcal Z_N(q_k;f)^k]
  =A_k C_\beta(q_k)^k w_0 N\log N
                         +B_{\beta,k}[f]N+o(N).
\end{equation}
In particular, $B_{\beta,k}[1]=B_{\beta,k}$. The conclusion applies
to every nonnegative bounded function of bounded variation, including
the indicator of a finite union of arcs. It also gives the $O(N)$
remainder in \eqref{eq:weighted-critical-moment} under
\eqref{eq:critical-dini-modulus}.
\end{corollary}

\begin{proof}[Proof of \cref{cor:weighted-critical-second-term}]
We control the change of the averaged weight near the diagonal by its
translation modulus. This is the only additional estimate needed to
repeat the subtraction in the proof of \cref{thm:critical-second-term}.
The generalized H\"older inequality gives
\[
 0\le w_f(r,\mathbf u)\le\|f\|_k^k.
\]
Replacing the $k$ translated factors one at a time and applying
H\"older's inequality to each difference gives
\begin{equation}\label{eq:critical-dini-weight-bound}
 |w_f(r,\mathbf u)-w_0|
 \le k\|f\|_k^{k-1}\omega_k(f,r),\qquad 0<r\le1.
\end{equation}
The estimate is uniform on the closed shape simplex, since
$0\le t_j(\mathbf u)\le1$. Translations are continuous in $L^k$,
so $w_f(r,\mathbf u)\to w_0$ uniformly in $\mathbf u$ as
$r\downarrow0$. By \eqref{eq:critical-chord-comparison},
\[
 |D_k(r,\mathbf u)w_f(r,\mathbf u)-w_0|
                  \le C_{k,f}\{r^2+\omega_k(f,r)\},\qquad 0<r\le1.
\]
The right side is integrable against $dr/r$ by
\eqref{eq:critical-dini-modulus}. On the remaining interval
$1\le r\le b(\mathbf u)$ the factors are uniformly bounded.
Since $V_k$ is integrable by \cref{lem:critical-shape-integrals},
this proves absolute convergence of the
circular integral in \eqref{eq:weighted-critical-B}.

The gap formula \eqref{eq:exact-gap-integration} is valid for measurable
nonnegative $f$ by Tonelli's theorem. We continue to use
$T_N(y,\mathbf u)=N^{-k}F_N(q_k\mathbf1,y\mathbf t(\mathbf u)/N)$,
the rescaled moment from the proof of \cref{thm:critical-second-term}.
Then
\[
 \frac1N\E_\beta[\mathcal Z_N(q_k;f)^k]
 =c_k\int_{\Delta_{k-2}}\int_0^{Nb(\mathbf u)}
       y^{k-2}T_N(y,\mathbf u)w_f(y/N,\mathbf u)
                                             \,dy\,d\mathbf u.
\]
As in the proof of \cref{thm:critical-second-term}, subtract
$C_kV_k(\mathbf u)D_k(y/N,\mathbf u)\mathbf1_{\{y\ge1\}}/y$,
now multiplied by $w_f(y/N,\mathbf u)$. This multiplier is bounded by
$\|f\|_k^k$ and, by \eqref{eq:critical-dini-weight-bound}, converges
to $w_0$ uniformly on bounded $y$ intervals and on the closed
shape simplex. The integrable bounds in
\eqref{eq:critical-scaled-global} and
\eqref{eq:critical-tail-envelope} therefore justify dominated
convergence exactly as before, including the regions with a small
internal gap. The integral of the difference after this subtraction
converges to $w_0 B_{\beta,k}^{\mathrm{mic}}$.

The baseline integral equals
\[
 c_k C_k\int_{\Delta_{k-2}}V_k(\mathbf u)
 \left\{w_0\log N+w_0\log b(\mathbf u)
       +\int_{1/N}^{b(\mathbf u)}
         \frac{D_k(r,\mathbf u)w_f(r,\mathbf u)-w_0}{r}\,dr
 \right\}\,d\mathbf u.
\]
The omitted integral over $(0,1/N)$ is bounded by
\[
 C_{k,f}\int_0^{1/N}\{r^2+\omega_k(f,r)\}\,\frac{dr}{r}=o(1).
\]
By \eqref{eq:critical-shape-normalization}, the coefficient of
$\log N$ is $A_kC_kw_0$. Combining this calculation with the preceding
limit proves \eqref{eq:weighted-critical-second-term}.

Finally, for bounded $f$ of bounded variation on the circle,
\[
 \|f(\cdot+h)-f\|_k^k
 \le (2\|f\|_\infty)^{k-1}\|f(\cdot+h)-f\|_1
 \le C_f|h|.
\]
Hence $\omega_k(f,r)\le C_f r^{1/k}$, which proves
\eqref{eq:critical-dini-modulus} for the stated examples.
\end{proof}

For the leading coefficient alone, an approximation argument removes
all regularity assumptions on the spatial weight beyond its natural
integrability.

\begin{proposition}
\label{prop:critical-lp-weights}
Let $p>1$, $\beta>0$, and $q=\sqrt{2\beta/p}$. For every complex
$f\in L^p(\T,\sigma)$,
\[
 \lim_{N\to\infty}
 \frac{\E|\mathcal Z_N(q;f)|^p}{N\log N}
 =\frac{(p-1)C_\beta(q)^p}{\Gamma(1-1/p)^p}
                                  \int_\T |f|^p\,d\sigma.
\]
\end{proposition}

\begin{proof}[Proof of \cref{prop:critical-lp-weights}]
We first bound the map from a spatial weight to its random integral.
The arc asymptotic and complex additivity then give the limit for
step functions, and this bound permits passage to $L^p$.
Write $\chi_N=|P_N|^q/M_N(q)$ for the normalized density and
$W_N=\int\chi_N\,d\sigma$ for its total mass. H\"older's inequality
with respect to the finite measure $\chi_N\,d\sigma$ gives
\[
 \left|\int g\chi_N\,d\sigma\right|^p
 \le W_N^{p-1}\int|g|^p\chi_N\,d\sigma.
\]
Rotation invariance and Tonelli's theorem imply
\begin{equation}\label{eq:critical-lp-operator-bound}
 \E\left|\int g\chi_N\,d\sigma\right|^p
 \le (\E W_N^p)\|g\|_p^p.
\end{equation}
Indeed $\E[W_N^{p-1}\chi_N(\theta)]$ is independent of $\theta$,
and its integral equals $\E W_N^p$. The unweighted limits in
\cref{thm:fractional-critical,thm:critical-moments} give
$\E W_N^p/\log N\to A_p$ for every real $p>1$.

Let $f=\sum_{j=1}^l c_j\mathbf1_{I_j}$, where the arcs form a finite
partition and $c_j\in\C$. By \eqref{eq:fractional-additivity},
\[
 \E\left|\int f\chi_N\,d\sigma\right|^p
 =\sum_{j=1}^l|c_j|^p\E W_N(I_j)^p+O(1).
\]
Apply \cref{cor:weighted-fractional-critical} to each arc indicator,
and use $M_N(q)^p/N\to C_\beta(q)^p$. We obtain
\[
 \frac1{\log N}\E\left|\int f\chi_N\,d\sigma\right|^p
 \longrightarrow A_p\sum_{j=1}^l|c_j|^p\sigma(I_j)
 =A_p\|f\|_p^p.
\]

For general complex $f\in L^p$, choose complex step functions $f_l$
on finite arc partitions such that $\|f_l-f\|_p\to0$.
Minkowski's inequality in the probability space and
\eqref{eq:critical-lp-operator-bound} give
\[
 \left|
 \left(\frac{\E|\int f\chi_N\,d\sigma|^p}{\log N}\right)^{1/p}
 -\left(\frac{\E|\int f_l\chi_N\,d\sigma|^p}{\log N}\right)^{1/p}
 \right|
 \le\left(\frac{\E W_N^p}{\log N}\right)^{1/p}\|f-f_l\|_p.
\]
First let $N\to\infty$ with $l$ fixed, then let $l\to\infty$.
The coefficient on the right converges to $A_p^{1/p}$, and
$\|f_l\|_p\to\|f\|_p$. Thus the normalized random integral has the
claimed limit. Finally multiply by
$M_N(q)^p/N\to C_\beta(q)^p$ to obtain the stated formula.
\end{proof}

\subsection{The transition through the critical moment}
\label{sec:critical-window}

The critical logarithm can be followed through a window of width
$1/\log N$ in the exponent. This gives one formula on both sides of the
moment threshold. We keep the integer moment order fixed; the result
below does not use analytic continuation in that order.

Fix $k\ge2$ and write $q_* =\sqrt{2\beta/k}$. For $s$ in a
neighborhood of one put
\begin{align*}
 q(s)&=q_*\sqrt{s},&
 V_{k,s}(\mathbf u)&=\prod_{i<j}(t_j(\mathbf u)-t_i(\mathbf u))^{-2s/k},\\
 S_k(s)&=c_k\int_{\Delta_{k-2}}V_{k,s}(\mathbf u)\,d\mathbf u,&
 J_{\beta,k}(s)&=C_\beta(q(s))^kS_k(s).
\end{align*}
Apply the general shape criterion of \cref{lem:critical-shape-integrals}
with $t=2s/k$. At $s=1$, its margins are at least $(k-2)/k>0$ when
$k>2$; for $k=2$ the shape space is a point. Hence the shape integrals
and their integrands multiplied by fixed powers of logarithms of gaps
are locally uniformly integrable for $s$ near one. Differentiation
under the integral is therefore valid there. In particular,
$S_k(1)=A_k$ and $J_{\beta,k}(1)=A_kC_\beta(q_*)^k$.
The shape factor also has a Selberg evaluation. For $k=2$ it is
$S_2(s)=1/\pi$. For $k\ge3$,
\begin{equation}\label{eq:window-shape-selberg}
 S_k(s)=\frac{c_k}{(k-2)!}
 \prod_{j=0}^{k-3}
 \frac{\Gamma(1-(j+2)s/k)^2\Gamma(1-(j+1)s/k)}
 {\Gamma(2-(k+j+1)s/k)\Gamma(1-s/k)}.
\end{equation}
Indeed $t_2,\ldots,t_{k-1}$ range over the ordered simplex in $(0,1)$;
the endpoint exponents are $-2s/k$ and the exponent between interior
points is also $-2s/k$. Selberg's formula with $k-2$ variables,
parameters $1-2s/k,1-2s/k,-s/k$, gives
\eqref{eq:window-shape-selberg}; see
\cite[equations~(1.1) and~(1.2)]{ForresterWarnaar2008}.
All integrability inequalities hold for $s$ sufficiently close to one.

Define the positive continuous function
\begin{equation}\label{eq:critical-window-function}
 \Psi_k(x)=e^x\int_0^1 e^{(k-1)xt}\,dt
 =\begin{cases}
 \displaystyle\frac{e^{kx}-e^x}{(k-1)x},&x\ne0,\\
 1,&x=0.
 \end{cases}
\end{equation}
Let $B^{\rm mic}_{\beta,k}$ be as in
\eqref{eq:critical-mic-constant}, and for a nonnegative
$f\in L^k(\T,\sigma)$ satisfying \eqref{eq:critical-dini-modulus} put
\[
 w_0=\int_\T f^k\,d\sigma,\qquad
 B^{\rm circ}_{\beta,k}[f]
       =B_{\beta,k}[f]-w_0B^{\rm mic}_{\beta,k}.
\]

\begin{theorem}
\label{thm:critical-window}
Let $\beta>0$, let $k\ge2$ be an integer, and let $f\ge0$ be a fixed
function in $L^k(\T,\sigma)$ satisfying
\eqref{eq:critical-dini-modulus}. Set
\[
 s_N(x)=1+\frac{x}{\log N},\qquad q_N(x)=q_*\sqrt{s_N(x)}.
\]
For every fixed $A<\infty$, uniformly for $|x|\le A$, as $N\to\infty$,
\begin{align}\label{eq:critical-window-two-term}
 \E_\beta\bigl[\mathcal Z_N(q_N(x);f)^k\bigr]
 ={}&N\log N\,w_0J_{\beta,k}(1)\Psi_k(x)\notag\\
 &+N\left\{e^{kx}w_0B^{\rm mic}_{\beta,k}
             +e^xB^{\rm circ}_{\beta,k}[f]
             +w_0x\Psi_k(x)J'_{\beta,k}(1)\right\}+o(N).
\end{align}
Equivalently, the first line may use $J_{\beta,k}(s_N(x))$ in place
of $J_{\beta,k}(1)$ if the term containing $J'_{\beta,k}(1)$ is
removed. At $x=0$ this is
\eqref{eq:weighted-critical-second-term}.
\end{theorem}

The microscopic part of \cref{thm:critical-window} requires control
as the exponent varies. \Cref{lem:window-charge-variation} provides
this control by comparing nearby circular Jacobi laws.

\begin{lemma}
\label{lem:window-charge-variation}
Let $q_N>0$ satisfy $q_N\to q_*>0$ and
$(q_N-q_*)\sqrt{\log N}\to0$. For a fixed integer $k\ge2$,
\[
 N^{-(kq_N)^2/(2\beta)}
 F_N\left(q_N\mathbf1,\frac{\mathbf y}{N}\right)
 \longrightarrow\mathcal H_{\beta,q_*\mathbf1}(\mathbf y)
\]
locally uniformly in $\mathbf y\in\mathbb R^k$.
\end{lemma}

\begin{proof}[Proof of \cref{lem:window-charge-variation}]
We compare the varying tilt with a fixed circular Jacobi law in
$L^r$, and then apply the locally uniform limit of the random analytic
functions under that fixed law.
Put $Q_*=kq_*$, $\Delta_N=k(q_N-q_*)$, and let $\mathbb Q_N$ be the
circular Jacobi law obtained by tilting by $|P_N(1)|^{Q_*}$.
The tilt with exponent $Q_*+\Delta_N$ has Radon--Nikodym derivative
\[
 L_N=|P_N(1)|^{\Delta_N}\frac{M_N(Q_*)}{M_N(Q_*+\Delta_N)}
\]
relative to $\mathbb Q_N$. For every fixed $r>1$,
\[
 \E_{\mathbb Q_N}L_N^r
 =\frac{M_N(Q_*+r\Delta_N)M_N(Q_*)^{r-1}}
          {M_N(Q_*+\Delta_N)^r}\longrightarrow1.
\]
To check the limit, use \eqref{eq:onepointasymp} uniformly near $Q_*$.
The coefficient of $\log N$ in the logarithm of the displayed ratio is
$r(r-1)\Delta_N^2/(2\beta)$; the logarithms of the constant factors
also tend to zero. Since $\E L_N=1$, the case $r=2$ gives convergence
to one in $L^2$. Bounds with larger fixed $r$ and interpolation then
give convergence in every fixed $L^r$.

Anchor the first coordinate at zero. Under $\mathbb Q_N$ let
$f_N(z)=e^{-iz/2}P_N(e^{iz/N})/P_N(1)$. The analytic circular Jacobi
limit used in \cref{thm:positive-merging} couples these functions so
that $f_N\to\xi^{\beta,Q_*/2}$ locally uniformly almost surely.
Choose $p>1$ such that $p(k-1)q_*<Q_*$. For all large $N$,
$p(k-1)q_N\le Q_*$, and \eqref{eq:cj-holder-bound} gives
\[
 \sup_{\mathbf y\in\mathbb R^k}
 \E_{\mathbb Q_N}
 \left[\left(\prod_{j=2}^k|f_N(y_j-y_1)|^{q_N}\right)^p\right]\le1.
\]
The products thus converge in $L^1$ along every convergent sequence
of real arguments. This includes limiting zeros: for large $N$, $q_N\ge q_*/2>0$,
so a factor whose base tends to zero also tends to zero after raising
it to the power $q_N$. H\"older's inequality and the $L^r$ convergence of
$L_N$ show that multiplying the products by $L_N$ does not change
the limit. The exact circular Jacobi change of measure now gives the
assertion, because
$M_N(kq_N)N^{-(kq_N)^2/(2\beta)}\to C_\beta(Q_*)$.
A compactness argument, as in \cref{thm:positive-merging}, upgrades
sequential convergence to the stated local uniform convergence.
\end{proof}

\begin{proof}[Proof of \cref{thm:critical-window}]
We use the same subtraction as in \cref{thm:critical-second-term},
now uniformly over the bounded window. The microscopic part converges
by \cref{lem:window-charge-variation}; the remaining radial integral
gives the transition function.
Write $s=s_N(x)$, $q=q_N(x)$ and $a=(k-1)(s-1)$. For this proof and
the proof of \cref{cor:general-drifting-critical-crossover}, put
\[
 C(s)=C_\beta(q(s))^k,\qquad
 T_{N,s}(y,\mathbf u)
 =N^{-ks}F_N\left(q(s)\mathbf1,\frac yN\mathbf t(\mathbf u)\right).
\]
The gap formula \eqref{eq:exact-gap-integration} yields
\begin{equation}\label{eq:window-exact-radial}
 \frac1N\E_\beta\mathcal Z_N(q;f)^k
 =e^{kx}c_k\int_{\Delta_{k-2}}\int_0^{Nb(\mathbf u)}
    y^{k-2}T_{N,s}(y,\mathbf u)w_f(y/N,\mathbf u)
                                      \,dy\,d\mathbf u.
\end{equation}

The exponents $q_N(x)$ lie in a fixed compact subset of $(0,\infty)$
for $|x|\le A$ and all large $N$. The uniform statement of
\cref{lem:moments-global-bound} therefore gives the clipped bound
with a constant independent of $x,N$. The separated Fisher--Hartwig
estimate is uniform on the same set by
\cref{thm:signed-characteristic-FH}.
Consequently, for $y\ge1$ in the integration domain,
\begin{align}\label{eq:window-tail-difference}
 &\left|y^{k-2}T_{N,s}(y,\mathbf u)
 -C(s)V_{k,s}(\mathbf u)y^{-1-a}D_k(y/N,\mathbf u)^s\right|
 \notag\\
 &\hspace{25mm}\le C V_{k,s}(\mathbf u)y^{-1-a}
                                      \min\{1,(yv)^{-\kappa}\},
\end{align}
with one $\kappa>0$ and one constant for all $|x|\le A$ and all
large $N$. On $yv\ge R$, this follows from
\cref{thm:signed-characteristic-FH}; on its complement, we use the
bound \eqref{eq:moments-global-bound} with distances replaced by their
maximum with $N^{-1}$. For
$0<y\le1$, the integrand without $w_f$ is bounded by $Cy^{k-2}$.

Choose $\epsilon,\gamma>0$ so small that
\[
 \int_{\Delta_{k-2}}V_{k,1+\epsilon}(\mathbf u)
                      v(\mathbf u)^{-\gamma}\,d\mathbf u<\infty,
 \qquad \gamma<\kappa.
\]
For $k=2$ the shape space is a point, so only $\gamma<\kappa$ is
needed. For all large $N$, $V_{k,s}\le V_{k,1+\epsilon}$.
Moreover $y^{-a}\le C_{A,k}$ when $1\le y\le Nb(\mathbf u)$.
After shape integration, \eqref{eq:window-tail-difference} is therefore
bounded by $Cy^{-1-\gamma}$. The multiplier $w_f$ is bounded.
\Cref{lem:window-charge-variation} and this integrable bound imply
that the integral in \eqref{eq:window-exact-radial}, after subtracting
\[
 C(s)V_{k,s}(\mathbf u)y^{-1-a}D_k(y/N,\mathbf u)^s
                         w_f(y/N,\mathbf u)\mathbf1_{\{y\ge1\}},
\]
converges to $w_0 B^{\rm mic}_{\beta,k}/c_k$. This convergence is
uniform for $|x|\le A$: every sequence witnessing failure has a
subsequence on which $x$ converges, and the same dominated convergence
argument applies on that subsequence.

It remains to integrate the subtracted term. Returning to $r=y/N$
shows that its contribution to the left side of
\eqref{eq:window-exact-radial} is
\begin{align*}
 e^x C(s)c_k\int V_{k,s}(\mathbf u)
 \bigg\{&w_0\log N\int_0^1e^{(k-1)xt}\,dt
       +w_0\frac{1-b(\mathbf u)^{-a}}a\\
       &+\int_{1/N}^{b(\mathbf u)}
        r^{-1-a}\bigl[D_k(r,\mathbf u)^sw_f(r,\mathbf u)-w_0\bigr]
             \,dr\bigg\}\,d\mathbf u.
\end{align*}
The quotient at $a=0$ means $\log b$. Uniformly near $s=1$,
\eqref{eq:critical-dini-weight-bound} bounds the square bracket by
$C_{k,f}\{r^2+\omega_k(f,r)\}$ for $0<r\le1$; it is bounded on
$1\le r\le b(\mathbf u)$. For $|x|\le A$ and $1/N\le r\le1$,
\[
 r^{-a}=\exp\left(\frac{(k-1)x}{\log N}\log\frac1r\right)
                       \le e^{(k-1)A}.
\]
Thus the Dini condition supplies an integrable radial majorant,
independent of $x,N$. The shape integrals are dominated by an
integrable $V_{k,1+\epsilon}$. Since $a\to0$, dominated convergence
applies to the last integral and to the shape integral. Uniformity
for $|x|\le A$ follows by the same sequential argument as for the
microscopic part. The last two terms therefore contribute
$e^xB^{\rm circ}_{\beta,k}[f]+o(1)$. The first contributes exactly
$w_0J_{\beta,k}(s)\Psi_k(x)\log N$.
Combining this calculation with the subtracted microscopic part proves
the equivalent form of \eqref{eq:critical-window-two-term}. Finally,
\[
 J_{\beta,k}(s_N(x))
 =J_{\beta,k}(1)+\frac{x}{\log N}J'_{\beta,k}(1)
                                  +O_A((\log N)^{-2}),
\]
uniformly for $|x|\le A$, proving the stated form.
\end{proof}

For an arbitrary drift through criticality, we retain a H\"older
assumption. Its positive power modulus absorbs the radial factor
$r^{-a}$ uniformly on a fixed neighborhood of exponents.

\begin{corollary}
\label{cor:general-drifting-critical-crossover}
Fix $\beta>0$, an integer $k\geq2$, and a nonnegative H\"older
function $f$. There is $\epsilon>0$ such that, uniformly for
$|s-1|\leq\epsilon$ and $N\geq2$,
\begin{equation}\label{eq:general-drifting-critical-crossover}
 \E_\beta\bigl[\mathcal Z_N(q(s);f)^k\bigr]
 =w_0J_{\beta,k}(s)N^s
       \frac{N^{(k-1)(s-1)}-1}{(k-1)(s-1)}
       +O_{\beta,k,f}\bigl(N^s+N^{ks-k+1}\bigr).
\end{equation}
At $s=1$ the quotient is interpreted as $\log N$.
If $f$ is not identically zero, then for every sequence $s_N\to1$,
\begin{equation}\label{eq:general-drifting-critical-relative}
 \E_\beta\bigl[\mathcal Z_N(q(s_N);f)^k\bigr]
 =w_0J_{\beta,k}(s_N)N^{s_N}
       \frac{N^{(k-1)(s_N-1)}-1}{(k-1)(s_N-1)}
       \left(1+O_{\beta,k,f}
                    \left(|s_N-1|+\frac1{\log N}\right)\right).
\end{equation}
\end{corollary}

\begin{proof}[Proof of \cref{cor:general-drifting-critical-crossover}]
We keep a fixed interval of exponents about the critical value and
bound both parts of the subtracted integral uniformly on that interval.
The H\"older modulus controls the additional radial power near zero.
Let $\alpha>0$ be a H\"older exponent of $f$, and put
$a=(k-1)(s-1)$. First fix a neighborhood of $s=1$ on which
$q(s)>0$. The uniform bound of \cref{lem:moments-global-bound} and
\cref{thm:signed-characteristic-FH} apply on a compact interval
inside that neighborhood, with a common exponent $\kappa>0$.
Choose $\gamma>0$ smaller than $\kappa$ and, when $k>2$, smaller
than the positive shape margin at $s=1$. Then decrease $\epsilon>0$
so that $[1-\epsilon,1+\epsilon]$ lies in that interval and
\begin{equation}\label{eq:drifting-shape-margin}
 \int_{\Delta_{k-2}}
      V_{k,1+\epsilon}(\mathbf u)v(\mathbf u)^{-\gamma}
                                      \,d\mathbf u<\infty,
 \qquad
 \gamma<\kappa,
 \qquad
 (k-1)\epsilon<\tfrac12\min\{\gamma,\alpha,2\}.
\end{equation}
When $k=2$, the shape integral equals one. For $k>2$, the general
criterion of \cref{lem:critical-shape-integrals} ensures that these
choices are possible.

Keep the notation $C(s)$ and $T_{N,s}$ introduced in the proof of
\cref{thm:critical-window}. Before imposing any relation between
$s$ and $N$, the exact radial formula is
\begin{equation}\label{eq:drifting-exact-radial}
 \E_\beta\bigl[\mathcal Z_N(q(s);f)^k\bigr]
 =N^{ks-k+1}c_k
   \int_{\Delta_{k-2}}\int_0^{Nb(\mathbf u)}
    y^{k-2}T_{N,s}(y,\mathbf u)w_f(y/N,\mathbf u)
                                        \,dy\,d\mathbf u.
\end{equation}
Apply the uniform statement of \cref{lem:moments-global-bound} with
$Q=k\max_{|s-1|\le\epsilon}q(s)$. Together with
\cref{thm:signed-characteristic-FH}, it proves
\eqref{eq:window-tail-difference} throughout $|s-1|\le\epsilon$,
with one constant and one $\kappa$, in the form
\[
 \left|y^{k-2}T_{N,s}(y,\mathbf u)
      -C(s)V_{k,s}(\mathbf u)y^{-1-a}
                    D_k(y/N,\mathbf u)^s\right|
 \leq C V_{k,s}(\mathbf u)y^{-1-a}
                          \min\{1,(yv)^{-\kappa}\},\qquad y\geq1.
\]
Since $\min\{1,z^{-\kappa}\}\leq z^{-\gamma}$ for every $z>0$,
shape integration and \eqref{eq:drifting-shape-margin} bound this
error by $Cy^{-1-a-\gamma}$. Its integral over $[1,\infty)$ is
uniformly bounded:
\[
 \int_1^\infty y^{-1-a-\gamma}\,dy
       =\frac1{a+\gamma}\leq\frac2\gamma.
\]
The omitted region $0<y<1$ also contributes a uniformly bounded
integral, because $T_{N,s}\leq C$ there and $k\geq2$.
The factor $w_f$ is bounded. It follows from
\eqref{eq:drifting-exact-radial} that
\begin{align}\label{eq:drifting-subtracted-radial}
 \E_\beta\bigl[\mathcal Z_N(q(s);f)^k\bigr]
 ={}&N^s C(s)c_k\int_{\Delta_{k-2}}V_{k,s}(\mathbf u)
       \int_{1/N}^{b(\mathbf u)}r^{-1-a}
                   D_k(r,\mathbf u)^sw_f(r,\mathbf u)
                                      \,dr\,d\mathbf u\notag\\
 &+O_{\beta,k,f}(N^{ks-k+1}).
\end{align}
In this change of variables, $y=Nr$ contributes $N^{-a}$ and
$ks-k+1-a=s$.

Uniformly on the integration domain and on the chosen interval of exponents,
\[
 D_k(r,\mathbf u)^sw_f(r,\mathbf u)-w_0
                         =O_f(r^\alpha+r^2).
\]
Thus replacing this product by $w_0$ in
\eqref{eq:drifting-subtracted-radial} has error at most $CN^s$.
Indeed $|a|<\min\{\alpha,2\}/2$, the shape integrals are uniformly
bounded, and the resulting radial powers are integrable at zero with
a fixed positive margin. Finally,
\[
 \int_{1/N}^{b(\mathbf u)}r^{-1-a}\,dr
 =\frac{N^a-1}{a}+\int_1^{b(\mathbf u)}r^{-1-a}\,dr.
\]
The last integral is bounded uniformly in $a$ and $\mathbf u$,
because $\pi\leq b(\mathbf u)<2\pi$. Multiplying the first term by
$C(s)c_k\int V_{k,s}=J_{\beta,k}(s)$ proves
\eqref{eq:general-drifting-critical-crossover}.

For the relative assertion, $w_0>0$ and $J_{\beta,k}(s)$ has a
positive lower bound on a sufficiently small fixed interval about one.
Put $L=\log N$. The quotient $(N^a-1)/a$ is positive. For $a\ne0$,
\[
 \frac{1+N^a}{(N^a-1)/a}
 =|a|\coth\left(\frac{|a|L}{2}\right)
 \leq |a|+\frac2L,
\]
where the inequality follows from $e^u-1\geq u$ for $u\geq0$.
At $a=0$ the ratio is $2/L$. Since $N^{ks-k+1}=N^sN^a$,
dividing the additive error in
\eqref{eq:general-drifting-critical-crossover} by its leading term
gives \eqref{eq:general-drifting-critical-relative}.
\end{proof}

\begin{remark}
\Cref{thm:critical-window} identifies a second term.
\Cref{cor:general-drifting-critical-crossover} instead gives a leading asymptotic for every drift
$s_N\to1$, including $(s_N-1)\log N\to+\infty$ or $-\infty$.
Its proof needs only a uniform bound on the subtracted microscopic
integral, so it imposes no additional rate of convergence on the
exponent. In both results the outer moment order $k$ is an integer;
these arguments give nothing for noninteger outer moments.
\end{remark}

\begin{proof}[Proof of \cref{thm:main-critical}]
For noninteger $p>1$, \cref{thm:fractional-critical} proves the
unweighted asymptotic. The integer case follows from
\cref{thm:critical-moments} with $f=1$.
\Cref{prop:critical-lp-weights} proves the leading asymptotic for
every real $p>1$ and every complex $L^p$ weight. The constant is obtained
from the normalized
critical coefficient $A_p$ using
$M_N(q)^p/N\to C_\beta(q)^p$ from \cref{prop:onepoint}.
For integer $p=k$, \cref{cor:weighted-critical-second-term} proves
the second assertion under exactly the Dini condition in the theorem,
and proves absolute convergence of \eqref{eq:weighted-critical-B}.
\end{proof}

\section{Mixed modulus and argument weights}
\label{sec:mixed-root-jump}

We extend the root formula to include a real exponential weight in the
argument at every prescribed point. This extension gives the fixed arc
counting formulas of \cref{cor:fixed-arc-counting,cor:fixed-arc-moderate-deviations}
and the single argument background formula of
\cref{cor:loc-root-jump-background}. The proofs of
\cref{thm:main-root,thm:main-critical} do not use argument weights.
The argument has three steps: we extend the Fourier correction to the
mixed recursion, identify its limit using the common real and imaginary
Gaussian field, and integrate the final coefficient using a bound on
the weighted phase density.

As in \Cref{sec:onepoint}, put $L_N(\vartheta)=\log P_N(e^{i\vartheta})$,
$X_N(\vartheta)=\Re L_N(\vartheta)=\log|P_N(e^{i\vartheta})|$ and
$J_N(\vartheta)=\Im L_N(\vartheta)$, using the radial boundary value of
the analytic logarithm normalized to vanish at the disk origin. For real
$\mathbf q,\mathbf s$, the mixed moment and its growth exponent are
\[
 \mathcal F_N^{\rm rj}(\mathbf q,\mathbf s;\boldsymbol\vartheta)
 =\E\exp\left\{\sum_jq_jX_N(\vartheta_j)+s_jJ_N(\vartheta_j)\right\},
 \qquad
 \mathfrak a(\mathbf q,\mathbf s)
 =\frac1{2\beta}\sum_j(q_j^2+s_j^2).
\]
The same expectation is holomorphic in complex $\mathbf q,\mathbf s$
when $\Re q_j>-1$. We retain the two parameters from
\cref{prop:onepoint}, so that $q_jX_N+s_jJ_N=a_jL_N+c_j\overline{L_N}$:
\[
 a_j=\frac{q_j-is_j}{2},\quad c_j=\frac{q_j+is_j}{2},\qquad
 C_\beta^{\rm rj}(q,s)=\mathcal C_\beta((q-is)/2,(q+is)/2).
\]
Write $M_N^{\rm rj}(q,s)=\mathcal F_N^{\rm rj}(q,s;0)$ for the
one point moment. The Morris identity of \cref{prop:onepoint} becomes
\begin{equation}\label{eq:mixed-onepoint-morris}
 M_N^{\rm rj}(q,s)
 =\prod_{k=0}^{N-1}
 \frac{\Gamma(1+bk)\Gamma(1+bk+q)}
 {\Gamma(1+bk+(q-is)/2)\Gamma(1+bk+(q+is)/2)}
 =C_\beta^{\rm rj}(q,s)N^{(q^2+s^2)/(2\beta)}(1+O(N^{-1})).
\end{equation}
For real $q>-1$ and $s\in\mathbb R$, the constants are positive.
On sufficiently small complex neighborhoods of compact real parameter
sets, the defining gamma series is analytic and nonzero, and the error
in \eqref{eq:mixed-onepoint-morris} is locally uniform.

\begin{theorem}
\label{thm:mixed-FH}
Let $K\subset(-1,\infty)^m\times\mathbb R^m$ be compact and let
$\delta>0$. Put
\begin{align*}
 \eta_K&=\min\left\{1,\inf_{(\mathbf q,\mathbf s)\in K}
                              \min_j(1+q_j)\right\},\\
 A_K&=\sup_{(\mathbf q,\mathbf s)\in K}\sum_j(-q_j)_+,
 \qquad M_K=2+(A_K-1)_+.
\end{align*}
For every
\[
 0<\kappa<\frac{\eta_K}{M_K+\eta_K}
                              \min\left\{\frac4\beta,1\right\},
\]
uniformly for $(\mathbf q,\mathbf s)\in K$ and
$\boldsymbol\vartheta\in\mathcal S_\delta$,
\begin{align}\label{eq:mixed-FH}
 \mathcal F_N^{\rm rj}(\mathbf q,\mathbf s;\boldsymbol\vartheta)
 ={}&N^{\mathfrak a(\mathbf q,\mathbf s)}
       \prod_j C_\beta^{\rm rj}(q_j,s_j)\notag\\
 &\times\prod_{j<k}
 (1-e^{i(\vartheta_j-\vartheta_k)})^{-2a_jc_k/\beta}
 (1-e^{i(\vartheta_k-\vartheta_j)})^{-2a_kc_j/\beta}
 \left(1+O(N^{-\kappa})\right).
\end{align}
All logarithms in the pair factors are radial limits from the disk.
The formula also holds locally uniformly on a sufficiently small complex
neighborhood of each such compact real parameter set, with a possibly
smaller positive power rate. For $m=1$,
\eqref{eq:mixed-onepoint-morris} gives the stronger error $O(N^{-1})$.
\end{theorem}

For real $q_j,s_j$ and
$u_{jk}=(\vartheta_j-\vartheta_k)\bmod2\pi\in(0,2\pi)$, the pair factor
in \eqref{eq:mixed-FH} is the positive quantity
\begin{equation}\label{eq:mixed-real-pair-proved}
 |e^{i\vartheta_j}-e^{i\vartheta_k}|^{-(q_jq_k+s_js_k)/\beta}
 \exp\left\{\frac{u_{jk}-\pi}{2\beta}(q_js_k-q_ks_j)\right\}.
\end{equation}
This follows from
$\log(1-e^{iu})=\log|1-e^{iu}|+i(u-\pi)/2$ and fixes the sign of the
argument interaction. The parameters $s_j$ in the real theorem are
Laplace parameters. Oscillatory jump weights correspond to imaginary $s_j$.
The mixed OPUC estimates in \cref{lem:mixed-corrector,thm:mixed-opuc}
identify the phase modes used in the correlation expansions of the
unpublished companion manuscript \cite{XIIIB}.

\Cref{thm:mixed-FH} contains \cref{thm:full-root-FH} as the case
$\mathbf s=0$: then $a_j=c_j=q_j/2$, the pair factors reduce to
$|e^{i\vartheta_j}-e^{i\vartheta_k}|^{-q_jq_k/\beta}$, the constants
$C_\beta^{\rm rj}(q_j,0)$ reduce to $C_\beta(q_j)$, and the admissible
rates are those of \eqref{eq:full-root-rate}. What is new is the
argument weight: arbitrary real Laplace parameters $s_j$ at every
point, on the same root domain $q_j>-1$ and with the same power rate,
together with uniformity on a complex neighborhood of the real
parameters. This local complex uniformity is what the counting
application of \Cref{sec:fixed-arc-counting} and the correlation
expansions in the same companion manuscript \cite{XIIIB} use. Approaching points are treated in
\Cref{sec:mixed-varying-separation}.

To prove \cref{thm:mixed-FH}, we separate the moment of the reversed
orthogonal polynomial from the contribution of its final coefficient.
For the rest of this section, we suppress the prescribed angles in the
OPUC weight, and suppress the fixed exponents when no confusion is
possible. The weight and its Fourier modes are
\[
 \mathcal B_n(\mathbf q,\mathbf s)
 =\exp\left\{\sum_jq_j\Re\log\Phi_n^*(e^{i\vartheta_j})
                   +s_j\Im\log\Phi_n^*(e^{i\vartheta_j})\right\},
 \qquad
 T_n=\mathcal B_n e^{i\boldsymbol\ell\cdot\boldsymbol\psi_n}.
\]
Here $\boldsymbol\ell\in\mathbb Z^m$, and $\boldsymbol\psi_n$
is the vector of phases from \eqref{eq:Tdef}. We use the analytic
logarithm, not its imaginary part reduced modulo $2\pi$:
$\Im\log\Phi_n^*(e^{i\vartheta})=((n+1)\vartheta-\psi_n(\vartheta))/2$.
In particular, the deterministic term $(n+1)\vartheta/2$ remains
part of the weight.

To integrate the last independent coefficient, define
\[
 h_{q,s}(t)=|1-e^{it}|^q
                    e^{s\operatorname{Arg}(1-e^{it})},\qquad
 W_{\mathbf q,\mathbf s}(\boldsymbol\varphi)
 =\int_\T\prod_jh_{q_j,s_j}(u+\varphi_j)\,\sigma(du).
\]
The branch agrees with the radial logarithm; values at the finitely
many singular points do not affect these integrals. The exact identity for the last coefficient is
\begin{equation}\label{eq:mixed-terminal-identity}
 \mathcal F_{n+1}^{\rm rj}
 =\E\bigl[\mathcal B_n
                W_{\mathbf q,\mathbf s}(\boldsymbol\psi_n)\bigr].
\end{equation}
To obtain \eqref{eq:mixed-terminal-identity}, take the analytic
logarithm in \eqref{eq:terminal-polynomial}, with normalization at zero,
and then pass to radial boundary values.
In particular, no choice of independent branches at the prescribed
points is made.

For real parameters, \cref{thm:mixed-opuc} shows that the normalized expectation
$n^{-\mathfrak a}\E\mathcal B_n$ converges to the explicit function
in \eqref{eq:mixed-opuc-limit}. It remains to show that
\[
 \frac{\E[\mathcal B_n W_{\mathbf q,\mathbf s}
                      (\boldsymbol\psi_n)]}{\E\mathcal B_n}
 \longrightarrow \prod_j c(q_j,s_j),
 \qquad c(q,s)=\int_\T h_{q,s}(t)\,\sigma(dt).
\]
The two constants combine through \eqref{eq:mixed-terminal-c}.
For the second limit, the phase estimate \eqref{eq:mixed-mode-bound}
controls each nonconstant Fourier mode of $W_{\mathbf q,\mathbf s}$, while
\cref{lem:mixed-terminal-decomposition} justifies approximation of the
possibly unbounded function $W_{\mathbf q,\mathbf s}$. The proof of
\cref{thm:mixed-FH} below combines these estimates and gives their rates.

\subsection{The mixed recursion and its holomorphic limit}

\begin{lemma}\label{lem:mixed-corrector}
The one step estimate \eqref{eq:onestep}, the corrected recursion
\eqref{eq:corrected-recursion}, and the full block estimate
\eqref{eq:unifiedblock} hold for these observables after the replacements
\begin{align}\label{eq:mixed-diagonal}
 D(\boldsymbol\varphi)
 &=\left(\sum_j(a_j-\ell_j)e^{i\varphi_j}\right)
   \left(\sum_j(c_j+\ell_j)e^{-i\varphi_j}\right),\notag\\
 D_0&=\sum_j\left\{\frac{q_j^2+s_j^2}{4}-\ell_j^2-is_j\ell_j\right\},
 \qquad
 \lambda=\mathfrak a-d(\boldsymbol\ell)
                      -\frac{2i}{\beta}\sum_js_j\ell_j.
\end{align}
For real parameters in compact sets and fixed separation, the normalized
zero sum nonzero modes satisfy, for every
$0<\gamma<\min\{4/\beta,1\}$,
\begin{equation}\label{eq:mixed-mode-bound}
 \left|\frac{\E T_n}{\E\mathcal B_n}\right|
 \le C(1+\|\boldsymbol\ell\|_1)^2 n^{-\gamma}.
\end{equation}
Moreover $n^{-\mathfrak a}\E\mathcal B_n$ has a locally uniform
holomorphic limit on
\[
 \mathcal U_{\rm rj}=
 \left\{(\mathbf q,\mathbf s)\in\mathbb C^{2m}:
 \Re q_j>-1-b,\quad
 h:=\frac{\|\Im\mathbf q\|^2+\|\Im\mathbf s\|^2}{2\beta}<1\right\}.
\]
On compact subsets the zero mode error is $O(n^{-1+h_0})$, where
$h_0$ is the supremum of $h$ on that compact set. The limit is positive
for real parameters.
\end{lemma}

\begin{proof}[Proof of \cref{lem:mixed-corrector}]
We first check the one step expansion and the positive envelope. We then
control the nonzero modes with the scalar recursion, and use the same
envelope to obtain holomorphic convergence of the zero mode.
The exact multiplier is
\[
 M_k(\boldsymbol\varphi)
 =\prod_j(1-\alpha_ke^{i\varphi_j})^{a_j-\ell_j}
          (1-\bar\alpha_ke^{-i\varphi_j})^{c_j+\ell_j}.
\]
Writing $L_j=\log|1-\alpha_ke^{i\varphi_j}|$ and
$A_j=\operatorname{Arg}(1-\alpha_ke^{i\varphi_j})$, the two analytic logarithms are
$L_j+iA_j$ and $L_j-iA_j$. Hence, also for complex parameters,
\[
 |M_k(\boldsymbol\varphi)|
 =\exp\left\{\sum_j
       \bigl((\Re q_j)L_j+(\Re s_j)A_j\bigr)\right\}.
\]
This is the positive mixed one step weight at
$(\Re\mathbf q,\Re\mathbf s)$, independently of the integer mode. The
argument of each factor lies in $[-\pi/2,\pi/2]$, so the annulus
estimate \eqref{eq:signed-annulus} changes only by a constant depending
on the parameter compact set. On the inner
disk, Taylor expansion and rotational averaging give
$\E M_k=1+v_kD+O(L^4(k+1)^{-2})$, with $L=1+\|\boldsymbol\ell\|_1$.
The signed first phase moment and the weighted second phase moment
satisfy the bounds in \eqref{eq:mixed-one-step} and \eqref{eq:weighted-phase-second}.
Consequently \eqref{eq:corrector-equation} has the same Fourier
solution, with coefficients $(a_j-\ell_j)(c_k+\ell_k)$.
Only the denominators $e^{i(\vartheta_j-\vartheta_k)}-1$ occur,
and the corrected remainder is $O(L^4/(k^2\delta))$ times the
positive envelope. Expansion of $(a_j-\ell_j)(c_j+\ell_j)$ gives
\eqref{eq:mixed-diagonal}, including its negative imaginary sign.

For zero mode and real parameters, positivity gives the same real
envelope and relative moment comparison as in
\cref{lem:negative-envelope,lem:relative-moment-comparison}, with
$a$ replaced by $\mathfrak a$. These proofs use the zero mode block
estimate and positivity, so no identification of the limiting constant
is needed. At a fixed initial index the absolute value of the argument of
$\Phi_n^*$ is bounded by $n\pi/2$: each of its zero factors has positive
real part in the disk. Thus the mixed weight is bounded by $e^{C_Kn}$
times the root weight. The sharp domain and finite index uniform
integrability follow from \cref{lem:negative-initial}; a box of nearby
real parameters also supplies dominated holomorphy.

For a nonzero real mixed mode, recall from the proof of
\cref{prop:quantitative-phase-mixing} that
$p_k(a)=\Gamma(k+\tau_\beta+a)/\Gamma(k+\tau_\beta)$ and
$\tau_\beta=1+1/b$. We use this normalization at $a=\mathfrak a$.
Choose $n_*=\lceil CL^2/\delta\rceil$ and
$h=\E\mathcal B_{n_*}/p_{n_*}(\mathfrak a)$ to normalize the positive
envelope at the starting index. Here $h$ is a positive
normalizing constant; it is distinct from the imaginary parameter
quantity in the definition of $\mathcal U_{\rm rj}$.
The scalar coefficient is
\[
 c_k=1-\frac{d+i\eta}{u_k},\qquad
 u_k=k+\tau_\beta+\mathfrak a,\qquad
 \eta=\frac2\beta\mathbf s\cdot\boldsymbol\ell.
\]
On a compact real parameter set, $\eta^2\le C_Kd$.
Choose the constant in $n_*$ so that $d/u_k\le1/2$ for all
$k\ge n_*$. Then
\[
 |c_k|=\left(1-\frac d{u_k}\right)
       \left(1+\frac{\eta^2}{(u_k-d)^2}\right)^{1/2}
 \le\left(1-\frac d{u_k}\right)
             \exp(2\eta^2/u_k^2).
\]
For $n_*\le r\le n$, \eqref{eq:phase-actual-propagator} gives
\begin{equation}\label{eq:mixed-actual-propagator}
 \prod_{k=r}^{n-1}|c_k|
 \le C(r/n)^d\exp\left(2\eta^2
                         \sum_{k=r}^{n-1}u_k^{-2}\right)
 \le C(r/n)^d\exp(C_Kd/r)\le C'(r/n)^d.
\end{equation}
The constants are independent of the mode: both $d/r$ and the
index shift cost are bounded at $r\ge n_*$. The normalized forcing
is still $O(L^4/(k^2\delta))$, while removal of the correction costs
$O(L^2/(n\delta))$. The split into $d>2$ and the finite set $d\le2$
in the proof of \cref{prop:quantitative-phase-mixing} therefore gives
$CL^2\mathfrak r_{4/\beta}(n\delta)$. In particular it proves
\eqref{eq:mixed-mode-bound} at fixed separation.

For the later continuation to complex parameters we also retain the
elementary estimate
\[
 \log\left|1+\frac{\lambda}{k+\tau_\beta}\right|
 \le\frac{\Re\lambda}{k+\tau_\beta}
                    +\frac{C_KL^4}{(k+1)^2},\qquad k\ge C_KL^2.
\]
It follows from the Taylor expansion of $\log(1+z)$ for $|z|\le1/2$.
Starting instead at $k\ge C_KL^4$ makes the sum of its quadratic
errors bounded uniformly in the mode. The complex parameter estimates
below continue to use this bound.

For complex parameters,
$|\mathcal B_n(\mathbf q,\mathbf s)|
=\mathcal B_n(\Re\mathbf q,\Re\mathbf s)$ and
\[
 \mathfrak a(\Re\mathbf q,\Re\mathbf s)-
            \Re\mathfrak a(\mathbf q,\mathbf s)=h.
\]
The mixed version of \eqref{eq:unifiedblock} therefore bounds the
increment of the normalized zero mode across a doubling block by
$O(n^{-1+h_0})$. The sum over $n,2n,4n,\ldots$ converges because
$h_0<1$. Local uniform convergence of the holomorphic finite moments
then proves the asserted holomorphy.
\end{proof}

\subsection{The conditional identity for mixed root and argument weights}
\label{sec:mixed-chaos-identity}

The real and imaginary parts of the common Gaussian analytic function
are correlated as fields. To identify mixed weights, we extend the
conditional identity of Lambert and Najnudel directly, using their
tail decomposition. A positive coefficient expansion replaces both
scalar comparison functions used in their proof.

For $q,s\in\mathbb R$, we distinguish the exponential weight of a
complex number, its value on the reversed polynomial, and its one point
mean:
\[
 W_{q,s}(h)=\exp\{q\Re h+s\Im h\},\qquad
 w_n^{q,s}(z)=W_{q,s}(\log\Phi_n^*(z)),\qquad
 \nu_n(q,s)=\E w_n^{q,s}(1).
\]
Every analytic logarithm here vanishes at the origin. For $q>-1$,
the one point OPUC product \eqref{eq:mixed-opuc-onepoint} gives
$0<\nu_n(q,s)<\infty$.
For each fixed boundary point the normalized weights form a martingale:
rotation invariance of $\alpha_n$ makes the conditional mean of
the next multiplier independent of the phase, and that mean is
$\nu_{n+1}(q,s)/\nu_n(q,s)$. We identify the limit of the resulting random measures with the chaos
of the common field in \cref{prop:mixed-chaos-identity}.
The covariance of the common field $\varphi$ in \eqref{eq:common-gaf}
shows that $q\Re\varphi+s\Im\varphi$ has the same law as
$\sqrt{q^2+s^2}\Re\varphi$. When $q^2+s^2<2\beta$, the radial
convergence \eqref{eq:gmc-radial} therefore defines a mean one chaos,
denoted by $U^{q,s}$. Its marginal law is that of the root chaos with
parameter $\sqrt{q^2+s^2}$. The different mixed chaoses below are all
constructed from the same $\varphi$.

\begin{proposition}
\label{prop:mixed-chaos-identity}
Suppose $q>-1$ and $q^2+s^2<2\beta$. Then, for every bounded continuous
$f$ and every $n\geq1$,
\begin{equation}\label{eq:mixed-chaos-conditional}
 \E[U^{q,s}(f)\mid\mathcal F_n]
 =\int_\T f(\theta)
       \frac{w_n^{q,s}(e^{i\theta})}{\nu_n(q,s)}\,\sigma(d\theta).
\end{equation}
Consequently these random measures converge almost surely weakly
to $U^{q,s}$, and their integrals against $f$ converge in $L^p$ whenever
$p\geq1$ and $p(q^2+s^2)<2\beta$.
\end{proposition}

The tail decomposition in \cref{lem:mixed-analytic-tail} separates the
first $n$ coefficients from the remaining independent coefficients.
Its rotation property will allow us to average the mixed weight over a
single phase in the proof of \cref{prop:mixed-chaos-identity}.

\begin{lemma}
\label{lem:mixed-analytic-tail}
Fix $n\ge0$. With the matrix $A(\alpha,z)$ in
\eqref{eq:szego-recursion}, put
\[
 T_{n,N}(z)=A(\alpha_{N-1},z)\cdots A(\alpha_n,z),\qquad
 F_{n,N}^{\pm}(z)=[T_{n,N}(z)]_{22}\pm[T_{n,N}(z)]_{21}.
\]
Almost surely, $F_{n,N}^{\pm}$ converge locally uniformly in $\D$
to nonvanishing analytic functions $F_n^\pm$, equal to one at zero.
For every fixed $z\in\D$ and $c>0$,
\[
 \E\exp\{c|\log F_n^\pm(z)|\}<\infty.
\]
Set $A_n=(F_n^++F_n^-)/2$ and $B_n=(F_n^+-F_n^-)/2$;
this notation is local to the tail decomposition.
Then $A_n(0)=1$, $B_n(0)=0$, $A_n$ has no zero, and
$|B_n(z)/A_n(z)|\le|z|<1$. For every fixed $z\in\D$ and $c>0$,
\[
 \E e^{c|\log A_n(z)|}<\infty,\qquad
 \E\exp\left\{c\sup_{|\lambda|\le1}
       |\log(A_n(z)+\lambda B_n(z))|\right\}<\infty.
\]
These analytic logarithms vanish at zero. Moreover,
\begin{equation}\label{eq:mixed-exact-tail}
 e^{\varphi(z)}=\Phi_n^*(z)
 \left(A_n(z)+\frac{\Phi_n(z)}{\Phi_n^*(z)}B_n(z)\right).
\end{equation}
A common rotation $\alpha_j\mapsto\lambda\alpha_j$ for $j\ge n$,
where $|\lambda|=1$, fixes $A_n$ and sends $B_n$ to $\lambda B_n$.
\end{lemma}

\begin{proof}[Proof of \cref{lem:mixed-analytic-tail}]
The analytic convergence and exponential moments are precisely
\cite[Lemma~5]{LN2026}, applied to the shifted coefficient sequence
and to its negative. The function $F_{n,N}^+$ is the reversed
polynomial started from $(1,1)$ at index $n$; replacing every tail
coefficient by its negative gives $F_{n,N}^-$.
The positive real part identity for their quotient gives
$|(F_n^+-F_n^-)/(F_n^++F_n^-)|<1$; see
\cite[Fact~2]{LN2026}. It also makes $A_n$ nonvanishing.
The quotient $B_n/A_n$ is analytic, maps $\D$ into itself and
vanishes at zero. Schwarz's lemma therefore gives
$|B_n(z)/A_n(z)|\le|z|$.

The remaining identities are finite matrix algebra. The recursion gives
\[
 \Phi_N^*=\Phi_n^*
 \left([T_{n,N}]_{22}
       +\frac{\Phi_n}{\Phi_n^*}[T_{n,N}]_{21}\right).
\]
Passing to the analytic limit proves \eqref{eq:mixed-exact-tail}.
For $D_\lambda=\operatorname{diag}(1,\lambda)$, direct multiplication gives
\[
 A(\lambda\alpha,z)=D_\lambda A(\alpha,z)D_\lambda^{-1},\qquad
 T_{n,N}\mapsto D_\lambda T_{n,N}D_\lambda^{-1}.
\]
Thus entry $(2,2)$ stays fixed and entry $(2,1)$ is multiplied by
$\lambda$. The same holds in the limit. All logarithms are analytic
and zero at the origin. The bound $|B_n(z)/A_n(z)|\le|z|$ gives
\[
 \sup_{|\lambda|\le1}|\log(1+\lambda B_n(z)/A_n(z))|
 \le-\log(1-|z|).
\]
Since $F_n^+=A_n+B_n$, it follows that
\[
 \sup_{|\lambda|\le1}|\log(A_n(z)+\lambda B_n(z))|
 \le |\log F_n^+(z)|-2\log(1-|z|).
\]
The analogous bound for $|\log A_n(z)|$ has only one copy of
$-\log(1-|z|)$. The asserted exponential moments therefore follow
from those of $\log F_n^+(z)$, for each fixed $n$ and $z\in\D$.
\end{proof}

\begin{proof}[Proof of \cref{prop:mixed-chaos-identity}]
We first compare the conditional radial weight with its boundary value
using the rotation of the independent tail. We then identify the
normalizer and pass to the radial chaos limit in $L^1$.
Fix $n,q,s$ for this proof. To lighten notation, write
$W=W_{q,s}$, $w_n=w_n^{q,s}$ and $\nu_n=\nu_n(q,s)$; when the
parameters are multiplied by a moment exponent, we keep them visible.
Write $a=(q-is)/2$. The phase average needed for the comparison is,
for $0\leq r<1$,
\[
 G(r)=\int_{\mathbb T}
       |(1+re^{it})^a|^2\,\sigma(dt)
       =\sum_{k=0}^\infty\left|\binom ak\right|^2r^{2k}.
\]
The power uses the analytic logarithm on the right half-plane.
Parseval gives the last equality. Its coefficients are nonnegative,
and, unless the series terminates,
$|\binom ak|^2=O_{q,s}(k^{-q-2})$. Thus $G$ is increasing,
$G(r)\geq1$, and, for every
$0<\eta<\min\{1,q+1\}$,
\begin{equation}\label{eq:mixed-radial-comparison}
 0\leq G(r)-G(\rho r)
 \leq C_{q,s,\eta}(1-\rho)^\eta G(r),
 \qquad 0\leq r,\rho\leq1.
\end{equation}
Indeed $1-\rho^{2k}\leq C_\eta k^\eta(1-\rho)^\eta$, and the
resulting series is summable because $\eta<q+1$.

Use \cref{lem:mixed-analytic-tail} and abbreviate $A=A_n$, $B=B_n$.
Write $\Phi_n/\Phi_n^*=\rho_n\lambda_n$, where $\rho_n$ is the
modulus and $|\lambda_n|=1$ away from zeros; at a zero the phase may
be chosen arbitrarily. The decomposition \eqref{eq:mixed-exact-tail} gives
\[
 e^{\varphi(z)}=\Phi_n^*(z)
       \bigl(A(z)+\rho_n(z)\lambda_n(z)B(z)\bigr).
\]
The logarithm of the second factor is
$\log A+\log(1+\rho_n\lambda_n B/A)$, by normalization at zero.

For fixed $z$, interpolate between the two tail factors by the
conditional mean below, where $\rho\in[0,1]$ and
$\E_n=\E[\,\cdot\mid\mathcal F_n]$:
\[
 \Theta(\rho;z)=
 \E_n W\bigl(\log A(z)+
               \log(1+\rho\lambda_n(z)B(z)/A(z))\bigr).
\]
By \cref{lem:mixed-analytic-tail}, a common rotation of the tail
fixes $A$ and rotates $B$. The joint coefficient law is invariant
under this rotation.
Averaging over $\lambda$ therefore gives
\[
 \Theta(\rho;z)=
  \E\left[W(\log A(z))
           G(\rho|B(z)/A(z)|)\right].
\]
The exponential moment assertion of \cref{lem:mixed-analytic-tail}
makes these expectations finite. At $\rho=1$, rotation invariance and
independence of the tail give a deterministic mean depending only on
$|z|$. Denote this mean by $T(|z|)$; the fixed index $n$ is suppressed:
\[
 \Theta(1;z)=T(|z|),\qquad
 T(r)=\E W(\log F_n^+(r)).
\]
Hence \eqref{eq:mixed-radial-comparison} gives the two-sided bound
\begin{equation}\label{eq:mixed-tail-sandwich}
 (1-C(1-\rho_n(z))^\eta)w_n(z)T(|z|)
 \leq\E_n W(\varphi(z))
 \leq w_n(z)T(|z|).
\end{equation}
The estimate extends over zeros of $\Phi_n$ by continuity.

For fixed $n$, the polynomial $\Phi_n^*$ has no zero on the closed
disk, and $|\Phi_n/\Phi_n^*|=1$ on its boundary. Continuity therefore
gives $\rho_n(ru)\to1$ uniformly in $u\in\mathbb T$ almost surely.
To pass to expectations, choose $t>1$ with $tq>-1$, which is possible
because $q>-1$. For each realization,
$w_n^{tq,ts}(z)=|\exp((tq-its)\log\Phi_n^*(z)/2)|^2$ is
subharmonic. Its circular means increase with radius. Rotational
invariance therefore gives
\[
 \sup_{0\leq r\leq1}\E[w_n(r)^t]
 \leq\nu_n(tq,ts)<\infty.
\]
Continuity and uniform integrability yield
$\E w_n(r)\to\nu_n$ and
$\E[w_n(r)(1-\rho_n(r))^\eta]\to0$.
Taking expectation in \eqref{eq:mixed-tail-sandwich} consequently proves
\begin{equation}\label{eq:mixed-tail-normalizer}
 \frac{\E W(\varphi(r))}{T(r)}
 \longrightarrow\nu_n.
\end{equation}

The radial chaos measure is
\[
 U_r^{q,s}(d\theta)=
 \frac{W(\varphi(re^{i\theta}))}
      {\E W(\varphi(r))}\,\sigma(d\theta).
\]
Its normalizer is $(1-r^2)^{-(q^2+s^2)/(2\beta)}$.
The Gaussian covariance in \eqref{eq:common-gaf} shows that
$q\Re\varphi+s\Im\varphi$ has the same marginal law as
$\sqrt{q^2+s^2}\Re\varphi$. Since $q^2+s^2<2\beta$, the radial
convergence \eqref{eq:gmc-radial} for that root parameter gives
$U_r^{q,s}(f)\to U^{q,s}(f)$ in $L^1$.
For $f\geq0$, integrate \eqref{eq:mixed-tail-sandwich}, divide by the
radial normalizer, and use \eqref{eq:mixed-tail-normalizer}.
The uniform boundary continuity of the fixed polynomial and of
$\rho_n$ shows that
\[
 \E_n U_r^{q,s}(f)\longrightarrow
   \frac1{\nu_n}\int f(\theta)w_n(e^{i\theta})
                                      \,\sigma(d\theta)
 \quad\hbox{almost surely}.
\]
On the other hand, conditional expectation is an $L^1$ contraction,
so the same left side converges in $L^1$ to $\E_nU^{q,s}(f)$.
Uniqueness of the limit in probability proves
\eqref{eq:mixed-chaos-conditional}. Taking differences treats general
bounded continuous $f$.

The field and its chaos are measurable with respect to
$\mathcal F_\infty$. Applying martingale convergence in
\eqref{eq:mixed-chaos-conditional} to a countable dense family in
$C(\T)$, including the constant function, proves almost sure weak
convergence of the measures. Finally, \eqref{eq:gmc-moment-input} and
the marginal law just identified give a finite $t$th moment of the
chaos mass whenever $t(q^2+s^2)<2\beta$. For the stated $p$, choose
$t>p$ satisfying this inequality. Conditional Jensen gives a uniform
$t$th moment bound for the conditional expectations, hence uniform
integrability of their $p$th powers and convergence in $L^p$.
\end{proof}

\begin{remark}
For finitely many $(q_j,s_j)$ in a sufficiently small real open
neighborhood of zero, the moment argument of
\cref{sec:gmc-identification} applies to these measures. All are built
from the same analytic field. Gaussian integration then gives the pair
factor
\[
 \prod_{j<k}(1-e^{i(\vartheta_j-\vartheta_k)})^{-2a_jc_k/\beta}
                 (1-e^{i(\vartheta_k-\vartheta_j)})^{-2a_kc_j/\beta},
 \qquad a_j=(q_j-is_j)/2,\quad c_j=(q_j+is_j)/2.
\]
The $2m$ real parameters $(\mathbf q,\mathbf s)$ vary independently
on this open set. Applying the one variable identity theorem
successively in each coordinate shows that these values determine
a holomorphic function of all $2m$ complex parameters. The root case
alone corresponds to $\mathbf s=0$ and does not provide such a
uniqueness set. We shall use the mixed identity to determine the
holomorphic OPUC limit, then integrate the last coefficient to obtain
the characteristic polynomial asymptotic.
\end{remark}

\subsection{Identification and transfer at all integrable root exponents}

To fix the normalization without introducing a singularity from the last coefficient,
put $x_k=1+bk$ and define
\begin{equation}\label{eq:mixed-opuc-onepoint}
 \nu_n(q,s)=\prod_{k=1}^n
 \frac{\Gamma(x_k)\Gamma(x_k+q)}
 {\Gamma(x_k+(q-is)/2)\Gamma(x_k+(q+is)/2)},\qquad
 D_\beta^{\rm rj}(q,s)=
 \lim_{n\to\infty}n^{-(q^2+s^2)/(2\beta)}\nu_n(q,s).
\end{equation}
The product is the one point OPUC moment. Its normalized tail converges
locally uniformly for $\Re q>-1-b$ and bounded complex $s$, by the
$O(k^{-2})$ remainder in the gamma expansion. Finite denominator
reciprocals are entire, so possible zeros cause no loss of holomorphy.
Where $\Re q>-1$, conditioning on the last coefficient gives
\begin{equation}\label{eq:mixed-terminal-c}
 c(q,s)=\frac{\Gamma(1+q)}
 {\Gamma(1+(q-is)/2)\Gamma(1+(q+is)/2)},\qquad
 C_\beta^{\rm rj}(q,s)=c(q,s)D_\beta^{\rm rj}(q,s).
\end{equation}
We define $D_\beta^{\rm rj}$ by the product limit in
\eqref{eq:mixed-opuc-onepoint}. This definition remains holomorphic
when $c(q,s)=0$, whereas the quotient $C_\beta^{\rm rj}(q,s)/c(q,s)$
would require a separate argument.

\begin{theorem}\label{thm:mixed-opuc}
On the tube $\mathcal U_{\rm rj}$ in \cref{lem:mixed-corrector},
\begin{equation}\label{eq:mixed-opuc-limit}
 n^{-\mathfrak a}\E\mathcal B_n
 =H(\mathbf q,\mathbf s;\boldsymbol\vartheta)+O(n^{-1+h_0}),
 \qquad
 H=\prod_jD_\beta^{\rm rj}(q_j,s_j)\,
       \mathcal K(\mathbf q,\mathbf s;\boldsymbol\vartheta),
\end{equation}
locally uniformly in the parameters and in separated configurations.
Here $\mathcal K$ is the product of the two pair factors in
\eqref{eq:mixed-FH}, and $h_0<1$ is the supremum of $h$ on the
parameter compact set. For every real compact set
$K\subset(-1-b,\infty)^m\times\R^m$, the stronger relative estimate
\begin{equation}\label{eq:mixed-opuc-relative}
 \E\mathcal B_n=n^{\mathfrak a}H
       \left(1+O_K((n\delta)^{-1})\right)
\end{equation}
holds when $n\delta\ge R_K$, with constants independent of the
configuration and separation. The function $H$ is nonzero on the tube
and positive for real parameters.
\end{theorem}

\begin{proof}[Proof of \cref{thm:mixed-opuc}]
We identify the limit first for small real parameters, extend that
identity by holomorphy, and then recover a relative estimate uniform
in the separation from the block comparison.
For small real $(q_j,s_j)$, choose a neighborhood on which the chaos
masses have moments of order strictly larger than $m$. The mixed
conditional identity \eqref{eq:mixed-chaos-conditional} and
H\"older's inequality then give convergence
of their products in $L^1$. Integrating against products of continuous
tests supported in disjoint arcs, Gaussian integration gives the pair
kernel $\mathcal K$. The covariance kernels are bounded on these
supports, so dominated convergence applies. Local uniform convergence
from \cref{lem:mixed-corrector} identifies the continuous limiting
OPUC density at each distinct configuration, exactly as in the proof
of \cref{thm:opuc-full}.

The proposed limit $H$ in \eqref{eq:mixed-opuc-limit} is holomorphic throughout $\mathcal U_{\rm rj}$.
In fact its gamma factors have no zeros there: for either sign,
\[
 \Re\left(1+b+\frac{q\pm is}{2}\right)
 >\frac{1+b}{2}-\sqrt b
 =\frac{(\sqrt b-1)^2}{2}\ge0,
\]
since $\Re q_j>-1-b$ and $|\Im s_j|^2\le2\beta h<2\beta$, that is,
$|\Im s_j|<2\sqrt b$.
Numerator arguments also have positive real parts. The small real
open box is a uniqueness set, by applying the one variable identity
theorem successively in all $2m$ coordinates. Connectedness of the
tube proves \eqref{eq:mixed-opuc-limit} everywhere on it.

For real parameters, apply the mixed full block estimate of
\cref{lem:mixed-corrector} at zero Fourier mode. Positivity gives,
for $n\le M\le2n$,
\[
 \left|\left(\frac Mn\right)^{-\mathfrak a}
       \frac{\E\mathcal B_M}{\E\mathcal B_n}-1\right|
 \le\frac{C_K}{n\delta}.
\]
For $n\delta$ sufficiently large, each ratio is bounded away from
zero, and its logarithm has error $O((n\delta)^{-1})$. Over successive
doubling blocks these errors sum to at most
$C\sum_{j\ge0}(2^jn\delta)^{-1}$. A final partial block treats any
$M\ge n$, so multiplication proves the same relative estimate. Let $M$ tend to infinity in the already identified positive
limit, then invert the ratio after increasing $R_K$. This gives
\eqref{eq:mixed-opuc-relative} without a separation dependent constant.
\end{proof}

\begin{lemma}
\label{lem:mixed-terminal-decomposition}
Let $K\subset(-1,\infty)^m\times\mathbb R^m$ be compact, with
$m\ge2$. For real parameters in $K$, the weighted phase probability
measure
\[
 \mu_n^{\mathbf q,\mathbf s}(B)
 =\frac{\E[\mathcal B_n\mathbf1_{\{\boldsymbol\psi_n\in B\}}]}
             {\E\mathcal B_n}
\]
is entirely absolutely continuous when $n\delta\ge R_K$, and
\[
 \mu_n^{\mathbf q,\mathbf s}=f_n\sigma^{\otimes m},\qquad
 \|f_n\|_\infty
 \le\frac{C_K\delta^{-m/2}}{\mathcal V(\boldsymbol\vartheta)}
 \le C_K\delta^{-m^2/2}.
\]
The constants are independent of the distinct prescribed angles.
In particular the density is uniformly bounded at fixed separation.
\end{lemma}

\begin{proof}[Proof of \cref{lem:mixed-terminal-decomposition}]
We adapt the adjoint propagation and the smoothing estimate from the
root case, then normalize the resulting density by
\eqref{eq:mixed-opuc-relative}. For the circle diffeomorphism
$g_\alpha$ used in
\cref{lem:adjoint-density-propagation}, the inverse identity holds as
an identity of logarithms:
\[
 \log(1-\alpha e^{ig_{-\alpha}(y)})
 =\log(1-|\alpha|^2)-\log(1+\alpha e^{iy}).
\]
Both nonreal logarithms have arguments in $(-\pi/2,\pi/2)$, which
verifies their branches. Including the inverse Jacobian therefore
changes the mixed parameters to
\[
 q'_j=-q_j-2,\qquad s'_j=-s_j,
\]
and shifts the radial parameter by $S=\sum_jq_j+m$, with prefactor
$t/(t+S)$ as in \eqref{eq:adjoint-dual-identity}.
The diagonal cancellation becomes
\[
 \sum_j\frac{(q_j+2)^2+s_j^2}{4}-S
       =\sum_j\frac{q_j^2+s_j^2}{4}.
\]
The two factors in the dual zero mode polynomial are complex
conjugates, so its off-diagonal part is real. The adjoint corrector
\eqref{eq:adjoint-corrector} and the positive supersolution from
\cref{lem:adjoint-density-propagation} consequently propagate bounded densities
by $C(n/k)^{\mathfrak a}$, uniformly when $k\delta\ge R_K$.
Indeed their $C^2$ norms are $O(\delta^{-1})$, the one step remainder
is $O(k^{-2}\delta^{-1})$, and its sum from $k$ onwards is
$O((k\delta)^{-1})$, as in the proof of
\cref{lem:adjoint-density-propagation}.

For the whole block kernel, the phase map and its polynomial
Jacobian are unchanged. On a compact set of real $s_j$, every
incremental argument factor and its reciprocal are bounded by a
constant, since $|\operatorname{Arg}(1-\alpha e^{ix})|<\pi/2$. Thus the Gaussian
envelope \eqref{eq:whole-weighted-Gaussian} holds for the mixed
weight on the entire coefficient polydisk. The same observation
applies to the dual parameters. The row and column estimates, and their
interpolation in \cref{lem:whole-block-Lp,thm:unconditional-phase-kernel},
give the full mixed kernel bound $Ck^{m/2}/\mathcal V$ after $r_m$
steps. This argument uses arbitrary real dual root exponents only in
conditional integrals at sufficiently large indices.

Choose $k=\lceil C_K/\delta\rceil$ and increase $R_K$ so that
$n\ge k+r_m$. Smooth the incoming measure of mass
$A_k=\E\mathcal B_k$ and propagate its density by the preceding
adjoint estimate. The resulting unnormalized density is bounded by
\[
 \frac{Ck^{m/2}}{\mathcal V}A_k(n/k)^{\mathfrak a}.
\]
The positive real mixed OPUC comparison
\eqref{eq:mixed-opuc-relative} gives
$A_k/A_n\le C(k/n)^{\mathfrak a}$, with a constant independent of
the separation. Division by $A_n$ cancels these powers.
Finally $k\le C_K/\delta$ and
$\mathcal V\ge\delta^{m(m-1)/2}$ prove the asserted bounds.
\end{proof}

\begin{proof}[Proof of \cref{thm:mixed-FH}]
We first transfer the real OPUC asymptotic through the last coefficient
by Fourier approximation. We then repeat that approximation in a small
complex neighborhood, where the additional factor $n^h$ can still be
absorbed in a decreasing error.
Suppose first that all parameters are real and $m\ge2$. The proof of
\cref{lem:terminal-translation} applies to the mixed factor as follows.
For every $0<\alpha<\eta_K$, split a translation difference into a
neighborhood of its singularity and the complementary arcs to obtain
\[
 \|h_{q,s}(\cdot+t)-h_{q,s}\|_1\le C\|t\|_{2\pi}^{\alpha}.
\]
The factors involving the argument are bounded on compact $s$ sets;
away from zero their derivatives are bounded by a constant times
$|t|^q$, whereas differentiating the root gives $O(|t|^{q-1})$.
The interval crossing the jump contributes $O(|t|^{q+1})$.
These observations also cover $q=0$, where the bounded jump costs
$O(|t|)$. Telescoping the product and using Fubini's theorem gives the same
$L^1$ translation estimate for $W_{\mathbf q,\mathbf s}$.

In this proof write $W=W_{\mathbf q,\mathbf s}$, with the parameters
fixed, and let $W_L$ be its product Fej\'er approximation of degree $L$
in each coordinate. It satisfies
\[
 \|W_L-W\|_1\le CL^{-\alpha},\qquad
 \widehat W_L(0)=\prod_jc(q_j,s_j).
\]
Simultaneous translation invariance leaves only zero sum modes.
The mixed form of \cref{lem:terminal-truncated-mass} bounds their
absolute Fourier mass by $CL^{(A_K-1)_++\zeta}$ for every $\zeta>0$.
Combining \cref{lem:mixed-terminal-decomposition} and
\eqref{eq:mixed-mode-bound}, for $\gamma<\min\{4/\beta,1\}$, gives
\[
 \left|\int W\,d\mu_n^{\mathbf q,\mathbf s}
                   -\prod_jc(q_j,s_j)\right|
 \le C\{L^{-\alpha}+n^{-\gamma}L^{M_K+\zeta}\}.
\]
Take $L$ of order $n^{\gamma/(M_K+\zeta+\alpha)}$, multiply by
\eqref{eq:mixed-opuc-limit}, and use
$c(q,s)D_\beta^{\rm rj}(q,s)=C_\beta^{\rm rj}(q,s)$.
Letting $\alpha,\gamma$ increase to their upper limits and $\zeta$
decrease to zero proves every stated real power rate. Replacing $n$ by $N=n+1$ costs $1+O(N^{-1})$.

The local complex argument is as follows; it specifies the scope of
the continuation.
Choose a compact complex neighborhood $K_{\C}$ whose real parts stay
in a slightly larger admissible real compact set, with
$\max_j|\Im s_j|\le\epsilon<1$ and $h\le h_0$.
For an integer mode, write $L_{\boldsymbol\ell}=1+\|\boldsymbol\ell\|_1$.
The effective damping in its normalized recursion satisfies
\[
 d_{\rm eff}:=\frac2\beta\left(\|\boldsymbol\ell\|^2-
                       \Im\mathbf s\cdot\boldsymbol\ell\right)
 \ge(1-\epsilon)d(\boldsymbol\ell),
\]
because $\|\boldsymbol\ell\|_1\le\|\boldsymbol\ell\|^2$.
We choose a fixed exponent
$0<\gamma<\min\{(1-\epsilon)4/\beta,1-h_0\}$.
For every nonzero mode with zero sum, $d_{\rm eff}>\gamma$.
The logarithmic propagator estimate in \cref{lem:mixed-corrector},
after normalization by the zero mode power, gives
\[
 |n^{-\mathfrak a}\E T_n|
       \le C L_{\boldsymbol\ell}^4 n^{-\gamma}.
\]
We justify the constant uniformly in the mode. Start the corrected
recursion at $n_*=\lceil C_KL_{\boldsymbol\ell}^4\rceil$.
The quadratic errors in the logarithmic propagator estimate then
have bounded sum. Thus the normalized propagator from $r$ to $n$
is bounded by $C(r/n)^{d_{\rm eff}}$ for $n_*\le r\le n$.
The real envelope bounds the normalized initial value by $Cn_*^h$;
the correction at $n_*$ is absorbed in this bound. Its contribution
at time $n$ is at most
\[
 Cn^{-\gamma}n_*^{h+\gamma}
       \le C L_{\boldsymbol\ell}^4n^{-\gamma},
\]
since $h+\gamma\le h_0+\gamma<1$.
The normalized forcing at index $k$ is
$O(L_{\boldsymbol\ell}^4k^{h-2})$. To sum it without losing
uniformity near a damping exponent equal to $1-h$, we use
\[
 \begin{split}
 n^{-d_{\rm eff}}\sum_{k=n_*}^n k^{d_{\rm eff}+h-2}
 &\le n^{-\gamma}\sum_{k=n_*}^n k^{\gamma+h-2}\\
 &\le n^{-\gamma}\sum_{k=1}^{\infty}k^{\gamma+h_0-2}
 \le C_{\gamma,h_0}n^{-\gamma}.
 \end{split}
\]
Removing the endpoint correction costs
$O(L_{\boldsymbol\ell}^2n^{h-1})$, which is absorbed by the same bound.
For $n<n_*$, the real envelope gives
$|n^{-\mathfrak a}\E T_n|\le Cn^h
 \le Cn_*^{h+\gamma}n^{-\gamma}$, so the mode estimate holds
in that range as well. All constants here may depend on the fixed
separation and the parameter compact set, but not on the mode.

The total variation of the normalized complex weighted phase measure
is bounded by $Cn^h$ times the real mixed phase law at
$(\Re\mathbf q,\Re\mathbf s)$.
By \cref{lem:mixed-terminal-decomposition}, this real phase law has
a full bounded density at fixed separation. The normalized complex
measure therefore has a density whose absolute value is at most
$Cn^{h_0}$.
The translation estimate for $W_{\mathbf q,\mathbf s}$ is uniform on this complex compact
set, with any $0<\alpha<\min\{1,\inf_{K_{\C},j}(1+\Re q_j)\}$.
Write $A_{\C}=\sup_{K_{\C}}\sum_j(-\Re q_j)_+$ for its total negative
root strength. For $m\ge2$, the strict inequality $A_{\C}<m$ permits
$\zeta>0$ such that $(A_{\C}-1)_++4+\zeta\le m+3$.
The complex Fourier mass bound in \cref{lem:terminal-truncated-mass},
which uses \eqref{eq:mixed-terminal-fourier-decay}, combines with
$L_{\boldsymbol\ell}^4\le C_mL^4$ for modes of degree at most $L$
to give the power $L^{m+3}$ below.
The same approximation consequently has error
\[
 C\{n^{h_0}L^{-\alpha}+n^{-\gamma}L^{m+3}\}.
\]
First take the neighborhood so small that $h_0(m+3)<\alpha\gamma$,
and then choose $r$ with $h_0/\alpha<r<\gamma/(m+3)$ and
$L=\lfloor n^r\rfloor$. Both powers decrease. Together with the
holomorphic OPUC limit \eqref{eq:mixed-opuc-limit}, this proves local uniform mixed asymptotics
with a positive power error. The limiting factors are nonzero near
the compact real parameter set, so the error can be written in relative
form. The case $m=1$ follows directly from
\eqref{eq:mixed-onepoint-morris}.
\end{proof}

\begin{remark}
Lemma~\ref{lem:mixed-terminal-decomposition} treats all integrable root
exponents in the mixed weight at fixed separation. On $\mathcal D_m$,
\cref{thm:mixed-uniform-FH} below gives constants independent of the
separation. For arbitrary signed roots,
\cref{thm:full-signed-polynomial} gives a quantified merging range.
\end{remark}

\subsection{Counting eigenvalues in a fixed arc}
\label{sec:fixed-arc-counting}

The mixed formula gives an explicit refinement of Gaussian counting
fluctuations. The central limit theorem was proved by Killip
\cite{Killip2008}; Feng, Tian and Wei \cite{FengTianWei2023} give a
Berry--Esseen estimate and uniform variance bounds, including arcs
whose lengths vary with $N$. For a fixed arc we obtain the constant
term in the variance and all limiting cumulants of higher order.
The formulas below were predicted through Fisher--Hartwig asymptotics
by Fyodorov and Le Doussal \cite{FyodorovLeDoussal2020} and by Smith,
Le Doussal, Majumdar and Schehr \cite[Section~3]{SmithEtAl2021}.
Thus this corollary proves their circular ensemble formulas; it does
not address their further predictions for general trapping potentials
or extrema of the counting function.

Let $I=(\vartheta_1,\vartheta_2)$ be a positively oriented arc of
fixed length $\Delta\in(0,2\pi)$, with lifts satisfying
$\vartheta_2-\vartheta_1=\Delta$. Denote the eigenvalue count, its
centered version and the endpoint chord length by
\[
 \mathcal N_N(I)=\#\{j:e^{i\theta_j}\in I\},\qquad
 \mathcal X_N(I)=\mathcal N_N(I)-\frac{N\Delta}{2\pi},\qquad
 d=2\sin(\Delta/2).
\]
The choice of open or closed endpoints does not affect these random
variables almost surely. For a random variable $Z$ with an exponential
moment near zero, write
$\kappa_j(Z)=\left.\frac{d^j}{dt^j}\log\E e^{tZ}\right|_{t=0}$.
Let $\psi^{(j)}$ denote the $j$th derivative of the digamma function,
so that $\psi^{(1)}=(\log\Gamma)''$.

\begin{corollary}
\label{cor:fixed-arc-counting}
For every $\beta>0$, locally uniformly for complex $t$ in a
neighborhood of zero,
\begin{equation}\label{eq:counting-transform}
 \E e^{t\mathcal X_N(I)}
 =(Nd)^{t^2/(\beta\pi^2)}
   \mathcal C_\beta\left(\frac{it}{2\pi},-\frac{it}{2\pi}\right)^2
       \bigl(1+O(N^{-\kappa})\bigr)
\end{equation}
for some $\kappa>0$. The formula also holds uniformly on each compact
set of real $t$, with a possibly different positive exponent. In particular,
\begin{align}
 \Var\mathcal N_N(I)
 &=\frac{2}{\beta\pi^2}\{\log(Nd)+\mathfrak c_\beta\}
                +O(N^{-\kappa}),\label{eq:counting-variance}\\
 \mathfrak c_\beta
 &=\gamma_{\rm E}+
   \sum_{k=0}^\infty\left\{
      \frac\beta2\psi^{(1)}\left(1+\frac{\beta k}{2}\right)
                  -\frac1{k+1}\right\}.\label{eq:counting-variance-constant}
\end{align}
For every integer $r\ge2$,
\begin{equation}\label{eq:counting-even-cumulants}
 \lim_{N\to\infty}\kappa_{2r}(\mathcal X_N(I))
 =\frac{4(-1)^{r+1}}{(2\pi)^{2r}}
       \sum_{k=0}^\infty
         \psi^{(2r-1)}\left(1+\frac{\beta k}{2}\right).
\end{equation}
All odd cumulants of order at least three tend to zero, and
$\kappa_1(\mathcal X_N(I))=0$ exactly. Each fixed higher cumulant
converges with a positive power error.
\end{corollary}

\begin{proof}[Proof of \cref{cor:fixed-arc-counting}]
For $u\in(0,2\pi)$,
$\Im\log(1-e^{iu})=(u-\pi)/2$, with the radial logarithm used here.
Its slope is $1/2$ between consecutive zeros, and it has a jump of
$-\pi$ at a zero. Summing over the eigenvalues gives the exact identity
\begin{equation}\label{eq:counting-argument-identity}
 J_N(\vartheta_2)-J_N(\vartheta_1)
 =\frac{N\Delta}{2}-\pi\mathcal N_N(I)
 =-\pi\mathcal X_N(I).
\end{equation}
Apply \cref{thm:mixed-FH} with
\[
 (q_1,q_2)=(0,0),\qquad (s_1,s_2)=(t/\pi,-t/\pi).
\]
The power of $N$ is $t^2/(\beta\pi^2)$ and the chord factor is
$d^{t^2/(\beta\pi^2)}$. Symmetry of $\mathcal C_\beta(a,c)$
in its two arguments makes the two one point constants equal.
This proves \eqref{eq:counting-transform}, including its local
complex uniformity.

The limiting factor is analytic and nonzero near zero. On a smaller
closed complex disk, the relative error tends uniformly to zero;
its analytic logarithm is therefore $O(N^{-\kappa})$ there.
Cauchy's formula gives the same bound for each fixed derivative at
zero. We can consequently differentiate the logarithm of
\eqref{eq:counting-transform} term by term.

Set $u=t/(2\pi)$ and recall $x_k=1+\beta k/2$.
The defining gamma series \eqref{eq:mixedconstantseries} gives
\begin{align*}
 2\log\mathcal C_\beta(iu,-iu)
 ={}&\frac{4u^2}{\beta}\gamma_{\rm E}\\
 &+2\sum_{k=0}^\infty\left\{
   2\log\Gamma(x_k)-\log\Gamma(x_k+iu)-\log\Gamma(x_k-iu)
                       -\frac{2u^2}{\beta(k+1)}\right\}.
\end{align*}
The coefficient of $t^2$ in this expression is
$\mathfrak c_\beta/(\beta\pi^2)$.
For $r\ge2$ the coefficient of $t^{2r}$ is
\[
 \frac{4(-1)^{r+1}}{(2r)!(2\pi)^{2r}}
            \sum_{k=0}^\infty\psi^{(2r-1)}(x_k).
\]
These series converge: the subtracted quadratic summand is
$O((k+1)^{-2})$, while
$\psi^{(2r-1)}(x_k)=O((k+1)^{1-2r})$ for $r\ge2$.
All odd coefficients vanish. Rotation invariance gives
$\E\mathcal N_N(I)=N\Delta/(2\pi)$ and completes the proof.
\end{proof}

For comparison with the classical cases, the series
\eqref{eq:counting-variance-constant} yields
\[
 \mathfrak c_1=\gamma_{\rm E}+1-\frac{\pi^2}{8},\qquad
 \mathfrak c_2=\gamma_{\rm E}+1,\qquad
 \mathfrak c_4=\gamma_{\rm E}+1+\log2+\frac{\pi^2}{8}.
\]
These evaluations follow from
$\psi^{(1)}(x)=\sum_{n\ge0}(x+n)^{-2}$ by separating the even and
odd indices in finite partial sums and then passing to the limit.
The constant denoted $c_\beta$ in \cite[equation~(56)]{SmithEtAl2021}
is $\mathfrak c_\beta+\log2$: that paper uses
$\log(N\sin(\Delta/2))$ in place of $\log(Nd)$.
At $\beta=2$, \eqref{eq:counting-even-cumulants} gives
$\lim\kappa_4=-3\zeta(3)/(2\pi^4)$, where $\zeta$ is the Riemann
zeta function.
For $\beta=2$, Gaussian moderate deviations for fixed arc counts
follow from the stronger lattice estimates in the unpublished manuscript
of Khodiakova, Leh\'ericy, M\'eliot and Nikeghbali
\cite[Corollary~4.12 and Section~8]{KLMN2026}.
These include individual counting probabilities, precise tails with
integer rounding, and conditional excess laws. The following
corollary, \cref{cor:fixed-arc-moderate-deviations}, extends the Gaussian
tail ratio to every $\beta>0$, with an explicit relative error. The proof uses the exponential
tilting argument of \cite{FMN2016}, with only local complex estimates,
and retains any atom at the threshold.

Separate the growing Gaussian variance parameter from the limiting
analytic factor by setting
\begin{equation}\label{eq:counting-moderate-normalization}
 v_N=\frac{2}{\beta\pi^2}\log(Nd),\qquad
 \Psi_\beta(z)=
 \mathcal C_\beta\left(\frac{iz}{2\pi},-\frac{iz}{2\pi}\right)^2.
\end{equation}
Then \eqref{eq:counting-transform} reads
\begin{equation}\label{eq:counting-local-mod-Gaussian}
 e^{-v_Nz^2/2}\E e^{z\mathcal X_N(I)}
 =\Psi_\beta(z)\bigl(1+O(N^{-\kappa})\bigr)
\end{equation}
locally uniformly for complex $z$ near zero. This is local
mod-Gaussian convergence with parameter $v_N$ and limiting function
$\Psi_\beta$: multiplication by $e^{-v_Nz^2/2}$ makes the transforms
converge locally uniformly on a complex neighborhood of zero. We take $N$ sufficiently
large that $Nd>1$. This function is even,
$\Psi_\beta(0)=1$, and
\begin{equation}\label{eq:counting-residue-quadratic}
 \log\Psi_\beta(z)
 =\frac{\mathfrak c_\beta}{\beta\pi^2}z^2+O(z^4),
 \qquad z\longrightarrow0.
\end{equation}
The variance constant in \eqref{eq:counting-variance} is thus included
in $\Psi_\beta$. At $\beta=2$, writing $t_N=\pi^{-2}\log N$, the
normalization in the unpublished manuscript
\cite[Section~8]{KLMN2026} is exactly the same:
\[
 v_N=t_N+\pi^{-2}\log d,\qquad
 \Phi_I(z)=d^{z^2/(2\pi^2)}\Psi_2(z),
\]
where $\Phi_I$ denotes the limiting factor for the arc $I$ in that manuscript.
Let $\Phi_{\rm G}$ be the standard Gaussian distribution function and
$\overline\Phi_{\rm G}$ its upper tail:
\[
 \Phi_{\rm G}(x)=\frac1{\sqrt{2\pi}}
                  \int_{-\infty}^{x}e^{-u^2/2}\dd u,
 \qquad \overline\Phi_{\rm G}(x)=1-\Phi_{\rm G}(x),
 \qquad x\in\R.
\]

\begin{corollary}[Moderate deviations for a fixed arc]
\label{cor:fixed-arc-moderate-deviations}
Fix $\beta>0$ and an arc $I$ of length $\Delta\in(0,2\pi)$.
Let $(b_N)_{N\ge1}$ be a sequence of nonnegative real numbers such
that $b_N=o(\sqrt{v_N})$. Then, uniformly for $0\le x\le b_N$ and
$\epsilon\in\{-1,1\}$,
\begin{equation}\label{eq:counting-moderate-ratio}
 \frac{\Prob\{\epsilon\mathcal X_N(I)\ge x\sqrt{v_N}\}}
      {\overline\Phi_{\rm G}(x)}
 =1+O\left(\frac{1+x}{\sqrt{v_N}}\right).
\end{equation}
The implied constant depends only on $\beta$ and $I$; the index after
which the estimate holds may also depend on the sequence $(b_N)$.
In particular, for every sequence $x_N\longrightarrow\infty$ with
$x_N=o(\sqrt{v_N})$ and each $\epsilon\in\{-1,1\}$,
\begin{equation}\label{eq:counting-moderate-equivalent}
 \Prob\{\epsilon\mathcal X_N(I)\ge x_N\sqrt{v_N}\}
 \sim\frac{e^{-x_N^2/2}}{x_N\sqrt{2\pi}}.
\end{equation}
\end{corollary}

\begin{proof}[Proof of \cref{cor:fixed-arc-moderate-deviations}]
We use the local analytic estimate to control exponentially tilted
laws, compare those laws with a Gaussian distribution, and then reverse
the change of measure. The integration at the threshold keeps any atom
there. We prove the upper tail first.
Throughout this proof $X_N$ denotes the centered count
$\mathcal X_N(I)$. Write its transform and its Gaussian normalization as
\[
 M_N(z)=\E e^{zX_N},\qquad
 A_N(z)=e^{-v_Nz^2/2}M_N(z).
\]
The random variable $X_N$ is bounded for each $N$, so $M_N$ and
$A_N$ are entire. Choose $\rho>0$ such that the closed disk
$\{|z|\le4\rho\}$ lies in the neighborhood of
\eqref{eq:counting-local-mod-Gaussian}, and $\Psi_\beta$ has no
zero on this disk. On this disk we may write
\[
 A_N(z)=\Psi_\beta(z)(1+R_N(z)),\qquad
 \sup_{|z|\le4\rho}|R_N(z)|\le C N^{-\kappa},
\]
where $R_N$ is analytic. Cauchy's formula on circles of radius
$\rho$ therefore gives
\begin{equation}\label{eq:counting-moderate-Cauchy}
 \sup_{|z|\le3\rho}
       \bigl(|R_N(z)|+|R_N'(z)|\bigr)
 \le C N^{-\kappa}.
\end{equation}
Since $\Psi_\beta$ and its derivative are bounded on the smaller
disk, and $\Psi_\beta(h)>0$ for real $h$ near zero, it follows,
for a sufficiently large index $N_0$, that
\begin{equation}\label{eq:counting-moderate-analytic-bounds}
 \sup_{|z|\le3\rho}\bigl(|A_N(z)|+|A_N'(z)|\bigr)\le C,
 \qquad \inf_{|h|\le\rho}A_N(h)\ge c>0,
 \qquad N\ge N_0.
\end{equation}
Here and below, $h$ is real, and the constants are independent of
$N$ and $h$ in the indicated ranges.

For $|h|\le\rho$, define a probability measure $\Prob_{N,h}$ by
\[
 \frac{\mathrm d\Prob_{N,h}}{\mathrm d\Prob}
 =\frac{e^{hX_N}}{M_N(h)},
\]
and let $F_{N,h}$ be the distribution function, under this measure,
of
\[
 W_{N,h}=\frac{X_N-v_Nh}{\sqrt{v_N}}.
\]
We first show that these distribution functions are uniformly close
to the standard Gaussian distribution function. Their characteristic
functions satisfy the exact identity
\begin{align*}
 \E_{N,h}e^{iuW_{N,h}}
 &=e^{-iuh\sqrt{v_N}}
       \frac{M_N(h+iu/\sqrt{v_N})}{M_N(h)}\\
 &=e^{-u^2/2}
       \frac{A_N(h+iu/\sqrt{v_N})}{A_N(h)}.
\end{align*}
If $|u|\le\rho\sqrt{v_N}$, the segment joining $h$ and
$h+iu/\sqrt{v_N}$ lies in $\{|z|\le2\rho\}$.
Integrating $A_N'$ along this segment and using
\eqref{eq:counting-moderate-analytic-bounds} gives
\[
 \left|\E_{N,h}e^{iuW_{N,h}}-e^{-u^2/2}\right|
 \le\frac{C|u|}{\sqrt{v_N}}e^{-u^2/2}.
\]
The Fourier smoothing inequality for a comparison distribution with
bounded density, recalled in
\cite[Section~2.4, especially Theorem~19]{FMN2019}, states here
that, for every $T>0$,
\[
 \sup_{y\in\R}|F_{N,h}(y)-\Phi_{\rm G}(y)|
 \le C\int_{-T}^{T}
    \frac{\left|\E_{N,h}e^{iuW_{N,h}}-e^{-u^2/2}\right|}
         {|u|}\dd u+\frac CT.
\]
The integrand at zero is interpreted by its limiting bound.
Taking $T=\rho\sqrt{v_N}$ yields
\begin{equation}\label{eq:counting-tilted-Gaussian-bound}
 \sup_{|h|\le\rho}\sup_{y\in\R}
      |F_{N,h}(y)-\Phi_{\rm G}(y)|
 \le\frac C{\sqrt{v_N}}.
\end{equation}
The center $v_Nh$ need not equal the tilted mean: the characteristic
function calculation proves precisely the estimate with this center.

Now let $0\le x\le\rho\sqrt{v_N}$ and take
$h=x/\sqrt{v_N}$. Reversing the change of measure gives
\begin{align}
 \Prob\{X_N\ge x\sqrt{v_N}\}
 &=M_N(h)e^{-x^2}
       \E_{N,h}\left[e^{-xW_{N,h}}
                         \one_{\{W_{N,h}\ge0\}}\right]\notag\\
 &=e^{-x^2/2}A_N(h)
       \int_{[0,\infty)}e^{-xw}\,\mathrm dF_{N,h}(w).
       \label{eq:counting-moderate-change-measure}
\end{align}
The closed endpoint in this integral retains a possible atom at zero.
For any distribution function $F$ and $x>0$, integration by parts
gives
\[
 \int_{[0,\infty)}e^{-xw}\,\mathrm dF(w)
 =-F(0-)+x\int_0^\infty e^{-xw}F(w)\dd w.
\]
Consequently, if
$\sup_{y\in\R}|F(y)-\Phi_{\rm G}(y)|\le D$, then
\begin{equation}\label{eq:counting-moderate-endpoint-bound}
 \left|\int_{[0,\infty)}e^{-xw}\,\mathrm dF(w)
       -\int_0^\infty e^{-xw}\,
                            \mathrm d\Phi_{\rm G}(w)\right|
 \le 2D.
\end{equation}
Indeed, the endpoint term has absolute value at most $D$, by
continuity of $\Phi_{\rm G}$, and the integral term is at most
$Dx\int_0^\infty e^{-xw}\dd w=D$.
For $x=0$, the same conclusion follows from
$\Prob\{W\ge0\}=1-F(0-)$, with the sharper bound $D$.

Completing the square gives
\begin{equation}\label{eq:counting-moderate-Gaussian-integral}
 J(x):=\int_0^\infty e^{-xw}\,\mathrm d\Phi_{\rm G}(w)
 =e^{x^2/2}\overline\Phi_{\rm G}(x),\qquad x\ge0.
\end{equation}
Moreover, restricting the integral to $0\le w\le(1+x)^{-1}$
shows directly that
\[
 J(x)=\frac1{\sqrt{2\pi}}\int_0^\infty e^{-xw-w^2/2}\dd w
 \ge\frac{e^{-3/2}}{\sqrt{2\pi}(1+x)},\qquad x\ge0.
\]
Applying \eqref{eq:counting-moderate-endpoint-bound} with
\eqref{eq:counting-tilted-Gaussian-bound}, and dividing by this
lower bound for $J(x)$, we obtain
\[
 \int_{[0,\infty)}e^{-xw}\,\mathrm dF_{N,h}(w)
 =J(x)\left(1+O\left(\frac{1+x}{\sqrt{v_N}}\right)\right).
\]
Substituting in \eqref{eq:counting-moderate-change-measure} yields
\begin{equation}\label{eq:counting-moderate-residue-form}
 \frac{\Prob\{X_N\ge x\sqrt{v_N}\}}
      {\overline\Phi_{\rm G}(x)}
 =\Psi_\beta\left(\frac{x}{\sqrt{v_N}}\right)
    \left(1+O\left(\frac{1+x}{\sqrt{v_N}}+N^{-\kappa}\right)\right),
 \qquad 0\le x\le\rho\sqrt{v_N}.
\end{equation}

For $0\le x\le b_N=o(\sqrt{v_N})$, the quadratic expansion
\eqref{eq:counting-residue-quadratic} implies
\[
 \Psi_\beta\left(\frac{x}{\sqrt{v_N}}\right)
 =1+O\left(\frac{x^2}{v_N}\right).
\]
For all sufficiently large $N$, we have $b_N\le\sqrt{v_N}$, so
$x^2/v_N\le x/\sqrt{v_N}$ throughout this range. Also
$N^{-\kappa}=O(v_N^{-1/2})$. These two observations give
\eqref{eq:counting-moderate-ratio} for $\epsilon=1$.
For $\epsilon=-1$, replace $X_N$ by $-X_N$. Its normalized
transform is $A_N(-z)$, and the evenness of $\Psi_\beta$ gives
the same estimates, with the same constants after taking their
maximum. Finally,
\[
 \overline\Phi_{\rm G}(x)
 =\frac{e^{-x^2/2}}{x\sqrt{2\pi}}\bigl(1+O(x^{-2})\bigr),
 \qquad x\longrightarrow\infty,
\]
which proves \eqref{eq:counting-moderate-equivalent}.
\end{proof}

The Gaussian approximation in
\eqref{eq:counting-tilted-Gaussian-bound} has error of order
$(\log N)^{-1/2}$, as in
\cite[Theorem~1.1]{FengTianWei2023}.
In \eqref{eq:counting-moderate-residue-form}, the difference
$\Psi_\beta(x/\sqrt{v_N})-1$ is smaller than the stated error
when $x=o(\sqrt{v_N})$.
The lattice local and precise tail estimates in the unpublished
manuscript \cite{KLMN2026} additionally use control of the tilted
characteristic function over
a full Fourier period. We do not establish that control for general
$\beta$ here. Our conclusion concerns fixed arcs; individual
counting probabilities and growing windows of the
$\mathsf{Sine}_\beta$ process are not treated.
The results for every $\beta>0$ in the same manuscript
\cite[Section~5]{KLMN2026} concern the different random variable $\log|P_N(1)|$.

\subsection{Mixed moments at approaching points}
\label{sec:mixed-varying-separation}

The real mixed weights retain the positivity needed for a uniform
comparison of OPUC moments. The Fourier coefficients of $h_{q,s}$ also
have the same decay as those of the corresponding root factor. These
two facts give a mixed formula above the microscopic scale on
$\mathcal D_m$, where the total negative root strength is less than one.

\begin{theorem}
\label{thm:mixed-uniform-FH}
Let $K\subset\mathcal D_m\times\mathbb R^m$ be compact and set
$\varepsilon_K=1-\sup_K\sum_j(-q_j)_+>0$. For every
\begin{equation}\label{eq:mixed-uniform-rate}
 0<\kappa<\varepsilon_K\min\{2/\beta,1/2\},
\end{equation}
there are $C,R<\infty$ such that \eqref{eq:mixed-FH} holds with
relative error $O((N\delta)^{-\kappa})$, uniformly for
$(\mathbf q,\mathbf s)\in K$, $\boldsymbol\vartheta\in\mathcal S_\delta$
and $N\delta\ge R$. The constants do not depend on $\delta$.
\end{theorem}

To prove the theorem, we combine a weighted absolute sum of Fourier
coefficients with a phase estimate uniform in the minimum separation.
These are \cref{lem:mixed-terminal-fourier,lem:mixed-relative-uniform-modes}.

\begin{lemma}
\label{lem:mixed-terminal-fourier}
For complex $q,s$ with $\Re q>-1$ and $k\in\mathbb Z$,
\begin{equation}\label{eq:mixed-terminal-fourier}
 \widehat h_{q,s}(k)
 =\frac{(-1)^k\Gamma(1+q)}
 {\Gamma(1+(q-is)/2-k)\Gamma(1+(q+is)/2+k)}.
\end{equation}
Uniformly on compact subsets of $\{\Re q>-1\}\times\C$,
\begin{equation}\label{eq:mixed-terminal-fourier-decay}
 |\widehat h_{q,s}(k)|\le C(1+|k|)^{-1-\Re q}.
\end{equation}
If $K\subset\mathcal D_m\times\mathbb R^m$ is compact and
$\varepsilon_K=1-\sup_K\sum_j(-q_j)_+$, then for every
$0<u<\varepsilon_K$,
\begin{equation}\label{eq:mixed-terminal-weighted-fourier}
 \sum_{\sum_j\ell_j=0}(1+\|\boldsymbol\ell\|_1)^u
     \prod_j|\widehat h_{q_j,s_j}(\ell_j)|\le C_{K,u}.
\end{equation}
\end{lemma}

\begin{proof}[Proof of \cref{lem:mixed-terminal-fourier}]
We first prove the coefficient formula for real $q>-1$ and real $s$.
On $0<t<2\pi$, the chosen argument is $(t-\pi)/2$, so
$h'_{q,s}/h_{q,s}=(q/2)\cot(t/2)+s/2$. The function
$\sin(t/2)h_{q,s}(t)$ vanishes at both endpoints and has an integrable
derivative. Integration by parts against $e^{-i(k-1/2)t}$ gives
\[
 (q+2-2k-is)\widehat h_{q,s}(k-1)
              +(q+2k+is)\widehat h_{q,s}(k)=0.
\]
The zero coefficient is $c(q,s)$ in \eqref{eq:mixed-terminal-c}.
Iterating the recurrence proves \eqref{eq:mixed-terminal-fourier}
for nonnegative indices. The identity
$h_{q,s}(-t)=h_{q,-s}(t)$ holds away from the singularity, with the
arguments interpreted modulo $2\pi$. It gives
$\widehat h_{q,s}(-k)=\widehat h_{q,-s}(k)$ and therefore proves the
formula for negative indices as well.

Both sides of \eqref{eq:mixed-terminal-fourier} are holomorphic in
$(q,s)$ on $\{\Re q>-1\}\times\C$. For the integral, this follows
from locally uniform domination at the endpoints: the argument factor
is bounded on compact parameter sets, and powers of
$|\log|1-e^{it}||$ remain integrable against
$|1-e^{it}|^{\Re q}$. For the gamma quotient, the numerator has no
pole on this domain and reciprocal gamma functions are entire.
For each real $s$, the identity theorem first extends the formula
in $q$; for each resulting complex $q$, it then extends the formula
in $s$. This also includes vanishing coefficients.

For sufficiently large positive $k$, uniformly on the chosen compact
set, reflection rewrites the coefficient as
\[
 -\frac{\Gamma(1+q)\sin(\pi(q-is)/2)}\pi
 \frac{\Gamma(k-(q-is)/2)}{\Gamma(k+1+(q+is)/2)}.
\]
The difference of the two gamma shifts is $-1-q$, so uniform gamma
ratio estimates give \eqref{eq:mixed-terminal-fourier-decay} for positive
indices. Negative indices follow by replacing $s$ by $-s$.
The finitely many remaining coefficients are bounded by continuity
on the compact parameter set. For real parameters in the weighted
assertion, \cref{lem:terminal-weighted-fourier} uses only this
coefficient bound and therefore proves
\eqref{eq:mixed-terminal-weighted-fourier} as well.
\end{proof}

\begin{lemma}
\label{lem:mixed-relative-uniform-modes}
For a compact real set $K\subset(-1,\infty)^m\times\mathbb R^m$
and every nonzero zero sum integer mode,
\begin{equation}\label{eq:mixed-relative-uniform-modes}
 \frac{|\E[\mathcal B_n e^{i\boldsymbol\ell\cdot\boldsymbol\psi_n}]|}
 {\E\mathcal B_n}
 \le C_K(1+\|\boldsymbol\ell\|_1)^2
                       \mathfrak r_{4/\beta}(n\delta)
\end{equation}
when $n\delta\ge1$. The constant is independent of the mode, the
configuration and $0<\delta\le1$. In particular, the right side may
be replaced by $C_{K,\gamma}(1+\|\boldsymbol\ell\|_1)^2
(n\delta)^{-\gamma}$ for every $0<\gamma<g_\beta$,
where $g_\beta=\min\{4/\beta,1\}$. The rate function
$\mathfrak r_{4/\beta}$ is defined in \eqref{eq:phase-rate-function}.
\end{lemma}

\begin{proof}[Proof of \cref{lem:mixed-relative-uniform-modes}]
The proof of \cref{lem:mixed-corrector} already retains the actual
mode exponent in \eqref{eq:mixed-actual-propagator}. It remains to
check that its normalization is uniform in the separation. Put
$A_k=\E\mathcal B_k$, $n_*=\lceil CL^2/\delta\rceil$ and
$h=A_{n_*}/p_{n_*}(\mathfrak a)$. The relative mixed OPUC estimate
\eqref{eq:mixed-opuc-relative} gives
\[
 C^{-1}\le\frac{A_k}{p_k(\mathfrak a)h}\le C,
                       \qquad k\ge n_*,
\]
with constants independent of $\delta$ and the mode. Here
$L=1+\|\boldsymbol\ell\|_1$ and $p_k$ is the gamma normalization
recalled in the proof of \cref{lem:mixed-corrector}. Thus the normalized
forcing and propagator in that proof have the uniform bounds
$CL^4/(k^2\delta)$ and $C(r/n)^d$. Summing as in
\eqref{eq:mixed-actual-propagator}, removing the correction, and using
the probability bound for $n<n_*$ prove
\eqref{eq:mixed-relative-uniform-modes}. No estimate for the initial
mass separately from its normalizer is required.
\end{proof}

\begin{proof}[Proof of \cref{thm:mixed-uniform-FH}]
The factors $e^{s_j\operatorname{Arg}(1-e^{it})}$ are uniformly bounded
above and below on a compact real parameter set. The H\"older argument
of \cref{lem:terminal-holder} therefore proves that
$W_{\mathbf q,\mathbf s}$ is bounded and continuous on
$\mathcal D_m\times\mathbb R^m$: use the same H\"older exponents and
continuity of translations in each of the resulting $L^{p_j}$ spaces.
Its Fourier coefficient at
$\boldsymbol\ell$ is
$\mathbf1_{\{\sum_j\ell_j=0\}}\prod_j\widehat h_{q_j,s_j}(\ell_j)$,
and \eqref{eq:mixed-terminal-weighted-fourier} makes the series
absolutely convergent.

Given $\kappa$ in \eqref{eq:mixed-uniform-rate}, choose
$0<u<\varepsilon_K$ and $0<\gamma<g_\beta$ so that
$\gamma u/2>\kappa$. Positivity and
\eqref{eq:mixed-relative-uniform-modes} bound a nonconstant normalized
phase mode by
\[
 \min\{1,CL^2(n\delta)^{-\gamma}\}
       \le CL^u(n\delta)^{-\gamma u/2},
 \qquad L=1+\|\boldsymbol\ell\|_1.
\]
Integrating the absolute Fourier series gives
\[
 \frac{\E[\mathcal B_nW_{\mathbf q,\mathbf s}
                                  (\boldsymbol\psi_n)]}
      {\E\mathcal B_n}
   =\prod_jc(q_j,s_j)+O((n\delta)^{-\gamma u/2}).
\]
The mean $\prod_jc(q_j,s_j)$ is bounded below by a positive constant on $K$.
Combine this identity with \eqref{eq:mixed-terminal-identity} and
\eqref{eq:mixed-opuc-relative}. The exponent $\gamma u/2$ exceeds
the prescribed $\kappa$, while the OPUC error has exponent one.
Replacing $n$ by $N=n+1$ changes the result by $1+O(N^{-1})$,
which is absorbed in the stated error.
\end{proof}

\begin{corollary}
\label{cor:mixed-subset-rate}
For a compact set of real mixed parameters with root coordinates
in $\mathcal D_m$, the error in \cref{thm:mixed-uniform-FH} allows
every exponent in \eqref{eq:positive-terminal-rate}, with $s_K$
computed from the root coordinates.
\end{corollary}
\begin{proof}[Proof of \cref{cor:mixed-subset-rate}]
The decay \eqref{eq:mixed-terminal-fourier-decay} is the same as
for roots. Thus the shell argument in
\cref{lem:terminal-weighted-fourier} gives every weight $s<s_K$.
Apply \eqref{eq:mixed-relative-uniform-modes}, the probability bound,
and the interpolation in the proof of \cref{cor:terminal-subset-rate}.
The mixed OPUC error is $O((N\delta)^{-1})$ and is absorbed in the
stated rate.
\end{proof}

\section{Analytic backgrounds and singular loop equations}
\label{sec:smooth}

Recall that $X_N(\vartheta)=\log|P_N(e^{i\vartheta})|$. We add the
factor $e^{\sum_{\ell=1}^N V(\theta_\ell)}$ to the singular moment.
The results of this section extend \cref{thm:main-root}; its proof,
the critical moment proof and the correlation expansions do not use
this extension.
Without singularities, the background asymptotic is the strong Szeg\H{o}
asymptotic of Johansson \cite{Johansson1988}. The additional term in
\cref{thm:analytic-background} describes the interaction of the
background with the prescribed roots.

In this section $H$ denotes the circular Hilbert transform. For
$V_k=\int V(x)e^{-ikx}\,d\sigma(x)$, define its quadratic energy and
analytic norm together with our sign convention for $H$:
\[
 H(e^{ikx})=-i\operatorname{sgn}(k)e^{ikx},\quad H1=0,\qquad
 \sigma^2(V)=\sum_{k\in\mathbb Z}|k|V_kV_{-k},
 \qquad
 \|V\|_\rho=\sum_{k\in\mathbb Z}e^{\rho|k|}|V_k|.
\]
Thus $H\cos x=\sin x$, $H\sin x=-\cos x$, and
$Hf(x)=(2\pi)^{-1}\operatorname{p.v.}\int_0^{2\pi}
\cot((x-y)/2)f(y)\,dy$ for smooth periodic $f$.
Every real analytic function has finite $\|V\|_\rho$ for some $\rho>0$.
A bounded analytic ball means a family of real functions $h$ satisfying
$\|h\|_\rho\le B$ for fixed $\rho,B>0$.

\begin{theorem}\label{thm:analytic-background}
Fix $\rho,B>0$, a compact $K\subset(-1,\infty)^m$ and $\delta>0$.
Uniformly for real $V$ with $\|V\|_\rho\le B$,
$\mathbf q\in K$ and $\boldsymbol\vartheta\in\mathcal S_\delta$,
\begin{align}\label{eq:loc-full-positive-fh}
 &\E_\beta
 e^{\sum_{\ell=1}^N V(\theta_\ell)+\sum_{j=1}^m q_jX_N(\vartheta_j)}
 \notag\\
 &\quad=F_N(\mathbf q,\boldsymbol\vartheta)
 \exp\left\{NV_0+\frac{\sigma^2(V)}\beta
       -\frac1\beta\sum_{j=1}^m q_j(V(\vartheta_j)-V_0)\right\}
       (1+O(N^{-\kappa})),
\end{align}
for every $\kappa$ in \eqref{eq:full-root-rate}.
For $m=1$, the relative error in \eqref{eq:loc-full-positive-fh}
is $O(N^{-1})$.
For $K\subset\mathcal D_m$, the separation condition may be replaced
by $N\delta(\boldsymbol\vartheta)\ge R$, with error
$O((N\delta(\boldsymbol\vartheta))^{-\kappa})$ for every $\kappa$
in \eqref{eq:intro-rate} and sufficiently large $R$.
The displaced and compatible formulas
\eqref{eq:loc-moving-formula} and \eqref{eq:loc-smooth-fh}
have relative error $O(N^{-1})$ on these sets.
In particular all assertions hold uniformly for real trigonometric
polynomials in compact coefficient sets of fixed degree.
\end{theorem}

For a compact $K\subset\mathcal D_m$, the relative error in
\eqref{eq:loc-full-positive-fh} at moving prescribed points tends to
zero when $N\delta\to\infty$, uniformly under the stated bounds on
$V$ and for $\mathbf q\in K$. The estimate gives a uniform bound
for $N\delta\ge R$, but does not give a vanishing relative error
when $N\delta$ stays fixed.

The proof first obtains moments growing at most logarithmically,
then a transport formula, and finally bounded moments by ratios of
the resulting partition functions. This order avoids using the
improved moment bounds to prove the formula from which they follow.

\subsection{A change of variables and preliminary moment bounds}

Throughout this section, $Q_{N,\mathbf q,\boldsymbol\vartheta;V}$ denotes
the law tilted by the root factors and the background. We write $Q$
when these parameters are fixed; its density is proportional to
\[
 \prod_{\ell<k}|e^{i\theta_\ell}-e^{i\theta_k}|^\beta
 \prod_{\ell,j}|e^{i\theta_\ell}-e^{i\vartheta_j}|^{q_j}
 e^{\sum_\ell V(\theta_\ell)}.
\]
The empirical measure, its centered version and the centered linear
statistic are denoted by $\mu_N=\sum_{\ell=1}^N\delta_{\theta_\ell}$,
$\widetilde\mu_N=\mu_N-N\sigma$ and
$\operatorname{Fl}_N(f)=\int f\,d\widetilde\mu_N$, respectively.

\begin{lemma}\label{lem:loc-pinned}
Suppose $q_j>-1$ and $V$ is real and smooth. If a real periodic
$C^2$ function $g$ vanishes at every $\vartheta_j$, then
\begin{align}\label{eq:loc-pinned-loop}
 0=\E_Q\Bigg[
 &\left(1-\frac\beta2\right)\int g'\,d\mu_N
 +\frac\beta4\iint K_g(x,y)\,d\mu_N(x)d\mu_N(y)
 +\int gV'\,d\mu_N\notag\\
 &+\frac12\sum_jq_j\int
       g(x)\cot\frac{x-\vartheta_j}{2}\,d\mu_N(x)\Bigg],
\end{align}
where $K_g(x,y)=(g(x)-g(y))\cot((x-y)/2)$ and
$K_g(x,x)=2g'(x)$. Every term is absolutely integrable; the
root terms use their continuous extensions at the prescribed points.
\end{lemma}

\begin{proof}[Proof of \cref{lem:loc-pinned}]
Change every particle variable by $T_u(x)=x+ug(x)$ for small real
$u$. This is an orientation preserving circle diffeomorphism fixing
the prescribed points. The fixed point condition makes the root
chord ratio extend continuously across each prescribed point.
Ratios of transformed to original chords are bounded above and below,
uniformly in their endpoints. Their
logarithmic derivatives are
$\tfrac12(g(x)-g(y))\cot((x-y)/2)$ and
$\tfrac12g(x)\cot((x-\vartheta_j)/2)$; both are bounded.
For fixed $N$, the differentiated integrand is therefore dominated
by a constant times the original integrable density, also for
negative root exponents. Differentiation gives the sum over distinct particle pairs. Adding its diagonal contributes
$(\beta/2)\int g'\,d\mu_N$, which is subtracted from the Jacobian
term and gives \eqref{eq:loc-pinned-loop}.
\end{proof}

\begin{lemma}\label{lem:loc-moment-transfer}
Fix $\rho,B>0$, a compact set $K\subset(-1,\infty)^m$ and a
separation $\delta_0>0$. For every $p\ge1$,
\begin{equation}\label{eq:loc-polylog-moments}
 \E_Q|\operatorname{Fl}_N(f)|^p\le C_p(1+\log N)^p,
\end{equation}
uniformly for $\mathbf q\in K$,
$\boldsymbol\vartheta\in\mathcal S_{\delta_0}$ and real functions
$V,f$ with $\|V\|_\rho\le B$ and $\|f\|_\rho\le B$.
If $K\subset\mathcal D_m$, the bound is uniform over all distinct
configurations. These conclusions include compact coefficient
sets of fixed degree. For $\|V\|_\rho\le B$ and
$0<\varepsilon<\rho/8$ they also give
\begin{equation}\label{eq:analytic-preliminary-moments}
 \|\operatorname{Fl}_N(e^{ikx})\|_{L^p(Q)}
 \le C_p e^{\varepsilon|k|}(1+\log N),\qquad k\ne0.
\end{equation}
\end{lemma}

\begin{proof}[Proof of \cref{lem:loc-moment-transfer}]
We first transfer a smooth exponential moment bound through the root
weights by H\"older's inequality. A polynomial bound in $N$ then gives
the logarithmic moment bound by tail integration. For configurations
with arbitrarily small gaps, the separate H\"older exponents below use
the defining condition of $\mathcal D_m$.
The smooth circular exponential moment bound
\cite[Theorem~1.2 and equation~(2.13)]{Lambert2021} gives
\begin{equation}\label{eq:loc-smooth-bound}
 \sup_N\E_\beta e^{\operatorname{Fl}_N(h)}<\infty
\end{equation}
uniformly on bounded analytic balls. Indeed, the required first
and third derivatives of $h$ and the third derivative of $Hh$
are bounded by the analytic norm. For fixed separation choose
$u>1$ close enough to one that $uK\subset(-1,\infty)^m$.
H\"older's inequality and the already proved \cref{thm:full-root-FH}
yield
\[
 \E_\beta e^{\sum_jq_jX_N(\vartheta_j)+\operatorname{Fl}_N(h)}
 \le F_N(u\mathbf q,\boldsymbol\vartheta)^{1/u}
       \left(\E_\beta e^{u'\operatorname{Fl}_N(h)}\right)^{1/u'}
 \le CN^A,\qquad u'=\frac{u}{u-1}.
\]
For uniformity on $\mathcal D_m$, the inequality
$\sum_j(-q_j)_+<1$ allows weights
$w_j>\max\{-q_j,0\}$ with $\sum_jw_j<1$. Set
$p_j=1/w_j$, $r=(1-\sum_jw_j)^{-1}$. Then
$p_j>1$, $r>1$, $p_jq_j>-1$ and
$\sum_jp_j^{-1}+r^{-1}=1$. These choices persist in a neighborhood
of $\mathbf q$, so a finite covering makes the bounds uniform on $K$.
The Morris estimate from \cref{prop:onepoint},
\begin{equation}\label{eq:loc-morris-bound}
 M_N(v)\le C_I N^{v^2/(2\beta)},\qquad
 v\in I\Subset(-1,\infty),
\end{equation}
and the generalized H\"older inequality give instead
\[
 \E_\beta e^{\sum_jq_jX_N(\vartheta_j)+\operatorname{Fl}_N(h)}
 \le\prod_jM_N(p_jq_j)^{1/p_j}
       \left(\E_\beta e^{r\operatorname{Fl}_N(h)}\right)^{1/r}
 \le CN^A
\]
at every distinct configuration. The denominator defining $Q$,
after removing $NV_0$, is at least one by Jensen's inequality:
$\E_\beta X_N(\vartheta_j)=0$ by rotation invariance and
$\int_0^{2\pi}\log|1-e^{ix}|\,dx=0$.
Apply the bound with $h=V\pm f$ to obtain
$\E_Qe^{\pm\operatorname{Fl}_N(f)}\le CN^A$.
Integrating
$Q(|\operatorname{Fl}_N(f)|>t)\le\min\{1,2CN^Ae^{-t}\}$
against $pt^{p-1}\,dt$ proves \eqref{eq:loc-polylog-moments}.
Finally the real tests
$e^{-\varepsilon|k|}\cos(kx)$ and
$e^{-\varepsilon|k|}\sin(kx)$ belong to one bounded analytic
ball of width $\varepsilon/2$. Apply the same estimate to these
tests and to $V$, viewed at that smaller width, to obtain
\eqref{eq:analytic-preliminary-moments}.
\end{proof}

\subsection{Transport of the particles and prescribed points}
\label{sec:loc-moving}

\begin{theorem}\label{thm:loc-moving-insertions}
Let $V$ be real with $\|V\|_\rho\le B$, let $c\in\R$, and set
$g=HV-c$. Define
\[
 \partial_t\phi_{N,t}(x)=-\frac{2}{\beta N}g(\phi_{N,t}(x)),
 \qquad \phi_{N,0}(x)=x,\qquad 0\le t\le1.
\]
For $\mathbf q\in(-1,\infty)^m$,
\begin{align}\label{eq:loc-moving-formula}
 &\frac{\E_\beta
  e^{\sum_\ell V(\theta_\ell)
       +\sum_jq_jX_N(\phi_{N,1}(\vartheta_j))}}
       {F_N(\mathbf q,\boldsymbol\vartheta)}
 \notag\\
 &\quad=\exp\left\{NV_0+\frac{\sigma^2(V)}\beta
        -\frac1\beta\sum_jq_j(V(\vartheta_j)-V_0)\right\}
       \left(1+O\left(\frac{(1+\log N)^2}{N}\right)\right).
\end{align}
The error is uniform on the fixed separation compact parameter
sets in \cref{thm:analytic-background}, with $c$ in a bounded real set.
For $K\subset\mathcal D_m$ it is uniform at every distinct
configuration. At fixed separation, or on
$K\subset\mathcal D_m$ with $N\delta\ge R$ for sufficiently large $R$,
\cref{cor:smooth-bootstrap} improves the error to $O(N^{-1})$.
\end{theorem}

\begin{proof}[Proof of \cref{thm:loc-moving-insertions}]
We move the particles and the prescribed points by the same flow,
compute the logarithmic derivative of the partition function, and
bound its remainder using \cref{lem:loc-moment-transfer}. The Fourier
expansion below also checks the summability required for analytic $V$.
Write $\vartheta_j(t)=\phi_{N,t}(\vartheta_j)$ and introduce the
interpolating partition function
$Z_N(t)=\E_\beta
 e^{t\operatorname{Fl}_N(V)+\sum_jq_jX_N(\vartheta_j(t))}$.
In this integral change variables
$\theta_\ell=\phi_{N,t}(x_\ell)$.
Both endpoints of every singular chord move under the same flow.
Its logarithmic derivative therefore contains a difference
$g(x)-g(y)$, which removes the apparent singularity.
The chord bounds used for \cref{lem:loc-pinned} justify
differentiation of the transformed integral for fixed $N$.
Returning to the current variables gives
\begin{equation}\label{eq:loc-moving-derivative}
 \frac d{dt}\log Z_N(t)
 =\E_{Q_t}\left[\operatorname{Fl}_N(V)
                      -\frac{2}{\beta N}B_{N,t}\right],
\end{equation}
where $Q_t=Q_{N,\mathbf q,\boldsymbol\vartheta(t);tV}$ is the
corresponding tilted law, and
\[
 \begin{split}
 B_{N,t}={}&\left(1-\frac\beta2\right)\int g'\,d\mu_N
       +\frac\beta4\iint K_g\,d\mu_Nd\mu_N
       +t\int gV'\,d\mu_N\\
 &+\frac12\sum_jq_j\int
       (g(x)-g(\vartheta_j(t)))
          \cot\frac{x-\vartheta_j(t)}2\,d\mu_N(x).
 \end{split}
\]
These terms come respectively from the Jacobian, pair interactions,
smooth potential and root interactions. The difference
$g(x)-g(\vartheta_j(t))$ makes each root integrand continuous at its
prescribed point, also when the corresponding root exponent is negative.

For $g=HV-c$, Fourier calculation gives
\begin{equation}\label{eq:loc-fourier-identities}
 \begin{gathered}
 \int K_g(x,y)\,d\sigma(y)=V(x)-V_0,\qquad
 \int gV'\,d\sigma=-\sigma^2(V),\\
 \int(g(x)-g(\vartheta))\cot\frac{x-\vartheta}{2}\,d\sigma(x)
       =V(\vartheta)-V_0.
 \end{gathered}
\end{equation}
For example,
$\int(e^{ikx}-e^{iky})\cot((x-y)/2)\,d\sigma(y)
=i\operatorname{sgn}(k)e^{ikx}=-H(e^{ikx})$ for $k\ne0$.
The second identity is Parseval and $\int V'\,d\sigma=0$;
the third is the same Hilbert identity at $\vartheta$.
Writing
$h_{j,t}(x)=(g(x)-g(\vartheta_j(t)))
\cot((x-\vartheta_j(t))/2)$, expansion of
$\mu_N=N\sigma+\widetilde\mu_N$ yields
\[
 B_{N,t}=\frac{\beta N}2\operatorname{Fl}_N(V)
 -tN\sigma^2(V)+\frac N2\sum_jq_j(V(\vartheta_j(t))-V_0)+R_{N,t},
\]
where the remainder contains the quadratic and linear fluctuations:
\[
 R_{N,t}=\frac\beta4\iint K_g\,d\widetilde\mu_Nd\widetilde\mu_N
  +\left(1-\frac\beta2\right)\operatorname{Fl}_N(g')
  +t\operatorname{Fl}_N(gV')
  +\frac12\sum_jq_j\operatorname{Fl}_N(h_{j,t}).
\]

We check absolute summability before using this expansion for
analytic $V$.
For $k>0$, $z=e^{ix}$ and $w=e^{iy}$,
\[
 (z^k-w^k)\cot\frac{x-y}{2}
       =i(z+w)\sum_{a=0}^{k-1}z^{k-1-a}w^a.
\]
There are at most $2k$ monomials of total Fourier degree $k$;
negative indices follow by conjugation. Consequently, for
$0<\varepsilon<\rho/8$,
\[
 \sum_{r,s}|(K_g)_{r,s}|e^{\varepsilon(|r|+|s|)}
 \le2\sum_{k\in\mathbb Z}|k|\,|V_k|e^{\varepsilon|k|}
 \le C_{\rho,\varepsilon}B.
\]
The same weighted absolute summability holds for $g'$, $gV'$
and $h_{j,t}$, uniformly in the prescribed points. By
\eqref{eq:analytic-preliminary-moments} and the Cauchy--Schwarz
inequality, the corresponding linear and quadratic series converge absolutely
in $L^1(Q_t)$ and
$\sup_{0\le t\le1}\E_{Q_t}|R_{N,t}|
\le C(1+\log N)^2$.
The flow distorts chords by $1+O(N^{-1})$, so it preserves
distinctness and fixed separation; the moment bounds apply throughout.
This also justifies the Fourier identities by absolute convergence.

Substitution in \eqref{eq:loc-moving-derivative} cancels the
centered statistic and gives
\[
 \frac d{dt}\log Z_N(t)
 =\frac{2t}\beta\sigma^2(V)
 -\frac1\beta\sum_jq_j(V(\vartheta_j(t))-V_0)
       +O(N^{-1}(1+\log N)^2).
\]
Since $\phi_{N,t}(x)=x+O(N^{-1})$ in $C^1$,
replacing $V(\vartheta_j(t))$ by $V(\vartheta_j)$ costs
$O(N^{-1})$. Integrate from $0$ to $1$, use
$Z_N(0)=F_N(\mathbf q,\boldsymbol\vartheta)$ and restore $NV_0$.
\end{proof}

\begin{proof}[Proof of \cref{thm:analytic-background},
apart from the improved transport error]
Apply \eqref{eq:loc-moving-formula} with $c=0$ at
$\widetilde\vartheta_j=\phi_{N,1}^{-1}(\vartheta_j)$, so that the flow
ends at the prescribed points $\vartheta_j$.
Differentiating the logarithm of a chord along the flow gives
a bounded multiple of
$(g(x)-g(y))\cot((x-y)/2)/N$. Hence, uniformly for distinct pairs,
\begin{equation}\label{eq:flow-chord-distortion}
 \frac{|e^{i\widetilde\vartheta_j}-e^{i\widetilde\vartheta_k}|}
 {|e^{i\vartheta_j}-e^{i\vartheta_k}|}=1+O(N^{-1}).
\end{equation}
The Fisher--Hartwig interaction and the background values change
by $O(N^{-1})$ in their logarithms.
At fixed separation \cref{thm:full-root-FH} therefore gives
$F_N(\mathbf q,\widetilde{\boldsymbol\vartheta})/
F_N(\mathbf q,\boldsymbol\vartheta)=1+O(N^{-\kappa})$.
On $\mathcal D_m$ use \cref{thm:signed-characteristic-FH};
the two minimum separations differ by $1+O(N^{-1})$.
Finally $N^{-1}(1+\log N)^2=O(N^{-\kappa})$ for $\kappa<1$
and $N^{-\kappa}\le(N\delta)^{-\kappa}$.
This proves \eqref{eq:loc-full-positive-fh} with all its uniformities.
\end{proof}

\subsection{Fixed prescribed points and a single argument weight}

\begin{theorem}\label{thm:loc-smooth-background}
If
\begin{equation}\label{eq:loc-compatible}
 HV(\vartheta_1)=\cdots=HV(\vartheta_m),
\end{equation}
then, under the fixed separation hypotheses of
\cref{thm:analytic-background},
\begin{align}\label{eq:loc-smooth-fh}
 &\frac{\E_\beta e^{\sum_\ell V(\theta_\ell)
                    +\sum_jq_jX_N(\vartheta_j)}}
        {F_N(\mathbf q,\boldsymbol\vartheta)}
 \notag\\
 &\quad=\exp\left\{NV_0+\frac{\sigma^2(V)}\beta
       -\frac1\beta\sum_jq_j(V(\vartheta_j)-V_0)\right\}
       \left(1+O(N^{-1}(1+\log N)^2)\right).
\end{align}
For one prescribed point compatibility is automatic.
Under these hypotheses, \cref{cor:smooth-bootstrap} improves the
relative error to $O(N^{-1})$.
\end{theorem}

\begin{proof}[Proof of \cref{thm:loc-smooth-background}]
Use \cref{thm:loc-moving-insertions} with
$c=HV(\vartheta_1)$. Its vector field vanishes at every prescribed
point, so these points remain fixed. With $R_{N,t}$ as defined in the
proof of \cref{thm:loc-moving-insertions}, the pinned equation
\eqref{eq:loc-pinned-loop} becomes
\begin{equation}\label{eq:loc-main-loop}
 0=\frac{\beta N}2\E_{Q_t}\operatorname{Fl}_N(V)
 -tN\sigma^2(V)+\frac N2\sum_jq_j(V(\vartheta_j)-V_0)
 +\E_{Q_t}R_{N,t},
\end{equation}
and therefore
\begin{equation}\label{eq:loc-smooth-response}
 \E_{Q_t}\operatorname{Fl}_N(V)
 =\frac{2t}\beta\sigma^2(V)
 -\frac1\beta\sum_jq_j(V(\vartheta_j)-V_0)
 +O(N^{-1}(1+\log N)^2).
\end{equation}
The logarithmic partition derivative is
$\E_{Q_t}\operatorname{Fl}_N(V)$. Integrating
\eqref{eq:loc-smooth-response} from $0$ to $1$ and restoring $NV_0$
proves \eqref{eq:loc-smooth-fh}.
\end{proof}

\begin{remark}\label{rem:loc-first-fourier}
For one root at zero, $V=0$ and $g(x)=\sin x$,
$K_g(x,y)=\cos x+\cos y$ and
$\sin x\cot(x/2)=1+\cos x$. Thus \eqref{eq:loc-pinned-loop} gives
\[
 \E_{Q_{N,q,0}}\sum_{\ell=1}^N\cos\theta_\ell
 =-\frac{qN}{2+\beta(N-1)+q},\qquad q>-1,
\]
whose limit $-q/\beta$ agrees with \eqref{eq:loc-smooth-response}
at $t=0$ and $V(x)=\cos x$.
\end{remark}

\begin{corollary}\label{cor:loc-root-jump-background}
For real $q>-1$, $s\in\mathbb R$ and a real trigonometric polynomial $V$,
\begin{align*}
 &\frac{\E_\beta
 e^{\sum_\ell V(\theta_\ell)+qX_N(\vartheta)+sJ_N(\vartheta)}}
 {M_N^{\rm rj}(q,s)}\\
 &\quad=\exp\left\{NV_0+\frac{\sigma^2(V)}\beta
       -\frac q\beta(V(\vartheta)-V_0)
       -\frac s\beta HV(\vartheta)\right\}
       \left(1+O(N^{-1}(1+\log N)^2)\right),
\end{align*}
uniformly on compact root and argument parameters and compact
coefficient sets of fixed degree. Under the same hypotheses,
\cref{cor:smooth-bootstrap} improves the relative error to $O(N^{-1})$.
\end{corollary}

\begin{proof}[Proof of \cref{cor:loc-root-jump-background}]
For $0<x-\vartheta<2\pi$ the argument summand is
$(\pi-(x-\vartheta))/2$, with derivative $-1/2$ and mean zero.
The proof of \cref{lem:loc-moment-transfer} applies with
$M_N^{\rm rj}(uq,us)^{1/u}$ for $u>1$ with $uq>-1$,
using \eqref{eq:mixed-onepoint-morris}; Jensen's lower bound is
unchanged. Take $g=HV-HV(\vartheta)$, whose flow fixes the jump
and preserves its two sides. Differentiation adds
$-(s/2)\int g\,d\mu_N$ to \eqref{eq:loc-pinned-loop}.
Its deterministic part is $sN HV(\vartheta)/2$ and its centered
part has expectation $O(1+\log N)$.
Solving for the mean in \eqref{eq:loc-main-loop} therefore adds
$-sHV(\vartheta)/\beta$ to \eqref{eq:loc-smooth-response}.
Integration proves the assertion.
\end{proof}

Here $s$ is a real Laplace parameter. For imaginary $s$ the weight
is complex, so the positivity and Jensen arguments require a
different justification; holomorphy at each finite $N$ alone does not supply it.

\subsection{Bounded tilted moments and the improved error}

\begin{corollary}\label{cor:smooth-bootstrap}
Under the fixed separation hypotheses of
\cref{thm:analytic-background}, for every $p\ge1$,
\begin{equation}\label{eq:bounded-tilted-moments}
 \sup_N\E_Q|\operatorname{Fl}_N(f)|^p<\infty
\end{equation}
uniformly for $f$ in a bounded analytic ball; in particular this
holds on compact coefficient sets of fixed degree.
The same conclusion holds for $K\subset\mathcal D_m$ and
$N\delta(\boldsymbol\vartheta)\ge R$, with sufficiently large $R$.
At fixed separation, the relative errors in
\eqref{eq:loc-moving-formula} and \eqref{eq:loc-smooth-fh} improve to $O(N^{-1})$. The same improvement holds for
$K\subset\mathcal D_m$ when $N\delta(\boldsymbol\vartheta)\ge R$,
with sufficiently large $R$.
The relative error in \cref{cor:loc-root-jump-background} also improves
to $O(N^{-1})$.
\end{corollary}

\begin{proof}[Proof of \cref{cor:smooth-bootstrap}]
We first bound tilted exponential moments by ratios of partition
functions. These bounds improve the remainder in the transport
identities from $O((1+\log N)^2)$ to $O(1)$.
The exponential moment of $\operatorname{Fl}_N(f)$ under $Q$ is
the ratio of the background partition functions at $V+f$ and $V$,
after removing $Nf_0$. The already proved
\eqref{eq:loc-full-positive-fh} cancels the root factor and bounds
the deterministic ratio. The factors multiplying the main terms in
\eqref{eq:loc-full-positive-fh} are bounded above and away from zero
for large $N$. Replacing $f$ by $-f$ gives
\[
 \E_Qe^{\operatorname{Fl}_N(f)}
   +\E_Qe^{-\operatorname{Fl}_N(f)}\le C.
\]
The exponential estimate in the proof of
\cref{lem:loc-moment-transfer} covers the finitely many remaining $N$
at fixed separation; in the moving regime choose
$R$ large at the outset. Exponential tail integration proves
\eqref{eq:bounded-tilted-moments}.

For analytic backgrounds apply this argument at width
$\varepsilon/2$, with $0<\varepsilon<\rho/8$, to the tests
$e^{-\varepsilon|k|}\cos(kx)$ and
$e^{-\varepsilon|k|}\sin(kx)$.
The resulting bound is
$\|\operatorname{Fl}_N(e^{ikx})\|_{L^p(Q)}
\le C_p e^{\varepsilon|k|}$.
The weighted Fourier sums in the proof of
\cref{thm:loc-moving-insertions} now give
$\sup_t\E_{Q_t}|R_{N,t}|\le C$. Thus
\eqref{eq:loc-moving-derivative} and \eqref{eq:loc-main-loop} have remainder
$O(N^{-1})$ after solving for the partition derivative,
which remains of that order after integration.
For the single root and argument weight with polynomial background,
use \cref{cor:loc-root-jump-background} in the same partition ratios,
with polynomial test functions of fixed degree.
The extra centered linear statistic is bounded by the resulting
moment estimates, which proves its improved error.
For $m=1$, compatibility in \cref{thm:loc-smooth-background} is
automatic, so the improved estimate gives
\eqref{eq:loc-full-positive-fh} with relative error $O(N^{-1})$.
This also proves the one insertion assertion of
\cref{thm:analytic-background} and completes its proof.
\end{proof}

\section{Further questions}
\label{sec:extensions}

The two principal theorems leave precise parts of the original
questions unresolved. We describe the estimates which are still needed.

\subsection{The full signed range at moving points}

For fixed separated points, \cref{thm:main-root} reaches every real
root exponent $q_j>-1$. \Cref{thm:moving-signed-two} gives a uniform
relative power error whenever $N\delta\to\infty$ on $\mathcal E_m$;
for two points this is the whole integrable range. The question
remains for three or more points with sum of the absolute values of the negative exponents at
least two. A sufficient estimate is
\[
 \sup_{\substack{\mathbf q\in K,\ n\delta\ge R}}
 \int W_{\mathbf q}(\boldsymbol\varphi)^p\,
                 d\mu_{n,\mathbf q}(\boldsymbol\varphi)<\infty
 \qquad\text{for some }p>1
\]
for every compact $K\subset(-1,\infty)^m$.
The pair density bound in \cref{prop:moving-pair-density} does not
imply this estimate: simultaneous
approach of several singularities with negative exponents requires corresponding
joint estimates. In the derivative argument the different phase
flows carry different positive Jacobian factors, so a pair estimate
cannot simply be iterated to give a full density bound uniform in
the separation.

\Cref{thm:analytic-background} also leaves a regularity question.
The transport remainder is controlled in an analytic Fourier norm.
The weaker regularity contemplated in \cite{FF2004} would require
a corresponding estimate under the singularly tilted law.
At $\beta=2$, Bradinoff and Duits
\cite[Definition~2.9 and Theorem~2.10]{BradinoffDuits2026}
give bounds uniform in the matrix size for higher coefficients in
logarithmic Toeplitz determinant expansions, under an admissibility
condition allowing slow Fourier decay. These estimates do not identify
the full Fisher--Hartwig constant.

For counting statistics, the lattice local and precise tail estimates
in the unpublished manuscript \cite{KLMN2026} at $\beta=2$ use control
of the tilted
characteristic function over a full Fourier period. At general $\beta$,
such control, together with the required uniformity in the window,
remains to be established before analogous individual counting
probabilities or growing windows of $\mathsf{Sine}_\beta$ are treated.

Extending the positive microscopic merging theorem
\cref{thm:positive-merging} to general mixed signs is another question:
when points coincide the sum of their exponents must be integrable, and
inverse products of the limiting random analytic function require
uniform integrability. Positive merging, which is sufficient for
the applications here, avoids this inverse product difficulty.

For genuine jump parameters $s_j=i\tau_j$, a natural further
boundary comes from the phase damping. More generally, at
$\mathbf s=\boldsymbol\sigma+i\boldsymbol\tau$ the real damping
relative to zero mode is
\[
 \frac2\beta\bigl(\|\boldsymbol\ell\|^2
                     -\boldsymbol\tau\cdot\boldsymbol\ell\bigr).
\]
It is positive for every nonzero integer vector with zero coordinate sum exactly when
$\max_j\tau_j-\min_j\tau_j<2$. Indeed, putting
$S=\sum_j(\ell_j)_+\ge1$, one has
$\|\boldsymbol\ell\|^2\ge2S$ and
$\boldsymbol\tau\cdot\boldsymbol\ell
 \le S(\max\tau-\min\tau)$; the converse follows from a mode
$\mathbf e_i-\mathbf e_j$ at the two extreme coordinates.
In the classical jump coordinates this is the width one condition
on the real parts of the Fisher--Hartwig parameters. The competing
representations at its boundary are classical
\cite{BasorTracy1991,DIK2011}. Positivity of these damping exponents
does not by itself extend our local complex theorem to that full
region: uniform estimates for the last coefficient and continuation across their
larger parameter domain are still needed.

\subsection{Weighted critical moments and the next term}

The leading critical coefficient is determined for every real $p>1$
and every complex spatial weight in $L^p$ by \cref{thm:main-critical}.
The mixed moment estimate for disjoint arcs in
\cref{prop:fractional-critical-additivity} identifies this local
coefficient even when $p\ge2$.

It remains to identify a finite term of order $N$ at noninteger
critical orders. The integer proof uses an integral near the diagonal;
the Palm argument at noninteger orders determines the leading
logarithm but does not provide its next term. For spatial weights, the result for the finite term in
\cref{thm:critical-second-term} concerns nonnegative weights satisfying
the stated Dini condition. We do not develop a formula for the finite term for
complex weights here.

There is an equivalent formulation through a growing integral under
a circular Jacobi law. The exact change of measure of
\cite[Proposition~2.3]{AN2026} gives, with $q=\sqrt{2\beta/p}$,
\[
 \E \mathcal Z_N(q)^p=M_N(pq)(2\pi N)^{1-p}
        \E_{{\rm CJ}(\beta,pq/2)} I_N^{p-1},\qquad
 I_N=\int_{-\pi N}^{\pi N}
      |\mathfrak q_N(e^{ix/N})|^q\,dx,
\]
where $\mathfrak q_N(z)=P_N(z)/P_N(1)$ under the rooted law.
Since $M_N(pq)\sim C_\beta(pq)N^p$ by \cref{prop:onepoint},
\cref{thm:main-critical} implies
\begin{equation}\label{eq:open-fractional-critical}
 \frac{\E_{{\rm CJ}(\beta,pq/2)}I_N^{p-1}}{\log N}
 \longrightarrow
 (2\pi)^{p-1}\frac{A_pC_\beta(q)^p}{C_\beta(pq)},
 \qquad A_p=\frac{p-1}{\Gamma(1-1/p)^p}.
\end{equation}
This growing interval conclusion follows from \cref{thm:main-critical};
compact convergence of the rooted random analytic functions alone
would not imply it.

For the Gaussian beta ensemble, Lambert and Paquette
\cite{LambertPaquette2025} give a bulk approximation of the
characteristic polynomial that retains both the microscopic
$\mathsf{Sine}_\beta$ fluctuations and the Gaussian fluctuations
on larger scales. Establishing the required joint moment estimates
would be a further step towards Fisher--Hartwig formulas in that setting.

\clearpage
\section*{Index of notation}
\addcontentsline{toc}{section}{Index of notation}

The table lists the symbols used in more than one section, with their
meaning and the place where they are defined. Notation confined to a
single proof is not listed.

\begin{longtable}{@{}>{\raggedright\arraybackslash}p{0.21\textwidth}>{\raggedright\arraybackslash}p{0.54\textwidth}>{\raggedright\arraybackslash}p{0.18\textwidth}@{}}
\toprule
Symbol & Meaning & Defined in\\
\midrule
\endfirsthead
\toprule
Symbol & Meaning & Defined in\\
\midrule
\endhead
\bottomrule
\endfoot
\multicolumn{3}{@{}l}{\emph{The ensemble and the main objects}}\\[2pt]
$\beta$, $b=\beta/2$ & inverse temperature and its half & \eqref{eq:ensemble}\\
$\T$, $\sigma$ & the circle $\R/(2\pi\Z)$ and normalized Haar measure $d\theta/(2\pi)$ & \cref{sec:intro-FH}\\
$\theta_1,\ldots,\theta_N$ & the particles (eigenangles) & \eqref{eq:ensemble}\\
$\vartheta_1,\ldots,\vartheta_m$ & the prescribed points & \eqref{eq:FN}\\
$P_N(z)$ & characteristic polynomial, normalized by $P_N(0)=1$ & \eqref{eq:PN}\\
$F_N(\mathbf q,\boldsymbol\vartheta)$ & joint moment $\E_\beta\prod_j|P_N(e^{i\vartheta_j})|^{q_j}$ & \eqref{eq:FN}\\
$M_N(q)$ & one point moment $F_N(q,0)=\E_\beta|P_N(1)|^q$ & \eqref{eq:FN}\\
$C_\beta(q)$ & one point constant, $M_N(q)\sim C_\beta(q)N^{q^2/(2\beta)}$ & \cref{prop:onepoint}, \eqref{eq:constantseries}\\
$\mathcal C_\beta(a,c)$, $\mathcal U_{\rm mix}$ & two variable one point constant and its domain; $C_\beta(q)=\mathcal C_\beta(q/2,q/2)$ & \cref{prop:onepoint}\\
$\mathcal A_\beta(\mathbf q,\boldsymbol\vartheta)$ & leading factor: one point constants and pair interactions & \eqref{eq:complex-root-leading-factor}\\
$\mathcal D_m$, $\mathcal E_m$ & total negative exponent below $1$ and below $2$, respectively; $q_j>-1$ is also required in $\mathcal E_m$ & \eqref{eq:intro-domain}\\
$\delta(\boldsymbol\vartheta)$, $\mathcal S_\delta$ & minimum chord distance capped at one; configurations with separation at least $\delta$ & \eqref{eq:intro-domain}\\
$A_K$, $\eta_K$, $g_\beta$, $M_K$ & total negative exponent on $K$, distance from the threshold, phase decay factor $\min\{4/\beta,1\}$, and $M_K=2+(A_K-1)_+$ & \cref{rem:admissible-exponents}, \cref{sec:full-negative-transfer}\\
$\varepsilon_K$, $\alpha_K$ & $1-A_K$ for $K\subset\mathcal D_m$, and the H\"older exponent $\varepsilon_K/(2m)$ & \eqref{eq:terminal-alpha}\\
$s_K$ & subset sum parameter giving the improved rate & \eqref{eq:terminal-subset-index}\\
$\mathcal Z_N(q;f)$, $\mathcal Z_N(q)$ & spatial integral $\int_\T f|P_N|^q\,d\sigma$ and its unweighted case & \cref{sec:intro-critical}\\
$A_p$, $d_p$ & critical coefficient $(p-1)/\Gamma(1-1/p)^p$ and $d_p=\Gamma(1-1/p)^{-p}$ & \cref{thm:main-critical}, \cref{sec:fractional-critical}\\
$q_k$, $A_k$ & critical exponent $\sqrt{2\beta/k}$ and coefficient at integer order $k$ & \eqref{eq:critical-charge-constant}\\
$B_{\beta,k}[f]$ & coefficient of $N$ in the weighted integer critical moment & \eqref{eq:weighted-critical-B}\\
\multicolumn{3}{@{}l}{\emph{Orthogonal polynomials and the coupling}}\\[2pt]
$\alpha_k$ & Verblunsky coefficients, independent with law \eqref{eq:coefficient-law} & \eqref{eq:coefficient-law}\\
$\mathcal F_n$, $\mathcal F_\infty$ & sigma-fields of the first $n$ coefficients, and of all of them & \cref{sec:coupling}\\
$\Phi_n$, $\Phi_n^*$ & monic orthogonal polynomial and its reversed polynomial & \eqref{eq:szego-recursion}\\
$B_n(z)$ & Blaschke product $z\Phi_n/\Phi_n^*$ & \cref{sec:coupling}\\
$\eta$, $t=e^{i\eta}$ & final independent angle and coefficient of modulus one & \cref{prop:coupled-ensemble}\\
$\psi_n(\theta)$ & lifted phase, $B_n(e^{i\theta})=e^{i\psi_n(\theta)}$ & \eqref{eq:phase-definition}\\
$\boldsymbol\psi_n$ & vector of phases at the prescribed points & \eqref{eq:terminal-expectation}\\
$Y_n(\mathbf q)$, $a(\mathbf q)$ & weight $\prod_j|\Phi_n^*(e^{i\vartheta_j})|^{q_j}$ and growth exponent $\frac1{2\beta}\sum_jq_j^2$ & \eqref{eq:root-weight}\\
$W_{\mathbf q}(\boldsymbol\varphi)$ & average over the final angle, $\int_\T\prod_j|1-e^{i(\eta+\varphi_j)}|^{q_j}\sigma(d\eta)$ & \eqref{eq:terminal-W}\\
$\nu_n(q)$, $D_\beta(q)$, $c(q)$ & OPUC one point moment, its constant, and the one point average $c(q)=C_\beta(q)/D_\beta(q)$ & \cref{lem:opuc-onepoint}\\
$T_n(\mathbf q,\boldsymbol\ell)$ & Fourier observable $Y_n(\mathbf q)e^{i\sum_j\ell_j\psi_n(\vartheta_j)}$ & \eqref{eq:Tdef}\\
$d(\boldsymbol\ell)$, $\lambda$, $\omega$, $L$ & damping exponent, net exponent $a(\mathbf q)-d(\boldsymbol\ell)$, rotation and size of the mode & \eqref{eq:exponents}\\
$\mu_{n,\mathbf q}$, $A_n(\mathbf q)$ & phase law weighted by $Y_n(\mathbf q)$, and $\E Y_n(\mathbf q)$ & \cref{sec:terminal}, \cref{sec:full-negative-transfer}\\
$f_{n,\mathbf q}$ & density of $\mu_{n,\mathbf q}$ & \eqref{eq:phase-law-decomposition}\\
$\mathcal V(\boldsymbol\vartheta)$ & Vandermonde product of the prescribed points & \cref{sec:full-negative-transfer}\\
\multicolumn{3}{@{}l}{\emph{Microscopic limit and critical moments}}\\[2pt]
$\mathcal{HP}_{\beta,\tau}$, $\xi^{\beta,\tau}$ & Hua--Pickrell process and its stochastic zeta function & \cref{sec:microscopic-limit}\\
$\mathcal H_{\beta,\mathbf q}(\mathbf x)$ & microscopic moment profile & \eqref{eq:merging-function}\\
$D_n$, $\mu_n$, $U_n$ & normalized OPUC density $|\Phi_n^*|^q/\nu_n(q)$, its measure and total mass & \cref{sec:fractional-critical}\\
$\chi_N$, $W_N$ & normalized density $|P_N|^q/M_N(q)$ and its total mass & \cref{sec:fractional-critical}\\
$\mu$, $\mu_\omega$ & limiting chaos measure and its realization on the coefficient sequence $\omega$ & \eqref{eq:fractional-chaos-martingale}, \cref{sec:fc-one-marked}\\
$\mathbb Q_{\boldsymbol\theta}$ & Palm distribution with marked points $\boldsymbol\theta$ & \cref{sec:fc-several}\\
$d_{ij}$, $K_t$, $K_\epsilon$ & chord distances, kernel $\prod_{i<j}d_{ij}^{-t}$, and clipped critical kernel $\prod_{i<j}(d_{ij}\vee\epsilon)^{-2/k}$ & \cref{sec:critical-reduction}\\
\multicolumn{3}{@{}l}{\emph{Mixed weights}}\\[2pt]
$L_N$, $X_N$, $J_N$ & $\log P_N$ on the circle, its real and imaginary parts & \cref{sec:onepoint}\\
$\mathcal F_N^{\rm rj}$, $\mathfrak a(\mathbf q,\mathbf s)$, $C^{\rm rj}_\beta(q,s)$ & mixed moment, its growth exponent and one point constant & \cref{sec:mixed-root-jump}\\
\end{longtable}


\begin{thebibliography}{99}

\bibitem{AlexanianKuznetsov2023}
S.~Alexanian and A.~Kuznetsov,
\emph{On the Barnes double gamma function},
Integral Transforms Spec. Funct. \textbf{34} (2023), no.~12, 891--914.
\href{https://doi.org/10.1080/10652469.2023.2238115}{doi:10.1080/10652469.2023.2238115}.
We use the normalization and sectorial formulas in
\href{https://arxiv.org/abs/2208.13876v2}{arXiv:2208.13876v2}.

\bibitem{Assiotis2022}
T.~Assiotis,
\emph{On the moments of the partition function of the C$\beta$E field},
J. Stat. Phys. \textbf{187} (2022), article~14.
\href{https://doi.org/10.1007/s10955-022-02906-6}{doi:10.1007/s10955-022-02906-6};
\href{https://arxiv.org/abs/2011.10323}{arXiv:2011.10323}.

\bibitem{AssiotisBKT2026}
T.~Assiotis,
\emph{BKT transition for $\mathsf{Sine}_\beta$},
preprint (2026),
\href{https://arxiv.org/abs/2610.13484}{arXiv:2610.13484},
version~1, 9~October~2026.

\bibitem{AssiotisKeating2021}
T.~Assiotis and J.~P.~Keating,
\emph{Moments of moments of characteristic polynomials of random
unitary matrices and lattice point counts},
Random Matrices Theory Appl. \textbf{10} (2021), no.~2, 2150019;
\href{https://doi.org/10.1142/S2010326321500192}{doi:10.1142/S2010326321500192};
\href{https://arxiv.org/abs/1905.06072v2}{arXiv:1905.06072v2}.
Theorem~1.1 is cited from this arXiv version.

\bibitem{AN2026}
T.~Assiotis and J.~Najnudel,
\emph{Moments of C$\beta$E field partition function, $\mathsf{Sine}_\beta$ correlations and stochastic zeta},
preprint (2026), \href{https://arxiv.org/abs/2602.08739}{arXiv:2602.08739}, version~1, 9~February~2026.

\bibitem{BaileyKeating2019}
E.~C.~Bailey and J.~P.~Keating,
\emph{On the moments of the moments of the characteristic polynomials
of random unitary matrices},
Comm. Math. Phys. \textbf{371} (2019), 689--726;
\href{https://arxiv.org/abs/1807.06605}{arXiv:1807.06605},
\href{https://doi.org/10.1007/s00220-019-03503-7}{doi:10.1007/s00220-019-03503-7}.

\bibitem{Basor1978}
E.~L.~Basor,
\emph{Asymptotic formulas for Toeplitz determinants},
Trans. Amer. Math. Soc. \textbf{239} (1978), 33--65.
\href{https://doi.org/10.1090/s0002-9947-1978-0493480-x}{doi:10.1090/s0002-9947-1978-0493480-x}.

\bibitem{BasorTracy1991}
E.~L.~Basor and C.~A.~Tracy,
\emph{The Fisher--Hartwig conjecture and generalizations},
Physica A \textbf{177} (1991), 167--173.
\href{https://doi.org/10.1016/0378-4371(91)90149-7}{doi:10.1016/0378-4371(91)90149-7}.

\bibitem{BWW2018}
N.~Berestycki, C.~Webb and M.~D.~Wong,
\emph{Random Hermitian matrices and Gaussian multiplicative chaos},
Probab. Theory Related Fields \textbf{172} (2018), 103--189.
\href{https://doi.org/10.1007/s00440-017-0806-9}{doi:10.1007/s00440-017-0806-9}.

\bibitem{BottcherSilbermann1985}
A.~B\"ottcher and B.~Silbermann,
\emph{Toeplitz matrices and determinants with Fisher--Hartwig symbols},
J. Funct. Anal. \textbf{63} (1985), 178--214.
\href{https://doi.org/10.1016/0022-1236(85)90085-0}{doi:10.1016/0022-1236(85)90085-0}.

\bibitem{BourgadeFalconet2022}
P.~Bourgade and H.~Falconet,
\emph{Liouville quantum gravity from random matrix dynamics},
preprint (2022),
\href{https://arxiv.org/abs/2206.03029v4}{arXiv:2206.03029v4}, revised 9~July~2025.

\bibitem{BHNY2008}
P.~Bourgade, C.~P.~Hughes, A.~Nikeghbali and M.~Yor,
\emph{The characteristic polynomial of a random unitary matrix:
a probabilistic approach},
Duke Math. J. \textbf{145} (2008), no.~1, 45--69;
\href{https://arxiv.org/abs/0706.0333}{arXiv:0706.0333},
\href{https://doi.org/10.1215/00127094-2008-046}{doi:10.1215/00127094-2008-046}.

\bibitem{BNN2013}
P.~Bourgade, J.~Najnudel, and A.~Nikeghbali,
\emph{A unitary extension of virtual permutations},
Int. Math. Res. Not. IMRN (2013), no.~18, 4101--4134;
\href{https://arxiv.org/abs/1102.2633}{arXiv:1102.2633},
\href{https://doi.org/10.1093/imrn/rns167}{doi:10.1093/imrn/rns167}.

\bibitem{BNR2010}
P.~Bourgade, A.~Nikeghbali, and A.~Rouault,
\emph{Circular Jacobi ensembles and deformed Verblunsky coefficients},
Int. Math. Res. Not. IMRN (2009), no.~23, 4357--4394;
\href{https://arxiv.org/abs/0804.4512}{arXiv:0804.4512}, version~3 (2010).
\href{https://doi.org/10.1093/imrn/rnp092}{doi:10.1093/imrn/rnp092}.

\bibitem{BradinoffDuits2026}
N.~Bradinoff and M.~Duits,
\emph{On Toeplitz determinants with slow Fourier decay},
preprint (2026),
\href{https://arxiv.org/abs/2608.13182v1}{arXiv:2608.13182v1},
13~August~2026.

\bibitem{Burkholder1973}
D.~L.~Burkholder,
\emph{Distribution function inequalities for martingales},
Ann. Probab. \textbf{1} (1973), no.~1, 19--42.
\href{https://doi.org/10.1214/aop/1176997023}{doi:10.1214/aop/1176997023}.

\bibitem{CarberyWright2001}
A.~Carbery and J.~Wright,
\emph{Distributional and $L^q$ norm inequalities for polynomials over
convex bodies in $\mathbb R^n$},
Math. Res. Lett. \textbf{8} (2001), 233--248.
\href{https://doi.org/10.4310/MRL.2001.v8.n3.a1}{doi:10.4310/MRL.2001.v8.n3.a1}.

\bibitem{Charlier2019}
C.~Charlier,
\emph{Asymptotics of Hankel determinants with a one-cut regular potential
and Fisher--Hartwig singularities},
Int. Math. Res. Not. IMRN (2019), no.~24, 7515--7576;
\href{https://doi.org/10.1093/imrn/rny009}{doi:10.1093/imrn/rny009};
\href{https://arxiv.org/abs/1706.03579}{arXiv:1706.03579}.

\bibitem{XIIIB}
R. Chhaibi, P.-L. M\'eliot, J. Najnudel and A. Nikeghbali, \emph{Large distance expansions of Sine beta correlation functions}, unpublished companion manuscript, 2026.

\bibitem{CN2026}
R.~Chhaibi and J.~Najnudel,
\emph{On the circle, Gaussian multiplicative chaos and Beta ensembles match exactly},
J. Eur. Math. Soc. \textbf{28} (2026), 1--77,
\href{https://doi.org/10.4171/JEMS/1700}{doi:10.4171/JEMS/1700};\newline\href{https://arxiv.org/abs/1904.00578}{arXiv:1904.00578}.

\bibitem{CNN2017}
R.~Chhaibi, J.~Najnudel, and A.~Nikeghbali,
\emph{The Circular Unitary Ensemble and the Riemann zeta function:
the microscopic landscape and a new approach to ratios},
Invent. Math. \textbf{207} (2017), 23--113;
\href{https://arxiv.org/abs/1410.1440}{arXiv:1410.1440},
\href{https://doi.org/10.1007/s00222-016-0669-1}{doi:10.1007/s00222-016-0669-1}.

\bibitem{ClaeysFahs2016}
T.~Claeys and B.~Fahs,
\emph{Random matrices with merging singularities and the Painlev\'e~V equation},
SIGMA \textbf{12} (2016), 031, 44~pp.;
\href{https://doi.org/10.3842/SIGMA.2016.031}{doi:10.3842/SIGMA.2016.031};
\href{https://arxiv.org/abs/1508.06734}{arXiv:1508.06734}.

\bibitem{ClaeysKrasovsky2015}
T.~Claeys and I.~Krasovsky,
\emph{Toeplitz determinants with merging singularities},
Duke Math. J. \textbf{164} (2015), no.~15, 2897--2987;
\href{https://doi.org/10.1215/00127094-3164897}{doi:10.1215/00127094-3164897}.
See also \href{https://arxiv.org/abs/1403.3639v4}{arXiv:1403.3639v4} (2022).

\bibitem{CFKRS2005}
J.~B.~Conrey, D.~W.~Farmer, J.~P.~Keating, M.~O.~Rubinstein and N.~C.~Snaith,
\emph{Integral moments of $L$-functions},
Proc. London Math. Soc. (3) \textbf{91} (2005), 33--104.
\href{https://doi.org/10.1112/S0024611504015175}{doi:10.1112/S0024611504015175}.

\bibitem{CostinLebowitz1995}
O.~Costin and J.~L.~Lebowitz,
\emph{Gaussian fluctuation in random matrices},
Phys. Rev. Lett. \textbf{75} (1995), 69--72.
\href{https://doi.org/10.1103/physrevlett.75.69}{doi:10.1103/physrevlett.75.69}.

\bibitem{DIK2011}
P.~Deift, A.~Its, and I.~Krasovsky,
\emph{Asymptotics of Toeplitz, Hankel, and Toeplitz+Hankel determinants with Fisher--Hartwig singularities},
Ann. of Math. (2) \textbf{174} (2011), no.~2, 1243--1299;
\href{https://arxiv.org/abs/0905.0443}{arXiv:0905.0443}, version~3 (2011).
\href{https://doi.org/10.4007/annals.2011.174.2.12}{doi:10.4007/annals.2011.174.2.12}.

\bibitem{DIK2013}
P.~Deift, A.~Its and I.~Krasovsky,
\emph{Toeplitz matrices and Toeplitz determinants under the impetus of the Ising model: some history and some recent results},
Comm. Pure Appl. Math. \textbf{66} (2013), 1360--1438.
\href{https://doi.org/10.1002/cpa.21467}{doi:10.1002/cpa.21467}.

\bibitem{DesrosiersLiu2014}
P.~Desrosiers and D.-Z.~Liu,
\emph{Asymptotics for products of characteristic polynomials in classical $\beta$-ensembles},
Constr. Approx. \textbf{39} (2014), 273--322.
\href{https://doi.org/10.1007/s00365-013-9206-2}{doi:10.1007/s00365-013-9206-2}.

\bibitem{Ehrhardt2001}
T.~Ehrhardt,
\emph{A status report on the asymptotic behavior of Toeplitz determinants with Fisher--Hartwig singularities},
in \emph{Recent Advances in Operator Theory},
A.~Dijksma, M.~A.~Kaashoek and A.~C.~M.~Ran (eds.),
Oper. Theory Adv. Appl. \textbf{124}, Birkh\"auser, Basel, 2001,
217--241.
\href{https://doi.org/10.1007/978-3-0348-8323-8_11}{doi:10.1007/978-3-0348-8323-8\_11}.

\bibitem{Fahs2021}
B.~Fahs,
\emph{Uniform asymptotics of Toeplitz determinants with
Fisher--Hartwig singularities},
Comm. Math. Phys. \textbf{383} (2021), 685--730;
\href{https://arxiv.org/abs/1909.07362}{arXiv:1909.07362}.
\href{https://doi.org/10.1007/s00220-021-03943-0}{doi:10.1007/s00220-021-03943-0}.

\bibitem{Fang2026}
W.~Fang,
\emph{Second-order fusion asymptotics for $\mathsf{Sine}_\beta$ correlation functions},
preprint (2026), \href{https://arxiv.org/abs/2608.23742}{arXiv:2608.23742}, version~1.

\bibitem{FangCritical2026}
W.~Fang,
\emph{Critical and subcritical fusion asymptotics for
$\mathsf{Sine}_\beta$ correlation functions},
preprint (2026),
\href{https://arxiv.org/abs/2609.07239}{arXiv:2609.07239},
version~1, 7~September~2026.

\bibitem{FangFourth2026}
W.~Fang,
\emph{Fourth-order fusion asymptotics for
$\mathsf{Sine}_\beta$ correlation functions},
preprint (2026),
\href{https://arxiv.org/abs/2609.11543}{arXiv:2609.11543},
version~1, 10~September~2026.

\bibitem{FangContinuity2026}
W.~Fang,
\emph{Joint continuity and Selberg--ODE equivalence for the
$\mathsf{Sine}_\beta$ pair correlation function},
preprint (2026),
\href{https://arxiv.org/abs/2609.17596}{arXiv:2609.17596},
version~1, 12~September~2026.

\bibitem{FengTianWei2023}
R.~Feng, G.~Tian and D.~Wei,
\emph{The Berry--Esseen theorem for circular $\beta$-ensemble},
Ann. Appl. Probab. \textbf{33} (2023), no.~6B, 5050--5070;
\href{https://arxiv.org/abs/1905.09448}{arXiv:1905.09448}, version~3.
\href{https://doi.org/10.1214/23-AAP1940}{doi:10.1214/23-AAP1940}.

\bibitem{FMN2016}
V.~F\'eray, P.-L.~M\'eliot and A.~Nikeghbali,
\emph{Mod-$\phi$ convergence: Normality zones and precise deviations},
SpringerBriefs in Probability and Mathematical Statistics,
Springer, Cham, 2016.
\href{https://doi.org/10.1007/978-3-319-46822-8}{doi:10.1007/978-3-319-46822-8};
\href{https://arxiv.org/abs/1304.2934}{arXiv:1304.2934}.

\bibitem{FMN2019}
V.~F\'eray, P.-L.~M\'eliot and A.~Nikeghbali,
\emph{Mod-$\phi$ convergence, II: Estimates on the speed of convergence},
in \emph{S\'eminaire de Probabilit\'es L},
Lecture Notes in Mathematics \textbf{2252},
Springer, Cham, 2019, 405--477.
\href{https://doi.org/10.1007/978-3-030-28535-7_15}{doi:10.1007/978-3-030-28535-7\_15};
\href{https://arxiv.org/abs/1705.10485v2}{arXiv:1705.10485v2}. The smoothing locator refers to this arXiv version.

\bibitem{FisherHartwig1968}
M.~E.~Fisher and R.~E.~Hartwig,
\emph{Toeplitz determinants: some applications, theorems, and conjectures},
Adv. Chem. Phys. \textbf{15} (1969), 333--353.
\href{https://doi.org/10.1002/9780470143605.ch18}{doi:10.1002/9780470143605.ch18}.

\bibitem{Forrester2010}
P.~J.~Forrester,
\emph{Log-Gases and Random Matrices},
London Mathematical Society Monographs Series, vol.~34,
Princeton University Press, Princeton, NJ, 2010.
\href{https://doi.org/10.1515/9781400835416}{doi:10.1515/9781400835416}.

\bibitem{FF2004}
P.~J.~Forrester and N.~E.~Frankel,
\emph{Applications and generalizations of Fisher--Hartwig asymptotics},
J. Math. Phys. \textbf{45} (2004), 2003--2028;
\href{https://arxiv.org/abs/math-ph/0401011v1}{arXiv:math-ph/0401011v1}.
The conjecture and equation locators refer to this arXiv version.
\href{https://doi.org/10.1063/1.1699484}{doi:10.1063/1.1699484}.

\bibitem{ForresterWarnaar2008}
P.~J.~Forrester and S.~O.~Warnaar,
\emph{The importance of the Selberg integral},
Bull. Amer. Math. Soc. (N.S.) \textbf{45} (2008), 489--534;
\href{https://arxiv.org/abs/0710.3981v1}{arXiv:0710.3981v1}.
The integral locators refer to this arXiv version.
\href{https://doi.org/10.1090/s0273-0979-08-01221-4}{doi:10.1090/s0273-0979-08-01221-4}.

\bibitem{FyodorovBouchaud2008}
Y.~V.~Fyodorov and J.-P.~Bouchaud,
\emph{Freezing and extreme-value statistics in a random energy model
with logarithmically correlated potential},
J. Phys. A: Math. Theor. \textbf{41} (2008), 372001;
\href{https://arxiv.org/abs/0805.0407}{arXiv:0805.0407},
\href{https://doi.org/10.1088/1751-8113/41/37/372001}{doi:10.1088/1751-8113/41/37/372001}.

\bibitem{FHK2012}
Y.~V.~Fyodorov, G.~A.~Hiary and J.~P.~Keating,
\emph{Freezing transition, characteristic polynomials of random matrices, and the Riemann zeta function},
Phys. Rev. Lett. \textbf{108} (2012), 170601.
\href{https://doi.org/10.1103/physrevlett.108.170601}{doi:10.1103/physrevlett.108.170601}.

\bibitem{FyodorovKeating2014}
Y.~V.~Fyodorov and J.~P.~Keating,
\emph{Freezing transitions and extreme values: random matrix theory,
$\zeta(1/2+it)$ and disordered landscapes},
Philos. Trans. R. Soc. A \textbf{372} (2014), 20120503;
\href{https://arxiv.org/abs/1211.6063}{arXiv:1211.6063}.
\href{https://doi.org/10.1098/rsta.2012.0503}{doi:10.1098/rsta.2012.0503}.

\bibitem{FyodorovLeDoussal2020}
Y.~V.~Fyodorov and P.~Le Doussal,
\emph{Statistics of extremes in eigenvalue-counting staircases},
Phys. Rev. Lett. \textbf{124} (2020), 210602;
\href{https://doi.org/10.1103/PhysRevLett.124.210602}{doi:10.1103/PhysRevLett.124.210602};
\href{https://arxiv.org/abs/2001.04135}{arXiv:2001.04135}.

\bibitem{HKO2001}
C.~P.~Hughes, J.~P.~Keating and N.~O'Connell,
\emph{On the characteristic polynomial of a random unitary matrix},
Comm. Math. Phys. \textbf{220} (2001), 429--451.
\href{https://doi.org/10.1007/s002200100453}{doi:10.1007/s002200100453}.

\bibitem{Johansson1988}
K.~Johansson,
\emph{On Szeg\H{o}'s asymptotic formula for Toeplitz determinants
and generalizations},
Bull. Sci. Math. (2) \textbf{112} (1988), no.~3, 257--304.

\bibitem{KeatingSnaith2000}
J.~P.~Keating and N.~C.~Snaith,
\emph{Random matrix theory and $\zeta(1/2+it)$},
Comm. Math. Phys. \textbf{214} (2000), no.~1, 57--89.
\href{https://doi.org/10.1007/s002200000261}{doi:10.1007/s002200000261}.

\bibitem{KW2022}
J.~P.~Keating and M.~D.~Wong,
\emph{On the critical--subcritical moments of moments of random
characteristic polynomials: a GMC perspective},
Comm. Math. Phys. \textbf{394} (2022), 1247--1301;
\href{https://arxiv.org/abs/2012.15851}{arXiv:2012.15851},
\href{https://doi.org/10.1007/s00220-022-04429-3}{doi:10.1007/s00220-022-04429-3}.

\bibitem{Keles2025}
A.~Keles,
\emph{Non-intersecting Brownian Motions and Gaussian Multiplicative Chaos},
preprint (2025),
\href{https://arxiv.org/abs/2508.11505v1}{arXiv:2508.11505v1},
15~August~2025.

\bibitem{KLMN2026}
M.~Khodiakova, T.~Leh\'ericy, P.-L.~M\'eliot and A.~Nikeghbali,
\emph{Precise deviations via Fourier local limit theorems:
lattice effects, circular ensembles and rooted graphs},
unpublished manuscript, 17 October 2026.

\bibitem{Killip2008}
R.~Killip,
\emph{Gaussian fluctuations for $\beta$ ensembles},
Int. Math. Res. Not. IMRN (2008), rnn007, 19 pp.;
\href{https://arxiv.org/abs/math/0703140}{arXiv:math/0703140},
\href{https://doi.org/10.1093/imrn/rnn007}{doi:10.1093/imrn/rnn007}.

\bibitem{KK2017}
R.~Killip and R.~Kozhan,
\emph{Matrix models and eigenvalue statistics for truncations of
classical ensembles of random unitary matrices},
Comm. Math. Phys. \textbf{349} (2017), 991--1027;
\href{https://arxiv.org/abs/1501.05160}{arXiv:1501.05160}.
\href{https://doi.org/10.1007/s00220-016-2658-z}{doi:10.1007/s00220-016-2658-z}.

\bibitem{KN2004}
R.~Killip and I.~Nenciu,
\emph{Matrix models for circular ensembles},
Int. Math. Res. Not. IMRN (2004), no.~50, 2665--2701;
\href{https://arxiv.org/abs/math/0410034}{arXiv:math/0410034}.
\href{https://doi.org/10.1155/S1073792804141597}{doi:10.1155/S1073792804141597}.

\bibitem{Lambert2021}
G.~Lambert,
\emph{Mesoscopic central limit theorem for the circular beta-ensembles and applications},
Electron. J. Probab. \textbf{26} (2021), paper no.~7;
\href{https://arxiv.org/abs/1902.06611}{arXiv:1902.06611},
\href{https://doi.org/10.1214/20-EJP559}{doi:10.1214/20-EJP559}.

\bibitem{LN2026}
G.~Lambert and J.~Najnudel,
\emph{Subcritical multiplicative chaos and the characteristic polynomial
of the C$\beta$E},
Probab. Theory Related Fields \textbf{196} (2026), 2103--2154;
\href{https://doi.org/10.1007/s00440-026-01501-0}{doi:10.1007/s00440-026-01501-0};
\href{https://arxiv.org/abs/2407.19817}{arXiv:2407.19817}.

\bibitem{LOS2018}
G.~Lambert, D.~Ostrovsky and N.~Simm,
\emph{Subcritical multiplicative chaos for regularized counting statistics from random matrix theory},
Comm. Math. Phys. \textbf{360} (2018), 1--54.
\href{https://doi.org/10.1007/s00220-018-3130-z}{doi:10.1007/s00220-018-3130-z}.

\bibitem{LambertPaquette2025}
G.~Lambert and E.~Paquette,
\emph{Bulk asymptotics of the Gaussian $\beta$-ensemble characteristic polynomial},
preprint (2025),
\href{https://arxiv.org/abs/2508.01458}{arXiv:2508.01458},
version~1, 2~August~2025.

\bibitem{LiValko2022}
Y.~Li and B.~Valk\'o,
\emph{Operator level limit of the circular Jacobi $\beta$-ensemble},
Random Matrices Theory Appl. \textbf{11} (2022), no.~4, 2250043;
\href{https://arxiv.org/abs/2108.11039}{arXiv:2108.11039}.
\href{https://doi.org/10.1142/s2010326322500435}{doi:10.1142/s2010326322500435}.

\bibitem{Liu2017}
D.-Z.~Liu,
\emph{Limits for circular Jacobi beta-ensembles},
J. Approx. Theory \textbf{215} (2017), 40--67;
\href{https://arxiv.org/abs/1408.0486}{arXiv:1408.0486},
\href{https://doi.org/10.1016/j.jat.2016.11.005}{doi:10.1016/j.jat.2016.11.005}.

\bibitem{Matsumoto2007}
S.~Matsumoto,
\emph{Moments of characteristic polynomials for compact symmetric spaces and Jack polynomials},
J. Phys. A \textbf{40} (2007), 13567--13586.
\href{https://doi.org/10.1088/1751-8113/40/45/006}{doi:10.1088/1751-8113/40/45/006}.

\bibitem{MeliotNikeghbali2022}
P.-L.~M\'eliot and A.~Nikeghbali,
\emph{On the precise deviations of the characteristic polynomial
of a random matrix},
Random Matrices Theory Appl. \textbf{15} (2026), no.~3, 2550019;
\href{https://doi.org/10.1142/S2010326325500194}{doi:10.1142/S2010326325500194};
\href{https://arxiv.org/abs/2202.08336}{arXiv:2202.08336}.

\bibitem{NSW2020}
M.~Nikula, E.~Saksman and C.~Webb,
\emph{Multiplicative chaos and the characteristic polynomial of the CUE:
the $L^1$-phase},
Trans. Amer. Math. Soc. \textbf{373} (2020), no.~6, 3905--3965.
\href{https://doi.org/10.1090/tran/8020}{doi:10.1090/tran/8020}.

\bibitem{Remy2020}
G.~Remy,
\emph{The Fyodorov--Bouchaud formula and Liouville conformal field theory},
Duke Math. J. \textbf{169} (2020), no.~1, 177--211;
\href{https://arxiv.org/abs/1710.06897}{arXiv:1710.06897},
\href{https://doi.org/10.1215/00127094-2019-0045}{doi:10.1215/00127094-2019-0045}.

\bibitem{RhodesVargas2019}
R.~Rhodes and V.~Vargas,
\emph{The tail expansion of Gaussian multiplicative chaos and the
Liouville reflection coefficient},
Ann. Probab. \textbf{47} (2019), no.~5, 3082--3107;
\href{https://arxiv.org/abs/1710.02096}{arXiv:1710.02096},
\href{https://doi.org/10.1214/18-AOP1333}{doi:10.1214/18-AOP1333}.

\bibitem{Ruijsenaars2000}
S.~N.~M.~Ruijsenaars,
\emph{On Barnes' multiple zeta and gamma functions},
Adv. Math. \textbf{156} (2000), 107--132.
\href{https://doi.org/10.1006/aima.2000.1946}{doi:10.1006/aima.2000.1946}.

\bibitem{Ryckman2008}
E.~Ryckman,
\emph{Linear statistics of point processes via orthogonal polynomials},
J. Stat. Phys. \textbf{132} (2008), 473--486;
\href{https://arxiv.org/abs/0805.3516}{arXiv:0805.3516},
\href{https://doi.org/10.1007/s10955-008-9564-5}{doi:10.1007/s10955-008-9564-5}.

\bibitem{ShenForrester2025}
B.-J.~Shen and P.~J.~Forrester,
\emph{Moments of characteristic polynomials for classical $\beta$ ensembles},
preprint (2025),
\href{https://arxiv.org/abs/2502.07142v2}{arXiv:2502.07142v2},
revised 29~June~2025.

\bibitem{Simon2005}
B.~Simon,
\emph{Orthogonal polynomials on the unit circle. Part 1: Classical theory},
American Mathematical Society Colloquium Publications, vol.~54,
American Mathematical Society, Providence, RI, 2005.
\href{https://doi.org/10.1090/coll054.1}{doi:10.1090/coll054.1}.

\bibitem{SmithEtAl2021}
N.~R.~Smith, P.~Le Doussal, S.~N.~Majumdar and G.~Schehr,
\emph{Full counting statistics for interacting trapped fermions},
SciPost Phys. \textbf{11} (2021), no.~6, 110;
\href{https://doi.org/10.21468/SciPostPhys.11.6.110}{doi:10.21468/SciPostPhys.11.6.110};
\href{https://arxiv.org/abs/2106.05014v2}{arXiv:2106.05014v2}.
Equation~(56) is cited from this arXiv version.

\bibitem{ValkoVirag2022}
B.~Valk\'o and B.~Vir\'ag,
\emph{The many faces of the stochastic zeta function},
Geom. Funct. Anal. \textbf{32} (2022), 1160--1231;
\href{https://arxiv.org/abs/2009.04670}{arXiv:2009.04670}.
\href{https://doi.org/10.1007/s00039-022-00613-8}{doi:10.1007/s00039-022-00613-8}.

\bibitem{Webb2015}
C.~Webb,
\emph{The characteristic polynomial of a random unitary matrix and
Gaussian multiplicative chaos---The $L^2$-phase},
Electron. J. Probab. \textbf{20} (2015), no.~104, 1--21.
\href{https://doi.org/10.1214/EJP.v20-4296}{doi:10.1214/EJP.v20-4296}.

\bibitem{Widom1973}
H.~Widom,
\emph{Toeplitz determinants with singular generating functions},
Amer. J. Math. \textbf{95} (1973), 333--383.
\href{https://doi.org/10.2307/2373789}{doi:10.2307/2373789}.

\bibitem{Wieand2002}
K.~Wieand,
\emph{Eigenvalue distributions of random unitary matrices},
Probab. Theory Related Fields \textbf{123} (2002), 202--224.
\href{https://doi.org/10.1007/s004400100186}{doi:10.1007/s004400100186}.

\bibitem{Wong2020}
M.~D.~Wong,
\emph{Universal tail profile of Gaussian multiplicative chaos},
Probab. Theory Related Fields \textbf{177} (2020), 711--746;
\href{https://arxiv.org/abs/1902.04054}{arXiv:1902.04054},
\href{https://doi.org/10.1007/s00440-020-00960-3}{doi:10.1007/s00440-020-00960-3}.

\end{thebibliography}
\end{document}